universite PARIS-SACLAY

# Forking in valued fields and related structures

*Déviation dans les corps valués et les structures apparentées*

**Thèse de doctorat de l'université Paris-Saclay**

École doctorale n° 574 : mathématiques Hadamard (EDMH)
Spécialité de doctorat : Mathématiques fondamentales
Graduate School : Mathématiques. Référent : Faculté des sciences d'Orsay

Thèse préparée dans l'unité de recherche **LMO UMR 8628, Laboratoire de mathématiques d'Orsay (Université Paris-Saclay, CNRS)**, sous la direction de **Elisabeth BOUSCAREN**, directrice de recherche émérite CNRS, la co-direction de **Silvain RIDEAU-KIKUCHI**, chargé de recherche CNRS (Ecole normale supérieure - Paris sciences et lettres)

**Thèse soutenue à Paris-Saclay, le 03 juillet 2024, par**

## Akash HOSSAIN

### Composition du jury

Membres du jury avec voix délibérative

| | |
|---|---|
| **Jean-Benoît BOST**<br>Professeur à l'université Paris-Saclay | Président |
| **H.Dugald MACPHERSON**<br>Professeur à l'université de Leeds | Rapporteur |
| **Amador MARTIN-PIZARRO**<br>Professeur à l'université de Fribourg | Rapporteur & Examinateur |
| **Sylvy ANSCOMBE**<br>Maîtresse de conférence à l'université Paris-Cité | Examinatrice |
| **Martin HILS**<br>Professeur à l'université de Münster | Examinateur |
| **Thomas SCANLON**<br>Professeur à l'université de Californie à Berkeley | Examinateur |

THESE DE DOCTORAT

ÉCOLE DOCTORALE
de mathématiques
Hadamard (EDMH)

**Titre :** Déviation dans les corps valués et les structures apparentées

**Résumé :** Cette thèse est une contribution à la théorie des modèles des corps valués. On étudie la déviation dans les corps valués, ainsi que certains de leurs réduits. On s'intéresse particulièrement aux corps pseudo-locaux, les ultraproduits de caractéristique résiduelle nulle des corps valués p-adiques.

Nous considérons d'abord les groupes des valeurs des corps valués qui nous intéressent, les groupes Abéliens ordonnés réguliers. Nous y établissons une description géométrique de la déviation, ainsi qu'une classification détaillée des extensions globales non-déviantes ou invariantes d'un type donné.

Nous démontrons ensuite des principes d'Ax-Kochen-Ershov pour la division et la déviation dans la théorie resplendissante des expansions de suites exactes courtes pures de structures Abéliennes, telles qu'étudiées dans l'article sur la distalité d'Aschenbrenner-Chernikov-Gehret-Ziegler. En particulier, nos résultats s'appliquent aux groupes des termes dominants des (expansions de) corps valués.

Pour finir, nous donnons diverses conditions suffisantes pour qu'un ensemble de paramètres soit une base d'extension dans un corps valué Hensélien de caractéristique résiduelle nulle. En particulier, nous démontrons que la déviation coïncide avec la division dans les corps pseudo-locaux de caractéristique résiduelle nulle.

Nous discutons aussi des résultats de Ealy-Haskell-Simon sur la déviation pour les extensions séparées de corps valués Henséliens de caractéristique résiduelle nulle. Nous contribuons à la question en démontrant que, dans le cas d'une extension Abhyankar, et avec quelques hypothèses supplémentaires, la non-déviation d'un type dans un corps pseudo-local implique l'existence d'une mesure de Keisler globale invariante dont le support contient ce type, à l'instar des corps pseudo-finis.

ÉCOLE DOCTORALE
de mathématiques
Hadamard (EDMH)

**Title :** Forking in valued fields and related structures

**Abstract :** This thesis is a contribution to the model theory of valued fields. We study forking in valued fields and some of their reducts. We focus particularly on pseudo-local fields, the ultraproducts of residue characteristic zero of the p-adic valued fields.

First, we look at the value groups of the valued fields we are interested in, the regular ordered Abelian groups. We establish for these ordered groups a geometric description of forking, as well as a full classification of the global extensions of a given type which are non-forking or invariant.

Then, we prove an Ax-Kochen-Ershov principle for forking and dividing in expansions of pure short exact sequences of Abelian structures, as studied by Aschenbrenner-Chernikov-Gehret-Ziegler in their article about distality. This setting applies in particular to the leading-term structure of (expansions of) valued fields.

Lastly, we give various sufficient conditions for a parameter set in a Henselian valued field of residue characteristic zero to be an extension base. In particular, we show that forking equals dividing in pseudo-local fields of residue characteristic zero.

Additionally, we discuss results by Ealy-Haskell-Simon on forking in separated extensions of Henselian valued fields of residue characteristic zero. We contribute to the question in the setting of Abhyankar extensions, where we show that, with some additional conditions, if a type in a pseudo-local field does not fork, then there exists some global invariant Keisler measure whose support contains that type. This behavior is well-known in pseudo-finite fields.

Mathématiques Orsay

1794

ENS | PSL

# Remerciements

Je tiens d'abord à remercier mon co-directeur et mentor Silvain RIDEAU-KIKUCHI pour son investissement et sa bienveillance. Silvain a lu la quasi-totalité de ce que j'ai produit durant mon doctorat, et m'a prodigué des conseils indispensables, autant sur les standards d'écriture que sur des subtilités mathématiques plus techniques. Mon étude des groupes ordonnés, qui aborde des questions très délicates, a nécessité plusieurs versions intermédiaires pour atteindre une qualité qui me satisfaisait, et Silvain a lu chacune d'entre elles avec la plus grande patience. Je le remercie également pour sa disponibilité, y compris pendant la période critique des confinements de première année de doctorat, où les contacts de chacun étaient très restreints.

Je remercie ma directrice de thèse Elisabeth BOUSCAREN, qui a abordé mon travail avec beaucoup de recul, me faisant ainsi des critiques très pertinentes qui m'ont aidé à améliorer la clarté de mes productions. Son large domaine d'expertise sur la théorie des modèles, et ses conseils expérimentés sur les codes du milieu académique, ont été précieux et utiles.

I am very grateful to Dugald MACPHERSON and Amador MARTIN-PIZARRO for reviewing my thesis. I greatly appreciate their suggestions and their comments.

I am thankful to Sylvy ANSCOMBE, Jean-Benoît BOST, Martin HILS and Tom SCANLON for participating in the jury. Je remercie tout particulièrement Martin pour son accueil chaleureux lors de ma visite à Münster, pendant laquelle j'ai développé les fondations de mes travaux sur les suites exactes courtes.

Durant ma thèse, j'ai eu le privilège de bénéficier d'un soutien de l'ANR Geomod, qui a financé ma participation à de nombreuses conférences, une école d'été, et ma visite à Münster.

J'adresse un grand merci aux personnels de l'Ecole doctorale de mathématiques Hadamard, du Laboratoire de mathématiques d'Orsay, de l'Université Paris-Saclay et du Département de mathématiques et applications de l'Ecole normale supérieure - Paris sciences et lettres ; pour leurs services administratifs, leur suivi personnalisé de ma formation, et leur organisation de mes missions d'enseignement. Je remercie entre autres Séverine SIMON, Marie-Christine MYOUPO, Mathilde ROUSSEAU, Clotilde D'EPENOUX, Olga JAMIN, Claire LE POULENNEC, Stéphane NONNENMACHER, Olivier SCHIFFMANN, Frédéric BOURGEOIS, Benjamin SCHRAEN, Hans RUGH et Filipa CAETANO.

Je salue mes camarades théoriciens des modèles de Paris, Pablo, Paulo, Xavier, et plus particulièrement Paul et Stefan qui se sont impliqués activement dans les groupes de travail que nous avons organisé. Le nombre de jeunes théoriciens des modèles de Paris évolue de manière très instable. J'étais essentiellement le seul au début de ma thèse, bien que je connaissais déjà Paul qui en avait encore pour deux ans avant de commencer la sienne. Le moment où les autres m'ont rejoint a été un point de rupture positif dans ma vie académique, d'autant plus qu'il fut précédé des périodes de confinement. Je suis content de savoir que beaucoup de jeunes camarades les rejoindront l'année prochaine.

I wish good luck to my fellow graduate students from Orsay, in particular my office comrade Shang, who will defend his thesis not long after me.

Pour finir, je serais bien ingrat de ne pas remercier mes proches qui m'ont toujours soutenu. Le fait d'avoir une famille investie et des amis loyaux est un privilège invisible qui a beaucoup plus d'influence sur notre succès que ce que l'on perçoit. J'ai une pensée particulière pour mon père qui nous a quittés il y a bientôt six ans.

# Contents

# Introduction

Ce travail s'inscrit dans le cadre de la théorie des modèles, et plus précisément la théorie des modèles des corps valués.

La théorie des modèles étudie les structures mathématiques à partir des propriétés des ensembles qui y sont définissables, au sens de la logique du premier ordre. Elle a établi depuis de nombreuses années des liens solides avec l'algèbre, la théorie des nombres et la géométrie algébrique, et y a eu de nombreuses applications. La théorie des modèles des corps valués est exemplaire de ce point de vue et particulièrement riche en résultats. De premiers travaux datant des années soixante ([AK65a], [AK65b], [AK66]) ont montré la puissance des outils modèle-théoriques classiques pour obtenir des résultats algébriques sur les corps valués algébriquement clos ou plus généralement Henséliens. Plus récemment à partir des travaux de Haskell-Hrushovski-Macpherson ([HHM05]), ce sont les outils plus sophistiqués développés par la théorie des modèles pure moderne (théorie des modèles géométrique et néo-stabilité) qui ont été utilisés pour étudier des questions de type différent (élimination des quotients imaginaires, par exemple) et ont ensuite donné lieu à de nouvelles applications (aux espaces de Berkovich [HL16], à l'intégration motivique [HK06]. . . ).

Dans les travaux exposés ici, nous portons une attention toute particulière à la notion modèle-théorique de *déviation*, une généralisation des diverses notions algébriques de dépendance classiquement étudiées en algèbre et en géométrie.

# Bref historique et motivations

## Théorie des modèles des corps valués

La série d'articles qui a éveillé les intérêts de la communauté pour la théorie des modèles des corps valués est sans doute celle d'Ax-Kochen ([AK65a], [AK65b], [AK66]). Ces textes établissent avant tout un résultat de nature arithmétique : il est démontré qu'étant donné une propriété du premier ordre $P$, $P$ est vraie dans les corps $p$-adiques $\mathbb{Q}_p$ pour $p$ suffisamment grand, si et seulement si $P$ est vraie dans les corps de séries de Laurent $\mathbb{F}_p((t))$ pour $p$ suffisamment grand. En résumé, les comportements asymptotiques au premier ordre des corps valués $p$-adiques et des $\mathbb{F}_p((t))$ sont exactement les mêmes. Ce résultat a permis en particulier de démontrer un affaiblissement asymptotique de la conjecture (fausse) d'Artin sur le fait que les corps $p$-adiques soient $C_2$ : étant donné que les $\mathbb{F}_p((t))$ sont $C_2$, si on fixe un degré $d$, pour $p$ suffisamment grand, tout polynôme de degré $d$ avec strictement plus que $d^2$ variables a un zéro dans $\mathbb{Q}_p$.

La propriété modèle-théorique remarquable sur laquelle ces applications à la théorie des nombres s'appuient est ce que l'on appelle aujourd'hui le théorème d'Ax-Kochen et Ershov, démontré indépendamment par Ershov dans ([Ers65]) :

**Théorème 1** ([AK65a], [AK65b], [AK66], [Ers65]). *Soient $M$ et $N$ deux corps valués Henséliens de caractéristique résiduelle nulle. Alors $M$ et $N$ sont élémentairement équivalents si et seulement si c'est le cas de leurs groupes des valeurs et de leurs corps des résidus respectifs.*

Ce résultat fondamental est une instance d'un principe plus général, que l'on appelle le *principe d'Ax-Kochen-Ershov* (AKE), selon lequel il suffit, pour comprendre une notion modèle-théorique dans un corps valué Hensélien de caractéristique résiduelle nulle $K$, de la comprendre dans son groupe des valeurs $\Gamma(K)$ et son corps des résidus $k(K)$. Ce principe consiste donc à réduire un problème sur les corps valués à une question plus simple, sur des structures moins riches et mieux connues. En l'occurence, le théorème d'Ax-Kochen-Ershov est un principe d'Ax-Kochen-Ershov pour *l'équivalence élémentaire*. En réalité, le détail de leurs preuves établit un principe d'Ax-Kochen-Ershov pour *l'élimination relative des quantificateurs*.

Il existe en fait un principe intermédiaire qui n'est pas explicite dans leur preuve originale, mais qui a été mis en évidence par Basarab ([Bas91]). Ce principe consiste à réduire une question sur $K$ à une question sur la structure du premier ordre composée de la suite exacte courte :

$$1 \longrightarrow k(K) \longrightarrow \mathrm{RV}(K) \longrightarrow \Gamma(K) \longrightarrow 0$$

où $\mathrm{RV}(K)$ est le *groupe des termes dominants* de $K$, le quotient de $K^*$ par la boule ouverte de rayon 0 autour de 1. Pour passer d'une question sur (l'expansion à toute la structure de $k(K)$ et $\Gamma(K)$ de) cette suite exacte courte à deux questions sur $k(K)$ et $\Gamma(K)$, on utilise essentiellement le fait que toute suite exacte courte est scindée dans une extension élémentaire, et que le membre au centre de la suite joue donc le rôle d'un produit des deux autres membres. Nous nous devons également de mentionner le travail contemporain d'Aschenbrenner-Chernikov-Gehret-Ziegler ([ACGZ22], Section 4) qui formalise et généralise ce principe à la *théorie respendissante des expansions de suites exactes courtes pures de structures Abéliennes.*

Les articles fondateurs d'Ax-Kochen et Ershov ont été suivis de très nombreux travaux de la communauté des théoriciens des modèles, et la théorie des modèles des corps valués s'est imposée comme une branche importante et fructueuse du sujet. D'autres instances très variées du principe d'Ax-Kochen-Ershov ont été établies depuis, pour un certain nombre de notions modèle-théoriques, et ont donné lieu à de nouvelles applications à d'autres domaines des mathématiques. Nous pouvons mentionner le livre de Haskell-Hrushovski-Macpherson ([HHM05]) sur les corps non-trivialement valués algébriquement clos (ACVF), sur lequel s'appuient Hrushovski-Loeser ([HL16]) dans leurs résultats sur les espaces de Berkovich. L'étude des ensembles définissables dans ACVF joue également un rôle important dans la contribution de Hrushvski-Kazhdan ([HK06]) à l'intégration motivique.

C'est dans ce courant que s'inscrit notre travail : sachant que l'étude des propriétés du premier ordre des corps valués à un fort potentiel d'être appliquée à d'autres problèmes complexes de théorie des nombres ou de géométrie algébrique, l'objectif premier de cette thèse est de mieux comprendre la théorie des modèles des corps valués.

Nous travaillerons toujours dans un cadre où les corps valués seront Henséliens et de caractéristique résiduelle nulle. Une grande partie des résultats de la litérature s'applique également aux corps Henséliens de caractéristique *mixte* (à l'instar de $\mathbb{Q}_p$, le corps résiduel est de caractéristique positive, et le domaine de la valuation est de caractéristique nulle), mais des hypothèses supplémentaires doivent être faites, et les preuves en caractéristique mixte sont souvent adaptées à partir de celles en caractéristique résiduelle nulle. La généralisation de nos résultats à la caractéristique mixte devra donc se faire dans un second temps. Quand aux corps valués de caractéristique pure positive (à l'instar de $\mathbb{F}_p((t))$), ils sont encore très peu compris en théorie des modèles, et de nombreuses questions élémentaires (telle qu'une axiomatisation raisonnable de la théorie de $\mathbb{F}_p((t))$) sont encore ouvertes. On peut en dire autant des corps valués non-Henséliens, leur théorie des modèle est très sauvage. Parmi les corps Henséliens de caractéristique résiduelle nulle, nous étudierons avec attention ceux qui ont joué un rôle important dans la litérature, c'est à dire ceux qui sont algébriquement clos ou *pseudo-locaux*, les ultraproduits de corps valués locaux. Les théoriciens des modèles qui lisent cette introduction auront peut-être deviné que le théorème d'Ax-Kochen fait intervenir les ultraproduits des $\mathbb{Q}_p$ et ceux des $\mathbb{F}_p((t))$. Les corps pseudo-locaux sont également au coeur des développements contemporains sur l'intégration motivique. Bien que la mise en avant de ces deux classes de corps valués ne soit pas explicite dans nos résultats, les hypothèses que l'on y prend sont conçues pour être au moins vraies dans les corps pseudo-locaux, tout en étant aussi générales que possible.

## Déviation

De nombreux travaux ont été réalisés en théorie des modèles des corps valués, mais la *déviation* est une notion fondamentale de théorie des modèles contemporaine qui, bien qu'étant indispensable pour une étude fine des ensembles définissables, est encore très peu comprise dans les corps valués. C'est sur cette notion que nous nous concentrons dans ce travail.

Une manière naturelle d'étudier les propriétés du premier ordre d'une structure algébrique est d'interpréter ses ensembles définissables comme des objets *géométriques*. Par exemple, on voudrait pouvoir parler de *dimension* d'un ensemble définissable, et de *points génériques* de cet ensemble, d'une manière qui interagit naturellement avec la théorie des modèles. Ces notions de dimension et de point générique prennent un sens naturel quand on travaille dans une structure algébrique usuelle, mais la théorie des modèles porte sur les structures du premier ordre *quelconques*, qui s'écrivent comme la donnée d'une collection d'ensembles (les *domaines*), et d'une collection d'applications et de relations entre ces domaines. C'est un défi assez compliqué que de faire de la géométrie dans un contexte aussi abstrait. Les définitions que l'on peut proposer devront nécessairement être très abstraites, et il faudrait aussi qu'elles unifient les notions intuitives (de dimension, etc...) des structures algébriques usuelles concrètes. Les développements modernes de la *théorie géométrique de la stabilité* apportent des solutions à ces problèmes. On y définit de diverses manières abstraites plusieurs notions modèle-théoriques qui permettent de formaliser nos intuitions géométriques dans n'importe quelle structure, et il est démontré que ces notions coïncident dans une certaine classe de structures, dites *stables*. La classe des structures stables, à l'instar des autres classes de complexité de la théorie des modèles introduites par Shelah, est définie comme la classe des structures dans lesquelles il est impossible de représenter, à l'aide de formules de la logique du premier ordre, certaines configurations combinatoires complexes. Dans notre cas, les configurations à éviter sont les ordres :

**Définition 2.** Une théorie $T$ est *stable* quand dans tout modèle $M$ de $T$, pour toute formule $\phi(x, y)$ sans paramètres, il n'existe aucune suite $(a_i b_i)_{i<\omega}$ d'éléments de $M$ telle que $\phi(a_i, b_j)$ soit vraie dans $M$ si et seulement si $i \leqslant j$.

Comme nous l'avons mentionné, les théories stables ont de très bonnes propriétés, notamment il y a une notion uniforme et robuste de dimension. Cependant, bien que certaines structures algébriques usuelles soient stables, comme les corps séparablement clos, les groupes Abéliens ou les espaces vectoriels, un bon nombre de structures fondamentales en mathématiques sont instables, notamment le corps des réels, n'importe quel groupe Abélien ordonné non-trivial, et donc n'importe quel corps non-trivialement valué. En fait, il découle facilement de la définition ci-dessus que n'importe quelle expansion d'un ordre total infini est instable. Pour contourner ce problème, il faut considérer d'autres classes de complexité introduites par Shelah, qui vont être plus grandes, mais aussi moins régulières. Sans donner toutes les définitions techniques, la propriété de l'ordre (OP), qui est la négation de la définition ci-dessus, peut être renforcée en la propriété de l'arbre (TP), la propriété d'indépendance (IP), ou ce qu'on appelle la seconde propriété de l'arbre ($TP_2$). La *théorie de la classification* de Shelah fournit d'autres propriétés qui définissent d'autres classes de complexité, mais nous nous en tiendrons à celles-là car elles sont pertinentes pour notre travail. Les classes que l'on a ainsi définies sont plus larges et contiennent bien plus de structures algébriques usuelles : par exemple, le corps des réels, les corps valués algébriquement clos et les corps valués $p$-adiques ne satisfont pas IP (ont dit qu'ils sont NIP), les corps algébriquement clos avec un automorphisme générique et les corps des résidus des corps pseudo-locaux sont NTP, tandis que les corps valués pseudo-locaux sont $NTP_2$. Il est à noter que les théories stables sont exactement celles dans l'intersection de NIP et de NTP, tandis que les théories $NTP_2$ contiennent strictement la réunion de NIP et de NTP.

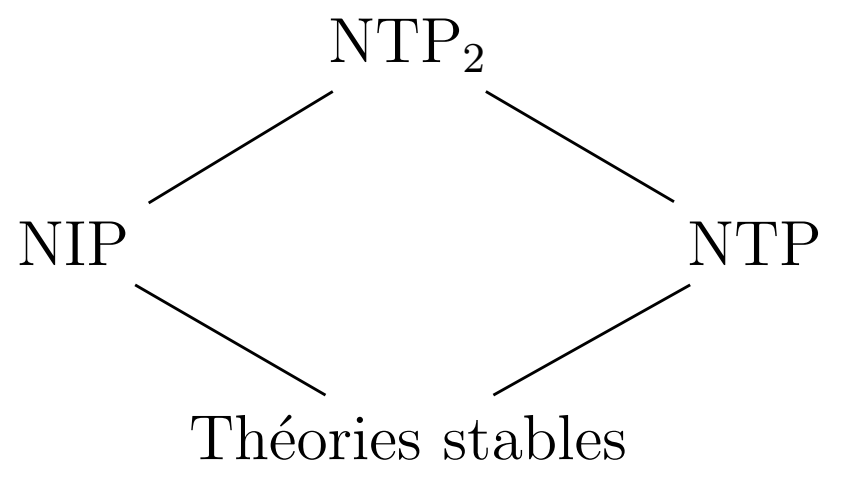

Bien évidemment, plus la classe de complexité que l'on considère est large, moins les notions que l'on définit sur cette classe ne se comportent bien. En particulier, il n'y a plus forcément de notion raisonnable de dimension. La déviation est une notion *d'indépendance* (ou plutôt de dépendance, c'est la non-déviation qui est une notion d'indépendance) qui, même à ce niveau de généralité, hérite d'un certain nombre des bonnes propriétés qu'elle possède dans les théories stables. C'est entre autres pour cette raison qu'elle joue un rôle important dans la théorie des modèles contemporaine.

La définition formelle de la déviation n'est pas très pertinente pour les non-experts, nous nous contenterons d'expliquer comment nous l'interprétons intuitivement. Le lecteur qui est interessé par la définition formelle peut se référer aux deux premières sections du chapitre 1. Dans une structure $M$ dans laquelle nous avons une notion de dimension, un phénomène géométrique notable est d'avoir une inclusion entre deux espaces de dimensions différentes (par exemple, une droite dans un plan). L'espace de dimension plus petite est vu comme étant d'un ordre de grandeur bien plus petit que l'autre, et l'on considérera que ses points ne sont pas assez *génériques* et satisfont une relation de *dépendance*. Etant donné un point $c$, et des sous-structures $A \leqslant B$ de $M$, on dit alors que $c$ est indépendant de $B$ au-dessus de $A$, noté $c \underset{A}{\downarrow} B$, quand la plus petite dimension d'un ensemble $A \cup B$-définissable contenant $c$, la dimension de $c$ sur $A \cup B$, est égale à la dimension de $c$ sur $A$. Cette notion informelle d'indépendance devient formelle dès que l'on fournit une manière de définir la dimension des ensembles définissables.

Si on prend l'exemple d'un corps algébriquement clos (qui est stable), les ensembles définissables sont exactement les ensembles Zariski-constructibles, et leur dimension au sens de la théorie des modèles coincide avec la dimension de Krull, donc la notion d'indépendance que l'on récupère n'est autre que la notion classique de disjonction algébrique de la théorie des corps. De même, dans les espaces vectoriels, la non-déviation correspond à l'indépendance linéaire.

On voit donc que la déviation admet des caractérisations algébriques naturelles dans des théories spécifiques telles que celle des corps algébriquement clos. L'idée qui guide notre travail est de pouvoir établir de telles caractérisations dans d'autres structures algébriques communes dans lesquelles nous avons déjà des intuitions géométriques, telles que les corps valués, et de voir à quel point ces intuitions peuvent être réconciliées avec la définition technique de déviation.

# Survol de la litérature et de nos contributions

## Synthèse

Commençons par un bref résumé avant de donner plus de détails. Notre objectif est d'établir des caractérisations algébriques de la déviation dans les corps valués. Une manière naturelle de procéder serait d'établir un principe de transfert à la Ax-Kochen-Ershov entre la déviation dans un corps valué, et la déviation dans son groupe des valeurs et son corps des résidus. Nous verrons page 14 qu'un tel principe ne peut être établi, à cause d'obstructions venant des extensions immédiates (les extensions de corps valués qui laissent inchangés le groupe des valeurs et le corps des résidus). Le groupe des termes dominants RV défini précédemment est une structure algébrique sous-jacente qui, pour des raisons techniques, encode essentiellement tout ce qui ne relève pas des extensions immédiates. Un plan plus raisonnable pour étudier la déviation serait donc :

1. Comprendre la déviation dans les extensions immédiates.

2. En parallèle, comprendre la déviation dans RV. Par exemple on pourrait :

    a Etablir un principe de transfert à la Ax-Kochen-Ershov entre la déviation dans RV, et la déviation dans le groupe des valeurs et le corps des résidus.

    b Comprendre la déviation dans le groupe des valeurs.

    c Comprendre la déviation dans le corps des résidus.

3. Recoller les points précédents en une caractérisation de la déviation dans le corps valué.

Les résultats de la litérature dont nous avons connaissance sont :

- L'étape 2.c : une caractérisation de la déviation dans les corps des résidus de ACVF, et des corps pseudo-locaux.
- Un principe de tranfert pour la déviation entre le corps valué et son groupe des valeurs et son corps des résidus, dans le cas particulier d'une extension séparée (voir la définition 3). Un tel principe de transfert a d'abord été établi dans ACVF ([HHM05], Théorème 13.7), avant d'être généralisé aux corps Henséliens de caractéristique résiduelle nulle ([EHS23], Corollaire 3.2 et [Vic21], Théorème 4.4).

Les contributions apportées par cette thèse sont :

- L'étape 2.b : une caractérisation de la déviation dans les groupes des valeurs de ACVF, des corps pseudo-locaux, et dans une classe de groupes Abéliens ordonnés un peu plus large (voir la fin de la page 12).
- L'étape 2.a. Le principe de tranfert est établi dans le contexte bien plus général des expansions de suites exactes courtes de structures Abéliennes (voir [ACGZ22], Section 4).
- L'étude d'un aspect important de la déviation : les bases d'extensions (voir la définition 4), dans un contexte très général qui prend en compte les extensions immédiates.

Nous établissons également dans la section 7.3 des liens entre la déviation et les mesures de Keisler pour les corps pseudo-locaux.

Avant d'énoncer formellement nos résultats, nous allons détailler davantage ce qui a été fait dans la litérature, et ce qui est fait dans cette thèse.

## Groupes Abéliens ordonnés

Selon le principe d'Ax-Kochen-Ershov, il est naturel de commencer par essayer de comprendre la déviation dans le groupe des valeurs et le corps des résidus de nos corps valués. En d'autres termes, ont voudrait comprendre la déviation dans les corps, et dans les groupes Abéliens ordonnés. Beaucoup de progrès ont été réalisés par la communauté en ce qui concerne la déviation dans les corps. Comme nous l'avons déjà mentionné, dans les corps algébriquement clos, qui sont les corps des résidus des modèles de ACVF, la déviation correspond à la dépendance algébrique classique. C'est également le cas dans les corps pseudo-finis (on peut le montrer grâce aux travaux de Chatzidakis-van den Dries-Macintyre [CvdDM92]), qui sont les corps des résidus des corps pseudo-locaux. Qu'en est-il de la déviation dans les groupes Abéliens ordonnés ? Etonnamment, il n'existe à ce jour aucun résultat (à notre connaisance) décrivant de manière complète et algébrique la déviation dans un groupe Abélien ordonné non-trivial, à l'exception d'un résultat de Haskell-Hrushovski-Macpherson sur les groupes divisibles dont nous allons discuter un peu plus bas (Théorème 5) qui nécessite des hypothèses assez fortes. Les seuls autres résultats sont des caractérisations abstraites de la déviation dans des cadres généraux (tels que les théories NIP) dans lesquels s'inscrivent les groupes Abéliens ordonnés. C'est d'autant plus étonnant que Dolich a su décrire dans ([Dol04]) la déviation dans les expansions o-minimales du corps des réels. On connaît donc la déviation dans le corps des réels, on la connait dans le groupe des réels (c'est l'indépendance linéaire classique des $\mathbb{Q}$-espaces vectoriels), mais on ne la connaît pas dans le groupe ordonné des réels.

Nous étudions cette question dans le chapitre 3 de ce travail, étude qui nous permet de donner une caractérisation géométrique simple de la déviation dans les groupes Abéliens ordonnés divisibles non-triviaux (DOAG). Une source de difficulté vient du fait que nous devons gérer des groupes Abéliens ordonnés dans lesquelles certains éléments sont infinitésimaux par rapport à d'autres. Il semble que le fait de travailler dans un corps permet à Dolich de contourner ce problème, et de parvenir à comprendre la déviation uniquement en termes de cellules dans une structure o-minimale, mais cela ne fonctionne pas dans notre contexte.

Sachant que le groupe des valeurs d'un corps pseudo-local est toujours une extension élémentaire du groupe ordonné $\mathbb{Z}$, nous étendons nos résultats du chapitre 3 sur DOAG (c'est-à-dire sur $\mathbb{Q}$) à tous les sous-groupes de $\mathbb{R}$. La théorie de ces groupes ordonnés est notée ROAG, et ses modèles sont les groupes Abéliens ordonnés dits *réguliers*. D'un point de vue modèle-théorique, le langage qui permet de facilement comprendre ROAG est le langage de Presburger, qui est l'enrichissement du langage des groupes ordonnés par des prédicats qui testent la divisibilité des éléments par les puissances des nombres premiers. Les groupes réguliers sont exactement ceux qui éliminent les quantificateurs dans ce langage.

## Extensions séparées et bases d'extension

On se demande quelle serait une notion d'indépendance géométrique naturelle dans un corps valué. Pour répondre à cette question, on peut s'inspirer de l'algèbre linéaire. Etant donné un uplet $c$ d'éléments d'un espace vectoriel, il satisfait une relation de dépendance quand on peut obtenir zéro à partir d'une combinaison linéaire non-triviale de $c$. On peut généraliser cette idée dans un espace vectoriel valué. La phénomène de dépendance non-trivial le plus simple qui peut survenir dans un espace vectoriel valué est quand on a deux points $x$, $y$ d'ordre de grandeur infini, qui sont infinitésimalement proches l'un de l'autre, en d'autres termes quand on a $\mathrm{val}(x-y) < \max(\mathrm{val}(x), \mathrm{val}(y))$. De manière plus générale :

**Définition 3.** La famille $c = (c_i)_i$ est *séparée* si pour toute combinaison linéaire $\sum_i \lambda_i c_i$, on a :

$$\mathrm{val}\left(\sum_i \lambda_i c_i\right) = \max_i \{\mathrm{val}(\lambda_i c_i)\}$$

Nous travaillons dans un contexte où la valeur de zéro est plus petite que toutes les autres. Certaines convention préféreront inverser cet ordre, auquel cas il faudra remplacer le max par un min.

Une extension $B \geqslant A$ de corps valués est *séparée* si tout sous-$A$-espace vectoriel de $B$ de dimension finie admet une base séparée par rapport à la structure de $A$-espace vectoriel valué de $B$.

Les extensions séparées de corps valués se comportent particulièrement bien. En réalité, les seuls résultats ([HHM05], Théorème 13.7, [EHS23], Corollaire 3.2, [Vic21], Théorème 4.4) sur la déviation dans les corps valués antérieurs à notre travail nécessitent l'hypothèse que l'on travaille avec des extensions séparées. Nous utilisons notamment les résultats de Ealy-Haskell-Simon ([EHS23]) pour obtenir dans le chapitre 7 de nouvelles descriptions de la déviation dans les corps pseudo-locaux, pour le type d'une extension qui est séparée sur les paramètres. Cependant, contrairement au fait que tout espace vectoriel admet une base, et comme le suggère le fait que l'on aie à définir les extensions séparées de corps valués, il existe des extensions de corps valués qui ne sont pas séparées. Une classe d'extensions de corps valués qui peut être considérée comme orthogonale aux extensions séparées est la classe des extensions *immédiates*, les extensions $C \geqslant A$ telles que $\Gamma(C) = \Gamma(A)$ et $k(C) = k(A)$. Les extensions immédiates de degré de transcendance positif ne sont jamais séparées, et la plupart des corps valués que l'on trouve dans la nature admettent une telle extension, et on ne peut pas les ignorer. L'existence des extensions immédiates rend d'ailleurs impossible le fait de décrire entièrement la déviation dans un corps valué par un principe d'Ax-Kochen-Ershov. En effet, si, au sein d'un modèle de ACVF, on a une extension de sous-corps $C \geqslant A$ qui est immédiate et de degré de transcendance positif, alors $\Gamma(C)$ (resp. $k(C)$) est indépendant de lui-même au-dessus de $\Gamma(A)$ (resp. $k(A)$) au sens de la déviation, mais $C$ n'est pas indépendant de lui-même au-dessus de $A$. Ce contre-exemple un peu dégénéré ne peut pas être ignoré, car il peut être utilisé pour construire des contre-exemples plus subtils. L'étude des extensions immédiates est un sujet très difficile, mais un résultat de [HMRC18] sur les extensions immédiates simples nous permet de bien comprendre la déviation en dimension un (quand l'ensemble de paramètres à gauche de la relation $\downarrow$ est une extension simple de celui à la base). Cela nous permet d'établir au chapitre 6 des résultats sur les bases d'extensions dans les corps valués :

**Définition 4.** Dans une structure du premier ordre $M$, un ensemble de paramètres $A$ est une *base d'extension* si dans toute extension élémentaire de $M$, tout ensemble de paramètres est indépendant de $A$ au-dessus de $A$ au sens de la déviation.

Le fait de travailler au-dessus d'une base d'extension est bien-sûr une condition absolument nécessaire pour espérer obtenir des résultats intéressants sur la déviation, les ensembles qui ne sont pas des bases d'extensions sont beaucoup trop dégénérés. Il se trouve que cette condition est également suffisante dans les théories NIP et $\mathrm{NTP}_2$ pour généraliser un grand nombre de résultats de la théorie de la stabilité. Les bases d'extensions jouent donc un role fondamental dans les corps valués qui nous intéressent. Il est à noter que tout ensemble est une base d'extension dans les théories NTP. On trouve ainsi dans le chapitre 6 des conditions suffisantes assez faibles pour qu'un ensemble de paramètres au sein d'un corps valué Hensélien de caractéristique résiduelle nulle soit une base d'extension. On montre en particulier que dans les corps pseudo-locaux, tout ensemble est une base d'extension, et donc la déviation satisfait un certain nombre de bonnes propriétés théoriques héritées des théories stables.

### Suites exactes courtes

Pour finir, la dernière contribution sur laquelle nous avons travaillé pendant la préparation de cette thèse concerne la déviation dans les expansions de suites exactes courtes pures de structures Abéliennes telles qu'étudiées dans ([ACGZ22], Section 4). Même si l'on est particulièrement intéressé par les suites exactes courtes pures de groupes Abéliens induites par les corps valués, le contexte dans lequel on se place est bien plus général. L'intérêt de ce cadre abstrait est qu'il concerne également les suites exactes courtes induites par des expansions de corps valués, notamment ceux avec des opérateurs (une dérivation, un automorphisme, ou même les deux) qui ont une place importante dans de nombreux travaux contemporains. Les résultats attendus devraient prendre la forme d'un principe d'Ax-Kochen-Ershov : on cherche à décrire la déviation dans la suite exacte courte en fonction de la déviation dans ses membres de droite et de gauche. Dans la description que l'on obtient, il reste des ensembles définissables (essentiellement des torseurs non-relevés) dont le comportement ne peut pas être décrit par un principe d'Ax-Kochen-Ershov. Modulo ces torseurs, on obtient donc un principe d'Ax-Kochen-Erhov partiel. Quand l'ensemble de base est un modèle, une sous-structure élémentaire du modèle ambiant, les obstructions à nos résultats attendus disparaissent, et nous obtenons un principe d'Ax-Kochen-Ershov aussi fort que l'on pouvait espérer.

## Présentation des résultats

### Extensions séparées (chapitre 7)

Le premier des résultats de la littérature sur la déviation dans les corps valués est dû à Haskell-Hrushovski-Macpherson. On énonce ici une variante légèrement différente qui peut se déduire facilement de l'énoncé original avec des arguments classiques :

**Théorème 5** ([HHM05], Théorème 13.7)**.** *Dans un modèle $M$ de* ACVF *dont la caractéristique résiduelle est nulle, soient $C \geqslant A \leqslant B$ des sous-corps tels que l'extension $C \geqslant A$ soit séparée. Supposons de plus que pour tout $\gamma$ dans* *la clôture divisible* $\mathrm{div}(\Gamma(C))$ *du groupe des valeurs de* $C$*, qui ne soit pas dans* $\mathrm{div}(\Gamma(A))$*, le groupe Abélien ordonné* $\mathrm{div}(\Gamma(A)) + \mathbb{Q}\gamma$ *a strictement plus de sous-groupes convexes que* $\mathrm{div}(\Gamma(A))$*. Alors les conditions suivantes sont équivalentes :*

- $C$ *est indépendant de* $B$ *au-dessus de* $A$ *au sens de la déviation.*
- $k(C)$ *est algébriquement disjoint de* $k(B)$ *au-dessus de* $k(A)$ *d'une part, et d'autre part tout intervalle fermé borné de* $\mathrm{div}(\Gamma(B))$ *qui intersecte* $\mathrm{div}(\Gamma(C))$ *intersecte déjà* $\mathrm{div}(\Gamma(A))$.

La seconde condition de ce théorème est non-seulement un principe d'Ax-Kochen-Ershov, mais elle ne fait intervenir que des phénomènes géométriques très concrets. C'est le genre de caractérisation de la déviation que nous cherchons à obtenir dans notre travail. Il est à noter que la condition de ce théorème sur les intervalles caractérise la déviation dans DOAG quand la base $A$ est un groupe dont la clôture divisible est *Archimède-complète*, c'est-à-dire quand toute extension propre de $\mathrm{div}(A)$ a strictement plus de sous-groupes convexes que $\mathrm{div}(A)$. Nous montrons au chapitre 3 que la même condition en termes d'intervalles caractérise la déviation au-dessus d'une base quelconque.

Le théorème de Haskell-Hrushovski-Macpherson a été généralisé aux corps Henséliens de caractéristique pure nulle par Ealy-Haskell-Simon ([EHS23]) et Vicaria ([Vic21]). Cependant, ces résultats concernent une relation de dépendance plus forte (qui induit donc une relation d'indépendance plus faible) que la déviation, qu'on appelle la division. On sait que dans ACVF, la déviation et la division coincident, ce qui est considéré comme une bonne propriété modèle-théorique, mais on ignore si c'est le cas dans tout corps valué. Enonçons le résultat de Ealy-Haskell-Simon :

**Théorème 6** ([EHS23], Corollaire 3.4)**.** *Soit* $M$ *un corps valué Hensélien de caractéristique résiduelle nulle. Supposons que* $[K(M)^* : (K(M)^*)^N]$ *soit fini pour tout* $N > 0$*. Soient* $C \geqslant A \leqslant B$ *des sous-corps de* $M$ *tels que* $B$ *et* $C$ *sont définissablement clos, et tels que l'extension* $C \geqslant A$ *est séparée. Supposons que pour tout* $N > 0$*, tout coset de* $\mathrm{RV}(M) / \mathrm{RV}(M)^N$ *se relève dans* $A$*, que* $k(C)$ *est une extension régulière de* $k(A)$*, et que* $\Gamma(C)$ *est une extension pure de* $\Gamma(A)$*. Alors les conditions suivantes sont équivalentes :*

- $C$ *est indépendant de* $B$ *au-dessus de* $A$ *au sens de la division.*

- *L'ensemble de paramètres $k(C)$ (resp. $\Gamma(C)$) est indépendant de $k(B)$ (resp. $\Gamma(B)$) au-dessus de $A$ au sens de la division.*

En réalité, l'énoncé de [EHS23] est plus faible, bien que la preuve correspond exactement à ce que nous avons écrit ci-dessus : le même énoncé est donné avec comme hypothèse supplémentaire que la déviation coïncide avec la division au-dessus de $A$. Vicaria a un résultat similaire dans un langage un peu plus riche, avec des hypothèses plus faibles (pas de finitude des indices).

Le fait que la seconde condition ne fasse pas intervenir de géométrie n'est pas un problème : le corps des résidus et le groupe des valeurs de $M$ sont quasiment arbitraires (on sait que l'indice de $(k(M)^*)^N$ dans $k(M)^*$ est fini pour tout $N > 0$, en dehors de cela $k(M)$ peut être très sauvage, et $\Gamma(M)$ peut vraiment être quelconque). On ne peut donc pas espérer décrire la déviation ou la division de manière purement géométrique, on ne sait même pas le faire dans tous les corps et les groupes Abéliens ordonnés. Le mieux que l'on puisse espérer est un principe d'Ax-Kochen-Ershov, et le théorème ci-dessus nous en fournit effectivement un. Ces résultats de Ealy-Haskell-Simon et Vicaria sont avant tout présentés comme corollaires de résultats portant sur la domination stable, qui est le sujet principal de leurs articles. Nous utilisons ces résultats de domination stable au chapitre 7 pour présenter une version plus simple et plus forte du théorème ci-dessus. En particulier, dans l'énoncé ci-dessus, la division est étudiée au-dessus de $A$ dans la structure $M$, alors que l'on voudrait idéalement se réduire à étudier la division au-dessus de $k(A)$ et $\Gamma(A)$ dans les structures $k(M)$ et $\Gamma(M)$. Nous donnons ci-dessous notre énoncé dans le cas d'un corps pseudo-local, bien que nous l'ayons prouvé dans un cadre plus général (voir les hypothèses 7.0.2) :

**Théorème 7** (Corollaire 7.1.4 et Proposition 7.2.3)**.** *Soit $M$ un corps valué pseudo-local de caractéristique résiduelle nulle. Soient $C \geqslant A \leqslant B$ des sous-corps de $M$ tels que l'extension $C \geqslant A$ est séparée. Supposons que pour tout $N > 0$, tout coset de $\mathrm{RV}(M) / \mathrm{RV}(M)^N$ se relève dans $A$. Supposons qu'il en est de même pour les cosets de $k(M)^* / (k(M)^*)^N$, et que $A$ contienne un uniformisateur (un élément de valeur $1$). Alors les conditions suivantes sont équivalentes :*

- *$C$ est indépendant de $B$ au-dessus de $A$ pour la division.*
- *Dans le réduit $k(M)$ (resp. $\Gamma(M)$), $k(C)$ (resp. $\Gamma(C)$) est indépendant de $k(B)$ (resp. $\Gamma(B)$) au-dessus de $k(A)$ (resp. $\Gamma(A)$) pour la division.*

*Si, en plus de cela, $k(C)$ est une extension de corps régulière de $k(A)$, et $\Gamma(C)$ est une extension de groupes Abéliens pure de $\Gamma(A)$, alors la même équivalence est vraie pour la déviation.*

En résumé, nous supposons en plus que $A$ contient des solutions de certaines équations, mais nous n'avons plus besoin de supposer que $B$ et $C$ soient définissablement clos, nous n'avons plus besoin de l'hypothèse de régularité de $k(C)$ pour la division, et nous prouvons le résultat analogue pour la déviation. De plus, on se réduit bel et bien à la déviation et la division dans les réduits $\Gamma(M)$ et $k(M)$, au lieu de $M$.

Nos résultats sur la déviation dans $\mathbb{Z}$ du chapitre 4 permettent de montrer que les conditions ci-dessus sont équivalentes aux conditions géométriques de Haskell-Hrushovski-Macpherson dans le Théorème 5. Nous préférons tout-de-même garder cet énoncé non-géométrique, car il reste vrai dans une classe de corps valués qui est un peu plus large que ceux qui sont pseudo-locaux (voir le début du chapitre 7). Les points de $A$ que l'on doit nommer pour relever tous ces coset ne forment qu'un ensemble au plus dénombrable. Le fait de devoir relever ainsi une collection bornée d'imaginaires est commun en théorie des modèles, et cela ne pose pas de problème notable, en dehors du fait que l'on exclut quelques cas trop dégénérés. La seule hypothèse un peu forte dont nous avons besoin pour la généralisation du Théorème 7 (plus précisément, cette hypothèse n'est nécessaire que pour la seconde partie de l'énoncé, qui concerne la déviation), est le fait que la théorie du corps des résidus soit algébriquemant bornée, ce qui est vrai dans les corps pseudo-locaux.

Etant donné un corps pseudo-local, la non-déviation dans le corps des résidus et le groupe des valeurs a toujours pour témoin une mesure de Keisler (une mesure de probabilité finiment additive sur les ensembles définissables) globale invariante qui étend le type en question (dans le groupes des valeurs, cette mesure est un type). Nous étudions donc dans la section 7.3 à quel point la déviation peut être mise en relation avec les mesures de Keisler. Pour des raison techniques, on se réduit à une sous-classe des extensions séparées, celles qui sont *Abhyankar*, essentiellement celles qui n'admettent aucune extension intermédiaire immédiate (voir les hypothèses 7.3.1). Cela nous amène à manipuler des notions d'indépendance modèle-théoriques faisant intervenir les mesures de Keisler. Bien que cette contribution reste modeste, nous développons dans la section 1.3 une théorie de l'indépendance par les mesures de Keisler, dont nous nous servons pour nos résultats de la section 7.3 sur les corps pseudo-locaux :

**Théorème 8** (Corollaire 7.3.14)**.** *Avec les mêmes hypothèses que le Théorème 7, supposons que l'extension $C \geqslant A$ soit Abhyankar. Alors les conditions suivantes sont équivalentes :*

- *$C$ est indépendant de $B$ au-dessus de $A$ pour la déviation.*
- *Il existe une mesure de Keisler globale $\mathrm{Aut}(M/A)$-invariante dont le support contient tous les ensembles $B$-définissables qui contiennent (une énumération donnée de) $C$.*

## Groupes Abéliens ordonnés divisibles (chapitre 3)

On peut commencer par décrire la déviation dans les réduits naturels de DOAG :

**Fait 9.** *Au sein d'un $\mathbb{Q}$-espace vectoriel, soient $A$, $B$, $C$ des ensembles de paramètres. Soient $A', B', C'$ les sous-$\mathbb{Q}$-espaces vectoriels respectivement engendrés par $A, A \cup B, A \cup C$. Alors $C$ est indépendant de $B$ au-dessus de $A$ au sens de la déviation si et seulement si $C' \cap B' = A'$.*

**Fait 10.** *Au sein d'un ordre total infini dense sans minimum ni maximum (*DLO*), soient $A$, $B$, $C$ des ensembles de paramètres. Alors $C$ est indépendant de $B$ au-dessus de $A$ au sens de la déviation si et seulement si tout intervalle fermé borné de $A \cup B$ qui intersecte $A \cup C$ intersecte déjà $A$.*

La caractérisation de la déviation dans DLO ci-dessus est formellement démontrée dans la sous-section 1.2.2. Il est important que le lecteur interessé par nos résultats sur DOAG lise cette sous-section, car les méthodes utilisées sont des versions très simples de ce que l'on fait au chapitre 3.

On devine que la déviation dans DOAG doit être un amalgame des deux notions d'indépendance ci-dessus. On montre dans notre Théorème 3.1.5 qu'en effet, la déviation en est l'amalgame le plus naturel :

**Théorème 11.** *Au sein d'un modèle de* DOAG*, soient $A$, $B$, $C$ des ensembles de paramètres. Soient $A', B', C'$ les clôtures divisibles des sous groupes respectivement engendrés par $A, A \cup B, A \cup C$. Alors $C$ est indépendant de $B$ au-dessus de $A$ au sens de la déviation si et seulement si tout intervalle fermé borné de $B'$ qui intersecte $C'$ intersecte déjà $A'$.*

En particulier, pour vérifier si $C$ est indépendant, il suffit de vérifier que chaque point du $\mathbb{Q}$-espace vectoriel engendré par $A \cup C$ est indépendant. Ce phénomène fort de réduction à la déviation en dimension un est bien connu dans les structures stables dites *monobasées*, et le fait que ce phénomène apparaisse dans DOAG soulève une question intéressante, la question de savoir si ce phénomène est causé par des raisons plus abstraites, si les techniques des théories monobasées et de la modularité locale peuvent être généralisées à des théories instables (o-minimales ?) telles que DOAG.

Le lecteur sera peut-être surpris par la complexité de la preuve de ce théorème, qui prend quand-même place dans une théorie bien comprise, et dont l'énoncé est si simple. Nous jugeons bon d'expliquer dans cette introduction quel cheminement nous a mené vers cette preuve.

La version en dimension un (quand l'extension $C \geqslant A$ est simple) de notre théorème est bien connue dans n'importe quelle théorie o-minimale. Sachant qu'en dimension supérieure, les seules propriétés élémentaires que l'on peut exprimer dans DOAG s'écrivent comme l'appartenance d'une combinaison $\mathbb{Q}$-linéaire des variables à des intervalles, une question naturelle que l'on peut se poser est : si l'on connaît les coupures respectives de quelques points, que peut-on dire de celles de leurs combinaisons linéaires ? Pour répondre à cette question, on voudrait commencer par classifier les coupures (on définit ici la *coupure* d'un point $c$ sur un sous-espace vectoriel $A$ comme l'intersection des intervalles de $A$ qui contiennent $c$) :

**Définition 12.** Fixons $M$ un modèle suffisamment large de DOAG, et $A$ un sous-$\mathbb{Q}$-espace vectoriel. Pour tout $c$ dans $M$, on définit la *classe* $\text{cl}(c/A)$ de $c$ sur $A$ de manière suivante :

- Supposons que la coupure de $c$ sur $A$ soit le coset d'un sous-groupe convexe $A$-type-définissable $G$. Alors $\text{cl}(c/A) = (1, G)$.
- Supposons qu'il existe $G$ un sous-groupe convexe $A$-type-définissable non-trivial tel que, si $H$ est le plus grand sous-groupe convexe $A$-$\vee$-définissable de $G$, on ait $c \in (A+G) \setminus (A+H)$. Alors la coupure de $c$ sur $A$ est l'une des deux composantes connexes d'un ensemble $A$-type-définissable de la forme $X \setminus Y$, avec $X$ un coset de $G$, et $Y$ un coset de $H$. Dans ce cas, on définit $\text{cl}(c/A) = (2, G)$.

Il se trouve que tout $c \in M$ vérifie exactement un des deux cas ci-dessus, et que $\text{cl}(c/A)$ est bien-définie. Il n'est pas difficile de vérifier que l'application $M \longrightarrow \mathrm{cl}(M/A)$ passe au quotient par $A$. L'observation étonnante que nous avons faite est que l'application induite ${}^{M}\!\diagup_{A} \longrightarrow \mathrm{cl}(M/A)$ est une *valuation de groupes* (pour un certain ordre sur les classes). Elle correspond à la valuation $\mathrm{val}_A^2$ de la définition 3.2.35. En raffinant notre analyse, nous avons découvert que diverses propriétés modèle-théoriques interagissaient de manière naturelle avec trois valuations de groupes non-conventionnelles sur des quotients. Par exemple, on peut remarquer que deux points dont les classes sur $A$ sont distinctes sont nécessairement *faiblement orthogonaux* au-dessus de $A$. Le problème initial de pouvoir deviner les coupures des combinaisons linéaires de points juste à partir des coupures de ces mêmes points a donc une solution quand on travaille sur des familles séparées (pour la valuation triviale sur $\mathbb{Q}$) par rapport à ces valuations. On montre que la condition d'indépendance géométrique de notre Théorème 11 implique que les uplets que l'on manipule sont toujours $A$-interdéfinissables avec une famille séparée, et c'est essentiellement ce qui nous permet de démontrer ce théorème.

De plus, la manière dont ces valuations non-conventionnelles interagissent avec la théorie des modèles nous permet d'établir une classification très fine et explicite de l'espace des extensions globales invariantes d'un type non-déviant dans DOAG. Cette classification est décrite dans la remarque 3.3.31.

Plutôt que de manipuler des valuations sur des quotients, nous préférerons manipuler des valuations que l'on autorise à s'annuler sur un sous-groupe non-trivial $A$, que l'on appelle des *$A$-valuations*. Les terminologies que l'on adopte sont définies à la section 2.1. Nous notons nos valuations de manière additive, et nous ordonnons les valeurs de manière à ce que la valeur de zéro soit plus petite que les autres, et non l'inverse. Nous préférons adopter ces conventions, étant donné que l'on est amené à manipuler de nombreuses valuations de groupes non-conventionnelles de nature purement géométrique.

## Groupes Abéliens ordonnés réguliers (chapitre 4)

Nous démontrons au chapitre 4 une description géométrique de la déviation pour les groupes Abéliens ordonnés qui éliminent les quantificateurs dans le langage de Presburger (ROAG). Il est à noter qu'un groupe ordonné est toujours sans-torsion, et que les groupes Abéliens sans-torsion éliminent les quantificateurs dans le langage de Presburger privé de l'ordre. Etant donné que ces groupes sont stables, la déviation est facile à calculer dans leur théorie. Il s'agit essentiellement d'un principe local-global :

**Fait 13.** *Au sein d'un groupe Abélien sans-torsion $M$, soient $A, B, C$ des ensembles de paramètres. Soient $A', B', C'$ les clôtures divisibles relatives à $M$ des sous-groupes engendrés respectivement par $A, A\cup B, A\cup C$. Alors $C$ est indépendant de $B$ au-dessus de $A$ si et seulement si les conditions suivantes sont vérifiées :*

- *(Indépendance modulo zéro) $C' \cap B' = A'$.*
- *(Indépendance modulo p) Pour tout nombre premier $p$, si $[M : pM]$ est infini, alors nous avons pour tout $N > 0$ :*

$$(C' + p^N M) \cap (B' + p^N M) = A' + p^N M$$

Nous ignorons si ce fait est prouvé explicitement quelque part, mais le fait est que l'on peut facilement en obtenir une preuve en ignorant tout ce qui relève de l'ordre dans les preuves du chapitre 4. Sachant que les théories en question sont de toute façon stables, nous allons admettre ce fait.

Un type complet dans ROAG s'écrit comme la réunion de deux types partiels particuliers (qui sont sans-quantificateurs par élimination). D'une part, on a un type partiel qui ne fait intervenir que l'ordre. Ce type partiel indique essentiellement le DOAG-type des variables dans la clôture divisible du modèle. D'autre part, on a un type partiel qui ne fait intervenir que les divisibilités par les puissances des nombres premiers. Ce type partiel indique essentiellement le type des variables dans le réduit stable. En choisissant un tuple $A$-interdéfinissable adéquat, on arrive à se placer dans une situation où l'on coupe le problème en deux problèmes indépendants, et où le type que l'on examine est non-déviant si et seulement si c'est le cas de chacun des deux types partiels. Il reste ensuite à caractériser géométriquement quand-est-ce que chacun des types partiels est non-déviant. Pour le type qui fait intervenir l'ordre, on utilise nos résultats sur DOAG, tandis que pour le type qui fait intervenir les divisibilités, on utilise le Fait 13 (on reprouve proprement tout ce dont on a besoin). On obtient ainsi le résultat suivant (Théorème 4.4.1) :

**Théorème 14.** *Au sein d'un groupe Abélien ordonné régulier $M$, soient $A, B, C$ des ensembles de paramètres. Soient $A', B', C'$ les clôtures divisibles relatives à $M$ des sous-groupes engendrés respectivement par $A, A \cup B, A \cup C$. Alors $C$ est indépendant de $B$ au-dessus de $A$ si et seulement si les conditions suivantes sont vérifiées :*

- *Tout intervalle fermé borné de $B'$ qui intersecte $C$ intersecte aussi la clôture divisible de $A$.*
- *Pour tout nombre premier $p$, si $[M : pM]$ est infini, alors nous avons pour tout $N > 0$ :*

$$(C' + p^N M) \cap (B' + p^N M) = A' + p^N M$$

En résumé, on intègre l'ordre à notre indépendance en ajoutant la description géométrique de l'indépendance de DLO à la condition zéro de notre principe local global. Notons que nous avons de nouveau une réduction à la déviation en dimension un comme dans les théories monobasées.

Nous montrons également que ces conditions sont équivalentes au fait que le type de $C$ sur $A \cup B$ admette une extension globale $\mathrm{acl}^{eq}(A)$-invariante, et que la déviation coincide avec la division dans ROAG. Nous caractérisons aussi quand ce type admet une extension globale $A$-invariante : il faut ajouter la condition que si $[M : pM]$ est fini, alors $C \subseteq A' + \bigcap_{N>0} p^N M$.

Pour finir, nous généralisons à ROAG dans la remarque 4.4.4 la classification des extensions globales non-déviantes d'un type que nous avions établie dans DOAG. Le fait d'ajouter des informations sur les divisibilités revient à faire un produit direct topologique avec un espace qui est naturellement homéomorphe à un produit Cartésien de fermés des $\mathbb{Z}_p$.

Nous terminons ce chapitre avec la section 4.5 qui est dédiée à divers exemples de comportements de la déviation dans plusieurs groupes Abéliens ordonnés.

## Suites exactes courtes et expansions (chapitre 5)

Décrivons plus en détail le contexte général dans lequel le chapitre 5 prend place. Au lieu de travailler avec des suites exactes courtes de groupes Abéliens, on travaille avec des suites de structures Abéliennes. Une *structure Abélienne* est une structure du premier ordre avec potentiellement plusieurs sortes, telles que chaque sorte est une expansion d'un groupe Abélien, chaque prédicat du langage est un sous-groupe (pour la structure de groupe produit sur son domaine), et chaque fonction du langage est un morphisme de groupes. Par exemple, un groupe Abélien, un module à gauche sur un anneau fixé, ou même un complexe de tels modules, est une structure Abélienne. On peut définir ce qu'est une suite exacte courte de structures Abéliennes ayant le même langage, d'une manière qui étend naturellement les notions classiques pour les groupes Abéliens. (Notons qu'un complexe de structures Abéliennes est toujours une structure Abélienne.) Il en est de même pour la notion de *pureté* d'une telle suite. En algèbre commutative classique, la pureté revient à dire que certains ensembles définissables dans le membre central de la suite admettent un point dans le membre de gauche. C'est la même chose pour les structures Abéliennes, il faut prendre une certaine classe d'ensembles définissables, ceux qui correspondent à ce que l'on appelle les *formules primitives positives*, plus communément appelées les *p.p-formules.* Nous étudions dans le chapitre 5 la déviation dans les structures du premier ordre qui sont des suites exactes courtes pures de structures Abéliennes :

$$0 \longrightarrow A \longrightarrow B \longrightarrow C \longrightarrow 0$$

auxquelles on a ajouté une expansion arbitraire (pas forcément interprétable) du membre de gauche $A^*$ et de celui de droite $C^*$. Par exemple, dans un corps valué $M$, la structure en question est la suite de groupes Abéliens:

$$0 \longrightarrow k(M)^* \longrightarrow \mathrm{RV}(M) \longrightarrow \Gamma(M) \longrightarrow 0$$

à laquelle on ajoute l'expansion du groupe $(k(M))^*$ au corps $k(M)$, et l'expansion du groupe $\Gamma(M)$ à sa structure de groupe ordonné.

Les résultats attendus sont des principes d'Ax-Kochen-Ershov où l'on réduirait la question de la déviation dans l'expansion $M$ :

$$0 \longrightarrow A^*(M) \longrightarrow B(M) \longrightarrow C^*(M) \longrightarrow 0$$

où $M$ est l'expansion d'une suite exacte courte pure, à la déviation dans les réduits respectifs $A^*(M)$ et $C^*(M)$. Soient $\Gamma \subseteq \Delta$ des ensembles de paramètres de $M^{eq}$, et $u$ un tuple de $M$. On cherche à caractériser quand $u$ est indépendant de $\Delta$ au-dessus de $\Gamma$ au sens de la déviation. On peut commencer par s'inspirer du résultat d'élimination des quantificateurs de ([ACGZ22], Corollaire 4.20). Ce corollaire peut être facilement utilisé pour montrer que le type de $u$ sur $\Delta$ s'écrit naturellement comme la réunion de deux types partiels $\mathrm{tp}_A(u/\Delta)$ et $\mathrm{tp}_C(u/\Delta)$, qui décrivent essentiellement les relations algébriques que $u$ a avec $A^*(\Delta)$ et $C^*(\Delta)$. Nous définissons formellement ces types partiels dans la définition 5.1.2. Le principe d'Ax-Kochen-Ershov naïf auquel on s'attend pour la déviation en fonction de ces deux types partiels échoue, comme nous le montrons dans la section 5.5. Il faut en fait définir un type partiel $p_0$ sur $\Gamma$ (qui dépend du type de $u$ sur $\Gamma$) satisfait par $u$ tel que, modulo de légères hypothèses ((H1) et (H2)), on ait le principe d'Ax-Kochen-Ershov suivant :

**Théorème 15** (Corollaire 5.4.6)**.** *Supposons (H1) et (H2). Soit $p_0$ le type partiel de la définition 5.2.1, et $\mathrm{tp}_A(u/\Delta), \mathrm{tp}_C(u/\Delta)$ ceux de la définition 5.1.2. Alors les conditions suivantes sont équivalentes :*

- *$u$ est indépendant de $\Delta$ au-dessus de $\Gamma$ au sens de la déviation.*
- *Aucun des deux types partiels $\mathrm{tp}_C(u/\Delta)$, $p_0 \cup \mathrm{tp}_A(u/\Delta)$ ne dévie sur $\Gamma$.*

Il est à noter qu'on étudie la déviation de types partiels toujours dans la suite exacte, et non dans les réduits.

Nous avons aussi une classification à la Ax-Kochen-Ershov de l'espace topologique des extensions non-déviantes d'un type en fonction de ces types partiels, modulo $p_0$ (voir Corollaire 5.4.5, et Corollaire 5.3.5 pour la division).

Pour finir, nous renforçons dans la section 5.6 nos hypothèses en supposant que $\Gamma$ est un modèle (on fait en réalité des suppositions plus faibles, voir la définition 5.6.1). La caractérisation de la déviation que nous obtenons dans le Corollaire 5.6.7 est alors très élégante :

**Théorème 16.** *Supposons (H2), et supposons que $\Gamma$ est un modèle. Soit $U$ la sous-structure de $M$ engendrée par $u$ et $\Gamma$. Alors les conditions suivantes sont équivalentes :*

- *$u$ est indépendant de $\Delta$ au-dessus de $\Gamma$ au sens de la déviation.*
- *Dans le réduit $A^*(M)$ (resp. $C^*(M)$), $A^*(U)$ (resp. $C^*(U)$) est indépendant de $A^*(\Delta)$ (resp. $C^*(\Delta)$) au-dessus de $A^*(\Gamma)$ (resp. au-dessus de $C^*(\Gamma)$) au sens de la déviation.*

Il faut faire attention néanmoins : le langage dans lequel on travaille ajoute quelques imaginaires, et la notion de sous-structure engendrée les prend en compte. Nous avons une caractérisation analogue de la division dans le Corollaire 5.6.5.

## Bases d'extension (chapitre 6)

Pour finir cette présentation exhaustive de nos résultats, parlons de notre travail du chapitre 6. Nous cherchons à montrer que des ensembles de paramètres dans les corps valués Henséliens de caractéristique résiduelle nulle sont des bases d'extension. A l'instar des suites exactes courtes du chapitre 5, la complexité modèle-théorique des structures que l'on manipule peut être arbitrairement mauvaise. En effet, il existe des corps (et donc des corps valués Henséliens) qui peuvent encoder des configurations combinatoires arbitrairement complexes. Par exemple, $\mathbb{Z}$ est définissable dans le corps $\mathbb{Q}$, et la théorie de la récursion nous dit que l'on peut essentiellement encoder la logique dans l'anneau $\mathbb{Z}$. Il est donc déraisonnable d'espérer décrire la déviation dans ces corps-là, mais un principe d'Ax-Kochen-Ershov est toujours un but réaliste.

**Définition 17.** On fixe $M$ un corps valué Hensélien dont la caractéristique résiduelle est nulle. Commençons par supposer que dans les réduits $k(M)$ et $\Gamma(M)$, tout ensemble est une base d'extension. Supposons également que les indices $[k(M)^* : (k(M)^*)^N]$ soient finis pour tout $N > 0$. Notons $K$ la sorte qui définit le domaine de la valuation.

On énonce la série de résultats que nous démontrons au chapitre 6, les Théorèmes 6.4.8, 6.4.9, 6.4.10, 6.4.11 :

**Théorème 18.** *Si la déviation coincide avec la dépendance algébrique dans $k(M)$, alors tout sous-ensemble de $K(M)$ est une base d'extension.*

*De plus, si $\Gamma(M)$ est élémentairement équivalent à $\mathbb{Z}$, alors tout sous-ensemble de $M^{eq}$ est une base d'extension.*

**Corollaire 19.** *Dans les corps pseudo-locaux, tout ensemble est une base d'extension, et la déviation coincide avec la division. Cela reste vrai si on ajoute des imaginaires, c'est à dire que tout ensemble est une base d'extension dans une structure bi-interprétable avec un corps pseudo-local.*

Ces résultats sont encourageants, car ils montrent que l'on peut espérer de bons comportements de la part de la déviation dans les corps pseudo-locaux. On voudrait raffiner un peu plus ces résultats, car les hypothèses du théorème ci-dessus sont un peu trop fortes. Pour se faire, un point technique de notre travail consiste à chercher des conditions suffisantes pour que certains imaginaires (les boules, les valeurs...) soient relevés dans nos ensembles de paramètres. Tout comme dans le chapitre 7, nous nommons donc une collection au plus dénombrable de constantes :

**Définition 20.** Soit $A$ un sous-corps de $M$. On suppose que pour tout $N > 0$, tout coset de ${k(M)^*}/{(k(M)^*)^N}$ se relève dans $A$. On suppose également que toute valeur de $\Gamma(\mathrm{dcl}^{eq}(\varnothing))$ se relève dans $A$.

**Théorème 21.** *Avec les hypothèses supplémentaires ci-dessus, si* $\Gamma(M)$ *est dense ou* $k(M)$ *est algébriquement borné, alors* $A$ *est une base d'extension.*

Ces conditions sont pour le coup très faibles. Avant de commencer les préliminaires, on remarque que nos résultats sur les bases d'extensions unifient les résultats de la litérature selon lesquels tout ensemble est une base d'extension dans ACVF (à cause de la C-minimalité) et dans les corps valués réels clos (à cause de l'o-minimalité faible). Nos résultats sur les corps pseudo-locaux (ils s'appliquent facilement à d'autres corps valués, comme par exemple les séries de Laurent sur $\mathbb{C}$), eux, sont bien nouveaux.

# General notations

Given a first-order structure $M$, and an infinite cardinal $\kappa$, $M$ is *$\kappa$-saturated* if every type with strictly less than $\kappa$ parameters has a realization in $M$. Likewise, $M$ is *strongly $\kappa$-homogeneous* if every partial elementary isomorphism of $M$ having domain strictly smaller than $\kappa$ extends to an automorphism of $M$. If $\mathcal{M}$ is another first-order structure on the same language as $M$, we write $M \equiv \mathcal{M}$ when $M$ and $\mathcal{M}$ are elementarily equivalent, and we write $M \preceq \mathcal{M}$ when $M$ is an elementary substructure of $\mathcal{M}$, in which case we say that $M$ is a *model* in the theory of $\mathcal{M}$. We write $M \simeq \mathcal{M}$ when $M$ and $\mathcal{M}$ are isomorphic. We also use the symbol $\simeq$ for homeomorphism of topological spaces.

Tuples are written with lowercase letters, without an overline. We sometimes keep the size of a tuple unspecified. We use concatenation to refer to the union of parameter sets and/or tuples in $M$. For instance, the parameter set $ABcd$ is the union of $A$, $B$, the set elements of $c$ and that of $d$. We may abuse notation and concatenate tuples of parameters or variables to refer to the tuple obtained by their actual concatenation, the contexts in which we use such concatenations should leave no ambiguity.

Definable sets are usually written $X, Y \dots$ We often identify formulas, definable sets, and clopen sets in the Stone space of types. Likewise, we identify type-definable sets with the corresponding closed sets of types. Given a definable set $X$ in $n$ variables, and a parameter set $A$, we write $X(A)$ for the set of tuples of $A^n$ which realize $X$. In a multi-sorted structure, we view (tuples of) sorts as particular instances of definable sets. In particular, if $s$ is a sort, and $A$ is a parameter set, then $s(A)$ refers to the elements of $A$ which are in the sort $s$. Given a tuple of sorted variables $x$, and a parameter set $A$, we write $S^x(A)$ for the Stone space of types in the variables $x$ with parameters in $A$. If $|x| = n$, we may abuse notations and write $S^n(A)$, or even $S(A)$ if the sorts of the variables and/or their number are unspecified.

We denote by $M^{eq}$ the expansion of $M$ by *imaginaries*, i.e. every quotient of $M^n$ by every definable equivalence relation, for every $n < \omega$. The sorts of $M^{eq}$ which are already sorts in $M$ are called *home sorts* or *real sorts*, and their elements are called *reals*. When $M$ only has a single home sort $s$, we may abuse notation and write $M$ to refer to the set $s(M)$.

Given a parameter set $A$, we write $\mathrm{dcl}(A)$ and $\mathrm{acl}(A)$ for the definable closure and the algebraic closure of $A$ in $M$. We write $\mathrm{dcl}^{eq}(A)$, $\mathrm{acl}^{eq}(A)$ for those closures in $M^{eq}$.

In a field $M$, given a subset $A \subseteq M$, we write $A^{\mathbf{alg}}$ for the field-theoretic algebraic closure of the field generated by $A$. If $M$ is a vector space over some field $k$, we write $\mathrm{Vec}_k(A)$ for the $k$-vector subspace of $M$ generated by $A$. We abuse notation and write $\mathrm{Vec}(A)$ when there is no ambiguity as to which field $k$ is involved. If $G$ is a torsion-free Abelian group, we write $\mathrm{div}(G)$ for the divisible closure of $G$, which can naturally be identified with a $\mathbb{Q}$-vector space.

In an Abelian group $G$ (written additively), given an element $c$ and a subgroup $H$, we write $c + H$ or $c$ mod $H$ for the coset of $H$ which contains $c$. Given another subgroup $L$ of $G$, we write $H + L$ for the subgroup generated by $H$ and $L$. Given $n < \omega$, we write $nG$ for the group of elements of $G$ divisible by $n$. If $G$ is written multiplicatively, we write $cH, HL, G^n$ instead of $c + H, H + L, nG$. In this thesis, the Abelian groups that we write multiplicatively are (subgroups of) groups of invertible elements of some commutative ring, and sometimes quotients of those groups, notably the group RV of leading terms of a valued field (see definition 2.2.4). The value group of a valued field can also be written as such a quotient, but we prefer to write it additively, with the value of zero being conventionally smaller than all the others.

The valued fields that we consider in this thesis are always Henselian valued fields of residue characteristic zero, with the exception of local fields, and section 2.2 where the setting is more abstract.

# Chapter 1

# Preliminaries on independence

The notions of forking and dividing were introduced by Shelah ([She82]) during the developments of stability theory. A crucial question was to find generic ways to extend types: given parameter sets $A \subseteq B$, and $p \in S(A)$, is there a natural way to choose a completion of $p$ in $S(B)$? What are the completions of $p$ which add as little relevant information as possible? There are several ways to choose such extensions (heirs, coheirs, definitions that involve the Morley rank…), and these choices usually coincide when the theory of the first-order structure at hand is stable. This gives a very natural notion of non-forking extensions of a type in stable theories. This notion is central in stability theory, in fact the class of stable theories can be characterized as the class of first-order theories for which there exists an abstract way to choose extensions of a type satisfying some axioms (see for instance [TZ12], Theorem 8.5.10). However, stable theories are quite rare, and, in the unstable world, the many ways to define non-forking extensions for stable theories no longer coincide, and they do not have the same behavior, so one has to choose which one best generalizes the behavior that we get in stable theories. In contemporary model theory, we tend to focus on the non-forking independence notion. On top of the fact that many axioms of stable forking still hold for forking in general, it turns out that this notion plays an important role in the new classes of theories which extend the class of stable theories:

- The class of simple theories has an axiomatic characterization in terms of forking which is very similar to that of stable theories (the celebrated Kim-Pillay theorem, see [KP97]).

- In the class of dependent theories (NIP), forking does not satisfy all of these axioms, but it coincides with other independence notions that are easier to understand (see Section 2 of [HP07]).
- Recent work of Chernikov and collaborators ([CK12], [Che14], [BYC14]) shows that forking also plays an important role in the study of $\mathrm{NTP}_2$, the class of theories without the tree property of the second kind.

For this reason, computing forking as a way to study an unstable theory has become a standard approach. This is what we do for concrete theories in this thesis, hence an introduction to forking is needed.

In this chapter, we first define forking and dividing in an arbitrary, potentially unstable theory. We state the basic properties that these two notions enjoy, with the intent to be able to formally manipulate those notions in any structure.

Then, we give easy examples on the behavior of forking in basic algebraic structures. The basic idea is that the non-forking extension of a type should not contain (as a set of formulas) much relevant data. By the usual Stone duality between properties and their contexts, a formula which gives too much data should be intuitively seen as a formula with few realizations. A way to define an independence notion would be to define which formulas are small in this regard (the set of small formulas should be an ideal of the Boolean algebra of definable sets), and expect that a non-forking extension of $p$ is a type which avoids every small formula.

Finally, a recent shift in the literature saw growing interest in generalizing various model-theoretic notions for types to Keisler measures. Not much investigation has been done as of today regarding generalizations of independence notions and non-forking extensions to Keisler measures, this field is still very open. Though we do not have a major breakthrough result about that in this thesis, we would like to suggest a formalism to carry out such a project. We discuss this at the end of this chapter, and we show minor abstract results. More results about Keisler measures in valued fields are present in chapter 7.

## 1.1 Forking and dividing

Let us first recall the definition of forking and dividing. We give examples in section 1.2.

In order to simplify the definitions, we fix $\kappa$ an arbitrarily large infinite cardinal, and $M$ a $\kappa$-saturated, strongly $\kappa$-homogeneous structure. A set $A$ is called *small* if $|A| < \kappa$. The first-order language of $M$ may very well be uncountable, but we assume that it is small.

### 1.1.1 Dividing

**Definition 1.1.1.** Let $A$ be a small subset of $M$, and $\phi(x, b)$ a definable subset of $M$ in a small number of variables. We say that $\phi(x, b)$ *divides* over $A$ when there exists $(b_i)_{i<\omega}$ an $A$-indiscernible sequence with $b_0 = b$, such that the partial type:

$$\{\phi(x, b_i) | i < \omega\}$$

is inconsistent. We shall say that a partial type $\Sigma$ *divides* over $A$ if $\Sigma$ implies some definable set which divides over $A$.

Let $B, C$ be small subsets of $M$. We write $C \mathop{\downarrow}^{\mathbf{d}}_A B$ if $tp(C/AB)$ does not divide over $A$ (for some/any enumeration of $C$).

This definition does not depend on the choice of $M$, one may freely replace $M$ by some arbitrary elementary extension which is sufficiently saturated and strongly homogeneous.

Note that if $C \not\mathop{\downarrow}^{\mathbf{d}}_A B$, then there exists $c_1 \dots c_n \in C$ such that $c_1 \dots c_n \not\mathop{\downarrow}^{\mathbf{d}}_A B$ (we call that property *left finite character*). Dividing also has right finite character, i.e. it satisfies the analogue of this property for $B$. The same cannot be said of $A$, because if $A' \subseteq A$, then $C \not\mathop{\downarrow}^{\mathbf{d}}_{A'} B$ does not necessarily imply $C \not\mathop{\downarrow}^{\mathbf{d}}_A B$, though the converse always holds. For instance, if $A = A'B$, then we always have $C \mathop{\downarrow}^{\mathbf{d}}_A B$, even if $C \not\mathop{\downarrow}^{\mathbf{d}}_{A'} B$.

One may also note that the class of $M$-definable sets which divide over $A$ is downward-closed: if $X$ divides over $A$, and $Y \subseteq X$, then $Y$ divides over $A$. Note that by definition, the empty set always divides over $A$.

*Remark* 1.1.2. An important property of dividing is that we always have $C \mathop{\downarrow}^{\mathbf{d}}_A A$ for every $A, C$. We say that every set is an *extension base* for dividing. This is a property which we look for in a well-behaved independence notion.

**Fact 1.1.3** ([TZ12], Lemma 7.1.4)**.** *Let $X$ be a definable set. Then the following are equivalent:*

1. *The definable set $X$ divides over $A$.*

2. *There exist $N < \omega$, and $(\sigma_n)_{n<\omega} \in \mathrm{Aut}(M/A)^\omega$ such that $\sigma_0 = id$, and, for all $P \in \mathcal{P}(\omega)$, if $|P| = N$, then $\bigcap_{n \in P} \sigma_n(X) = \varnothing$.*

**Definition 1.1.4.** We call a witness of the second condition of fact 1.1.3 a *witness for division* of $X$ over $A$.

This combinatorial characterization of dividing can be seen as a way to state that a definable set is "small": it is not large enough to have a significant overlap with its conjugates.

**Fact 1.1.5** ([TZ12], Corollary 7.1.5)**.** *The following are equivalent:*

- $C \mathop{\downarrow}^{\mathbf{d}}_{A} B$
- *For every $A$-indiscernible sequence $I$ with $B$ the first term of $I$, there exists $J \equiv_{AB} I$ such that $J$ is $AC$-indiscernible.*

The standard way to prove that some partial type divides is to find a witness for division. Fact 1.1.5 helps us to prove that a given complete type does not divide. Another tool that is useful to us in that regard is left-transitivity, which we may call "transitivity" in this thesis, as we never use right-transitivity:

**Fact 1.1.6** ([TZ12], Proposition 7.1.6)**.** *If $C \mathop{\downarrow}^{\mathbf{d}}_{A} B$ and $D \mathop{\downarrow}^{\mathbf{d}}_{AC} B$, then we have $CD \mathop{\downarrow}^{\mathbf{d}}_{A} B$.*

The converse does not hold in general, and we give a counterexample in subsection 4.5.2.

Let us now prove a few technical properties of dividing which we use in chapter 6:

**Proposition 1.1.7.** *If $f$ is an $A$-definable function, and $Y$ is a definable subset of the image of $f$, then $Y$ divides over $A$ if and only if $f^{-1}(Y)$ does.*

*Proof.* Let $(\sigma_n)_n \in \mathrm{Aut}(M/A)^\omega$. As $f$ is $A$-definable, we have the equality:

$$\sigma_n(f^{-1}(Y)) = f^{-1}(\sigma_n(Y))$$

Now, if $P \subseteq \omega$, then:

$$\bigcap_{n \in P} \sigma_n(f^{-1}(Y)) = \bigcap_{n \in P} f^{-1}(\sigma_n(Y)) = f^{-1}\left(\bigcap_{n \in P} \sigma_n(Y)\right)$$

As $Y$ and its $A$-conjugates are subsets of the image of $f$, such an intersection is empty if and only if $\bigcap_{n \in P} \sigma_n(Y) = \varnothing$. As a result, $(\sigma_n)_n$ is a witness for division of $Y$ over $A$ if and only it is also a witness for division of $f^{-1}(Y)$ over $A$, which concludes the proof. □

**Lemma 1.1.8.** *Let $A$ be a small subset of $M$, and let $Z$, $Z'$ be $A$-definable sets, let $R \subseteq Z \times Z'$ be an $A$-definable set, and $m < \omega$. Suppose that, for all $x \in Z$, the set $R_x = \{y \in Z' | R(x,y)\}$ has exactly $m$ elements.*

*Then, if $Y \subseteq Z'$ is an $M$-definable set that divides over $A$, then the definable set $X = \{x \in Z | \exists y \in Y \; R(x,y)\}$ divides over $A$.*

*Proof.* By contraposition, suppose $X$ does not divide over $A$. Let $N < \omega$, and $(\sigma_n)_{n<\omega} \in \mathrm{Aut}(M/A)^\omega$ such that $\sigma_0 = id$. As $X$ does not divide over $A$, there must exist $P$ a subset of $\omega$ of size $Nm$ such that $\bigcap_{n \in P} \sigma_n(X) \neq \varnothing$. Let $\alpha$ be in that intersection. As $\alpha \in Z$, let $\beta_1 \dots \beta_m$ be the elements of $R_\alpha$. For each $n \in P$, let $1 \leqslant f(n) \leqslant m$ such that $\beta_{f(n)} \in \sigma_n(Y)$. By the pigeonhole principle, there exists $1 \leqslant q \leqslant m$ such that $q$ has at least $N$ antecedents by $f$. Let $Q = \{n \in P | f(n) = q\}$. Then $Q$ is a subset of $\omega$ of size at least $N$ for which $\beta_q \in \bigcap_{n \in Q} \sigma_n(Y)$, hence $N, (\sigma_n)_n$ is not a witness for division of $Y$ over $A$. This holds for all $N, (\sigma_n)_n$, so $Y$ does not divide over $A$, concluding the proof. □

**Lemma 1.1.9.** *Let $A$, $C$ be small subsets of $M$, and $X$ a definable subset of $M$ in a small number of variables that divides over $A$. Let $\lambda$ be the (small) cardinal of the set of formulas with parameters in $AC$. Suppose $(2^\lambda)^+ < \kappa$ (an assumption that can always be made as $\kappa$ is arbitrarily large, and we can replace $M$ by a sufficiently large elementary extension). Then there exists $C' \equiv_A C$, $C' \subseteq M$, such that $X$ divides over $AC'$.*

*Proof.* Let $b \in M$ such that $X$ is $b$-definable. By a standard compactness argument, from a witness for division of $X$ over $A$, we can find a sequence $(\sigma_i)_{i<(2^\lambda)^+} \in \mathrm{Aut}(M/A)^{(2^\lambda)^+}$, $N < \omega$, such that, for all $P \subseteq (2^\lambda)^+$, if $|P| = N$, then $\bigcap_{i \in P} \sigma_i(X) = \varnothing$. Let us look at the map $f : i \longmapsto tp(\sigma_i(b)/AC)$. The image of $f$ is smaller than $2^\lambda$, so we can apply the pigeonhole principle to find $i_0 < i_1 < i_2 \dots$ such that, for all $n < \omega$, $f(i_{n+1}) = f(i_n)$. Let $C' = \sigma_{i_0}^{-1}(C) \equiv_A C$. For each $n < \omega$, $\sigma_{i_0}^{-1} \circ \sigma_{i_n}(b)$ is a $\sigma_{i_0}^{-1}(AC) = AC'$-conjugate of $b$, so we can use strong homogeneity to find $\tau_n \in \mathrm{Aut}(M/AC')$ ($\tau_0 = id$) such that $\tau_n(b) = \sigma_{i_o}^{-1} \circ \sigma_{i_n}(b)$, i.e. $\tau_n(X) = \sigma_{i_o}^{-1} \circ \sigma_{i_n}(X)$. Then $N$, $(\tau_n)_{n<\omega}$ is a witness for division of $X$ over $AC'$. □

The main problem with dividing is that this independence notion is too weak, i.e. there are definable sets which do not divide over $A$, despite the fact that they should. For instance, one may find a $B$-definable set $X$ which does not divide over $A$, but such that for all $c \vDash X$, we have $c \not\downarrow^{\mathbf{d}}_A B$. We would like to strengthen the notion of dividing to avoid that.

*Remark* 1.1.10. By compactness, the following are equivalent:

- For all $c \vDash X$, we have $c \not\downarrow^{\mathbf{d}}_A B$.
- $X$ may be written as a finite union of $B$-definable sets which divide over $A$.

As a result, we would like the class of definable sets which give too much data with respect to our independence notion to be closed under finite union. Just as for dividing, it should also be downward-closed, thus it should be an ideal. All this motivates the definition of forking.

### 1.1.2 Forking

**Definition 1.1.11.** We say that $X$ *forks* over $A$ if there exists a finite family $X_1, \dots, X_n$ of definable subsets of $M$ such that $X = \bigcup_i X_i$, and $X_i$ divides over $A$ for all $i$. In other words, the class of $M$-definable sets which fork over $A$ is the ideal of the Boolean algebra of $M$-definable sets which is generated by those that divide over $A$.

We shall say that a partial type $\Sigma$ *forks* over $A$ if $\Sigma$ implies a definable set which forks over $A$.

Let $B, C$ be small subsets of $M$. We write $C \downarrow^{\mathbf{f}}_A B$ if $tp(C/AB)$ does not fork over $A$.

This definition may also be taken in $M^{eq}$, for the natural correspondence between $\mathrm{Aut}(M/A)$ and $\mathrm{Aut}(M^{eq}/A)$ (when $A \subseteq M$) implies that we have $C \downarrow^{\mathbf{f}}_A B$ in $M$ if and only if $C \downarrow^{\mathbf{f}}_A B$ in $M^{eq}$.

*Remark* 1.1.12. For all $\sigma \in \mathrm{Aut}(M)$, we have $C \mathop{\downarrow}^{\mathbf{f}}_{A} B$ if and only if we also have $\sigma(C) \mathop{\downarrow}^{\mathbf{f}}_{\sigma(A)} \sigma(B)$, and $C \mathop{\downarrow}^{\mathbf{d}}_{A} B$ if and only if $\sigma(C) \mathop{\downarrow}^{\mathbf{d}}_{\sigma(A)} \sigma(B)$

Note that dividing always implies forking. However, the folklore example of cyclical orders ([TZ12], Exercise 7.1.6) shows that the converse does not hold in general.

Just as dividing, forking has left and right finite character.

An important subtlety of definition 1.1.11 is that if $X$ is $B$-definable and forks over $A$, then the witnesses $X_i$ for forking are $M$-definable sets which may not be $B$-definable. However, by replacing the $X_i$ by well-chosen $AB$-conjugates, we may assume that they can always be definable over some sufficiently saturated model containing $AB$. In particular:

*Remark* 1.1.13 ([TZ12], Exercise 7.1.3). If $B$ is an $|AC|^+$-saturated elementary substructure of $M$ that contains $A$, then we have $C \mathop{\downarrow}^{\mathbf{f}}_{A} B$ if and only if $C \mathop{\downarrow}^{\mathbf{d}}_{A} B$.

*Remark* 1.1.14. It easily follows from the definitions that any definable set which has a point in $A$ does not fork (hence does not divide) over $A$. In particular, any type which is finitely satisfiable in $A$ does not fork over $A$.

The non-forking independence notion is a strengthening of non-dividing which fixes the problem that we previously noted:

*Remark* 1.1.15. Let $X$ be a $B$-definable set. Then, by compactness, the following are equivalent:

- For all $c \models X$, we have $c \mathop{\not\downarrow}^{\mathbf{f}}_{A} B$
- $X$ forks over $A$.

**Corollary 1.1.16.** *If $p$ is a partial type over $AB$ which does not fork over $A$, then there exists a complete type $q \in S(AB)$ extending $p$ such that $q$ does not fork over $A$.*

**Corollary 1.1.17.** *The following are equivalent:*

- $C \mathop{\downarrow}^{\mathbf{f}}_{A} B$
- *In some elementary extension of $M$, there exists $D \equiv_{AB} C$ such that* $D \mathop{\downarrow}^{\mathbf{d}}_{A} M$

The fact that $M$ is sufficiently saturated is important here.

*Proof.* Extend by corollary 1.1.16 the partial type $\mathrm{tp}(C/AB)$ to a complete type $q$ over $M$ which does not fork over $A$. Choose $D$ a realization of $q$. Then $D \equiv_{AB} C$, and $D \mathop{\smile\!\!\!\!\mid}^{\mathbf{f}}_{A} M$. This concludes the proof by remark 1.1.13, as $M$ is sufficiently saturated. □

This extension property does not hold for dividing in general. Forking can be defined as the sub-relation of $\mathop{\smile\!\!\!\!\mid}^{\mathbf{d}}$ satisfying this extension property, that is "maximal" in some sense (see [Cas11], Proposition 12.14).

However, by going from dividing to forking, we lose a crucial property: not all sets are extension bases anymore. Therefore, the study of extension bases that we carry out in chapter 6 is the research of a setting where forking keeps the good properties that hold for dividing, where we can get the best of both worlds. One property that we always keep is transitivity:

**Proposition 1.1.18.** *Suppose $C \mathop{\smile\!\!\!\!\mid}^{\mathbf{f}}_{A} B$ and $D \mathop{\smile\!\!\!\!\mid}^{\mathbf{f}}_{AC} B$. Then we have $CD \mathop{\smile\!\!\!\!\mid}^{\mathbf{f}}_{A} B$.*

*Proof.* Let $M'$ be some $|M|^+$-saturated, strongly $|M|^+$-homogeneous elementary extension of $M$. By corollary 1.1.17 and strong homogeneity, let $\sigma \in \mathrm{Aut}(M'/AB)$ be such that $\sigma(C) \mathop{\smile\!\!\!\!\mid}^{\mathbf{d}}_{A} M$. By remark 1.1.12, we have:

$$C \mathop{\smile\!\!\!\!\mid}^{\mathbf{d}}_{A} \sigma^{-1}(M)$$

By corollary 1.1.16, $\mathrm{tp}(D/ABC)$ has a global extension to $\sigma^{-1}(M)$ which does not fork over $AC$, thus there exists $\tau \in \mathrm{Aut}(M'/ABC)$ such that:

$$\tau(D) \mathop{\smile\!\!\!\!\mid}^{\mathbf{f}}_{AC} \sigma^{-1}(M)$$

By transitivity for dividing (fact 1.1.6), we have $C\tau(D) \mathop{\smile\!\!\!\!\mid}^{\mathbf{d}}_{A} \sigma^{-1}(M)$. By applying $\sigma$, we have $\sigma \circ \tau(CD) \mathop{\smile\!\!\!\!\mid}^{\mathbf{d}}_{A} M$. We conclude by corollary 1.1.17. □

The counterexample in 4.5.2 shows that the converse is false in general.

**Fact 1.1.19.** *Let $A$, $B$, $C$ be small subsets of $M$. Then:*

1. *Any partial type divides/forks over $A$ if and only if it divides/forks over $\mathrm{acl}(A)$. (This is because every $A$-indiscernible sequence is $\mathrm{acl}(A)$-indiscernible.)*

2. *If $B$ is an $|AC|^+$-saturated, strongly $|AC|^+$-homogeneous elementary substructure of $M$, $p$ is a complete type over $B$ in a finite number of variables, and the orbit of $p$ under* $\mathrm{Aut}(B/A)$ *is smaller than* $|A|^+$ *(in particular, if this orbit is a point), then $p$ does not divide over $A$. (See for instance [TZ12], exercises 7.1.3 and 7.1.5.)*

Let us show that forking preserves algebraic closure. The following facts are already known, but they were not proved explicitly in the literature for a long time, until our article on extension bases ([Hos23a]), and recent work from Conant-Kruckman ([CK23]). Even though forking preserves algebraic closure, Conant and Kruckman show counterexamples where dividing does not.

**Proposition 1.1.20.** *Let $A$, $B$, $C$, $D$ be small subsets of $M$ such that* $\mathrm{acl}(AC) \subseteq \mathrm{acl}(AD)$. *If $D \mathop{\smile\!\!\!\!\!\!\!\;|}^{\mathbf{f}}_A B$, then $C \mathop{\smile\!\!\!\!\!\!\!\;|}^{\mathbf{f}}_A B$.*

*Proof.* Let $c$ be a finite tuple from $C$, $Y$ an $AB$-definable set that contains $c$, $p < \omega$, and $(Y_i)_{0 \leqslant i \leqslant p}$ some $M$-definable sets such that $Y \subseteq \bigcup_i Y_i$. It suffices to show that $Y$ does not fork over $A$. Let $\phi(x, y)$ be a formula with parameters in $A$ such that $\phi(D, y)$ is the least finite (say, of size $m < \omega$) $AD$-definable set that contains $c$. Let $\psi(x)$ be the $A$-definable formula stating that there exists exactly $m$ distinct points $y_1..y_m$ for which $\phi(x, y_i)$ holds for all $i$. Let $X$ be the $AB$-definable set defined by the following formula with free variable $x$:

$$\psi(x) \wedge (\exists y \ \ (\phi(x, y) \wedge y \in Y))$$

Define similarly $X_i$ from $Y_i$. By hypothesis, as $D \in X$, $X$ does not fork over $A$. We clearly have $X \subseteq \bigcup_i X_i$, so there exists $j$ such that $X_j$ does not divide over $A$. By lemma 1.1.8, $Y_j$ does not divide over $A$. This holds for every choice of a finite family $(Y_i)_i$, so $Y$ does not fork over $A$. □

**Corollary 1.1.21.** *We have $C \mathop{\smile\!\!\!\!\!\!\!\;|}^{\mathbf{f}}_A B$ if and only if* $\mathrm{acl}(AC) \mathop{\smile\!\!\!\!\!\!\!\;|}^{\mathbf{f}}_A B$.

Corollary 1.1.21 is not to be underrated, it is extremely useful for anyone (us in particular) that wishes to compute the non-forking independence relation $\mathop{\smile\!\!\!\!\!\!\!\;|}^{\mathbf{f}}$ in a concrete theory. Instead of working with $C$ (or a finite tuple from $C$ by left finite character), one may choose to work instead with a tuple which is $A$-interalgebraic with $C$, and whose type is easier to understand. The variant of this statement that we use the most is the following:

**Corollary 1.1.22.** *Suppose we have $\overline{C}$ a set of finite tuples from $C$ such that for every finite tuple $d$ of $C$, there exists $c$ an element of $\overline{C}$ such that* $\mathrm{acl}(Ad) \subseteq \mathrm{acl}(Ac)$*. Then $C \underset{A}{\perp^{\mathbf{f}}} B$ if and only if, for every $c \in \overline{C}$, we have $c \underset{A}{\perp^{\mathbf{f}}} B$.*

For instance, a direct application is the following:

**Corollary 1.1.23.** *Suppose $M$ expands a field. Let $c = (c_i)_i$ be some maximal (field-theoretically)* $\mathrm{acl}(A)$*-algebraically independent tuple from* $\mathrm{acl}(AC)$*. Then $C \underset{A}{\perp^{\mathbf{f}}} B$ if and only if, for every finite subtuple $d$ of $c$, we have $d \underset{A}{\perp^{\mathbf{f}}} B$.*

Working with a finite tuple of algebraically independent elements is arguably easier than an arbitrary infinite tuple. Note that in this instance, one can choose $c$ with coordinates in $C$ itself.

We should specify that acl denotes here (and throughout the whole thesis) model-theoretic algebraic closure. We recall that we denote the field-theoretic algebraic closure of a field $k$ by $k^{\mathbf{alg}}$. We will recall this notation whenever it is relevant. We conclude this section by studying how forking interacts with acl on the right:

**Proposition 1.1.24.** *If $A$, $B$, $C$ are small subsets of $M$, and $C \underset{A}{\not\perp^{\mathbf{f}}} \mathrm{acl}(AB)$, then $C \underset{A}{\not\perp^{\mathbf{f}}} B$.*

Note that the other direction is an easy consequence of the definition of forking.

*Proof.* Let $X$ be an $\mathrm{acl}(AB)$-definable set containing $C$ that forks over $A$. Let $X_1 \dots X_n$ witness that $X$ forks over $A$, i.e. $X \subseteq \bigcup_i X_i$ and each $X_i$ divides over $A$. Now, $X$ has finitely many $AB$-conjugates, so there exist $N < \omega$ and $\sigma_1 \dots \sigma_N \in \mathrm{Aut}(M/AB)$ such that $\{\sigma_j(X) | 1 \leqslant j \leqslant N\}$ is the orbit of $X$ under the action of $\mathrm{Aut}(M/AB)$. The union $Y = \bigcup_j \sigma_j(X)$ is an $AB$-definable set containing $C$. As each $\sigma_j$ pointwise-fixes $A$, each $\sigma_j(X_i)$ divides over $A$ (apply $\sigma_j$ to a witness for division). As $Y \subseteq \bigcup_{ij} \sigma_j(X_i)$, we have $C \underset{A}{\not\perp^{\mathbf{f}}} B$, which concludes the proof. □

**Corollary 1.1.25.** *If $A$, $B$, $C$ are small subsets of $M$, then $C \underset{A}{\perp^{\mathbf{f}}} B$ if and only if $C \underset{\mathrm{acl}(A)}{\perp^{\mathbf{f}}} B$.*

*Proof.* Suppose $C \underset{A}{\not\downarrow^{\mathbf{f}}} B$. Then, by fact 1.1.19, $\operatorname{tp}(C/AB)$ forks over $\operatorname{acl}(A)$, hence so does $\operatorname{tp}(C/\operatorname{acl}(A)B)$, thus $C \underset{\operatorname{acl}(A)}{\not\downarrow^{\mathbf{f}}} B$.

Conversely, suppose $C \underset{\operatorname{acl}(A)}{\not\downarrow^{\mathbf{f}}} B$. Then, we have $C \underset{\operatorname{acl}(A)}{\not\downarrow^{\mathbf{f}}} \operatorname{acl}(AB)$, thus $C \underset{A}{\not\downarrow^{\mathbf{f}}} \operatorname{acl}(AB)$. By proposition 1.1.24, we have $C \underset{A}{\not\downarrow^{\mathbf{f}}} B$. □

## 1.2 Examples

In this chapter, we compute forking in common first-order theories, with proofs that are more direct than using the usual machinery from stability theory.

### 1.2.1 Introductory easy examples

We compute forking in very simple theories, as a way to show in easy examples what the general approach should be when one wants to compute forking in a specific (potentially unstable) theory.

The most simple theory is INFSET, the theory of infinite first-order structures on the empty language. Given small parameter sets $A$, $B$, $C$ in this theory, when do we have $C \underset{A}{\downarrow^{\mathbf{f}}} B$? First of all, we use left finite character and corollary 1.1.22 to restrict to the case where $C$ has a specific form that is easier to understand:

**Lemma 1.2.1.** *The class $\overline{C}$ of finite tuples from $C \smallsetminus A$ with pairwise-distinct elements satisfies the hypothesis of corollary 1.1.22.*

Now let $c = (c_1, \dots, c_n) \in \overline{C}$, and let $p = \operatorname{tp}(c/A)$. We want to describe which extensions of $p$ in $S(AB)$ do not fork over $A$. In order to turn $p$ into a complete type over $AB$, we just need to specify which points of the tuple are in $B$. For each coordinate of the tuple, one either has to choose an element of $B \smallsetminus A$, or choose that the point does not belong to $B$. Choosing an arbitrary element of $B \smallsetminus A$ does not seem very natural, and indeed:

**Lemma 1.2.2.** *In any theory, if $b \notin \operatorname{acl}(A)$, then the formula $x = b$ divides over $A$.*

*Proof.* A witness for division is $N = 2$, and $(\sigma_n)_n$ a sequence of automorphisms such that $(\sigma_n(b))_n$ is a family of pairwise-distinct $A$-conjugates of $b$. □

The only candidate left for a non-forking extension of $p$ is the one according to which no point is in $B$:

**Lemma 1.2.3.** *The following complete type in $S(AB)$ does not fork over $A$:*

$$p(x_1 \dots x_n) \cup \bigcup_i \{x_i \neq b | b \in B \setminus A\}$$

*Proof.* This type admits a global $\mathrm{Aut}(M/A)$-invariant extension:

$$p(x_1 \dots x_n) \cup \bigcup_i \{x_i \neq m | m \in M \setminus A\}$$

and we conclude by the second item of fact 1.1.19. □

It follows that $c \mathop{\downarrow}^{\mathbf{f}}_A B$ if and only if no $c_i$ belongs to $B \setminus A$. As a result, we have $C \mathop{\downarrow}^{\mathbf{f}}_A B$ if and only if $C \cap B \subseteq A$. More generally, we have:

**Proposition 1.2.4.** *In any theory, if* $\mathrm{acl}(AC) \cap \mathrm{acl}(AB) \neq \mathrm{acl}(A)$*, then* $C \mathop{\not\downarrow}^{\mathbf{d}}_A B$*.*

*Proof.* This follows immediately from lemma 1.2.2, and the contraposite of corollary 1.1.22. □

Let us show that the same thing can be done in the theory $k-\mathrm{VEC}$ of infinite $k$-vector spaces on a fixed field $k$.

Suppose now $M$ is a vector space over some field $k$.

**Lemma 1.2.5.** *The class $\overline{C}$ of $k$-linearly free finite tuples $c$ from $C$, such that* $\mathrm{Vec}(c) \cap \mathrm{Vec}(A) = \{0\}$*, satisfies the hypothesis of corollary 1.1.22.*

**Definition 1.2.6.** Let $R$ be some commutative ring. Suppose $M$ is an expansion of some $R$-module, and let $f$ be some $n$-ary $\varnothing$-definable function in this structure. We write $f \in \mathrm{LC}^n(R)$ when there exist $\lambda_1 \dots \lambda_n \in R$ such that $f$ coincides with the function $x_1 \dots x_n \longmapsto \sum_i \lambda_i \cdot x_i$. We write $\mathrm{LC}(R)$ when $n$ is implicit.

Let $c = (c_1, \ldots, c_n) \in \overline{C}$, $p = \mathrm{tp}(c/A)$. Informally, to turn $p$ into a complete type over $AB$, one must know which linear combination of the tuple belongs to $\mathrm{Vec}(AB)$. This amounts to choose $V$ a vector subspace of $\mathrm{LC}(k)$ for the space of such linear combinations, and an embedding $V \longrightarrow \mathrm{Vec}(AB)$ for their respective values. Just as in INFSET, choosing a non-trivial embedding in to $\mathrm{Vec}(AB)$ is not natural, and it leads to a type that forks over $A$ by lemma 1.2.2. The extension of $p$ that remains is the one induced by the trivial embedding:

**Lemma 1.2.7.** *The following type does not fork over $A$:*

$$p(x) \cup \{f(x) \neq b | f \in \mathrm{LC}(k), f \neq 0, b \in \mathrm{Vec}(AB)\}$$

*Proof.* We find a global $\mathrm{Aut}(M/A)$-invariant extension:

$$p(x) \cup \{f(x) \neq m | f \in \mathrm{LC}(k), f \neq 0, m \in M\}$$

and we conclude by the second item of fact 1.1.19. □

**Corollary 1.2.8.** *Let $k$ be some field. Then, in the theory of infinite $k$-vector spaces, we have $C \mathop{\downarrow}^{\mathbf{f}}_{A} B$ if and only if* $\mathrm{Vec}_k(AC) \cap \mathrm{Vec}_k(AB) = \mathrm{Vec}_k(A)$.

In general, it is better to work with imaginaries, as there are examples of forking due to non-trivial intersection of the algebraic closures involving $\mathrm{acl}^{eq}$ instead of acl.

*Example* 1.2.9. In some elementary extension of the Abelian group $\mathbf{Z}^{\oplus\omega}$, let $A = \{0\}$, $B = \{(1,0,0\ldots)\}$, $C = \{(1,2,0,0\ldots)\}$.

Let us show that $\mathrm{acl}(AB) \cap \mathrm{acl}(AC) = \mathrm{acl}(A)$. For each permutation $\sigma$ of $\omega$, one can build an automorphism of $\mathbf{Z}^{\oplus\omega}$ permuting the coordinates of its elements according to $\sigma$. For each $x \notin A$, we can easily find infinitely many permutations $(\sigma_n)_n$ such that the images of $x$ are pairwise-distinct $A$-conjugates of $x$, thus $A = \mathrm{acl}(A)$. Likewise, we can easily show that $\mathrm{acl}(AB) = \mathbb{Z} \cdot (1,0,0\ldots)$, and $\mathrm{acl}(AC) = \mathbb{Z} \cdot (1,2,0,0\ldots)$, in particular $\mathrm{acl}(AB) \cap \mathrm{acl}(AC) = \mathrm{acl}(A)$.

However, $\mathrm{acl}^{eq}(AB) \cap \mathrm{acl}^{eq}(AC) \neq \mathrm{acl}^{eq}(A)$. To see this, let $G$ be the $A$-definable group of elements divisible by 2. Let $\alpha \in M^{eq}$ be some code of $(1,0,0,\ldots) + G$. Then we clearly have $\alpha \in \mathrm{acl}^{eq}(AB) \cap \mathrm{acl}^{eq}(AC)$. However, by exchanging the 0-th coordinate with the $n$-th coordinate for $n < \omega$, we get infinitely many $A$-conjugates of $\alpha$, thus $\alpha \notin \mathrm{acl}^{eq}(A)$. It follows from proposition 1.2.4 that $C \mathop{\not\downarrow}^{\mathbf{d}}_{A} B$ (in $M$ and $M^{eq}$, it does not matter !).

Though this may already be known in the literature, in a torsion-free Abelian group, we have $C \mathop{\downarrow}^{\mathbf{f}}_{A} B$ if and only the phenomenon from example 1.2.9 does not happen, that is when for all prime $l$ for which $[M : lM]$ is infinite, if $0 < N < \omega$, then $(C' + l^N M) \cap (B' + l^N M) = A' + l^N M$, with $A', B', C'$ the respective relative divisible closures of the subgroup generated by $A, B, C$. We will not give a proof since it is more involved, but a direct one can be reconstructed from the arguments in chapter 4, basically by removing everything involving the total order.

### 1.2.2 Dense total orders

In all the previous examples, the non-forking independence relation is in some sense trivial, it coincides with naive set-theoretic independence of the respective algebraic closures. The most simple example of a theory when forking has a non-trivial behavior is perhaps in totally ordered sets.

**Notation** In the theory of some total order $<$, a set is considered an interval (be it closed, open, half-open, possibly with infinite bounds) when its lower bound is smaller or equal to its upper bound. For instance, if we consider, say an interval $]a, b[$, then it is implicit that $a \leqslant b$.

*Example* 1.2.10. Suppose $M$ is an elementary extension of $(\mathbb{Q}, <) \vDash \mathrm{DLO}$. Let $A = \varnothing$, $B = \{0, 2\}$, $C = \{1\}$. Then $\mathrm{acl}(AC) \cap \mathrm{acl}(AB) = \mathrm{acl}(A)$, but $C \mathop{\not\downarrow}^{\mathbf{d}}_{A} B$. Indeed, if we define $I_n = [4n, 4n+2]$, then the $(I_n)_{n<\omega}$ are pairwise-distinct and $A$-conjugates, thus they divide over $A$, and $1 \in I_0$.

In this example, the definable set that forks is an interval which is disjoint from some conjugates. We can see with the same reasoning that such intervals always divide:

**Lemma 1.2.11.** *In the expansion of some total order, let $I$ be an interval with non-empty interior such that its bounds are $A$-conjugates. Then $I$ divides over $A$.*

**Corollary 1.2.12.** *Suppose that the theory of $M$ is o-minimal, and let $I = [b_1, b_2]$ be an interval which is closed, bounded and disjoint from* $\mathrm{dcl}(A)$ *(in particular, $b_1$ and $b_2$ are $A$-conjugates). Then $I$ divides over $A$.*

Note that we only consider bounded intervals: the bounds are not in $\{\pm\infty\}$, for they must be $A$-conjugates.

We have a basic example of a definable set which forks in DLO. Does it get more complicated than that ? We show that this is not the case for DLO. The readers that are interested in chapter 3 should read the proof, as it is a very simplified version of what we do in DOAG. First let us define the basic notion of cuts. We see cuts as type-definable sets which can be written as intersections of intervals, as opposed to sets of upper/lower bounds of said intervals. This non-standard formalism is much more natural for our work in chapters 3 and 4.

**Definition 1.2.13.** Suppose $M$ expands a total order. Let $c \in M$. The *cut* of $c$ over $A$ is the $A$-type-definable set $\mathrm{ct}(c/A)$, defined as the intersection of every interval (closed, open, and half-open) containing $c$ with bounds in $A \cup \{\pm\infty\}$ (a singleton is a closed interval).

Let $P$ be some partition of $M$ into convex sets. Then there is a natural total order on $P$, the only total order that makes the projection $M \longrightarrow P$ an order-preserving map.

The set of cuts over $A$ is clearly a partition of $M$ into convex sets, therefore the above paragraph applies.

Let us also write $x > A$ when $\forall a \in A,\ x > a$. If $\emptyset \neq A' \subseteq A$, define $\mathrm{ct}_>(A'/A)$ (resp. $\mathrm{ct}_<(A'/A)$) as the cut over $A$ that corresponds to the elements that are $> A'$ (resp. $< A'$), and strictly smaller (resp. larger) than any element of $A$ that is $> A'$ (resp. $< A'$). More generally, if $X$ is some (type/$\vee$)-definable subset of $M$, we write $\mathrm{ct}_<(X/A) = \mathrm{ct}_<(X(A)/A)$, and similarly for $\mathrm{ct}_>(X/A)$, where $X(A)$ is the set of elements of $A$ which belong to $X$.

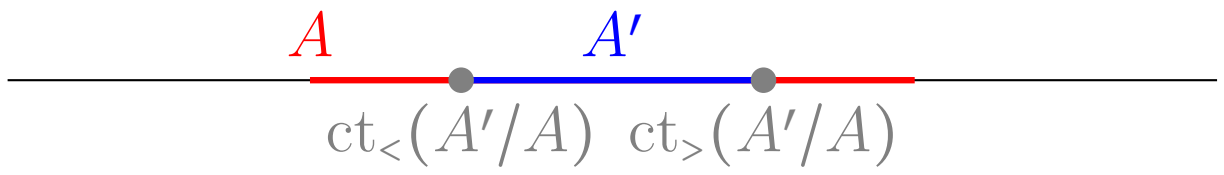

This definition will be used for abstract expansions of total orders in chapter 3.

**Proposition 1.2.14.** *Suppose the theory of $M$ is o-minimal. Then, for every singleton $c \in M$, the following are equivalent:*

1. $c \mathop{\downarrow}^{\mathbf{f}}_{A} B$
2. $\mathrm{tp}(c/AB)$ *has a global* $\mathrm{Aut}(M/A)$*-invariant extension.*

3. *Every closed interval with bounds in* $\operatorname{dcl}(AB)$ *containing* $c$ *has a point in* $\operatorname{dcl}(A)$.

*Proof.* $1 \implies 3$ follows from corollary 1.2.12, and $2 \implies 1$ always holds (we recall that $M$ is $\kappa$-saturated and strongly $\kappa$-homogeneous for some infinite $\kappa > |AB|$.

For $3 \implies 2$, assume condition 3 holds. Then the elements of $\operatorname{dcl}(AB)$ having the same cut as $c$ must all be either strictly smaller or strictly larger than $c$. As a result, either:

$$\operatorname{ct}_<(\operatorname{ct}(c/\operatorname{dcl}(A))/M)$$

or

$$\operatorname{ct}_>(\operatorname{ct}(c/\operatorname{dcl}(A))/M)$$

corresponds to a global (and complete by o-minimality) unary type that extends $\operatorname{tp}(c/AB)$. As $\operatorname{ct}(c/\operatorname{dcl}(A))$ is $\operatorname{Aut}(M/A)$-invariant, so are both those types, and we get the proposition. □

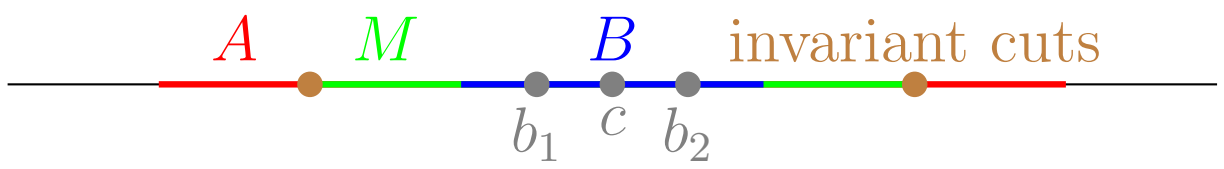

*Remark* 1.2.15. By the above proof, given a singleton $c \notin \operatorname{dcl}(A)$, $\operatorname{tp}(c/A)$ has exactly two global non-forking extensions: the one that we call left-generic, the only one that contains the formulas $x < d$ for every $d \equiv_A c$, and the one that we call right-generic containing $x > d$.

We showed that forking for singletons is entirely controlled by the intervals we found in corollary 1.2.12. Of course this is not true for arbitrary tuples when the geometry is too complicated, such as in RCF. However, the geometry in DLO is trivial (in the sense of trichotomy theory).

The technique that we use to compute forking in DLO is a simplified prototype of what we do in chapter 3: we split our tuple into smaller subtuples (the *blocks*), we find a non-forking extension of each block, and we glue them into a non-forking extension of the whole tuple. A setting where it is very easy to glue types is when they are orthogonal:

**Definition 1.2.16.** Let $p_1 \dots p_n \in S(A)$. We say that the family $(p_i)_i$ is *weakly orthogonal* if $\bigcup_i p_i(x_i)$ is a complete type over $A$ in the variables $x_1 \dots x_n$.

**Lemma 1.2.17.** *Let $c = (c_i)_i$, $d = (d_j)_j$ be two tuples of a model of* DLO*. Suppose we have* $\mathrm{ct}(c_i/A) \neq \mathrm{ct}(d_j/A)$ *for all* $i, j$*. Then* $\mathrm{tp}(c/A)$ *and* $\mathrm{tp}(d/A)$ *are weakly orthogonal.*

*Proof.* Recall that DLO eliminates quantifiers in the language $\{<\}$. An atomic formula in this language is either of the form $t = t'$ or $t < t'$, with $t$, $t'$ either a parameter or a single variable. As a result, an atomic formula with parameters in $A$ which involves both the tuples of variables $(x_i)_i$ and $(y_j)_j$ is either of the form $x_i = y_j$, $x_i < y_j$ or $x_i > y_j$ for some $i$, $j$. However, as $\mathrm{ct}(c_i/A) \neq \mathrm{ct}(d_j/A)$, there must exist $a \in A$ such that either $c_i < a < d_j$ or $c_i > a > d_j$. Either way, $\mathrm{tp}(c/A) \cup \mathrm{tp}(d/A)$ either implies $x_i < y_j$ or $x_i > y_j$. It follows that this partial type is complete, concluding the proof. □

**Lemma 1.2.18.** *In any theory, if* $\mathrm{tp}(c/AB)$ *and* $\mathrm{tp}(d/AB)$ *are weakly orthogonal,* $c \mathop{\perp\!\!\!\!\perp}^{\mathbf{f}}_A B$*, and* $d \mathop{\perp\!\!\!\!\perp}^{\mathbf{f}}_A B$*, then we have* $cd \mathop{\perp\!\!\!\!\perp}^{\mathbf{f}}_A B$*.*

*Proof.* By weak orthogonality, $\mathrm{tp}(c/AB)$ is complete in $S(ABd)$, hence $c \mathop{\perp\!\!\!\!\perp}^{\mathbf{f}}_A Bd$, which implies $c \mathop{\perp\!\!\!\!\perp}^{\mathbf{f}}_{Ad} B$. We conclude by proposition 1.1.18. □

**Corollary 1.2.19.** *In* DLO*, let $c = (c_{ij})_{ij}$ be a tuple enumerated such that $c_{ij} \equiv_{AB} c_{kl}$ if and only if $i = k$ for all $i, j, k, l$. Then we have $c \mathop{\perp\!\!\!\!\perp}^{\mathbf{f}}_A B$ if and only if, for every $i$, we have $(c_{ij})_j \mathop{\perp\!\!\!\!\perp}^{\mathbf{f}}_A B$.*

The blocks are the $((c_{ij})_j)_i$, and we may glue a non-forking extension of their types by union as proved in lemma 1.2.18. Now, the following provides a non-forking extension for each block:

**Lemma 1.2.20.** *In* DLO*, let $p$ be a unary type over $A$, and let $c = (c_j)_j$ be a tuple of realizations of $p$. Let $\phi$ be the parameter-free formula which isolates* $\mathrm{tp}(c/\varnothing)$*. Then the partial type:*

$$q(x) : \ \{\phi\} \cup \bigcup_j p(x_j)$$

*is complete, hence coincides with* $\mathrm{tp}(c/A)$*.*

*Proof.* Let $\psi$ be an atomic formula with parameters in $A$. If $\psi$ involves two distinct variables, then $\psi$ is parameter-free, thus it (resp. its negation) is satisfied by $c$ if and only if it (resp. its negation) is implied by $\phi$. If $\psi$ only involves one variable $x_j$, then the same applies with $p(x_j)$ instead of $\phi$. We conclude that $q$ is complete by quantifier elimination. □

**Lemma 1.2.21.** *In* DLO*, suppose* $p \in S(AB)$ *is a unary type which does not fork over* $A$*. Let* $c = (c_j)_j$ *be a tuple of realizations of* $p$*. Then* $c \mathop{\downarrow}^{\mathbf{f}}_{A} B$*.*

*Proof.* Let $\phi$ be the parameter-free formula which isolates $\operatorname{tp}(c/\varnothing)$. By proposition 1.2.14 (the implication $1 \implies 2$), let $q$ be a global $\operatorname{Aut}(M/A)$-invariant type extending $p$. One may show by a standard compactness argument that the partial type:

$$r(x): \{\phi\} \cup \bigcup_j q(x_j)$$

is consistent. This type is a complete global extension of $\operatorname{tp}(c/AB)$ by lemma 1.2.20, and it is $\operatorname{Aut}(M/A)$-invariant by definition, concluding the proof. □

**Corollary 1.2.22.** *In* DLO*, the following are equivalent:*

- $C \mathop{\downarrow}^{\mathbf{f}}_{A} B$
- *Every closed bounded* $AB$*-definable interval having a point in* $C$ *has a point in* $A$*.*

*Proof.* The top-to-bottom direction follows from lemma 1.2.11 by contraposition. Suppose the second condition holds. Let $c$ be a finite tuple from $C$. Enumerate $c = (c_{ij})_{ij}$ just as in corollary 1.2.19. We have $d \mathop{\downarrow}^{\mathbf{f}}_{A} B$ for each singleton $d$ from $c$ by proposition 1.2.14. It follows that each block is independent: $(c_{ij})_j \mathop{\downarrow}^{\mathbf{f}}_{A} B$ for every $i$, thus $c \mathop{\downarrow}^{\mathbf{f}}_{A} B$ by corollary 1.2.19, concluding the proof. □

We will show in chapter 3 that a similar description of forking can be done in DOAG, though the proof is much more technical.

### 1.2.3 Fields

In the previous subsection, we studied the behavior of forking in various basic first-order structures, such as structures with a trivial geometry, vector spaces, Abelian groups, and totally ordered sets. In these structures, nonforking refines, and sometimes coincides with various independence notions of geometric nature, such as trivial intersection, linear independence, or some variant of finite satisfiability for intervals. Fields are also examples of common first-order structures, and we have an independence notion of geometric nature in these structures, namely algebraic independence:

**Definition 1.2.23.** Suppose $M$ expands a field. We write $C \downarrow^{\mathbf{alg}}_A B$ if for any finite tuple $c$ from $C$, we have the equality of field-theoretic transcendence degrees $\mathrm{trdeg}(c/A') = \mathrm{trdeg}(c/B')$, with $A'$, $B'$ the respective subfields of $M$ generated by $A$, $AB$.

This independence notion is symmetric (see for instance [Lan05], Section 8.3): we have $C \downarrow^{\mathbf{alg}}_A B$ if and only if $B \downarrow^{\mathbf{alg}}_A C$.

It is well-known that forking refines algebraic independence (over the definable closure of the base) in (expansions of) infinite fields. It is not done very explicitly in the literature, and we wrote an explicit proof during the preparation of this thesis, with a rather technical argument. We later found out that Johnson gave in ([Joh20], Section 2) an elegant proof of this fact in the particular case of some folklore topological fields. It turns out that the argument of Johnson can be extended to any expansion of an infinite field, and this is essentially what we do in this subsection. We try to use as few stability-theoretic facts as possible.

**Assumptions 1.2.24.** Suppose $M$ is a multi-sorted first-order structure, and $k$ is a sort of $M$ such that $M$ expands an infinite field with domain $k$.

In one statement of this subsection, we make further assumptions:

**Assumptions 1.2.25.** On top of assumptions 1.2.24, suppose that the only home sort of $M$ is $k$, and suppose that $M$ is solely defined on the language of rings (i.e. $M$ is a pure field). Moreover, suppose $M$ is a completion of the theory ACF of algebraically closed fields.

**Definition 1.2.26.** With assumptions 1.2.24, we say that $M$ is *algebraically bounded* whenever $k(\mathrm{acl}(A)) = k(A)^{\mathbf{alg}} \cap k(M)$ for any substructure $A \leqslant M$.

Let us start by showing that non-forking coincides with algebraic independence in ACF. For more complicated theories of expansions of fields where non-forking coincides with algebraic independence (over the algebraic closure of the base), such as ACFA or pseudo-finite fields, the standard way to compute forking is using the Kim-Pillay theorem ([KP98]). One does not need to work that much to do it in ACF. The only non-trivial fact from geometric stability theory that we need is the following:

**Fact 1.2.27.** *In a stable theory, the converse of proposition 1.1.18 holds, i.e. we have* $CD \downarrow^{\mathbf{f}}_A B$ *if and only if* $C \downarrow^{\mathbf{f}}_A B$ *and* $D \downarrow^{\mathbf{f}}_{AC} B$*, for any small parameter sets* $A$, $B$, $C$, $D$.

We state other facts which are much more trivial:

**Fact 1.2.28.**

- *In a strongly minimal theory, if* $c$ *is a singleton, and* $A$, $B$ *are small parameter sets, then we have* $c \mathop{\downarrow}^{\mathbf{f}}_{A} B$ *if and only if either* $c \in \operatorname{acl}(A)$, *or* $c \notin \operatorname{acl}(AB)$.
- *Strongly minimal theories are stable.*
- *Any completion of ACF is strongly minimal, and ACF eliminates quantifiers in the language of rings.*
- *Any completion of ACF is algebraically bounded.*

**Corollary 1.2.29.** *With assumptions 1.2.25, we have* $\downarrow^{\mathbf{f}} = \downarrow^{\mathbf{alg}}$.

*Proof.* Let $c = c_1 \dots c_n$ be a finite tuple, $A$, $B$ be parameter sets, and $A'$, $B'$ be the respective subfields generated by $A$, $AB$. By corollary 1.1.25 and proposition 1.1.24, we have $c \mathop{\downarrow}^{\mathbf{f}}_{A} B$ if and only if $c \mathop{\downarrow}^{\mathbf{f}}_{A'} B'$. By fact 1.2.27, this is equivalent to having $c_i \mathop{\downarrow}^{\mathbf{f}}_{A'(c_{<i})} B'$ for all $i$. By fact 1.2.28, this is equivalent to having $c_i \notin \operatorname{acl}(B'(c_{<i})) \setminus \operatorname{acl}(A'(c_{<i}))$ for all $i$, which, by algebraic boundedness, is equivalent to having $c_i \notin (B'(c_{<i}))^{\mathbf{alg}} \setminus (A'(c_{<i}))^{\mathbf{alg}}$, which is equivalent to $\operatorname{trdeg}(c_i/A'(c_{<i})) = \operatorname{trdeg}(c_i/B'(c_{<i}))$. By additivity of the transcendence degree, this is equivalent to having $c \mathop{\downarrow}^{\mathbf{alg}}_{A} B$. We conclude by left finite character of forking that $\downarrow^{\mathbf{f}} = \downarrow^{\mathbf{alg}}$. □

In order to build an indiscernible sequence witnessing dividing, we need to admit a fact from algebra about the field of definition of a variety:

**Fact 1.2.30.** *With assumptions 1.2.24, let* $I$ *be a proper radical ideal of* $k(M)[x]$ *(*$x$ *is a finite tuple of indeterminates). Then the class of subfields* $K$ *of* $k(M)$ *for which* $I$ *is generated by* $I \cap K[x]$ *as a* $k(M)$*-vector space has a least element with respect to inclusion,* $DF(I)$ *(the field of definition of* $I$*). Moreover, for any field-automorphism* $\sigma$ *of* $k(M)$*, we have* $\sigma(I) = I$ *if and only if* $\sigma \in \operatorname{Aut}(k(M)/DF(I))$.

**Corollary 1.2.31.** *With assumptions 1.2.24, let* $V$ *be a* $k(M)$*-Zariski-closed set, i.e. a definable set which is the zero-locus of finitely many polynomial equations over* $k(M)$*. Then* $V$ *is* $k(\operatorname{dcl}^{eq}(\ulcorner V \urcorner))$*-Zariski-closed.*

The standard notation $\ulcorner V \urcorner$ refers to the canonical parameter corresponding to $V$ in $M^{eq}$. Note that we are working in an arbitrary expansion of a field here.

*Proof.* Every automorphism of $M^{eq}$ fixing $\ulcorner V \urcorner$ must globally fix its ideal $I$ in $k(M)$, thus its restriction to $k(M)$ is a field automorphism of $k(M)$ which pointwise-fixes $DF(I)$. It follows that $DF(I) \subseteq \mathrm{dcl}(\ulcorner V \urcorner)$, concluding the proof. $\square$

We recall that with assumptions 1.2.24, the definable closure is relative to the expanded structure $M$, and it may be much larger than what one could define using only the language of rings.

**Lemma 1.2.32.** *With assumptions 1.2.24, let $A \subseteq k(M)$ be such that $A = k(\mathrm{dcl}^{eq}(A))$, let $b$ be some tuple of $k$, and let $E$ be a small parameter subset of $k(M)$. Then there exists $b' \equiv_A b$ such that $b' \mathop{\downarrow}^{\mathbf{alg}}_A E$.*

Extension for forking in algebraically closed fields is not enough to prove this lemma, as the elementary equivalence of the statement is not in the language of rings, but in an expansion of some field.

*Proof.* We can freely assume that $E$ is a field containing $A$. Let $d$ be some maximal $A$-algebraically independent subfamily of $b$. Let us show by compactness that there exists an $E$-algebraically independent tuple $d'$ such that $d' \equiv_A d$. Let $X$ be an $A$-definable set containing $d$. Let $I$ be the set of all polynomials $P$ in $k(M)[x]$ such that $M \vDash x \in X \implies P(x) = 0$. Note that, if $(P_i)_i$ is a finite family of polynomials, then the formula $\bigvee_i P_i(x) = 0$ is equivalent to $\left(\prod_i P_i(x)\right) = 0$. Then $I$ is a radical proper ideal, the ideal of the $k(M)$-Zariski closure $V$ of $X$. As $X$ is $A$-definable, $X$ (and hence $I$, $V$) is stabilized by $\mathrm{Aut}(M/A)$. It follows that $\ulcorner V \urcorner \in \mathrm{dcl}^{eq}(A)$, and therefore, by corollary 1.2.31, $I$ is generated as a $k(M)$-vector space by $I \cap A[x]$. However, as $d \in X$ is $A$-algebraically independent, we have $I \cap A[x] = \{0\}$, thus $I = \{0\}$. It follows by compactness that $X$ has an $E$-algebraically independent realization. Again, by compactness, there does exist $d' \equiv_A d$ such that $d'$ is $E$-algebraically independent.

Choose $\sigma \in \mathrm{Aut}(M/A)$ such that $\sigma(d) = d'$, and set $b' = \sigma(b)$. We have: $\mathrm{trdeg}(b/A(d)) = 0$ and $b \equiv_A b'$, therefore we have: $\mathrm{trdeg}(b'/A(d')) = 0$ thus

$\mathrm{trdeg}(b'/E(d')) = 0$. It follows that

$$\mathrm{trdeg}(b'/E) = \mathrm{trdeg}(d'/E) = \mathrm{trdeg}(d'/A) = \mathrm{trdeg}(b'/A)$$

the second equality follows from the fact that $d'$ is $E$-algebraically independent. This concludes the proof. □

Once we have our indiscernible sequence, our argument to show dividing uses the following lemma:

**Lemma 1.2.33.** *With assumptions 1.2.24, let $A \subseteq k(M)$ be such that $A = k(\mathrm{dcl}^{eq}(A))$, and let $c = c_1 \dots c_n$ be a finite tuple of $k$. Suppose $(B_i)_{i<\omega}$ is an $A$-indiscernible sequence (in the expanded structure $M$) of subfields of $k(M)$, such that: $B_i \downarrow^{\mathbf{alg}}_A A((B_j)_{j<i})$ for all $i$. Then there must exist some $i$ such that $c \downarrow^{\mathbf{alg}}_A B_i$.*

*Proof.* We proceed by induction on the size of $c$. The statement is trivial if $|c| = 0$.

Suppose by contradiction that we have $c \not\downarrow^{\mathbf{alg}}_A B_i$ for all $i$. Then the following map is well-defined:

$$f: i \longmapsto \min\{j \leqslant |c| \ : \ 1 = \mathrm{trdeg}(c_j/A(c_{<j})) \neq \mathrm{trdeg}(c_j/B_i(c_{<j})) = 0\}$$

By the pigeonhole principle, there exists $j \leqslant |c|$ such that the fiber of $f$ at $j$ is infinite. Let $(D_i)_{i<\omega}$ be the strictly increasing re-enumeration of $(B_i)_{f(i)=j}$. Let $I = (D_{2i}D_{2i+1})_{i<\omega}$, which is $A$-indiscernible. Let us show that $c_{<j} \not\downarrow^{\mathbf{alg}}_A D_{2i}D_{2i+1}$ for all $i$, contradicting the induction hypothesis.

Choose $i$. Let $F_0$ be the field-theoretic algebraic closure of the field generated by $AD_{2i}c_{<j}$, $F_1$ be that generated by $AD_{2i+1}c_{<j}$, and $A' = (A(c_{<j}))^{\mathbf{alg}}$. By definition of $j$, and by additivity of the transcendence degree, we have $c_{<j} \downarrow^{\mathbf{alg}}_A D_{2i}$, $c_{<j} \downarrow^{\mathbf{alg}}_A D_{2i+1}$ and $c_j \in F_0 \setminus A'$, $c_j \in F_1 \setminus A'$. Let $b_0$ (resp. $b_1$) be some finite tuple of $D_{2i}$ (resp. $D_{2i+1}$) such that $c_j \in (A'(b_0))^{\mathbf{alg}}$, and $c_j \in (A'(b_1))^{\mathbf{alg}}$. By hypothesis on the $B_i$, we have $b_0 \downarrow^{\mathbf{alg}}_A b_1$, thus $\mathrm{trdeg}(b_0b_1/A) = \mathrm{trdeg}(b_0/A) + \mathrm{trdeg}(b_1/A)$. Let $F = (A(b_0)(c_<j))^{\mathbf{alg}} \cap (A(b_1)(c_<j))^{\mathbf{alg}}$. As $c_j \in F \setminus A'$, and $A'$ is algebraically closed, we have $\mathrm{trdeg}(F/A') > 0$. By symmetry of algebraic independence, we have $D_{2i} \downarrow^{\mathbf{alg}}_A c_{<j}$, thus $\mathrm{trdeg}(b_0/A') =$

$\operatorname{trdeg}(b_0/A)$, and likewise for $b_1$. Now, the computation:

$$\begin{aligned} \operatorname{trdeg}(b_0b_1/A) &= \operatorname{trdeg}(b_0/A) + \operatorname{trdeg}(b_1/A) \\ &= \operatorname{trdeg}(b_0/A') + \operatorname{trdeg}(b_1/A') \\ &= \operatorname{trdeg}(b_0/F) + \operatorname{trdeg}(F/A') + \operatorname{trdeg}(b_1/F) + \operatorname{trdeg}(F/A') \\ &\geqslant \operatorname{trdeg}(b_0b_1/F) + 2\operatorname{trdeg}(F/A') \\ &= \operatorname{trdeg}(b_0b_1/A') + \operatorname{trdeg}(F/A') \text{ (as } F^{\mathbf{alg}} \leqslant (A'(b_0b_1))^{\mathbf{alg}}) \\ &> \operatorname{trdeg}(b_0b_1/A') \end{aligned}$$

yields that $D_{2i}D_{2i+1}\not\mathop{\downarrow}^{\mathbf{alg}}_{A} c_{<j}$. We conclude by symmetry. □

The above proof can be seen as an algebraic proof of Kim's lemma in ACF (though we work in the general context here), which Johnson uses in their paper [Joh20]. We now have enough tools to prove the result that we aimed for:

**Proposition 1.2.34.** *With assumptions 1.2.24, let $A \subseteq k(M)$ be such that $A = k(\operatorname{dcl}^{eq}(A))$, and let $B, C$ be some small parameter subsets of $k(M)$. If we have $C \mathop{\downarrow}^{\mathbf{d}}_{A} B$, then we have $C \mathop{\downarrow}^{\mathbf{alg}}_{A} B$.*

*Proof.* By contraposition, suppose $C\not\mathop{\downarrow}^{\mathbf{alg}}_{A}B$. By lemma 1.2.32, we can build inductively a sequence $(B_i)_{i<\omega}$ such that $B_0 = B$, $B_{i+1} \equiv_A B_i$, and:

$$B_i \mathop{\downarrow}^{\mathbf{alg}}_{A} \left( \bigcup_{j<i} B_j \right)$$

We may assume by a standard Ramsey argument that $I = (B_i)_i$ is $A$-indiscernible.

Let $\iota$ be an algebraic embedding of fields $k(M) \longrightarrow k(M)^{\mathbf{alg}}$. We clearly have $\iota(C)\not\mathop{\downarrow}^{\mathbf{alg}}_{\iota(A)}\iota(B)$. By quantifier elimination in ACF, there exists a formula without quantifiers $\phi(x, \iota(AB))$ from the language of rings, and a finite tuple $c$ from $C$, such that $k(M)^{\mathbf{alg}} \vDash \phi(\iota(c), \iota(AB))$, and $\phi(x, \iota(AB))$ forks over $\iota(A)$. Naturally, we also have $M \vDash \phi(c, AB)$.

It suffices to show that the following partial type:

$$p = \{\phi(x, AB_i) | i < \omega\}$$

is inconsistent, for $\phi(x, AB)$ would divide over $A$ by fact 1.1.3, and we could conclude that $C\not\mathop{\downarrow}^{\mathbf{d}}_{A}B$.

Suppose by contradiction that we have $c'$ a realization of $p$. Let $i < \omega$. By quantifier elimination in ACF, as $B_i \equiv_A B$, we have $\iota(B_i) \equiv_{\iota(A)} \iota(B)$, in particular $\phi(x, \iota(AB_i))$ forks over $\iota(A)$ by remark 1.1.12. As $c'$ realizes $p$, and $\phi$ is quantifier-free, we have $M^{\mathbf{alg}} \vDash \phi(\iota(c'), \iota(AB_i))$, thus $\iota(c') \not\perp^{\mathbf{f}}_{\iota(A} \iota(B_i)$, therefore $\iota(c') \not\perp^{\mathbf{alg}}_{\iota(A)} \iota(B_i)$ by corollary 1.2.29. Naturally, this implies $c' \not\perp^{\mathbf{alg}}_{A} B_i$. This holds for all $i < \omega$, which contradicts lemma 1.2.33, concluding the proof. □

## 1.3 Keisler measures

Keisler measures were introduced by Keisler in [Kei87], as a way to generalize the notion of a type, in order to obtain new results about forking (in terms of Morley rank in a stable theory) in dependent theories. Using Keisler measures to study forking is relevant for our work, because there is a well-known interaction between Keisler measures and forking in pseudo-finite fields (hence in the ultraproducts of the $p$-adic fields if we look at their residue fields). Indeed, ([CvdDM92], Section 4) gives a notion of measure on definable sets in pseudo-finite fields, which corresponds to a Keisler measure which interacts well with algebraic independence, i.e. non-forking in pseudo-finite fields. We actually use this pseudo-finite measure to build in section 7.3 a Keisler measure witnessing non-forking in those valued fields.

While they were the initial motivation behind the study of Keisler measures, the relations between Keisler measure and forking were not investigated much further, though Keisler measures themselves gained recent increased interest in the literature for other settings (continuous logic, amenable groups...). Following the spirit of [HP07], we see Keisler measures as a way to generalize the notion of an extension of a type, and thus build new independence notions, which are usually stronger, easier to define than non-forking independence, and more tame than their counterpart for types.

In this section, we build a general formalism about those various notions, we give a survey about what was done in the literature regarding independence (very little that is, the subject is wide open), and we prove minor results. As in [Kei87] and [HP07], many statements about Keisler measures hold in the specific setting of dependent theories.

### 1.3.1 Definitions

**Definition 1.3.1.** Fix $M$ a large, sufficiently saturated first-order structure. A *Keisler measure* on the tuple of variables $x$ is a finitely additive probability measure on the Boolean algebra $\mathcal{B}^x(M)$ of $M$-definable sets on the variables $x$ (every measure we consider is global). Such a measure can be extended to (and may be identified with) a unique $\sigma$-additive regular Borel probability measure with universe the space $S^x(M)$ of global types. This construction is done in details in ([Sim15], Section 7.1).

We denote by $\mathrm{KM}^x$ the space of all Keisler measures on $x$, which we see as a closed subspace of the (compact) space of maps from $\mathcal{B}^x(M)$ to $[0,1]$, with the product topology induced by that on $[0,1]$. The space $S^x(M)$ of actual global types corresponds naturally to a closed subspace of $\mathrm{KM}^x$, by identifying a global type with its Dirac measure.

We will often need to prove that some measures give probability one/non-zero to specific closed sets, so we should point out the following:

**Fact 1.3.2.** *The measure of a closed set $F$ coincides with the infimum of the measures of the clopen sets containing $F$.*

*In particular, some partial type (seen as a closed set) has measure one (resp. non-zero) if and only if all its formulas have measure one (resp. the infimum of the measures is non-zero).*

Our informal idea for what should be an abstract independence relation $c \mathop{\downarrow}_A B$ is: $\mathrm{tp}(c/AB)$ admits a "tame" global extension. From this idea we define various independence relations involving Keisler measures. The classes of measures that we consider tame, as analogues of non-forking global types, will mostly be spaces of measures that are invariant under the action of some subgroup of $\mathrm{Aut}(M/A)$ (induced by the action on the set of $M$-definable sets defined as $\sigma \cdot X = \sigma(X)$). Instead of considering measures on $S^x(M)$, we may also want to consider measures on some closed subspace:

**Proposition 1.3.3.** *Let $G$ be a subgroup of $\mathrm{Aut}(M)$, and $F$ a $G$-invariant closed subspace of $S^x(M)$. Then there is a natural 1-1 correspondence between the $G$-invariant Keisler measures giving $F$ probability one, and the $G$-invariant finitely additive probability measures on the Boolean algebra associated to $F$.*

*Proof.* Let $\mathcal{B}_F$ be the Boolean algebra associated to $F$, and let $\rho: X \longmapsto X \cap F$ be the natural morphism of Boolean algebras $\mathcal{B}(M) \longrightarrow \mathcal{B}_F$. Then $\rho$ is a surjective morphism of kernel $I$, the ideal of definable sets disjoint from $F$. If $\mu$ is a Keisler measure giving $F$ probability one, then $\mu(Y) = 0$ for every $Y \in I$, thus $\mu(X + Y) = \mu(X)$ for every $X \in \mathcal{B}(M), Y \in I$. As a result, $\mu$ factors through a unique map $\overline{\mu}: \mathcal{B}_F \longrightarrow [0,1]$. It is easy to see that $\overline{\mu}$ is a finitely additive probability measure.

The correspondence $\mu \longmapsto \overline{\mu}$ is a bijection, for its inverse is clearly $\nu \longmapsto \nu \circ \rho$. □

Let us define how we generalize the notion of extension of a type:

**Definition 1.3.4.** Let $\mu \in \mathrm{KM}^x$, and let $p$ be a partial type over $M$. We define the *support* $\mathrm{supp}(\mu)$ of $\mu$ as the least closed subspace of $S^x(M)$ having measure 1. Equivalently, by compactness:

$$p \notin \mathrm{supp}(\mu) \Longleftrightarrow p(x) \vDash x \in X, \mu(X) = 0 \text{ for some definable set } X.$$

We say that $\mu$ is *supported* on some partial type $\pi$ if $\pi$ does not imply a definable set having measure zero.

We say that $\mu$ is *concentrated* on $\pi$ (which we identify with a closed subspace of $S^x(M)$) if $\mu(\pi) = 1$.

The motivation behind those two definitions is that $\mu$ may be supported on $\pi$, while $\mu(\pi) = 0$, an odd behavior which cannot occur if $\mu$ is concentrated on $\pi$.

For a concrete example, consider in RCF the Keisler measure $\mu$ induced by the Lebesgue measure on $[0,1]$: the measure of an interval is the Lebesgue measure of its standard part (we fix a substructure isomorphic to $\mathbb{R}$), and the measure of an arbitrary definable set is deduced from the measures of intervals by o-minimality. If $c$ is a positive infinitesimal element, then $\mu$ is supported on $\mathrm{tp}(c/\mathbb{R})$, however $\mu(\mathrm{tp}(c/\mathbb{R})) = 0$.

Now that we generalized the notion of extension, the new independence notions that we get can be defined as follows:

**Definition 1.3.5.** Let $\mathfrak{T}$ (which stands for Tame) be a map from the class of small subsets of $M$ to the powerset of the disjoint union: $\bigcup_{N<\omega} \mathrm{KM}^N$. (The idea is that $\mathfrak{T}(A)$ is the set of measures that are tame with respect to the small set $A$. It would make sense to assume additionally that $\mathfrak{T}(A) \cap \mathrm{KM}^x$ is closed in $\mathrm{KM}^x$ for all $A$.)

We define the ternary relation $\mathop{\smile\!\!\!\!\!\!\mid}^{\mathfrak{I}\text{-KM}}$ as $C \mathop{\smile\!\!\!\!\!\!\mid}_A^{\mathfrak{I}\text{-KM}} B$ if and only if for every finite tuple $c$ from $C$, there exists $\mu \in \mathrm{KM}^{|c|}$ such that $\mu$ is supported on $\mathrm{tp}(c/AB)$, and $\mu \in \mathfrak{I}(A)$.

We define $\mathop{\smile\!\!\!\!\!\!\mid}^{\mathfrak{I}}$ as above, with the additional condition that $\mu$ is (the Dirac measure on) a global type.

Note that $\mathop{\smile\!\!\!\!\!\!\mid}^{\mathfrak{I}}$ refines $\mathop{\smile\!\!\!\!\!\!\mid}^{\mathfrak{I}\text{-KM}}$: $C \mathop{\smile\!\!\!\!\!\!\mid}_A^{\mathfrak{I}} B \Longrightarrow C \mathop{\smile\!\!\!\!\!\!\mid}_A^{\mathfrak{I}\text{-KM}} B$.

**Definition 1.3.6.** We will focus on the following possibilities for $\mathfrak{I}$ (yielding in each case two independence relations as in definition 1.3.5):

**(invariance) inv** : $A \longmapsto$ the set of $\mathrm{Aut}(M/A)$-invariant Keisler measures.

**(bounded orbit) bo** : $A \longmapsto$ the set of Keisler measures with small orbit under the action of $\mathrm{Aut}(M/A)$.

**(Shelah-invariance) Sh** : $A \longmapsto$ the set of $\mathrm{Aut}(M^{eq}/\mathrm{acl}^{eq}(A))$-invariant Keisler measures.

We will also consider, in this section, independence notions induced by the following maps to establish relations with [HP07] and other works:

**(Kim-Pillay-invariance) KP** : $A \longmapsto$ the set of $\mathrm{Aut}(M/bdd(A))$-invariant Keisler measures ($bdd(A)$ is defined in page 5 of [HP07]).

**(Lascar-invariance) L** : $A \longmapsto$ the set of $\mathrm{Aut}f_A(M)$-invariant Keisler measures (see [HKP20], Definition 1.1,).

**(coheir-independence) fs** : $A \longmapsto$ the set of Keisler measures which are *finitely satisfiable* over $A$, i.e. the set of measures $\mu$ such that if $X$ has no point in $A$, then $\mu(X) = 0$.

For instance, we have $c \mathop{\smile\!\!\!\!\!\!\mid}_A^{\mathbf{inv}} B$ if and only if $\mathrm{tp}(c/AB)$ admits a global $\mathrm{Aut}(M/A)$-invariant extension, while $c \mathop{\smile\!\!\!\!\!\!\mid}_A^{\mathbf{inv}\text{-KM}} B$ means that there exists some $\mathrm{Aut}(M/A)$-invariant Keisler measure which is supported on $\mathrm{tp}(c/AB)$.

As $\mathrm{Aut}(M/A)$ is the largest group considered in these various definitions, $\mathop{\smile\!\!\!\!\!\!\mid}^{\mathbf{inv}}$ is the strongest independence notion involving groups of automorphisms. It is also clear that $\mathop{\smile\!\!\!\!\!\!\mid}^{\mathbf{inv}}$ is stronger than $\mathop{\smile\!\!\!\!\!\!\mid}^{\mathbf{bo}}$, and likewise for their KM variants. Let us show that $\mathop{\smile\!\!\!\!\!\!\mid}^{\mathbf{fs}}$ is the strongest independence notion that we defined:

**Proposition 1.3.7.** *We always have* $\mathbf{fs}(A) \subseteq \mathbf{inv}(A)$*, i.e. every Keisler measure which is finitely satisfiable over* $A$ *must be* $\mathrm{Aut}(M/A)$*-invariant.*

*Proof.* Let $\mu$ be a Keisler measure which is finitely satisfiable over $A$. Let $X$ be some definable set, and $\sigma \in \mathrm{Aut}(M/A)$. Then $X \smallsetminus \sigma(X)$ and $\sigma(X) \smallsetminus X$ do not have points in $A$, therefore $\mu(x \smallsetminus \sigma(X)) = 0 = \mu(\sigma(X) \smallsetminus X)$, which implies that $\mu(X) = \mu(X \cap \sigma(X)) = \mu(\sigma(X))$. It follows that $\mu$ is $\mathrm{Aut}(M/A)$-invariant. □

**Corollary 1.3.8.** *We have* $\downarrow^{\mathbf{fs}} \subseteq \downarrow^{\mathbf{inv}}$*, and* $\downarrow^{\mathbf{fs}-\mathrm{KM}} \subseteq \downarrow^{\mathbf{inv}-\mathrm{KM}}$*.*

### 1.3.2 Forking for measures

Non-dividing independence does not involve global types, but non-forking independence can be recovered using definition 1.3.5. There is a very natural way to define the notion of a measure which does not fork over $A$. Hrushovski and Pillay consider those measures in [HP07], with applications to amenable groups. However, the notion of independence induced by non-forking measures is, as we show in this subsection, equivalent to the usual notion of forking. This is also true of finitely satisfiable measures.

**Lemma 1.3.9.** *Let* $\pi(x)$ *be some partial type over* $M$*, and* $\mu \in \mathrm{KM}^x$*. Then* $\mu$ *is supported on* $\pi$ *if and only if there exists some complete global type* $p \in \mathrm{supp}(\mu)$ *such that* $p(x) \vDash \pi(x)$*.*

*Proof.* The right-to-left direction is trivial.

Conversely, suppose $\mu$ is supported on $\pi$. It suffices to show that the partial type $\pi(x) \cup \{x \notin X | \mu(X) = 0\}$ is consistent. If not, then by compactness, there must exist $M$-definable sets $X, Y_1, \dots Y_n$ such that $\pi(x) \vDash x \in X$, $\mu(Y_i) = 0$, and $X \subseteq \bigcup_i Y_i$. This implies that $\mu(X) = 0$, contradicting the fact that $\mu$ is supported on $\pi$, hence proved. □

**Proposition 1.3.10.** *Let* $\mu \in \mathrm{KM}^x$*. Then the following are equivalent:*

1. *No global type in* $\mathrm{supp}(\mu)$ *forks over* $A$*.*

2. *For all* $M$*-definable sets* $X$*, if* $X$ *divides over* $A$*, then* $\mu(X) = 0$*.*

*Proof.* Let us show $2 \implies 1$ by contraposition. Suppose we have $p \in \mathrm{supp}(\mu)$ such that $p$ forks over $A$. By remark 1.1.13, $p$ divides over $A$, thus there must exist $X$ such that $p(x) \models x \in X$, and $X$ divides over $A$. As $p \in \text{supp}(\mu)$, we have $\mu(X) > 0$, thus 2 fails.

Conversely, suppose $\mu(X) > 0$ for some $M$-definable set $X$ which divides over $A$. By lemma 1.3.9, there exists $p \in \text{supp}(\mu)$ such that $p(x) \models x \in X$. As $X$ divides over $A$, $p$ forks over $A$, and 1 fails. $\square$

**Definition 1.3.11.** If the conditions of proposition 1.3.10 hold, then we say that $\mu$ *does not fork over* $A$.

This gives the independence relations induced by $\mathbf{f}: A \longmapsto$ the set of Keisler measures which do not fork over $A$.

This definition is consistent with the classical definition of forking.

Note also that the definition of non-forking measures in [HP07] is condition 2 of proposition 1.3.10.

**Proposition 1.3.12.** *We have* $\downarrow^{\mathbf{f}} = \downarrow^{\mathbf{f}-\mathrm{KM}}$.

*Proof.* The left-to-right inclusion trivially holds for any notion of tameness.

Suppose $c \downarrow^{\mathbf{f}-\mathrm{KM}}_A B$ for some finite tuple $c$. Let $\mu \in \mathrm{KM}^{|c|}$ be supported on $\mathrm{tp}(c/AB)$, such that $\mu$ does not fork over $A$. By lemma 1.3.9, there exists $p$ a global extension of $\mathrm{tp}(c/AB)$ which is in $\text{supp}(\mu)$. By definition, $p$ does not fork over $A$, thus $c \downarrow^{\mathbf{f}}_A B$. $\square$

Similar arguments show the analogue of proposition 1.3.12 for finite satisfiability:

**Proposition 1.3.13.** *Let $\mu$ be a Keisler measure. Then $\mu$ is finitely satisfiable over some small set $A$ if and only if this is the case of every type in its support.*

*Proof.* If $\mu$ is not finitely satisfiable, then we have $\mu(X) > 0$ for some definable set $X$ having no point in $A$. By lemma 1.3.9, there exists $p \in \text{supp}(\mu)$ such that $p(x) \models x \in X$, thus $p$ is not finitely satisfiable over $A$.

The converse of the statement is trivial. $\square$

**Corollary 1.3.14.** *We have* $\downarrow^{\mathbf{fs}} = \downarrow^{\mathbf{fs}-\mathrm{KM}}$.

### 1.3.3 Invariance under various groups

Among the independence relations that we defined, we showed that adding Keisler measures to the setting did not add anything valuable to the study of $\mathop{\smile\!\!\!\!\!\!|}^{\mathbf{f}}$ and $\mathop{\smile\!\!\!\!\!\!|}^{\mathbf{fs}}$. We defined $\mathop{\smile\!\!\!\!\!\!|}^{\mathbf{bo}}$ because we need it in our work on regular ordered Abelian groups, but we do not know if its KM-variant is of interest, let alone refines non-forking. The relations that are left are $\mathop{\smile\!\!\!\!\!\!|}^{\mathbf{inv}}, \mathop{\smile\!\!\!\!\!\!|}^{\mathbf{Sh}}, \mathop{\smile\!\!\!\!\!\!|}^{\mathbf{KP}}, \mathop{\smile\!\!\!\!\!\!|}^{\mathbf{L}}$, which all involve measures which are invariant under the action of some subgroup of $\mathrm{Aut}(M/A)$. There is actually one direction of proposition 1.3.10 that still holds in this setting:

**Proposition 1.3.15.** *Let $G$ be some subgroup of $\mathrm{Aut}(M)$, and $\mu \in \mathrm{KM}^x$. If every type in $\mathrm{supp}(\mu)$ is $G$-invariant, then $\mu$ is $G$-invariant.*

*Proof.* Recall that in a Boolean algebra, the sum (a.k.a. the XOR or the symmetrical difference) of $X$ and $Y$ is defined as $X + Y = (X \vee Y) \smallsetminus (X \wedge Y) = (X \smallsetminus Y) \vee (Y \smallsetminus X)$ (with $X \smallsetminus Y = X \wedge \overline{Y} = X \wedge (1 + Y)$).

Let $X$ be some $M$-definable set, and $\sigma \in G$. Suppose by contradiction that $\mu(X + \sigma(X)) > 0$. Then, by lemma 1.3.9, there must exist $p \in \mathrm{supp}(\mu)$ such that $p(x) \vDash x \in (X \smallsetminus \sigma(X)) \cup (\sigma(X) \smallsetminus X)$. As $p$ is complete, we have either $p(x) \vDash x \in X \wedge x \notin \sigma(X)$ or $p(x) \vDash x \in \sigma(X) \wedge x \notin X$. Either way, $p$ is not $G$-invariant, which contradicts our hypothesis, thus we proved by contradiction that $\mu(X + \sigma(X)) = 0$.

In particular, we have $\mu(X \smallsetminus \sigma(X)) = 0 = \mu(\sigma(X) \smallsetminus X)$, thus $\mu(X) = \mu(X \cap \sigma(X)) = \mu(\sigma(X))$, and we conclude that $\mu$ is $G$-invariant. ☐

**Corollary 1.3.16.** *If every type which does not fork over $A$ is $G$-invariant, then every measure which does not fork over $A$ is $G$-invariant.*

This corollary greatly simplifies the proofs of ([HP07], Propositions 4.3 and 4.5) relating forking, Lascar-invariance and Kim-Pillay-invariance for measures in NIP theories.

*Proof.* Let $\mu$ be a Keisler measure which does not fork over $A$. Then all the global types in $\mathrm{supp}(\mu)$ do not fork over $A$ by definition. By hypothesis, they are all $G$-invariant, and we conclude by proposition 1.3.15. ☐

Let us prove that the converse of proposition 1.3.15 fails (otherwise, the study of independence relations involving Keisler measures would be totally pointless):

*Example* 1.3.17. The following example (which is a stable theory) illustrates several things:

- The converse of proposition 1.3.15 fails.
- This is an easy example where we have the inequalities $\mathop{\smile\!\!\!\!\!\!\;|}^{\mathbf{Sh}} \neq \mathop{\smile\!\!\!\!\!\!\;|}^{\mathbf{inv}} \neq \mathop{\smile\!\!\!\!\!\!\;|}^{\mathbf{inv}-\mathrm{KM}}$.
- The relation $\mathop{\smile\!\!\!\!\!\!\;|}^{\mathbf{inv}}$ is quite wild in general, for not every set is an extension base for this relation in this example. However, $\mathop{\smile\!\!\!\!\!\!\;|}^{\mathbf{inv}-\mathrm{KM}}$ and $\mathop{\smile\!\!\!\!\!\!\;|}^{\mathbf{Sh}}$ behave nicely here, in fact they coincide with non-forking.

Work in the language $\{E\}$, in the theory stating that $E$ is an equivalence relation with two classes, both infinite. In this theory, we have $\mathrm{acl}(\varnothing) = \varnothing$ and there is only one unary complete type over $\varnothing$. Let $\alpha, \beta \in \mathrm{acl}^{eq}(\varnothing)$ be the two $E$-classes. There are two unary global $\mathrm{Aut}(M^{eq}/\mathrm{acl}^{eq}(\varnothing))$-invariant types:

$$p_\alpha(x): \ \{x \in \alpha\} \cup \{x \neq d | d \in M\}$$
$$p_\beta(x): \ \{x \in \beta\} \cup \{x \neq d | d \in M\}$$

The other global unary types divide over $\varnothing$ by lemma 1.2.2, thus they are not $\mathrm{Aut}(M^{eq}/\mathrm{acl}^{eq}(\varnothing))$-invariant (they fork over $\varnothing$, so we use item 2 of fact 1.1.19). As $\alpha$ and $\beta$ are $\varnothing$-conjugates, neither $p_\alpha$ nor $p_\beta$ is $\mathrm{Aut}(M/\varnothing)$-invariant. In particular, there does not exist any unary global $\mathrm{Aut}(M/\varnothing)$-invariant type, thus $c \mathop{\not\smile\!\!\!\!\!\!\;|}^{\mathbf{inv}}_{\varnothing} \varnothing$ for any singleton (and hence any tuple) $c$. However, the average $\mu = \frac{1}{2} \cdot p_\alpha + \frac{1}{2} \cdot p_\beta$ is $\mathrm{Aut}(M/\varnothing)$-invariant. As there is only one unary type over $\varnothing$, it must belong to $\mathrm{supp}(\mu)$, thus we have $c \mathop{\smile\!\!\!\!\!\!\;|}^{\mathbf{inv}-\mathrm{KM}}_{\varnothing} \varnothing$ for any singleton $c$.

Moreover, there are only two unary types over $\mathrm{acl}^{eq}(\varnothing)$: the one isolated by $x \in \alpha$, which is consistent with $p_\alpha$, and the one isolated by $x \in \beta$, which is consistent with $p_\beta$. As a result, we have $c \mathop{\smile\!\!\!\!\!\!\;|}^{\mathbf{Sh}}_{\varnothing} \varnothing$ for any singleton $c$.

To conclude, note that $\mu$ is a counterexample to the converse of the previous proposition, for it is $\mathrm{Aut}(M/\varnothing)$-invariant, but no type in its support (in fact, no complete global type at all) is $\mathrm{Aut}(M/\varnothing)$-invariant.

This example suggests that the variants for measures of our independence notions behave better than those for types. The useful tool that we have from corollary 1.1.22 for computing forking turns out to work for those independence notions as well:

**Proposition 1.3.18.** *Let $G$ be some subgroup of* $\mathrm{Aut}(M/A)$*. Suppose we have* $\mathrm{acl}(Ac) \subseteq \mathrm{acl}(Ad)$*, and there exists a $G$-invariant Keisler measure which is supported on* $\mathrm{tp}(d/AB)$*. Then there exists some $G$-invariant Keisler measure which is supported on* $\mathrm{tp}(c/AB)$*.*

*Proof.* Let $\mu$ be a $G$-invariant Keisler measure, supported on $\mathrm{tp}(d/AB)$. We have to build $\nu$ a $G$-invariant Keisler measure which is supported on $\mathrm{tp}(c/AB)$.

Let $\phi(x,y)$ be some formula with parameters in $A$ such that the definable set $\phi(x,d)$ is finite (say, $|\phi(x,d)| = N < \omega$), and isolates the type of $c$ over $Ad$. For each formula $\chi(x,y)$, and $i < \omega$, let $\exists_{=i}x\chi(x,y)$ be the formula stating that there exists exactly $i$ pairwise-distinct tuples $x_1 \dots x_i$ such that $\chi(x_j,y)$ holds for all $j$. By replacing $\mu$ by $\psi(y) \longmapsto \dfrac{\mu(\psi(y) \cap \exists_{=N}x\phi(x,y))}{\mu(\exists_{=N}x\phi(x,y))}$ which is also $G$-invariant, and supported on $\mathrm{tp}(d/AB)$, we may assume $\mu(\exists_{=N}x\phi(x,y)) = 1$. The informal idea behind how we define $\nu$ is that, in order to pick randomly $x$, we first pick randomly $y$ according to $\mu$, then we pick randomly $x$ among the $N$ possibilities according to the uniform measure. Formally, we define:

$$\nu : \ \psi(x) \longmapsto \sum_{i=0}^{N} \frac{i}{N}\mu(\exists_{=i}x[\psi(x) \wedge \phi(x,y)])$$

Let us show that $\nu$ is a Keisler measure. If $\psi(x)$ holds for all $x$, then $\mu(\exists_{=N}x[\psi(x) \vee \phi(x,y)]) = 1$, and $\mu(\exists_{=i}x[\psi(x) \vee \phi(x,y)]) = 0$ for all $i \neq N$, thus $\nu(\psi(x)) = 1$. Suppose the definable sets $\psi(x)$ and $\theta(x)$ are disjoint. Let $\lambda_{i,j} = \mu([\exists_{=i}x(\psi(x) \wedge \phi(x,y))] \wedge [\exists_{=j}x(\theta(x) \wedge \phi(x,y))])$. Then we can see that:

$$\mu(\exists_{=i}x[\psi(x) \wedge \phi(x,y)]) = \sum_{j=0}^{N-i} \lambda_{i,j}$$

$$\mu(\exists_{=j}x[\theta(x) \wedge \phi(x,y)]) = \sum_{i=0}^{N-j} \lambda_{i,j}$$

$$\mu(\exists_{=k}x[(\psi(x) \vee \theta(x)) \wedge \phi(x,y)]) = \sum_{i=0}^{k} \lambda_{i,k-i}$$

therefore, we can compute:

$$\begin{aligned}
\nu(\psi(x) \vee \theta(x)) &= \sum_{k=0}^{N} \frac{k}{N} \sum_{i=0}^{k} \lambda_{i,k-i} \\
&= \sum_{i=0}^{N} \sum_{j=0}^{N-i} \frac{i+j}{N} \lambda_{i,j} \\
&= \sum_{i=0}^{N} \sum_{j=0}^{N-i} \frac{i}{N} \lambda_{i,j} + \sum_{i=0}^{N} \sum_{j=0}^{N-i} \frac{j}{N} \lambda_{i,j} \\
&= \sum_{i=0}^{N} \sum_{j=0}^{N-i} \frac{i}{N} \lambda_{i,j} + \sum_{j=0}^{N} \sum_{i=0}^{N-j} \frac{j}{N} \lambda_{i,j} \\
&= \nu(\psi(x)) + \nu(\theta(x))
\end{aligned}$$

It follows that $\nu$ is a Keisler measure, and it is $G$-invariant by definition. It remains to show that it is supported on $\mathrm{tp}(c/AB)$. Let $X$ be an $AB$-definable set such that $\mathrm{tp}(c/AB) \models x \in X$. Let $Y$ be the definable set $\phi(x,d) \wedge x \in X$, and let $i = |Y|$. As $c \in Y$, we have $i > 0$. As $\mu$ is supported on $\mathrm{tp}(d/AB)$, we have $\mu(\exists_{=i} x[x \in X \wedge \phi(x,y)]) > 0$. As a result, we have:

$$\nu(X) \geqslant \frac{i}{N} \mu(\exists_{=i} x[x \in X \wedge \phi(x,y)]) > 0$$

this concludes the proof. □

We conclude this chapter by drawing the Hasse diagram of the various independence notions that we defined. Before doing that, there is a last inclusion that needs to be proved:

**Proposition 1.3.19.** *We have* $\mathop{\smile\!\!\!\!\mid}^{\mathbf{L}\text{-KM}} \subseteq \mathop{\smile\!\!\!\!\mid}^{\mathbf{f}}$.

*Proof.* It suffices to prove that every measure which is $A$-Lascar-invariant does not fork over $A$. By contraposition, suppose we have $\mu$ a Keisler measure, and $\phi(x,b)$ a formula which divides over $A$, such that $\mu(\phi(x,b)) > 0$. Let $(b_i)_{i<\omega}$ be an $A$-indiscernible sequence such that $\{\phi(x,b_i) | i < \omega\}$ is inconsistent. Let $N$ be the largest integer for which:

$$\mu\left(\bigcap_{i<N} \phi(x,b_i)\right) > 0$$

Define $\lambda > 0$ as the left-handside of the above equation. For each integer $j < \omega$, let $X_j = \bigcap_{i<N} \phi(x, b_{jN+i})$. Suppose by contradiction that $\mu$ is $A$-Lascar-invariant. In an $A$-indiscernible sequence, two tuples of the same order-type are $A$-Lascar-conjugates, thus we have $\mu(X_i) = \mu(X_j)$ for all $i$, $j < \omega$. By maximality of $N$, we have $\mu(X_i \cap X_j) = 0$ whenever $i \neq j$. As a result, for $m < \omega$, we have:

$$\mu\left(\bigcup_{j<m} X_j\right) = m\lambda$$

if we choose $m$ large enough, we have $m\lambda > 1$, a contradiction. $\square$

We have enough data to draw the diagram in figure 1.1.

The upper relations are the weakest. Some of these inclusions might be equalities in general, and there may be non-trivial additional inclusions not appearing in the diagram. This leaves a great number of open question, but one that we find particularly interesting is whether the red relations coincide in general or not. The independence notion $\mathop{\smile\!\!\!\!\!\!|}^{\mathbf{bo}}$ can be shown to be weaker than $\mathop{\smile\!\!\!\!\!\!|}^{\mathbf{Sh}}$ and stronger than forking, but apart from that, we do not know exactly where it fits in the picture.

Before we begin the next chapter, we state the few results that we have from the literature:

**Fact 1.3.20.**

**([HP07], Proposition 2.11)** *In NIP theories, we have* $\mathop{\smile\!\!\!\!\!\!|}^{\mathbf{f}} = \mathop{\smile\!\!\!\!\!\!|}^{\mathbf{KP}}$.

**([HP07], Proposition 4.7)** *In a NIP theory, we have* $C \mathop{\smile\!\!\!\!\!\!|}^{\mathbf{f}}_A A$ *if and only if we have* $C \mathop{\smile\!\!\!\!\!\!|}^{\mathbf{inv}\text{-KM}}_A A$ *(this is a statement about extension bases). Moreover, if* $C \mathop{\smile\!\!\!\!\!\!|}^{\mathbf{inv}\text{-KM}}_A A$, *then a witness* $\mu$ *can be chosen such that* $\mu$ *is* ***concentrated*** *on* $\mathrm{tp}(C/A)$, *and not just supported.*

**([CHK+23], Proposition 2.12)** *We have* $\mathop{\smile\!\!\!\!\!\!|}^{\mathbf{f}} \neq \mathop{\smile\!\!\!\!\!\!|}^{\mathbf{inv}\text{-KM}}$ *in general. The counterexample is in a simple theory, and Pillay-Stonestrom recently found a counterexample in a NIP theory in ([PS23], subsection 4.3).*

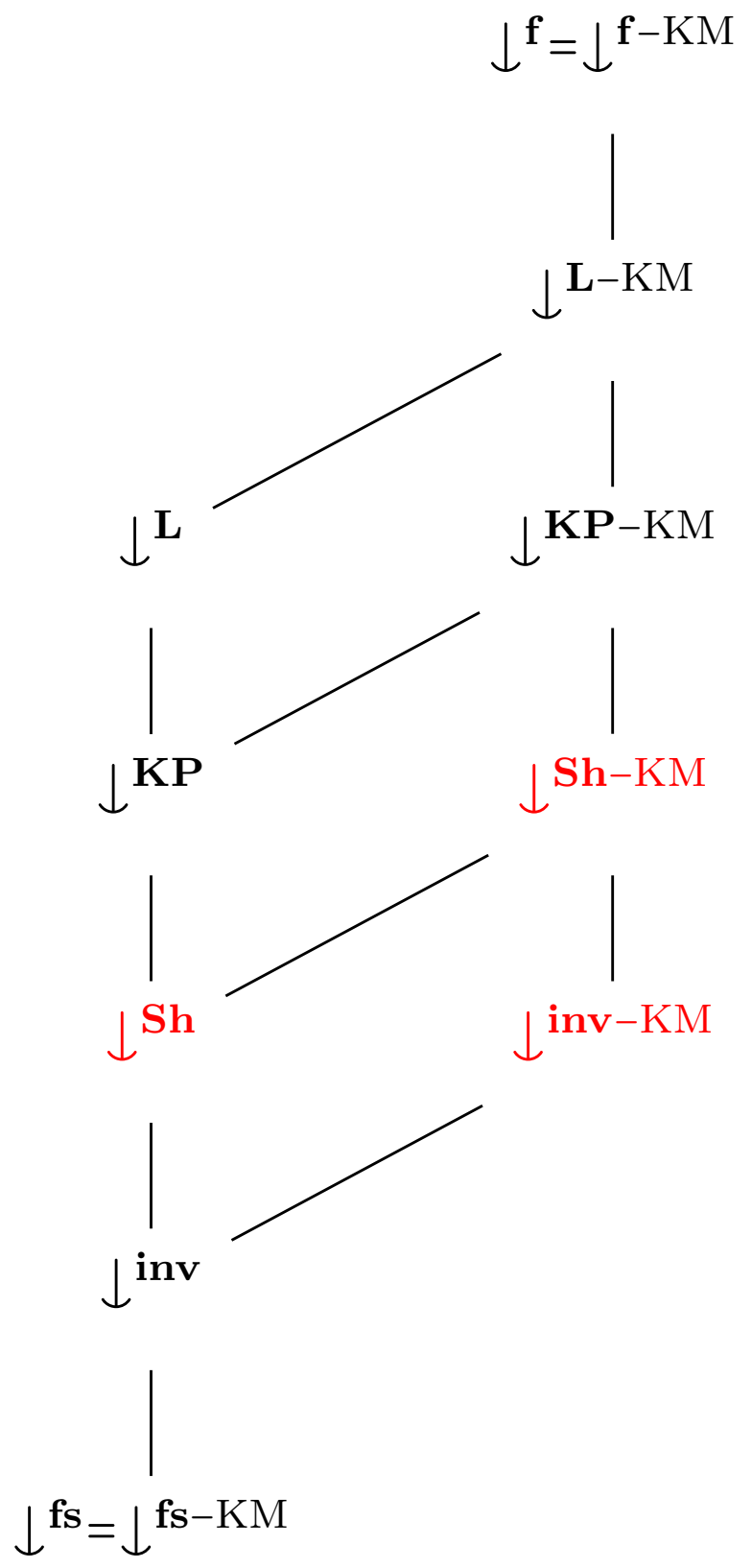

Figure 1.1: Diagram of independence notions

# Chapter 2

# An introduction to valuations

In this chapter, we define the mathematical structures that we manipulate in this thesis, and we try to build a geometric intuition for them. In particular, we want to have an informal idea of what independence should look like in this setting.

The valuations that one usually encounters in standard mathematics appear in three general settings, which are of course closely related:

- As the degree of polynomials, or, more generally, as the degree of the leading term of some series with transfinite support.
- As the convex valuation in ordered Abelian groups: $\mathrm{val}(x) < \mathrm{val}(y)$ whenever $x$ is infinitesimal compared to $y$. In case the structure is an ordered field, we see the elements having the same value as 1 as the *finite* elements, those of smaller value the *infinitesimal* elements, and those of greater value the *infinite* elements.
- In algebraic number theory, for each prime $p$, the $p$-adic valuation on $\mathbb{Q}$ measures the degree of divisibility by $p$ of a number. This valuation extends to the $p$-adic field $\mathbb{Q}_p$, which plays a fundamental role in number theory, where it is used to establish local-global principles.

Either way, one may see that a valuation measures the "order of magnitude" of some elements. Intuitively, the statement $\mathrm{val}(x - y) < \mathrm{val}(x - z)$ may be interpreted as "$y$ is infinitely closer to $x$ than $z$ is". So a valuation is a geometric notion, even though the geometric phenomena involving valuations differ from what happens in classical Archimedean geometry.

## 2.1 Valued groups

**Definition 2.1.1.** Let $G$ be an Abelian group, and $\Gamma$ a totally ordered set. A map $\mathrm{val}: G \longrightarrow \Gamma$ is a *semi-valuation* when it is onto, and when it satisfies the ultrametric inequality:

$$\forall x, y \in G \ \mathrm{val}(x - y) \leqslant \max(\mathrm{val}(x), \mathrm{val}(y))$$

*Remark* 2.1.2. Let $\mathrm{val}: G \longrightarrow \Gamma$ be a semi-valuation, and let $\Delta$ be a non-empty initial segment of $\Gamma$. Then, one may easily show with the axioms that $H = \{x \in G | \mathrm{val}(x) \in \Delta\}$ is a subgroup of $G$. In particular no element can have value smaller than that of 0, i.e. $\mathrm{val}(0)$ is the least element of $\Gamma$.

**Definition 2.1.3.** Let $G$ be an Abelian group, let $H$ be a subgroup of $G$, and $\Gamma$ a totally ordered set with a least element $-\infty$. A map $\mathrm{val}: G \longrightarrow \Gamma$ is an *$H$-valuation* if it is a semi-valuation such that for all $x \in G$:

$$\mathrm{val}(x) = -\infty \iff x \in H$$

a *valuation* over $G$ is a $\{0\}$-valuation.

If $\mathrm{val}'$ is another $H$-valuation of $G$, we say that $\mathrm{val}'$ *refines* val if there exists an order-preserving map $f: \mathrm{val}'(G) \longrightarrow \mathrm{val}(G)$ such that $\mathrm{val} = \mathrm{val}' \circ f$.

Equivalently, an $H$-valuation can be seen as a valuation over $G/H$. Note that any semi-valuation val is a $\mathrm{val}^{-1}(\mathrm{val}(0))$-valuation.

*Remark* 2.1.4. Let us remark for the readers unfamiliar with valuations that if $\mathrm{val}(x) < \mathrm{val}(y)$, then we have $\mathrm{val}(x+y) = \mathrm{val}(y)$. This comes from the fact that the subgroup $H$ of elements having value strictly smaller than $\mathrm{val}(y)$ contains $x$ and does not contain $y$, thus does not contain $x + y$. It follows that $\mathrm{val}(x+y) \geqslant \mathrm{val}(y)$, and the other direction is the ultrametric inequality up to dealing with subtraction.

The geometric meaning of this fact is that if we have a triangle with vertices three elements $x, y, z$ of $G$, such that the triangle is not equilateral with respect to the valuation (i.e. $\mathrm{val}(x - y), \mathrm{val}(x - z), \mathrm{val}(y - z)$ are not all equal), then the triangle is automatically isosceles, its two larger sides have the same value. If we relate this to our informal intuition on valuations, this translates to: "if $y$ is infinitely closer to $x$ than $z$ is, then the order of magnitude of the distance between $x$ and $z$ is the same as that between $y$ and $z$".

Speaking of distances, let us introduce more geometric objects:

**Definition 2.1.5.** A *ball* is a subset of $G$ of the form:

$$X = \{x \in G | \text{val}(x - c) R \gamma\}$$

with $c \in G$, $\gamma \in \Gamma$, and $R$ is either $<$ (in that case $X$ is *open*) or $\leqslant$ (in which case $X$ is *closed*). If $X$ is defined as above, then its *radius* is $\text{rad}(X) = \gamma$. Note that a closed ball can sometimes be written as an open ball with a different radius, thus the definition of the radius depends on the choice of the representation of the ball at hand. We allow the radius to be $+\infty$, in which case the induced open/closed ball is $G$. The ball is *non-trivial* if it can be written as a ball with finite (distinct from $\{\pm\infty\}$) radius.

Note that the trivial balls are $G$ (radius $+\infty$), $\varnothing$ (open of radius $-\infty$), and, if $-\infty$ has no successor in $\Gamma$, the singletons (closed of radius $-\infty$). If $-\infty$ has a successor $m$ in $\Gamma$, then the singletons are non-trivial, for they may be written as the open balls of radius $m$.

**Definition 2.1.6.** We equip $G$ with the topology having basic open sets the open balls.

*Remark* 2.1.7. This topology is compatible with the group structure on $G$. Note that, by remark 2.1.2, the open and closed balls are cosets of a chain of subgroups of $G$ (for the non-empty initial segments of $\Gamma$ form a chain). As a result, they are the elements of a chain of partitions of $G$: the non-empty balls of fixed radius and fixed shape (open or closed) form a partition of $G$. In particular, each open ball is clopen, for its complement is covered by the other open balls of same radius. With a similar reasoning, each non-trivial closed ball is clopen. Note that the topology is non-discrete if and only if the singletons are not open, that is when they are not non-trivial balls, i.e. when $-\infty$ does not have a successor.

The balls may be seen as the basic objects with which we build definable sets, the same way the intervals are the basic objects in, say, an o-minimal structure. With that analogy, what would correspond to the cuts? In a totally ordered set, an o-minimal 1-type (which we call a *cut* in this report) may be seen as the intersection of every interval containing some point in an extension. The corresponding notion in valued groups would be the intersection of every Boolean combination of balls containing some point in an extension. It is formally an intersection of (possibly infinitely many) crowns, a crown being a ball from which we removed a ball that can be empty. In other words, it is an intersection $I_G$ of balls (thus an intersection of cosets of a chain of subgroups), from which we removed the union $R_G$ of every sub-ball. Just as in remark 1.2.15, we may define notions of bottom/top generic extensions from a group $G$ to an extension $H$:

- The partial type given by $I_H \smallsetminus R_H$ is the top generic extension of $I_G \smallsetminus R_G$ if and only if $I_H = I_G$.
- $I_H \smallsetminus R_H$ is the bottom generic extension of $I_G \smallsetminus R_G$ if and only if the intersection of every ball in $H$ containing $R_G$ is $I_H$.

This gives us an intuition of what independence should look like in a structure with a valuation. However there are two issues:

- It may be possible that an $A$-definable ball has no point in $\operatorname{acl}(A)$, and thus cannot be defined with a standard formula $\operatorname{val}(x-c)R\gamma$ with parameters in $A$. This can be a problem when dealing with forking, as we do not have the freedom to change the parameter sets to add points. This leads us (in chapter 6) to finding sufficient conditions on parameter sets $A$ for every $A$-def ball, or every imaginary from a more abstract class, to have a point in A. It turns out the conditions we found do hold in the valued fields we are interested in (and, in fact, in a reasonably large class of valued fields), but they may fail in a very general setting.
- Just as in DLO or DOAG, once we understand forking in dimension one, the next step is to work in arbitrary dimension, possibly by gluing unary types into $n$-ary types. Such a thing cannot really be done in a field, the geometry is too rich. As a consequence, we do not have a reasonable description in, say, (transcendental) dimension 2, of what forking should look like. The problem basically comes from immediate extensions, which we do not understand well in higher dimension.

Let us end this small model-theoretic interlude and go back to algebra. The next lemma gives us a canonical way to build an $H$-valuation from a preorder over $G$ satisfying certain conditions:

**Lemma 2.1.8.** *Let $P$ be a preorder over $G$. Let $\sim$ be the associated equivalence relation over $G$, $\pi$ the quotient map, and $<$ the associated order on $G/\sim$. Suppose $<$ is total, and we have:*

$$\pi(x) < \pi(y) \Longrightarrow \pi(x+y) = \pi(y)$$

*for all $x, y \in G$ (the same relation as in remark 2.1.4). Then $\pi$ is a semi-valuation, thus $H = \pi^{-1}(\pi(0))$ is a subgroup of $G$, and $\pi$ is an $H$-valuation.*

*Proof.* Let $x, y \in G$. Then we cannot have $\pi(x - y) > \max(\pi(x), \pi(y))$, for we would have $\pi(x) = \pi(x - y + y) = \pi(x - y)$. □

We use this to build valuations in chapter 3.

## 2.2 Valued fields

An interesting example of a preorder on an Abelian group is the divisibility relation on a commutative unital ring $R$: $x \leqslant y$ if and only if $y$ divides $x$. One may ask when this preorder satisfies the key property $x < y \implies x + y \sim y$ from lemma 2.1.8:

**Proposition 2.2.1.** *The following are equivalent:*

1. *For all $x$, $y \in R$, if $y$ strictly divides $x$, then $R \cdot y = R \cdot (x + y)$.*
2. *For all $x \in R$, if $x$ is not invertible, then $1 + x$ is invertible.*
3. *$R$ is a local ring.*

*Proof.* Let us show $1 \implies 2$. Suppose $x$ is not invertible. Then 1 strictly divides $x$, thus $R = R \cdot (1 + x)$, thus $1 + x$ is invertible.

Let us show $2 \implies 3$ by contraposition. If 3 fails, then there must exist $x, y$ non-invertible elements of $R$ such that $x + y$ is invertible. Let $z \in R$ be such that $xz + yz = 1$. As $y$ is not invertible, neither is $-yz = xz - 1$. However, $1 - yz = xz$ is not invertible either, thus 2 fails.

Let us show $3 \implies 1$. Suppose $y$ strictly divides $x$. Then $x = yz$ with $z$ not invertible. As $R$ is a local ring, $z + 1$ is invertible, thus $R \cdot y = R \cdot y(z + 1) = R \cdot (x + y)$. □

However, the order induced by the divisibility relation is not always total in a local ring. It is total if and only if the principal ideals of $R$ form a chain. Then the set $Z$ of zero divisors must be an additive subgroup of $R$, and the divisibility relation induces a $Z$-valuation on $R$, for which it is clear that the fibers are the orbits of $R$ under the action of $R^*$, the group of invertible elements in $R$. A valued field is what we obtain when $R$ is additionally an integral domain, by naturally extending the valuation to the fraction field of $R$:

**Definition 2.2.2.** A *valued field* is a field $K$ with a valuation $\mathrm{val}: K \longrightarrow \Gamma$, equipped with the structure of an Abelian group on $\Gamma^* = \Gamma \smallsetminus \{-\infty\}$ such that:

- $(\Gamma^*, +, <)$ is an ordered Abelian group, i.e. every translation is an order-preserving map.
- val restricts to a group homomorphism $K^* \longrightarrow \Gamma^*$.

The order on $\Gamma$ is usually reversed in the literature, i.e. val(0) is the largest elements of $\Gamma$, and so on... This is a number-theoretic point of view, from which 0 is very large and 1 is very small. For our work, the geometric point of view is much more natural, and 0 is considered very small, while 1 is very large. Moreover, as we are using non-standard group valuations where the geometric point of view prevails, we think the best compromise is to keep consistent notations, and use the geometric notation for every valuation, even though it is less used in valued fields.

### 2.2.1 The algebraic structures at play

*Remark* 2.2.3. In a valued field, the closed ball $\mathcal{O}$ of radius 0 around 0 is a subring of $K$, called its *valuation ring*. We have $K = \mathrm{Frac}(\mathcal{O})$, and $\Gamma^*$ is naturally isomorphic to ${K^*}/{\mathcal{O}^*}$. For all $x$, $y$ in $K$, we have $\mathrm{val}(x) \leqslant \mathrm{val}(y)$ if and only if the $\mathcal{O}$-submodule of $K$ generated by $x$ is smaller than that generated by $y$. In particular, $\mathcal{O}$ is a local ring, the principal ideals of $\mathcal{O}$ form a chain, and the restriction of val to $\mathcal{O}$ coincides with the valuation on $\mathcal{O}$ induced by its divisibility predicate.

Conversely, given a *valuation ring* $\mathcal{O}$, i.e. an integral domain whose principal ideals form a chain, $\mathcal{O}$ is a local ring, and the map ($x \longmapsto$ the $\mathcal{O}$-submodule of $\mathrm{Frac}(\mathcal{O})$ generated by $x$) is a field valuation on $\mathrm{Frac}(\mathcal{O})$ whose valuation ring is $\mathcal{O}$. The restriction of the valuation to $\mathcal{O}$ coincides with the valuation from lemma 2.1.8 induced by the divisibility preorder.

Note that, as in remark 2.1.2, the preimage by val of a non-empty initial segment of $\Gamma$ is an $\mathcal{O}$-submodule of $K$. In particular, if the initial segment in question is included in $\mathrm{val}(\mathcal{O})$, then the preimage is an ideal of $\mathcal{O}$. Conversely, every ideal of $\mathcal{O}$ is a union of fibers of val, which must form an initial segment of $\mathrm{val}(\mathcal{O})$. In particular, the maximal ideal $\mathfrak{M}$ of $\mathcal{O}$ must be the preimage of the largest proper initial segment of $\mathrm{val}(\mathcal{O})$, namely that of $\mathrm{val}(\mathcal{O}) \smallsetminus \{0\}$, i.e. $\mathfrak{M}$ is the open ball of radius 0 around 0.

**Definition 2.2.4.** Let $K$ be a valued field. We call $\Gamma^*$ the *value group* of $K$, and $k = {}^{\mathcal{O}}\!/\!_{\mathfrak{M}}$ its *residue field.* We denote by res the projection $\mathcal{O} \longrightarrow k$. We have a short exact sequence of Abelian groups:

$$1 \longrightarrow \mathcal{O}^* \longrightarrow K^* \longrightarrow \Gamma^* \longrightarrow 0$$

If we quotient by $1 + \mathfrak{M} \leqslant \ker(\mathrm{val})$ (the group of elements infinitesimally close to 1), we obtain the short exact sequence:

$$1 \longrightarrow {}^{\mathcal{O}^*}\!/\!_{1+\mathfrak{M}} = k^* \longrightarrow {}^{K^*}\!/\!_{1+\mathfrak{M}} \longrightarrow \Gamma^* \longrightarrow 0$$

This leads us to define $\mathrm{RV}^* = {}^{K^*}\!/\!_{1+\mathfrak{M}}$ the *group of leading terms* of $K$. We add an element 0 to $\mathrm{RV}^*$, and we extend the canonical surjection to $\mathrm{rv} : K \longrightarrow \mathrm{RV} = \mathrm{RV}^* \cup \{0\}$ by sending 0 to 0. This way, val factors through rv, and RV is the quotient of $K$ by the equivalence relation:

$$x \sim y \Longleftrightarrow x = 0 = y \text{ or } \mathrm{val}(x - y) < \mathrm{val}(x)$$

*Example* 2.2.5. Let $k$ be a field, and $G$ be an ordered Abelian group. We denote by $k((t^G))$ the set of every series with coefficients in $k$, exponents in $G$, and anti-well-ordered support (i.e. every non-empty subset of the support has a largest element). The coordinate-wise sum and Cauchy product are well-defined over this set, making $k((t^G))$ a field. The function val mapping any series of this field to the largest element of its support is a field valuation, whose value group is $G$, and whose residue field is $k$ (the residue of an element of $\mathcal{O}$ is its 0-th coefficient). The class of a series in the group of leading terms is entirely determined by its leading term, i.e. the exponent and the coefficient of its term of largest exponent.

The field $k((t^G))$ is called a *Hahn field.*

### 2.2.2 Angular components

**Definition 2.2.6.** An *angular component* ac on a valued field is a group homomorphism $K^* \longrightarrow k^*$ extending res on $\mathcal{O}^*$.

Such a map clearly factors through rv, and the induced group homomorphism $\mathrm{RV}^* \longrightarrow k^*$ makes the following short exact sequence split:

$$1 \longrightarrow k^* \longrightarrow \mathrm{RV}^* \longrightarrow \Gamma^* \longrightarrow 0$$

**Definition 2.2.7.** A *cross-section* on $K$ is a section of val: a group homomorphism $s : \Gamma^* \longrightarrow K^*$ such that val: $\mathrm{val} \circ s = \mathrm{id}$.

A cross-section always induces an angular component, by setting $\mathrm{ac}(x) = \mathrm{res}\left(\dfrac{x}{s(\mathrm{val}(x))}\right)$.

Note that not every valued field admits an angular component. However, Hahn fields always admit a cross-section, by setting the cross-section of some value as the unitary monomial having exponent that value. The corresponding angular component maps a Hahn series to the coefficient of its leading term.

### 2.2.3 Separatedness

We describe here how the notion of a totally independent family in a vector space generalizes to valuations.

**Definition 2.2.8.** Let $(K, \nu)$ be a valued field with value group $\Gamma$. Let $V$ be a $K$-vector space, $W$ a $K$-vector subspace, and $\Delta$ a totally ordered set. Then a *$K$-linear $W$-valuation on $V$ of value set* $\Delta$ is a $W$-valuation $\mathrm{val} : V \longrightarrow \Delta$, together with an action of $\Gamma^*$ on $\Delta$, such that the following conditions hold for all $\gamma, \gamma' \in \Gamma^*$, $\delta, \delta' \in \Delta \setminus \{-\infty\}$, $\lambda \in K^*$, $v \in V$:

- The stabilizer of $-\infty$ is $\Gamma^*$.
- If $\gamma < \gamma'$, then $\gamma\delta < \gamma'\delta$.
- If $\delta < \delta'$, then $\gamma\delta < \gamma\delta'$.
- $\mathrm{val}(\lambda \cdot v) = \nu(\lambda)\mathrm{val}(v)$.

With this structure, $V$ is called a *$W$-valued $K$-vector space.* A $K$-linear valuation on $V$ is a $K$-linear $\{0\}$-valuation, and a valued $K$-vector space is a $\{0\}$-valued $K$-vector space. We usually abuse notations and call val the valuations of both the field, and the vector space.

In particular, any torsion-free divisible valued group is a valued $\mathbb{Q}$-vector space with respect to the trivial valuation on $\mathbb{Q}$. Another example of a valued vector space is in an extension of valued fields $K \leqslant L$, the valuation on the $K$-vector space $L$ is $K$-linear.

Recall the definition of $\mathrm{LC}(K)$ from definition 1.2.6. In a $K$-vector space $V$, given some family $v = (v_i)_{i \leqslant n}$ of $V$, the algebraic behavior of the $v_i$ only depends on which are the $f \in \mathrm{LC}(K)$ for which $f(v) = 0$. The family $v$ is $K$-linearly independent when the trivial map $f = 0$ is the only solution. More generally, we can speak of families $v$ which are $K$-linearly independent over some $K$-vector subspace $A$, when for every non-trivial $f \in \mathrm{LC}(K)$, we have $f(v) \notin A$. Equivalently, one non-conventional way to express independence over $A$ is to say that whether $f(v) \in A$ only depends on the support of $f$, i.e. the indices of its non-zero coordinates in the basis of $\mathrm{LC}(K)$ given by the projections:

*Remark* 2.2.9. The following conditions are equivalent:

1. For all $f, g \in \mathrm{LC}(K)$, if $f$ and $g$ have the same support, then $f(v) \in A$ if and only if $g(v) \in A$.
2. The elements of $v$ which are not in $A$ form a $K$-linearly independent family over $A$.

Now, if $V$ is equipped with a $K$-linear $W$-valuation, then the algebraic behavior of $v$ depends on the value of $f(v)$ for each $f \in \mathrm{LC}(K)$. When $K$ is not trivially valued, and the terms of $v$ are not all in $W$, this cannot possibly only depend on the support of $f$, for changing the value of one term of the support of $f$ can change the value of $f(v)$. However, $\mathrm{val}(f(v))$ may only depend on the values of the terms of $f$:

**Lemma 2.2.10.** *Suppose the residue field of $K$ is not $\mathbb{F}_2$. Then the following are equivalent:*

*a. If $(\lambda_i)_i, (\mu_i)_i \in K^n$, and $\mathrm{val}(\lambda_i) = \mathrm{val}(\mu_i)$ for all $i$, then:*

$$\mathrm{val}\left(\sum_i \lambda_i \cdot v_i\right) = \mathrm{val}\left(\sum_i \mu_i \cdot v_i\right)$$

*b. For every $(\lambda_i)_i \in K^n$, we have:*

$$\mathrm{val}\left(\sum_i \lambda_i \cdot v_i\right) = \max\left\{\mathrm{val}(\lambda_i \cdot v_i) | i \leqslant n\right\}$$

*Proof.* The direction $b \implies a$ is trivial. Let us prove the other direction by contraposition.

Suppose $b$ fails, say, $(\lambda_i)_i$ is a counterexample. Let $j$ be such that $\mathrm{val}(\lambda_j \cdot v_j)$ is maximal, and let $u = \sum_i \lambda_i \cdot v_i$. We have $\mathrm{val}(u) < \mathrm{val}(\lambda_j \cdot v_j)$. Let $\alpha \in \mathcal{O}$ be such that $\mathrm{res}(\alpha)$ is not 0 or 1. Then we have $\mathrm{val}(\lambda_j) = \mathrm{val}(\lambda_j \alpha) = \mathrm{val}(\lambda_j(1-\alpha))$. As $\mathrm{val}(u) < \mathrm{val}(\lambda_j \cdot v_j) = \mathrm{val}(\alpha\lambda_j \cdot v_j)$, we have $\mathrm{val}(u - \alpha\lambda_j \cdot v_j) = \mathrm{val}(\lambda_j v_j) > \mathrm{val}(u)$, and $a$ fails by replacing $\lambda_j$ by $\lambda_j(1-\alpha)$. □

**Definition 2.2.11.** We say that $v = (v_i)_i$ is val-*separated* if condition 2 from lemma 2.2.10 holds, and if no $v_i$ is in $W$. Whenever val is implicit, we say that the family is *separated*.

This definition naturally strengthens that of a linearly independent family. However, not every valued vector space is generated by a separated basis. For the remainder of this section, we focus on separatedness for extensions of valued fields.

**Definition 2.2.12.** An extension $L \geqslant K$ of valued fields is *separated* when every finite-dimensional $K$-vector subspace of $L$ admits a separated basis.

Let us describe the different phenomenon that one may encounter in an extension of valued fields:

**Fact 2.2.13.** *Let $K \leqslant L$ be an extension of valued fields, and $a \in L$ a singleton. Let $\Gamma(K), \Gamma(L)$ be the respective value groups of $K, L$, and $k(K)$, $k(L)$ their residue fields. Let $\mathcal{B}$ be the chain of every ball of $L$ having a point in $K$, a radius in $\Gamma(K)$, and which contains $a$. Then exactly one of the following four conditions hold:*

1. *We have $a \in K$. In that case, the least element of $\mathcal{B}$ is the singleton $\{a\}$.*

2. *The residue field of $K(a)$ is a proper extension of $k(K)$. In that case, the least element of $\mathcal{B}$ is a non-trivial closed ball which is not a singleton.*

3. *The value group of $K(a)$ is a proper extension of $\Gamma(K)$. In that case, $\mathcal{B}$ does not admit a least element which is a non-trivial closed ball, and its intersection has a point in $K$.*

4. *None of the above. In that case, $\mathcal{B}$ does not admit a least element which is a non-trivial closed ball, and its intersection is disjoint from $K$.*

**Definition 2.2.14.** If condition 2 of fact 2.2.13 holds, then we say that $a$ is *residual* over $K$. If condition 3 holds, then $a$ is *ramified*, and it is *immediate* if condition 4 holds. By convention, if condition 1 holds, then $a$ is residual, ramified and immediate at the same time.

The extension $K \leqslant L$ is *residual* (resp. *ramified*, *immediate*) if for every intermediate extension $F$, for every $a \in L$, $a$ is residual (resp. ramified, immediate) over $F$.

A *maximal* valued field is a valued field which does not admit any proper immediate extension.

There is also a similar case disjunction for transcendental extensions. The following is a variant of results by Abhyankar ([Abh56], Section 1):

**Fact 2.2.15.** *Let $F$ be an algebraically closed valued field. Note that the value group of $F$ is divisible and torsion-free, thus a $\mathbb{Q}$-vector space. Let $K \leqslant L$ be subfields of $F$ such that* $\operatorname{trdeg}(L/K)$ *is finite. Let $a_1 \dots a_n$ be a $K$-transcendence base of $L$. Then, the number of times when $a_i$ is ramified over $(K(a_{<i}))^{\mathbf{alg}}$ coincides with the dimension of the quotient $\Gamma(L^{\mathbf{alg}}) \big/ \Gamma(K^{\mathbf{alg}})$ as a $\mathbb{Q}$-vector space (and hence does not depend on the choice of the $a_i$). Likewise, the number of times when $a_i$ is residual over $(K(a_{<i}))^{\mathbf{alg}}$ coincides with the transcendence degree of $k(L)$ over $k(K)$.*

*Conversely, if $(a_i)_i$, $(b_j)_j$ are families of elements of $L$ such that:*

- $(\operatorname{val}(a_i))_i$ *is $\mathbb{Q}$-linearly independent over (the divisible closure of)* $\Gamma(K)$.
- $b_j \in \mathcal{O}$ *for all $j$, and* $(\operatorname{res}(b_j))_j$ *is algebraically independent over $k(K)$.*

*Then the $a_i$ and the $b_j$ form a $K$-algebraically independent family.*

The issue with separated extensions of valued fields is the following:

**Fact 2.2.16.** *A proper immediate extension of valued fields is never separated. However, any extension of a maximal valued field is separated.*

This shows that maximal valued fields are particularly nice. However, valued fields are rarely maximal, and maximality is not an elementary property. Maximality is a form of completeness, for a valued field $K$ is maximal if and only if the intersection of any chain of non-empty balls of $K$ has a point in $K$ (as the only case from the disjunction in fact 2.2.13 where the intersection does not have a point in $K$ is clearly the immediate case). A natural way to express completeness with first-order formulas is as follows:

**Fact 2.2.17.** *For any valued field $K$ of residue characteristic zero, the following are equivalent:*

- *For every polynomial $f \in \mathcal{O}[X]$ in one variable, if $\alpha \in \mathcal{O}$ is such that $\mathrm{val}(f(\alpha)) < 0$ and $\mathrm{val}(f'(\alpha)) = 0$, then there exists $\beta \in K$ such that $\mathrm{val}(\alpha - \beta) < 0$, and $f(\beta) = 0$.*
- *No proper immediate extension of $K$ is algebraic.*

**Definition 2.2.18.** A *Henselian valued field* (of arbitrary residue characteristic) is a valued field satisfying the first condition from fact 2.2.17.

The class of Henselian valued fields is clearly elementary. Most of the literature on the model theory of valued fields is focused on Henselian valued fields, and this thesis is no exception. We give a more detailed survey of the model theory of these structures in section 2.3. Despite the fact that the class of maximal valued fields is not elementary, we have the following inclusion:

**Fact 2.2.19.** *Every maximal valued field is Henselian.*

### 2.2.4 Local fields

**Definition 2.2.20.** A *non-Archimedean local field* (*local field* for short, since we only consider valued fields here) is a valued field which is locally compact and non-discrete for the topology induced by the valuation.

We state a classification of local fields in fact 2.3.10. While this classification is not an easy result, we establish in this subsection some easy properties of local fields which, from the point of view of model theory, already suffice to show results on their ultraproducts. Local fields play an important role in contemporary mathematics, notably in number theory and algebraic geometry. The study of local fields is the main motivation behind the celebrated article of Ax-Kochen ([AK65a], [AK65b], [AK66]), independently Ershov ([Ers65]), where they prove the Ax-Kochen-Ershov principle for Henselian valued fields.

*Remark* 2.2.21. Let us start by noting that in a valued field with non-discrete topology, the non-trivial open balls are pairwise-homeomorphic via the affine transformations. The same can be said of the non-trivial closed balls.

**Proposition 2.2.22.** *In a local field, the residue field is finite.*

*Proof.* Note that the non-trivial closed balls form a base for the topology. Any point has a compact neighbourhood, which can be assumed to be a non-trivial closed ball, which is homeomorphic to the non-trivial closed ball $\mathcal{O}$. It follows that $\mathcal{O}$ is compact. However, $\mathcal{O}$ is covered by the cosets of $\mathfrak{M}$, which are open balls. By compactness, there are only finitely many cosets, concluding the proof. □

**Proposition 2.2.23.** *In a local field, the value group is isomorphic to* $(\mathbb{Z}, +, <)$.

*Proof.* The value group is non-trivial as the topology is non-discrete. It suffices to show that every non-empty subset of $\Gamma$ which has a lower bound has a least element. Let $P$ be a non-empty subset of $\Gamma$ having a lower bound. Let $F$ be the complement of the union of every closed ball around 0 whose radius is a lower bound of $P$. Let $r \in P$, and let $X$ be the closed ball of radius $r$ around 0. Then $X$, and hence $X \cap F$, is a compact subspace of $K$. The intersection of (the chain of) every closed ball around 0 whose radius is in $P$ is disjoint from $F \cap X$, therefore, by compactness, some element $Y$ of this chain is disjoint from $F \cap X$. The radius of $Y$ must be a lower bound of $P$, but it is in $P$, concluding the proof. □

**Proposition 2.2.24.** *Any local field is maximal, and hence Henselian.*

*Proof.* We use the fact that a valued field is maximal if and only if it is spherically complete, i.e. any chain of balls has a non-empty intersection.

Let $\mathcal{B}$ be a chain of balls. If $\mathcal{B} = \varnothing$, then 0 (and, in fact, any point) is in the intersection. If $X \in \mathcal{B}$, then $X$ is compact, and the intersection of $\mathcal{B}$ may be written as the intersection of a chain of non-empty closed subsets of $X$. By compactness, this intersection has a point, hence the valued field is maximal. □

## 2.3 Model theory of Henselian valued fields

### 2.3.1 Notations

**Definition 2.3.1.** In model theory, there are several natural choices for a language in which to consider valued fields:

- The language $\mathcal{L}_{\mathrm{div}}$, with one sort for the domain of the valuation. This language is the enrichment of the language of rings with the binary predicate $\mathrm{val}(x) \leqslant \mathrm{val}(y)$.

- The language $\mathcal{L}_{\mathrm{val}}$ with two sorts: $K$ for the domain of the valuation, and $\Gamma$ for the value set. We have the ring language on $K$, the language of ordered groups plus the constant $-\infty$ on $\Gamma$, and a function symbol for the valuation $K \longrightarrow \Gamma$.

- The language $\mathcal{L}_{k,\Gamma}$, which is an enrichment of $\mathcal{L}_{\mathrm{val}}$. It adds a third sort $k$ for the residue field with (a copy of) the ring language of $k$, and a function symbol for the map:

$$\begin{array}{rccl} r: & K \times K & \longrightarrow & k \\ & (x,y) & \longmapsto & \begin{cases} \mathrm{res}\left(\dfrac{x}{y}\right) \text{ if } \mathrm{val}(y) \geqslant \mathrm{val}(x) > -\infty \\ 0 \text{ otherwise .} \end{cases} \end{array}$$

- The language $\mathcal{L}_{\mathrm{RV}}$ with two sorts: $K$ with the ring language, RV for the group of leading terms with the group language, a constant symbol for $0 = \mathrm{rv}(0)$, the map $\mathrm{rv}: K \longrightarrow \mathrm{RV}$, and the following ternary predicate on RV:

$$\oplus(x,y,z) \ : \ \exists x', y' \in K \ \ (\mathrm{rv}(x') = x \wedge \mathrm{rv}(y') = y \wedge \mathrm{rv}(x'+y') = z)$$

- The expansion $\mathcal{L}^{eq}$ of all/any of those languages by imaginaries.

These choices are all pairwise-bi-interpretable. We should also mention the geometric language introduced in [HHM05] in order to eliminate imaginaries.

For the majority of this thesis, we work implicitly in $\mathcal{L}^{eq}$, with parameter sets that are either generated by a subfield of $K$, or, on the contrary, arbitrary $\mathcal{L}^{eq}$-substructures. Given a first-order structure $M$ in those languages, and $S$ a (potentially imaginary) sort (such as $K$, $\Gamma$, RV...), we denote by $S(M)$ the reduct of $M$ to $S$.

*Remark* 2.3.2. If $M$ is a valued field, and $A \subseteq M^{eq}$, there may exist an element of $k(A)$ that cannot be pulled back to an element of $K(A)$. In other words, $k(A)$ and $res(K(A))$ may not coincide, these two notations are not aliases. Likewise, we may have $\mathrm{val}(K(A)) \neq \Gamma(A)$.

The notations $k$, $K$, $\Gamma$ refer to definable sets. When we deal with an actual field or group, we usually call them $k'$, $G$.

**Definition 2.3.3.** By condition 2 of fact 2.2.17, the class of Henselian valued fields is clearly elementary in those languages. Let $\mathrm{HVF}_{0,0}$ be the theory of Henselian non-trivially valued fields of residue characteristic zero. We assume the valuation to be non-trivial in order to have algebraic boundedness, as in corollary 6.1.2. Any Ax-Kochen-Ershov principle that one may think of trivially holds in trivially valued fields anyway.

**Definition 2.3.4.** We use valued fields with angular components or cross-sections in chapter 6. We see them as expansions of $\mathcal{L}_{k,\Gamma}$, $\mathcal{L}_{\mathrm{RV}}$ or $\mathcal{L}^{eq}$ by the additional function symbols $\mathrm{ac}: K \longrightarrow \mathrm{RV} \longrightarrow k$, and possibly $s: \Gamma \longrightarrow K$. We denote by $\mathcal{L}_{Pas}$ the enrichment of $\mathcal{L}_{k,\Gamma}$ to the function ac, and we denote $\mathcal{L}_{cs}$ the enrichment of $\mathcal{L}_{Pas}$ to the function $s$. We assume that in a $\mathcal{L}_{cs}$-valued field, the angular component is the one corresponding to the cross-section.

Recall that a valued field may not admit an angular component, and, even if it did, this structure would not be interpretable in the languages of valued fields.

**Fact 2.3.5** ([Che76], Corollary 28)**.** *If $M \vDash \mathrm{HVF}_{0,0}$ is $\aleph_1$-saturated, then it admits a cross-section, hence it can be expanded to a $\mathcal{L}_{cs}$-structure.*

### 2.3.2 Results from early literature

The most fundamental result about Henselian valued fields is the Ax-Kochen-Ershov (AKE) theorem:

**Theorem 2.3.6** ([AK65a], [AK65b], [AK66], [Ers65])**.** *Suppose $M$ and $N$ are two models of $\mathrm{HVF}_{0,0}$ such that $k(M) \equiv k(N)$, and $\Gamma(M) \equiv \Gamma(N)$. Then we have $M \equiv N$.*

Note that the naive version of the AKE principle for forking fails in general:

*Example* 2.3.7. Let $M$ be a model of $\mathrm{HVF}_{0,0}$ such that there exists, in an elementary extension of $M$, a proper immediate extension $B$ of $M$. We clearly have $k(B) \mathop{\downarrow}^{\mathbf{f}}_{k(M)} k(B)$ and $\Gamma(B) \mathop{\downarrow}^{\mathbf{f}}_{\Gamma(M)} \Gamma(B)$, however we have $B \not\mathop{\downarrow}^{\mathbf{d}}_{M} B$.

The first application of the AKE principle was to characterise the following class of residue characteristic zero valued fields:

**Definition 2.3.8.** We denote by $\mathrm{PL}_0$ the theory of valued fields of residue characteristic zero satisfying the sentences that hold in every local field. The models of this theory are the *pseudo-local fields of residue characteristic zero* (as our work only takes place in residue characteristic zero, we call them *pseudo-local fields* for short).

*Remark* 2.3.9. By proposition 2.2.24, $\mathrm{PL}_0$ extends $\mathrm{HVF}_{0,0}$. Moreover, by proposition 2.2.22 and proposition 2.2.23, the value group of a model of $\mathrm{PL}_0$ is an elementary extension of $(\mathbb{Z}, +, <)$, and its residue field is pseudo-finite.

These are easy necessary condition to be a model of $\mathrm{PL}_0$. It turns out that they are sufficient:

**Fact 2.3.10** ([Wei95], Theorem 5 and Theorem 8)**.** *Up to isomorphism, the local fields are exactly the finite extensions of $\mathbb{Q}_p$ with $p$ prime, with the topology induced by the $p$-adic valuation, and the valued fields of Laurent series over finite fields, i.e. of the form $\mathbb{F}_q((t^{\mathbb{Z}}))$.*

**Corollary 2.3.11.** *The models of* $\mathrm{PL}_0$ *are exactly the Henselian valued fields of residue characteristic zero with pseudo-finite residue field, whose value groups are elementary extensions of* $(\mathbb{Z}, +, <)$*.*

*Proof.* Suppose $M$ is such a field. Then there exists a non-principal ultrafilter $\mathcal{U}$ on the set of prime numbers such that $k(M)$ is elementarily equivalent to the ultraproduct of the $\mathbb{F}_p$ according to $\mathcal{U}$. By the AKE principle, $M$ is elementarily equivalent to the ultraproduct of the $\mathbb{Q}_p$ according to $\mathcal{U}$, which is a model of $\mathrm{PL}_0$.

Conversely, one may use standard model-theoretic arguments to prove that any model of $\mathrm{PL}_0$ is elementarily equivalent to some non-principal ultraproduct of local valued fields, which implies that its residue field is pseudo-finite, and its value group is elementarily equivalent to $\mathbb{Z}$, therefore an elementary extension of $\mathbb{Z}$, as $\mathbb{Z} = \mathrm{dcl}(\varnothing)$ is always an elementary substructure of any ordered group elementarily equivalent to $\mathbb{Z}$. □

Another notable class of Henselian valued fields plays an important role in the literature:

**Definition 2.3.12.** We denote by $\mathrm{ACVF}_{0,0}$ the theory of non-trivially valued algebraically closed fields of residue characteristic zero.

By condition 1 from fact 2.2.17, any algebraically closed valued field is Henselian. By the AKE principle, $\mathrm{ACVF}_{0,0}$ is exactly the theory of Henselian valued fields with non-trivial divisible value group (i.e. $\Gamma \vDash \mathrm{DOAG}$) and algebraically closed residue field of characteristic zero (i.e. $k \vDash \mathrm{ACF}_0$).

Just as DOAG and $\mathrm{ACF}_0$ naturally play the role of "canonical closures", the model-completions of the respective theories of ordered Abelian groups and fields of characteristic zero, we expect $\mathrm{ACVF}_{0,0}$ to play the same role for Henselian non-trivially valued fields of residue characteristic zero. In fact, we expect it to eliminate quantifiers. This is shown to be true in many languages, in a series of results initiated by Robinson ([Rob56]):

**Fact 2.3.13** ([Rob56], [Del82], [Wei84], [Hol95])**.** *The theory* $\mathrm{ACVF}_{0,0}$ *is complete and eliminates quantifiers in* $\mathcal{L}_{\mathrm{div}}$*,* $\mathcal{L}_{\mathrm{val}}$*,* $\mathcal{L}_{k,\Gamma}$*. Moreover, every unary definable subset of* $K$ *in this theory is a Boolean combination of balls.*

As a result, the informal notion of genericity from remark 2.1.7 is especially relevant for this theory, though the results that we establish for forking in general Henselian valued fields were already established in ([HP07], Corollary 2.14) in the particular case of ACVF (every set in a $C$-minimal theory is an extension base). Quantifier elimination results in ACVF suggest that there is a general notion of quantifier-free definable sets in Henselian valued fields. In particular, the language we consider for this notion does not matter:

**Proposition 2.3.14.** *Modulo the theory of valued fields of residue characteristic zero, every quantifier-free formula from one of the languages* $\mathcal{L}_{\mathrm{div}}$*,* $\mathcal{L}_{\mathrm{val}}$*,* $\mathcal{L}_{k,\Gamma}$ *with free variables in the sort* $K$ *is equivalent to a quantifier-free formula in the two other languages.*

*Proof.* Let $\phi(x)$ be such a formula. Let $\psi_1(x), \psi_2(x)$ be formulas of the two other languages which are equivalent to $\phi(x)$. Let $\varphi_i(x)$ be a quantifier-free formula which is equivalent to $\psi_i(x)$ modulo $\mathrm{ACVF}_{0,0}$. Let us show that $\phi(x)$ is equivalent to $\varphi_i(x)$ in any valued field of residue characteristic zero.

Let $M$ be such a valued field, and $a$ a tuple from $M$. Choose an embedding of valued fields $\iota$ from $M$ to some model $N$ of $\mathrm{ACVF}_{0,0}$ ($\iota$ is not unique if $M$ is not Henselian, but it does not matter). Then $\iota$ is clearly an embedding with respect to the languages $\mathcal{L}_{\mathrm{div}}, \mathcal{L}_{\mathrm{val}}, \mathcal{L}_{k,\Gamma}$, thus we have:

$$M \vDash \phi(a) \iff N \vDash \phi(\iota(a)), M \vDash \varphi_i(a) \iff N \vDash \varphi_i(\iota(a))$$

it follows that $M \vDash \phi(a)$ if and only if $M \vDash \varphi_i(a)$. This holds for any $a$, and any $M$, concluding the proof. □

# Chapter 3

# Divisible ordered Abelian groups

This chapter corresponds to Sections 1, 2, 3, 4 of our preprint ([Hos23b]).

## 3.1 Forking in DOAG

### 3.1.1 Cut-independence

For many common elementary classes of algebraic structures, we consider the model completion, when it exists, as the theory which "captures the algebraic part" of the class. For this completion, we usually have an algebraic intuition of what independence should be, and relating forking with this algebraic independence notion would be the starting point to compute forking in some theory from the class. See some concrete examples in the next page.

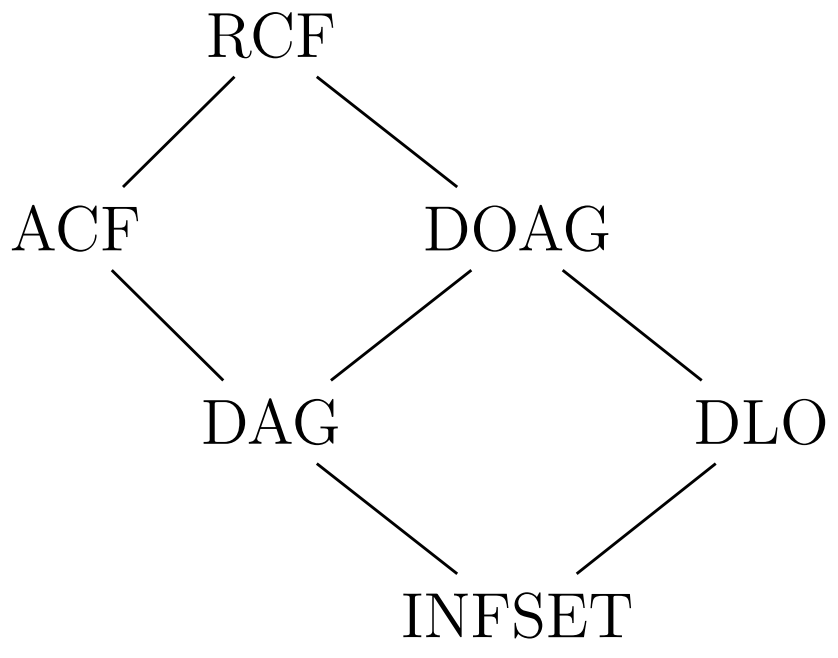

| Elementary class | model completion | Natural notions of independence |
|---|---|---|
| INFSET | INFSET | $C \cap B \subseteq A$ |
| Torsion-free Abelian Groups | Non-trivial divisible Abelian groups (DAG) | $\mathrm{Vec}_{\mathbb{Q}}(AC) \cap \mathrm{Vec}_{\mathbb{Q}}(AB) = \mathrm{Vec}_{\mathbb{Q}}(A)$ |
| Total orders | DLO | See subsection 1.2.2 |
| Fields | ACF | $\mathop{\smile\!\!\!\!\!\!\!\;\mid}^{\mathbf{alg}}$ |
| Ordered Abelian groups | DOAG | ? |
| Valued fields | ACVF | ? |
| Ordered fields | RCF | See [Dol04] |

The theory ACVF is much more complicated than the others, so let us ignore this issue for the moment. We have a diagram of theories, and we know forking for all of them but DOAG. This is very odd, as DOAG is not even the most complicated theory from the diagram. There is undoubtedly a gap to be filled, and this is what the results from this chapter achieve.

Recall that we defined at the end of subsection 1.3.1 several independence relations $\mathop{\smile\!\!\!\!\!\!|}^{\mathbf{inv}}, \mathop{\smile\!\!\!\!\!\!|}^{\mathbf{bo}}, \mathop{\smile\!\!\!\!\!\!|}^{\mathbf{Sh}}$ involving the existence of some tame global extension of some type. Recall that the global extension in question has to:

- be invariant over the base in the case of $\mathop{\smile\!\!\!\!\!\!|}^{\mathbf{inv}}$.
- be invariant over the imaginaries which are algebraic over the base in the case of $\mathop{\smile\!\!\!\!\!\!|}^{\mathbf{Sh}}$.
- have a bounded orbit under the action of the automorphisms which are invariant over the base in the case of $\mathop{\smile\!\!\!\!\!\!|}^{\mathbf{bo}}$.

*Remark* 3.1.1. Note that in DOAG, if $A$ is not included in $\{0\}$, then $\mathrm{dcl}(A)$ is a model. In an o-minimal theory (in fact in any NIP theory), $\mathop{\smile\!\!\!\!\!\!|}^{\mathbf{f}}$, $\mathop{\smile\!\!\!\!\!\!|}^{\mathbf{d}}$ and $\mathop{\smile\!\!\!\!\!\!|}^{\mathbf{inv}}$ all coincide over models, and one can easily check by hand that it is also the case in DOAG for $A \subseteq \{0\}$. More precisely, the fact that $\mathop{\smile\!\!\!\!\!\!|}^{\mathbf{f}}=\mathop{\smile\!\!\!\!\!\!|}^{\mathbf{d}}$ follows from ([CK12], Theorem 1.1), while $\mathop{\smile\!\!\!\!\!\!|}^{\mathbf{f}}=\mathop{\smile\!\!\!\!\!\!|}^{\mathbf{inv}}$ over models trivially follows from the fact that non-forking coincides with the independence notion given by Lascar-invariance (see for instance [HP07], Proposition 2.1). As a result, the abstract equality $\mathop{\smile\!\!\!\!\!\!|}^{\mathbf{f}}=\mathop{\smile\!\!\!\!\!\!|}^{\mathbf{inv}}$ is already well-known in DOAG, though this equality does not help in any way to relate forking with more concrete geometric phenomena.

We show in chapter 4 that $\mathop{\smile\!\!\!\!\!\!|}^{\mathbf{f}}$, $\mathop{\smile\!\!\!\!\!\!|}^{\mathbf{d}}$, $\mathop{\smile\!\!\!\!\!\!|}^{\mathbf{bo}}$ and $\mathop{\smile\!\!\!\!\!\!|}^{\mathbf{Sh}}$ all coincide in regular groups.

**Assumptions 3.1.2.** We work in DOAG, in the language of ordered $\mathbb{Q}$-vector spaces. In this language, every substructure is a definably closed $\mathbb{Q}$-vector subspace. The independence notions we are looking at are insensitive to definable closure (i.e. $C \mathop{\smile\!\!\!\!\!\!|}_A B$ if and only if $\mathrm{dcl}(AC) \mathop{\smile\!\!\!\!\!\!|}_{\mathrm{dcl}(A)} \mathrm{dcl}(AB)$), so we fix $C \geqslant A \leqslant B$ $\mathbb{Q}$-vector subspaces of $M$.

It is well-known that DOAG is complete, has quantifier elimination and is o-minimal in this language. Moreover, for $c_1 \dots c_n, d_1 \dots d_n \in M$, we have:

$$(c_i)_i \equiv_A (d_i)_i \iff \forall f \in \mathrm{LC}^n(\mathbb{Q}),\ \mathrm{ct}\,(f(c)/A) = \mathrm{ct}\,(f(d)/A)$$

where $\mathrm{LC}(\mathbb{Q})$ is defined in definition 1.2.6.

The independence notions from the other theories in the previous diagram actually give us an intuition on what forking could be in DOAG. Indeed, DOAG is to DLO what DAG is to INFSET. When going from INFSET to DAG, the only thing that changes for forking is that we replace the parameter sets with their respective algebraic closures. Our intuition about base change suggests that this may be all we need to do to obtain forking in DOAG from forking in DLO. This motivates the definition of what we call *cut-independence*:

**Definition 3.1.3.** Define $C \mathop{\downfree}^{\mathbf{cut}}_A B$ if the following equivalent conditions hold:

- $\forall c_1 \dots c_n \in C,\ \forall f \in \mathrm{LC}^n(\mathbb{Q}),\ (f(c)) \mathop{\downfree}^{\mathbf{inv}}_A B$
- $\forall c \in C,\ c \mathop{\downfree}^{\mathbf{inv}}_A B$.
- Every closed interval with bounds in $B$ that has a point in $C$ already has a point in $A$.

For $c_1 \dots c_n \in M$, we define:

$$c_1 \dots c_n \mathop{\downfree}^{\mathbf{cut}}_A B \iff (A + \mathbb{Q}c_1 + \dots + \mathbb{Q}c_n) \mathop{\downfree}^{\mathbf{cut}}_A B$$

The last item of definition 3.1.3 is a purely geometric description which depends only on $A$, $B$, $C$, hence it would make for a very satisfactory description of forking if we had $\downfree^{\mathbf{f}} = \downfree^{\mathbf{cut}}$.

**Lemma 3.1.4.** *We have* $\downfree^{\mathbf{inv}} \subseteq \downfree^{\mathbf{d}} \subseteq \downfree^{\mathbf{cut}}$ *in* DOAG.

*Proof.* See corollary 1.2.12, and fact 1.1.19. □

As announced, it turns out that this nice, but weak independence notion $\downfree^{\mathbf{cut}}$ is not weaker than $\downfree^{\mathbf{inv}}$ in DOAG. This is the fundamental result of this chapter:

**Theorem 3.1.5.** *In* DOAG*, we have* $\mathop{\downarrow}^{\mathbf{inv}} = \mathop{\downarrow}^{\mathbf{cut}}$.

By lemma 3.1.4, we just need to prove that $\mathop{\downarrow}^{\mathbf{cut}} \subseteq \mathop{\downarrow}^{\mathbf{inv}}$.

**Assumptions 3.1.6.** On top of assumptions 3.1.2, we assume $C \mathop{\downarrow}^{\mathbf{cut}}_{A} B$, and we fix $c_1 \dots c_n \in C$.

Our goal is to build a global $\mathrm{Aut}(M/A)$-invariant extension of the type $\mathrm{tp}(c_1 \dots c_n/B)$. This is achieved in this chapter, and we deal with ROAG and examples in the next chapter.

*Remark* 3.1.7. If $c$ is a singleton from $C \smallsetminus A$, then we have $c \mathop{\downarrow}^{\mathbf{cut}}_{A} B$ if and only if at least one of the following conditions holds:

- Every bounded $B$-definable closed interval included in $\mathrm{ct}(c/A)$ (viewed as an $A$-type-definable set) has all its points smaller than $c$.

- Every bounded $B$-definable closed interval included in $\mathrm{ct}(c/A)$ has all its points larger than $c$.

**Definition 3.1.8.** If the first condition of remark 3.1.7 holds, then we say that $c$ *leans right* with respect to $A$ and $B$. If the second condition holds, then $c$ *leans left*. The two conditions hold at the same time if and only if $\mathrm{ct}(c/\mathrm{dcl}(A))$ has no point in $\mathrm{dcl}(AB)$.

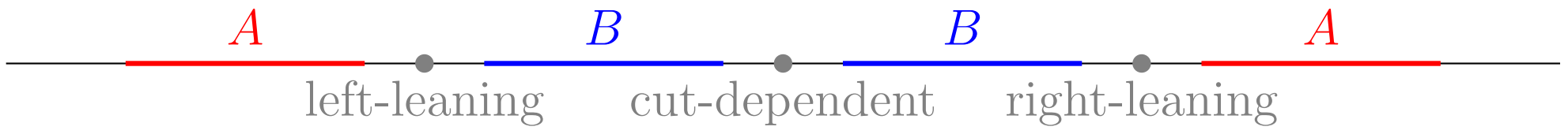

It is easy to see that $c$ leans right (resp. left) with respect to $A$ and $B$ if and only if there exists $A' \subseteq A$ such that $\mathrm{ct}(c/B) = \mathrm{ct}_{>}(A'/B)$ (resp. $\mathrm{ct}_{<}(A'/B)$), as defined in definition 1.2.13.

This definition also makes sense (and will be used) in more abstract total orders that do not come from ordered groups.

### 3.1.2 The Archimedean valuation

The content of this subsection is well-known in the literature, and most of the proofs are omitted.

**Assumptions 3.1.9.** In this subsection, we fix $G$ an ordered Abelian group.

**Proposition 3.1.10.** *Let $x, y \in G$. The following conditions are equivalent:*

1. $x \in \bigcup_{n<\omega} [-n|y|, n|y|]$.

2. *For every convex subgroup $H$ of $G$, if $y \in H$, then $x \in H$.*

3. *The convex subgroup of $G$ generated by $x$ is included in the one generated by $y$.*

**Definition 3.1.11.** These equivalent conditions define a relation on $G$ which is clearly a preorder that satisfies the necessary conditions of lemma 2.1.8, with $\pi^{-1}(\pi(0)) = \{0\}$. This gives a valuation over $G$, which we call $\Delta$, the *Archimedean valuation*. For $x \in G$, $\Delta^{-1}(\Delta(x))$ is called the *Archimedean class* of $x$. By abuse of notation, we may identify the Archimedean value of an element with its Archimedean class.

The idea behind the statement $\Delta(x) < \Delta(y)$ is that $x$ is "infinitesimal" compared to $y$, or that $y$ is "infinitely greater" than $x$.

Equivalently:

$$\begin{aligned} \Delta(x) < \Delta(y) \quad &\iff \quad x \in \bigcap_{n>0} \left]-\frac{1}{n}|y|, \frac{1}{n}|y|\right[ \\ &\iff \quad |y| \in \bigcap_{n>0} ]n|x|, +\infty[ \end{aligned}$$

Let us draw a picture (not scaled correctly) which intuitively shows how the Archimedean classes look like. Suppose $G$ has four Archimedean classes $\Delta(0) < \delta_1 < \delta_2 < \delta_3$. Then $G$ looks like this:

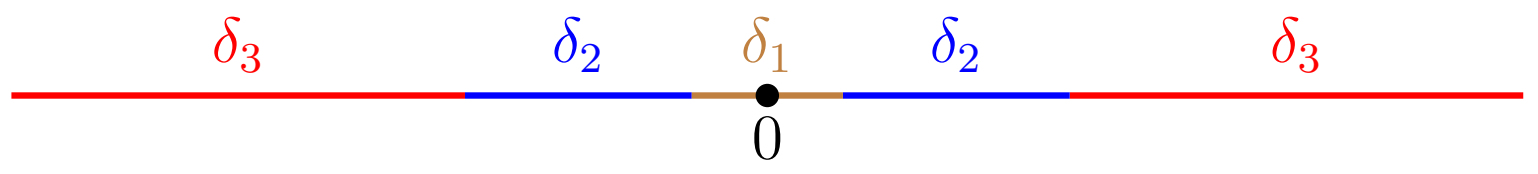

There is no standard terminology in the modern literature when it comes to the Archimedean valuation and related notions. Some use the keyword "convex" to refer to all those notions, others use the very obscure names "*i*-extension" or "*i*-completeness". We draw our inspiration from older sources: we choose to use the terminology of [Gra56], which, on top of being a good introductory paper, sets up notations which we find more intuitive and easier to work with.

*Remark* 3.1.12. A convex subgroup $H \leqslant G$ is equal to the union of the convex subgroups of $G$ generated by each $x \in H$. So $H$ is entirely determined by its direct image $\Delta(H)$, which is an initial segment of $\Delta(G)$. We can note that the set $\mathcal{C}(G)$ of convex subgroups of $G$ is totally ordered by inclusion, and that $H \longmapsto \Delta(H)$ is an order isomorphism between $\mathcal{C}(G)$ and the set of non-empty initial segments of $\Delta(G)$ ordered by inclusion.

We can also note that the cosets of some convex subgroup $H \leqslant G$ are all convex, because the translations are strictly increasing. So the quotient of $G$ by some convex subgroup $H$ is naturally endowed with a total order, as defined in definition 1.2.13.

**Definition 3.1.13.** Let $H$ be a subgroup of $G$, and $\delta \in \Delta(G)$. We denote by $H_{<\delta} = \{x \in H | \Delta(x) < \delta\}$. We introduce a similar notation for the conditions $\leqslant \delta$, $= \delta$... More generally, if $P \subseteq \Delta(G)$, we define $H_{\in P} = \{x \in H | \Delta(x) \in P\}$.

### 3.1.3 Outline of the proof

Although $\mathop{\downarrow}^{\mathbf{cut}}$ and $\mathop{\downarrow}^{\mathbf{inv}}$ are very combinatorial notions, we have to set up a rather technical algebraic machinery to prove that they coincide. Let us explain in this subsection what we do conceptually.

**Assumptions 3.1.14.** On top of assumptions 3.1.6, we assume that $M$ is $|B|^+$-saturated and strongly $|B|^+$-homogeneous.

We wish to show $c_1 \dots c_n \mathop{\downarrow}^{\mathbf{inv}}_A B$. For this, just as in subsection 1.2.2, we partition our family $c$ into smaller subfamilies, that we call "blocks", and we build a global $\mathrm{Aut}(M/A)$-invariant extension of the type of each block. We get a family of global $\mathrm{Aut}(M/A)$-invariant types, that we call (in definition 3.1.20) a "block extension", and what is left for us is to "glue" those types into a global $\mathrm{Aut}(M/A)$-invariant extension of the full type of $c$.

Subsection 1.2.2, involves weak orthogonality. There is a link between orthogonality and convex subgroups that appears in ([Men20], Proposition 4.5). We state a variant:

**Lemma 3.1.15.** *Let $G$ be some $\mathrm{Aut}(M/A)$-invariant convex subgroup of $M$, and $c, d \in M$. If $c \in (G + A) \not\ni d$, then $\mathrm{tp}(c/A)$ and $\mathrm{tp}(d/A)$ are weakly orthogonal.*

*Proof.* Suppose by contradiction that they are not. Then $\mathrm{tp}(d/A)$ is not complete in $S(Ac)$, i.e. $\mathrm{ct}(d/A)$ has a point in $A + \mathbb{Q} \cdot c$, and this point is in $A+G$. By strong homogeneity, there exists an automorhpism $\sigma \in \mathrm{Aut}(M/A)$ such that $\sigma(d) \in A + G$. As $G$ is $\mathrm{Aut}(M/A)$-invariant, we have $d \in A + G$, a contradiction. □

Note that lemma 3.1.15 holds in the more general context of a relatively divisible subgroup $G$ in the expansion of an Abelian group which is in some sense *algebraically bounded*, that is when the definable closure coincides with the relative divisible closure of the generated subgroup.

*Example* 3.1.16. Suppose $A$ is the lexicographical product:

$$\mathbb{Q} \times_{lex} \mathbb{Q} \times_{lex} \mathbb{Q} \leqslant \mathbb{R} \times_{lex} \mathbb{R} \times_{lex} \mathbb{R}$$

let $c_1 = (0, 0, \sqrt{2})$, $c_2 = (\sqrt{2}, 0, 0)$, $c_3 \in M$ such that:

$$c_3 \in \bigcap_{N>0} \left](0, N, 0), \left(\frac{1}{N}, 0, 0\right)\right[$$

let $G$ be the convex subgroup $\bigcap_{N>0} \left]\left(-\frac{1}{N}, 0, 0\right), \left(\frac{1}{N}, 0, 0\right)\right[$. Then we can show by lemma 3.1.15 that $\mathrm{tp}(c_2/A)$ is weakly orthogonal to both $\mathrm{tp}(c_1/A)$ and $\mathrm{tp}(c_3/A)$ using $G$.

Actually, we need not choose $G$ to show that the type of $c_1$ is weakly orthogonal to that of $c_2$, we may choose instead, for instance, the convex subgroup $\bigcap_{N>0} \left]\left(0, -\frac{1}{N}, 0\right), \left(0, \frac{1}{N}, 0\right)\right[$. In fact, this works for any $\mathrm{Aut}(M/A)$-invariant convex subgroup of $G$ which properly extends the convex subgroup $\bigcap_{N>0} \left]\left(0, 0, -\frac{1}{N}\right), \left(0, 0, \frac{1}{N}\right)\right[$.

The same cannot be said for $c_3$, $G$ seems to be the only valid witness for orthogonality. Indeed, $c_2$ belongs to $H = \bigcup_{N<\omega} [(-N, 0, 0), (N, 0, 0)]$, which is the least $\mathrm{Aut}(M/A)$-invariant convex subgroup properly extending $G$, whereas $c_3 \notin A + H'$, with $H' = \bigcup_{N<\omega} [(0, -N, 0), (0, N, 0)]$ the largest $\mathrm{Aut}(M/A)$-invariant proper convex subgroup of $G$ (note that $H'$ witnesses orthogonality between $c_3$ and $c_1$).

We conclude that the type of $c_2$ is weakly orthogonal to the two other types, but for very different reasons. As it is harder to prove weak orthogonality between $c_2$ and $c_3$, $c_2$ is, in some sense that will be made formal in the next section, more related to $c_3$ than it is to $c_1$.

Contrary to DLO, there seems to be "several layers" of weak orthogonality in DOAG. We understand them by defining several ways to partition our tuple, that refine each other. The finer the partition, the easier it becomes to find a block extension, the harder it is to "glue" the elements of the block extension.

In model theory, the standard way to "glue" global invariant types is via the tensor product:

**Definition 3.1.17.** Let $p, q$ be $\mathrm{Aut}(M/A)$-invariant global types. We define the *tensor product* of $p$ by $q$ to be the following $\mathrm{Aut}(M/A)$-invariant complete global type:

$$p(x) \otimes q(y) = \left\{\varphi(x,y,m) | m \in M,\ \exists a \vDash q_{|Am},\ p(x) \vDash \varphi(x,a,m)\right\}$$

The tensor product is associative, but not commutative in general. However, two types that are weakly orthogonal must necessarily commute.

**Fact 3.1.18.** *In* DOAG *(in fact in any distal theory, see [Sim13], Proposition 2.17), two global* $\mathrm{Aut}(M/A)$*-invariant types that commute must be weakly orthogonal.*

**Lemma 3.1.19.** *Let $F_1, F_2$ be closed subspaces of $S(A)$. Suppose $p_1$ and $p_2$ are weakly orthogonal for all $p_i \in F_i$. Then the set $F_1 \times F_2$, seen as a topological subspace of $S(A)$, is exactly the topological direct product $F_1 \times_{Top} F_2$, via the homeomorphism $h: p_1(x_1) \cup p_2(x_2) \longmapsto (p_1(x_1), p_2(x_2))$.*

*Proof.* The map $h$ is clearly a continuous bijection between compact separated spaces, thus it is a homeomorphism. ☐

The way we partition our tuple is with fibers of some ad-hoc valuations:

**Definition 3.1.20.** Let val be a $B$-valuation over $M$, as defined in definition 2.1.3. Given our family $c$, the val-*blocks* of $c$ are defined as the maximal subtuples of $c$ of elements of equal value. They form a partition of $c$.

By abuse of notation, if $c_i = (c_{ij})_j$ is a val-block of $c$, we define $\mathrm{val}(c_i) = \mathrm{val}(c_{ij})$, which does not depend on the choice of $j$.

A *weak* val-*block extension* of $c$ is a family $(p_i)_i$ of global $\mathrm{Aut}(M/A)$-invariant extensions of the types over $B$ of each val-block of $c$. Such a block extension is *strong* when $\bigcup_i p_i$ is consistent with $\mathrm{tp}(c/B)$.

In this paper, whenever we use the terminology of block extensions, it will always be for global $\mathrm{Aut}(M/A)$-invariant extensions of types over $B$, the parameter sets $A$, $B$ will not change.

In particular, if $\mathrm{val}'$ refines val, then the $\mathrm{val}'$-blocks of some family form a finer partition than its val-blocks.

The notion of separatedness from subsection 2.2.3 is very important in this chapter, for it is an algebraic way to state that the val-blocks of a tuple are independent from each other. The induced valued field will always be the trivially valued field $\mathbb{Q}$ here. One can note that $c$ is val-separated if and only if each of its val-blocks is. Moreover, a finite family that is separated with respect to some $A$-valuation must be a lift of some $\mathbb{Q}$-free family of $M/A$, i.e. a family in $M$ which maps via the canonical surjection to a $\mathbb{Q}$-free family in $M/A$. In particular, such a family is $\mathbb{Q}$-free.

*Remark* 3.1.21. The family $c$ is clearly $A$-interdefinable with the lift of a $\mathbb{Q}$-free family from $C/A$. Moreover, if $d$ is a tuple that is $A$-interdefinable with $c$, then one can easily show that we have $c \mathop{\downarrow}_A B$ if and only if $d \mathop{\downarrow}_A B$ for $\downarrow \in \{\downarrow^{\mathbf{cut}}, \downarrow^{\mathbf{inv}}\}$.

**Assumptions 3.1.22.** By remark 3.1.21, on top of assumptions 3.1.14, we may assume that $c$ is a lift of a $\mathbb{Q}$-free family from $C/A$.

*Remark* 3.1.23. Now, if $d$ is another lift of a $\mathbb{Q}$-free family from $C/A$, then $c$ and $d$ are $A$-interdefinable if and only if they have the same size $n$, and there exists $f \in GL_n(\mathbb{Q})$ such that the $i$-th component of $d$ is an $A$-translate of the $i$-th component of $f(c)$ for each $i$.

We just defined several ways to state that the blocks of some tuple are "independent" from each other and thus "easy to glue". The fact that the types of the blocks are weakly orthogonal is a model-theoretic way to define this idea, whereas the notion of separatedness with respect to some valuation is more algebraic. The core of our proofs relies on building valuations over $M$ for which these notions interact in a meaningful way.

More precisely, we build $D$-valuations $\mathrm{val}^i_D$ on $M$, with $i \in \{1, 2, 3\}$, and with $D \leqslant M$, such that:

(C1) $\mathrm{val}^3_D$ refines $\mathrm{val}^2_D$, which refines $\mathrm{val}^1_D$.

(C2) $\mathrm{val}^3_{\{0\}} = \Delta$ and $\mathrm{val}^2_{\{0\}}$ is trivial (it is refined by any valuation).

(C3) if $c \mathop{\downarrow}^{\mathbf{cut}}_A B$, then $c$ is $A$-interdefinable with a family $d$ that is simultaneously separated with respect to $\mathrm{val}^3_A$ and $\mathrm{val}^3_B$.

(C4) If $d$ is a $\mathrm{val}^3_B$-separated family, and $d \mathop{\downarrow}^{\mathbf{cut}}_A B$, then $d$ admits a strong $\mathrm{val}^3_B$-block extension.

(C5) Given $d$ a $\mathrm{val}^3_B$-separated family, and $(p_i)_i$ a strong $\mathrm{val}^3_B$-block extension of $d$, some tensor product of the $(p_i)_i$ is a global $\mathrm{Aut}(M/A)$-invariant extension of $\mathrm{tp}(d/B)$.

(C6) Suppose $d$ is $\mathrm{val}^2_D$-separated, and let $(b_i)_i$ be its $\mathrm{val}^2_D$-blocks. Then $(\mathrm{tp}(b_i/D))_i$ is weakly orthogonal. In particular, any weak $\mathrm{val}^2_B$-block extension is strong.

(C7) Given $d$ a $\mathrm{val}^1_B$-separated family, and $(p_i)_i$ a weak (thus strong by the previous point) $\mathrm{val}^1_B$-block extension of $d$, $(p_i)_i$ is weakly orthogonal, and thus they commute.

Now, if we could build such valuations, then it would follow that $c \mathop{\downarrow}^{\mathbf{inv}}_A B$, which would finish the proof of theorem 3.1.5. Indeed, if $d$ was given by (C3), then by remark 3.1.21 it would be enough to show that $d \mathop{\downarrow}^{\mathbf{inv}}_A B$. Now, again by remark 3.1.21, we would have $d \mathop{\downarrow}^{\mathbf{cut}}_A B$, so (C4) would give us a strong $\mathrm{val}^3_B$-block extension of $d$, and (C5) would give us a way to glue it into a witness of $d \mathop{\downarrow}^{\mathbf{inv}}_A B$, which would conclude the proof. Moreover, (C6), and (C7) would give us a fine understanding of the global $\mathrm{Aut}(M/A)$-invariant extensions of the type of $d$. In fact, we establish in section 3.3 an exhaustive classification of the Stone space of all these global $\mathrm{Aut}(M/A)$-invariant extensions.

Last but far from least, what remains for us to do is to actually define explicitly those valuations $\mathrm{val}^i_D$, and to show that all those nice properties hold.

### 3.1.4 Relation to Dolich-independence

We make short comments on how our result on DOAG relates to the work of Dolich ([Dol04]), their Section 8 in particular. Dolich defines in their paper a geometric independence notion in terms of “halfway-definable cells”, which we call $\mathop{\downarrow}^{\mathrm{Dolich}}$. It is shown that it coincides with $\mathop{\downarrow}^{\mathbf{f}}$ in any o-minimal expansion of RCF. To motivate their results, Dolich also gives five axioms for independence notions (allegedly from unpublished notes of Shelah), four of which are always satisfied by non-forking in any theory, and the fifth, called "Chain Condition", is a weakening of the independence theorem. They also define an independence notion called "non-1-dividing", with a combinatorial definition which strengthens that of non-forking ; let us write it $\mathop{\smile\hskip-0.9em\downarrow}^{\mathbf{1\text{-}d}}$. They claim that, in case some independence relation satisfies the five axioms in a given theory, non-1-dividing is the weakest relation satisfying those axioms, and they show that $\mathop{\smile\hskip-0.9em\downarrow}^{\mathrm{Dolich}}$ satisfies the five axioms in any o-minimal theory. In particular, we have $\mathop{\smile\hskip-0.9em\downarrow}^{\mathrm{Dolich}} \subseteq \mathop{\smile\hskip-0.9em\downarrow}^{\mathbf{1\text{-}d}} \subseteq \mathop{\smile\hskip-0.9em\downarrow}^{\mathbf{f}}$ in any o-minimal theory, and those inclusions are equalities in expansions of RCF.

Note that, if we generalize the independence notion $C \mathop{\smile\hskip-0.9em\downarrow}^{\mathbf{cut}}_A B$ to any o-minimal theory by the definition: "any closed bounded interval with bounds in $\mathrm{dcl}(AB)$ having a point in $\mathrm{dcl}(AC)$ already has a point in $\mathrm{dcl}(A)$", then it follows from corollary 1.2.12 that $\mathop{\smile\hskip-0.9em\downarrow}^{\mathbf{d}} \subseteq \mathop{\smile\hskip-0.9em\downarrow}^{\mathbf{cut}}$, hence in particular $\mathop{\smile\hskip-0.9em\downarrow}^{\mathrm{Dolich}} \subseteq \mathop{\smile\hskip-0.9em\downarrow}^{\mathbf{cut}}$. This inclusion is strict in general, for in RCF, one may find an example where $C \mathop{\smile\hskip-0.9em\downarrow}^{\mathbf{cut}}_A B$ holds, while $C \mathop{\smile\hskip-0.9em\downarrow}^{\mathbf{alg}}_A B$ fails, which implies $C \not\mathop{\smile\hskip-0.9em\downarrow}^{\mathrm{Dolich}}_A B$. So the main difference between the geometric independence notion introduced by Dolich and ours is that ours is easily shown to be weaker than non-forking in the general case, and the difficulties come when we prove that it is actually as strong as non-forking in DOAG, while $\mathop{\smile\hskip-0.9em\downarrow}^{\mathrm{Dolich}}$ is clearly stronger than non-forking, and it is difficult to prove the other direction in RCF.

While we do not know whether $\mathop{\smile\hskip-0.9em\downarrow}^{\mathrm{Dolich}} = \mathop{\smile\hskip-0.9em\downarrow}^{\mathbf{f}}$ in DOAG (it would be a strengthening of our result), we remark that $\mathop{\smile\hskip-0.9em\downarrow}^{\mathbf{1\text{-}d}} = \mathop{\smile\hskip-0.9em\downarrow}^{\mathbf{f}}$ in any theory (DOAG in particular) where $\mathop{\smile\hskip-0.9em\downarrow}^{\mathbf{f}} = \mathop{\smile\hskip-0.9em\downarrow}^{\mathbf{inv}}$. Indeed, one may easily show that $\mathop{\smile\hskip-0.9em\downarrow}^{\mathbf{inv}}$ satisfies the Chain Condition, using the fact that it satisfies Extension. If $\mathop{\smile\hskip-0.9em\downarrow}^{\mathbf{inv}} = \mathop{\smile\hskip-0.9em\downarrow}^{\mathbf{f}}$, then it follows that $\mathop{\smile\hskip-0.9em\downarrow}^{\mathbf{f}}$ satisfies the five axioms. By maximality of $\mathop{\smile\hskip-0.9em\downarrow}^{\mathbf{1\text{-}d}}$, we have $\mathop{\smile\hskip-0.9em\downarrow}^{\mathbf{f}} \subseteq \mathop{\smile\hskip-0.9em\downarrow}^{\mathbf{1\text{-}d}}$, and the other inclusion always holds.

## 3.2 How to glue types via ad-hoc valuations

In this section, we define the valuations $\mathrm{val}^i_A$ ($i \in \{1,2,3\}$), and we show that they satisfy the gluing properties (C5), (C6), (C7). We also give a classification of the cuts over some parameter set.

In this section, we work with assumptions 3.1.14.

### 3.2.1 Basic definitions and classification of the cuts

**Lemma 3.2.1.** *Let $c$, $d \in M$, and $a \in A$, such that* $\mathrm{ct}(c/A) = \mathrm{ct}(d/A)$. *Then* $\mathrm{ct}(c + a/A) = \mathrm{ct}(d + a/A)$.

*Proof.* If not, then there exists $a' \in A$ which lies strictly between $c + a$ and $d + a$.

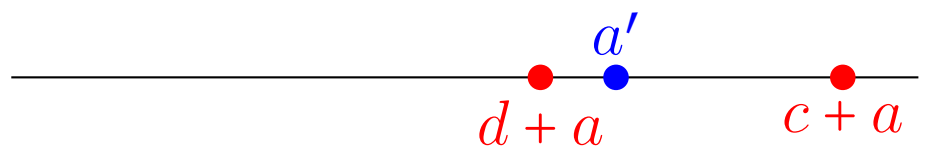

then $a' - a \in A$ lies strictly between $c$ and $d$, contradicting the hypothesis. □

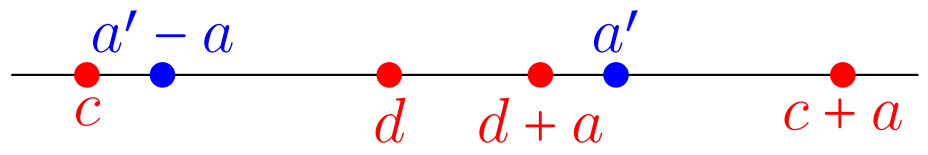

We hope the figures make the proofs easier to understand. However, we do not want them to be misleading, thus we would like to say that they may not cover all the possible cases. For instance, if $a$ was negative, then the correct picture would be reversed.

**Definition 3.2.2.** By lemma 3.2.1 (and by saturation), $A$ acts by translation over the set of all the cuts of $M$ over $A$. For $c \in M$, denote by $\mathrm{Stab}(c/A)$ the stabilizer of $\mathrm{ct}(c/A)$.

Note that such a stabilizer is only a subgroup of $A$, which is clearly convex in $A$, but not in $M$.

**Definition 3.2.3.** Let $d \in M$. If $d \notin A$, then we define:

$$G(d/A) = \cap\{]-|a|, |a|[\ :\ a \in A \setminus \mathrm{Stab}(d/A)\}$$

else we define $G(d/A) = \{0\}$. Either way, we also define:

$$H(d/A) = \cup\{[-|a|, |a|]\ :\ a \in \mathrm{Stab}(d/A)\}$$

We view $G(d/A)$ as an $A$-type-definable set, and $H(d/A)$ as an $A$-$\vee$-definable set. They have the same points in $A$, however they do not have the same points in $M$ when $d \notin A$. They are $A$-(type/$\vee$)-definable convex subgroups.

*Example* 3.2.4. In example 3.1.16, $G(c_2/A) = G(c_3/A) = G$, and $H(c_2/A) = H(c_3/A) = H'$. However, $G(c_1/A)$ is the group of elements which are infinitesimal with respect to $A$, and it is distinct from $G$. As for $H(c_1/A)$, it is trivial.

Likewise, for arbitrary $A$, if $\Delta(0) < \Delta(d) < \Delta(A \smallsetminus \{0\})$ (i.e. $d$ is infinitesimal with respect to $A$), then $H(d/A)$ is trivial, and $G(d/A)$ is the type-definable group of elements which are infinitesimal with respect to $A$.

If $\Delta(d) > \Delta(A)$, then $H(d/A)$ is the convex subgroup generated by $A$, whereas $G(d/A)$ is the whole group.

*Remark* 3.2.5. We always have $G(d/A) \geqslant H(d/A)$. By definition, $H(d/A)$ is the convex subgroup generated by $\mathrm{Stab}(d/A)$, i.e. the least $A$-$\vee$-definable convex subgroup containing $\mathrm{Stab}(d/A)$. In case $d \notin A$, $G(d/A)$ is the largest $A$-type-definable convex subgroup disjoint from $A \smallsetminus \mathrm{Stab}(d/A)$. For all $d$, $d' \in M$, we have in fact :

$$G(d/A) = G(d'/A) \Longrightarrow H(d/A) = H(d'/A)$$

the only case where the implication is not an equivalence being when one point is in $A$, and the other is infinitesimal with respect to $A$. We also have:

$$\mathrm{Stab}(d/A) = \mathrm{Stab}(d'/A) \Longleftrightarrow H(d/A) = H(d'/A).$$

**Proposition 3.2.6.** *Let $d, d' \in M$. If $G(d/A) < G(d'/A)$, then we have $G(d/A) < H(d'/A)$.*

*Proof.* By definition of $G(d/A)$, there must exist $a \in A \smallsetminus \mathrm{Stab}(d/A)$ such that $a \in G(d'/A)$. By remark 3.2.5, $a \in H(d'/A) \smallsetminus G(d/A)$. □

Let us now build $\mathrm{val}^1_A$.

**Lemma 3.2.7.** *Let $d_1, d_2 \in M$. Then $G(d_1 + d_2/A) \leqslant \max(G(d_1/A), G(d_2/A))$ with respect to inclusion.*

*Proof.* Suppose by contradiction $G(d_1 + d_2/A) > \max(G(d_1/A), G(d_2/A))$. Then, there must exist $a \in \mathrm{Stab}(d_1 + d_2/A)$ such that $a \notin \mathrm{Stab}(d_i/A)$. As $\mathrm{Stab}(d_i/A)$ is a convex subgroup of $A$, we also have $\frac{a}{2} \notin \mathrm{Stab}(d_i/A)$, so there exists $a_i \in A$ which lies strictly between $d_i$ and $d_i + \frac{a}{2}$.

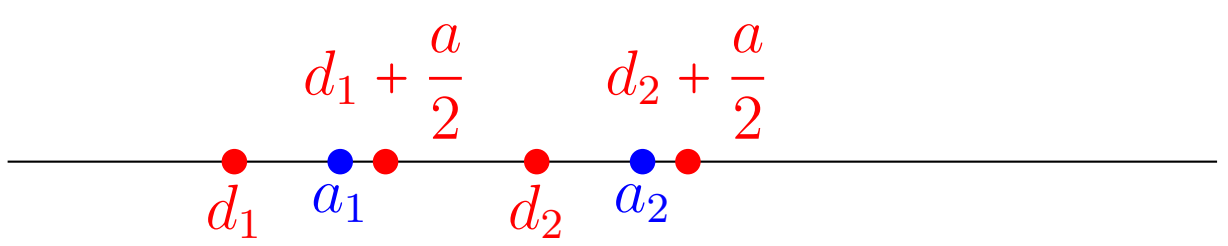

then $a_1 + a_2$ lies strictly between $d_1 + d_2$ and $d_1 + d_2 + a$, contradicting the fact that $a \in \mathrm{Stab}(d_1 + d_2/A)$. □

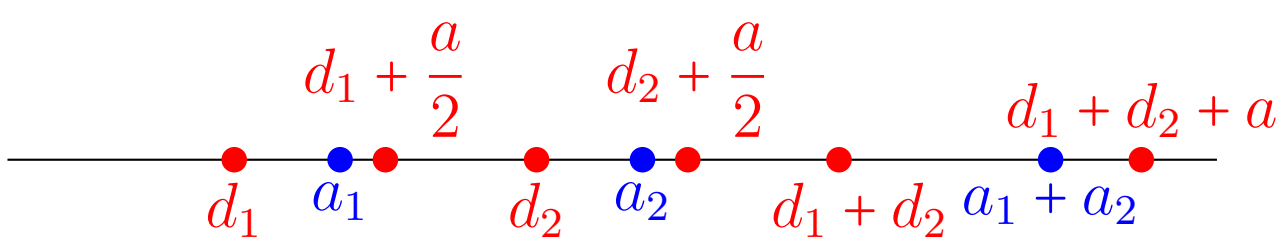

**Definition 3.2.8.** By lemma 3.2.7, the map $x \longmapsto G(x/A)$ is an $A$-valuation, so we set $\mathrm{val}^1_A(x) = G(x/A)$.

Note that, by definition, we always have $G(d/B) \leqslant G(d/A)$ for all $d \in M$. Furthermore:

**Proposition 3.2.9.** *For all $d \in M$, if $d \mathop{\downarrow}^{\mathbf{cut}}_{A} B$, then $H(d/B) \geqslant H(d/A)$.*

*Proof.* Suppose $H(d/B) < H(d/A)$. Then there exists $a \in \mathrm{Stab}(d/A)$ such that $a \notin \mathrm{Stab}(d/B)$. Thus, we can find $b \in B$ which lies strictly between $d$ and $d + a$.

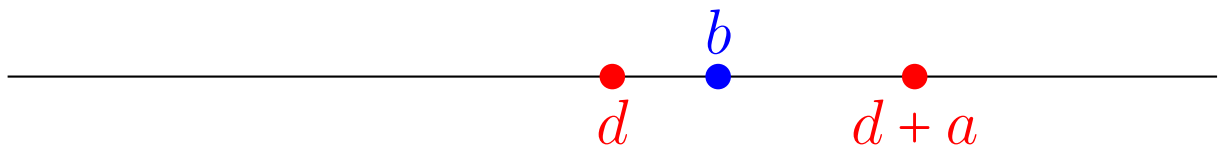

on one hand, $b - a$ lies strictly between $d$ and $d - a$.

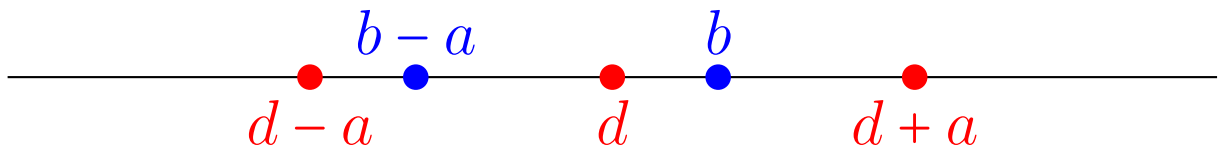

As $\mathrm{Stab}(d/A)$ is a group, we have on the other hand $-a \in \mathrm{Stab}(d/A)$, i.e. $\mathrm{ct}(d + a/A) = \mathrm{ct}(d/A) = \mathrm{ct}(d - a/A)$. In particular, no point from $A$ lies in the closed interval with bounds $b$ and $b - a$. On the other hand, $d$ lies strictly between $b$ and $b - a$, which implies $d \mathop{\not\downarrow}^{\mathbf{cut}}_{A} B$. □

In the cut-independent setting, one may intuitively see $G(c/A)$ as a distance that approximates the position of $c$ "from the top", while $H(c/A)$ approximates $c$ "from the bottom". In fact, $G(c/A)$ is the analogue of the radius of the elements of $I_A$ from remark 2.1.7. As for $H(c/A)$, it does not really correspond to $R_A$, but to $R_M$ for the global $\mathrm{Aut}(M/A)$-invariant extension. Then, when we go to a larger parameter set $B$, we get a finer approximation.

**Proposition 3.2.10.** *Let $a \in A$, and $d \in M$. Then the following conditions are equivalent:*

- $d - a \in G(d/A)$ *and* $d \notin A$.
- $\Delta(d-a) \notin \Delta(A)$.

*Proof.* Suppose $d-a \in G(d/A)$. If we had $\Delta(d-a) \in \Delta(A)$, then there would be $a' \in A$ such that $\Delta(d-a) = \Delta(a')$. This would imply $a' \in (G(d/A) \cap A) = \mathrm{Stab}(d/A)$. Since, for some $n$, $d - n|a'| \leqslant a \leqslant d + n|a'|$, as $na' \in \mathrm{Stab}(d/A)$, it follows that $\mathrm{ct}(d/A) = \mathrm{ct}(a/A)$, thus $d = a$, proving the first direction.

Conversely, suppose $\Delta(d-a) \notin \Delta(A)$. By definition of $G(d/A)$, in order to show that $a \in (d \bmod G(d/A))$, it suffices to show that $a \in ]d - |a'|, d + |a'|[$ (i.e. $|d-a| \leqslant |a'|$) for every $a' \in A \setminus \mathrm{Stab}(d/A)$. Suppose by contradiction this fails for some $a'$. Then $\Delta(d-a) \geqslant \Delta(a')$, and this inequality is strict by the hypothesis. As $a' \notin \mathrm{Stab}(d/A)$, let $a'' \in A$ which lies strictly between $d$ and $d + a'$.

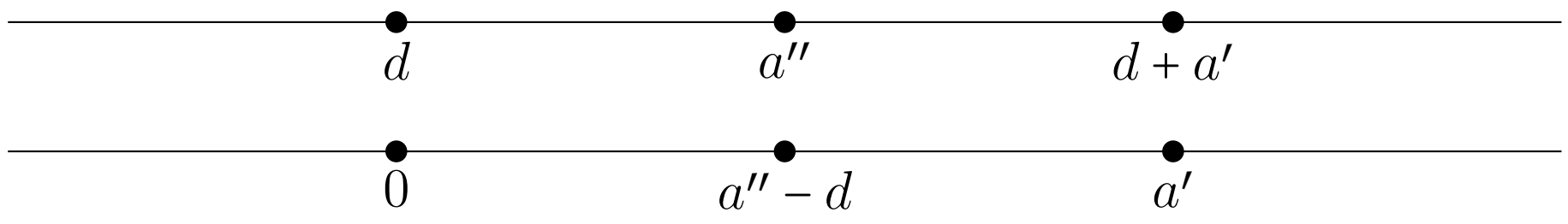

Then $|d - a''| \leqslant |a'|$, thus $\Delta(d - a'') < \Delta(d-a)$. We then apply remark 2.1.4 to get $\Delta(d-a) = \Delta(d - a'' - d + a) = \Delta(a - a'') \in \Delta(A)$, a contradiction. □

**Definition 3.2.11.** If $a$ satisfies the conditions of proposition 3.2.10, then we say that $a$ is a *ramifier of $d$ over $A$*. We write $\mathrm{ram}(d/A)$ for the set of those ramifiers, and we say that $d$ is *ramified over $A$* if this set is non-empty. We say that $d$ is *Archimedean over $A$* whenever it is not ramified over $A$.

In example 3.1.16, $c_1$ and $c_2$ are Archimedean over $A$, $c_3$ is ramified over $A$, and 0 is a ramifier.

*Remark* 3.2.12. If $d \in M$, then $d$ is ramified over $A$ if and only if it satisfies the following equivalent conditions:

1. The coset $d \bmod G(d/A)$ has a point in $A$, and $d \notin A$.

2. $\Delta(A+\mathbb{Q}d) \supsetneq \Delta(A)$.

Moreover, as $A \cap G(d/A) = \mathrm{Stab}(d/A)$, $\mathrm{ram}(d/A)$ is a coset from the quotient ${}^{A}\!/\!{}_{\mathrm{Stab}(d/A)}$, thus all its elements lie in the same coset $\mod H(d/A)$. It turns out that $d$ does not belong to that coset:

**Lemma 3.2.13.** *Let $d \in M \setminus A$ and $a \in A$. Then $d - a \notin H(d/A)$.*

*Proof.* Suppose by contradiction that there exists $a' \in \mathrm{Stab}(d/A)$ such that $|d-a| \leqslant |a'|$. Then $a \in \mathrm{ram}(d/A)$, thus $\Delta(d-a) \notin \Delta(A)$, and in particular $\Delta(d-a) < \Delta(a')$. By remark 2.1.4, $\Delta(d-(a+a')) = \Delta(a') \in \Delta(A)$, thus $a+a' \notin \mathrm{ram}(d/A)$. This contradicts the fact that $a \in \mathrm{ram}(d/A)$, and the fact that $a' \in \mathrm{Stab}(d/A)$. ☐

**Corollary 3.2.14.** *Let $d \in M$ be ramified over $A$, as well as two elements $a$ and $a'$ in $\mathrm{ram}(d/A)$. Then $\Delta(d-a) = \Delta(d-a')$.*

*Proof.* We have $d - a \notin H(d/A) \ni a - a'$, thus $\Delta(d-a) > \Delta(a-a')$, and by remark 2.1.4 we have $\Delta(d-a) = \Delta(d-a+a-a') = \Delta(d-a')$. ☐

**Definition 3.2.15.** Let $d \in M$ be ramified over $A$. We define $\delta(d/A)$ as $\Delta(d-a)$ for some $a \in \mathrm{ram}(d/A)$. This definition does not depend on the choice of the ramifier $a$ by corollary 3.2.14.

Note that $\delta(d/A)$ is the unique Archimedean class of $\Delta(A+\mathbb{Q}\cdot d) \setminus \Delta(A)$.

*Remark* 3.2.16. Let $d$ be ramified over $A$, and $a \in \mathrm{ram}(d/A)$. Then, as $G(d/A)$ and $H(d/A)$ have the same points in $A$, one can easily compute $\mathrm{ct}(\delta/\Delta(A))$. A convex subgroup is induced by its Archimedean classes, so by definition of the $\vee$-definable convex subgroup $H = H(d/A)$, for any model $N$ containing $A$, $\Delta(H(N))$ is the least initial segment of $\Delta(N)$ containing $\Delta(\mathrm{Stab}(d/A))$. Now, as $\delta \in \Delta(G(d/A)) \setminus \Delta(H)$, we clearly have $\mathrm{ct}(\delta/\Delta(A)) = \mathrm{ct}_{>}(\Delta(\mathrm{Stab}(d/A))/\Delta(A))$.

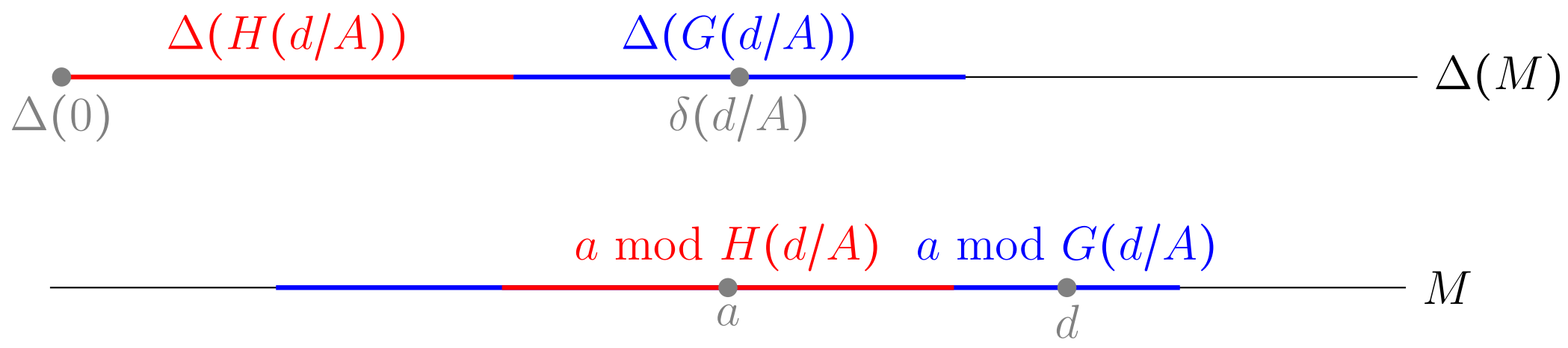

Many notions that we manipulate behave very differently with ramified and Archimedean points, it will often lead to case disjunctions. For instance, the following statements show us that cuts do not look the same for Archimedean and ramified points:

**Lemma 3.2.17.** *Let $d \in M$. Then* $\mathrm{ct}(d/A) \subseteq (d \bmod G(d/A))$ *(as type-definable sets).*

*Proof.* By contraposition, let $d' \in M$ such that $d - d' \notin G(d/A)$. Then, by definition, there exists $a \in A \setminus \mathrm{Stab}(d/A)$ such that $|d - d'| \geqslant |a|$. Let $a' \in A$ which lies strictly between $d$ and $d + a$. As either $d + a$ or $d - a$ lies strictly between $d$ and $d'$, either $a'$ or $a' - a$ lies strictly between $d$ and $d'$, thus $\mathrm{ct}(d/A) \neq \mathrm{ct}(d'/A)$. ☐

**Proposition 3.2.18.** *Let $d \in M$ be Archimedean over $A$. Then the $A$-type-definable set* $\mathrm{ct}(d/A)$ *coincides with* $d \bmod G(d/A)$.

*Proof.* Suppose $d \notin A$, else the proposition is trivial.

We have the inclusion $d \bmod G(d/A) \subseteq \mathrm{ct}(d/A)$, as $d \bmod G(d/A)$ is a convex set containing $d$ and disjoint from $A$. The other direction follows from lemma 3.2.17. ☐

**Proposition 3.2.19.** *Let $d \in M$ be ramified over $A$. Then* $\mathrm{ct}(d/A)$ *is an $A$-translate of one of the two connected components of $G(d/A) \setminus H(d/A)$.*

These translates are exactly the two blue segments in the last figure of remark 3.2.16. The connected components of some set will always refer to its maximal convex subsets.

*Proof.* Let $a \in \mathrm{ram}(d/A)$. By lemma 3.2.17, we have:

$$\mathrm{ct}(d/A) \subseteq d \bmod G(d/A) = a \bmod G(d/A)$$

By lemma 3.2.13, $d \notin (a \bmod H(d/A))$, which is the convex closure of the set $\mathrm{ram}(d/A)$, which is a subset of $A$, thus $a \bmod H(d/A)$ must be disjoint from $\mathrm{ct}(d/A)$. It follows that $\mathrm{ct}(d/A)$ is included in the connected component of $(a + G(d/A)) \setminus (a + H(d/A))$ containing $d$.

Let $\mathcal{C}$ be that connected component, and $d' \in \mathcal{C}$. There does not exist $a' \in A$ lying between $d$ and $d'$, for such $a'$ would be in $(a + G(d/A)) \cap A = \mathrm{ram}(d/A) \subseteq a + H(d/A)$, and $a'$ would also belong to $\mathcal{C}$ as $\mathcal{C}$ is convex, contradicting the fact that $\mathcal{C}$ is disjoint from $a + H(d/A)$. ☐

**Corollary 3.2.20.** *With the hypothesis from proposition 3.2.19, if $a$ is in* $\mathrm{ram}(d/A)$*, and $d'$ is another point from $M$, then $d' \equiv_A d$ if and only if the following conditions hold:*

- $\mathrm{ct}(\Delta(d'-a)/\Delta(A)) = \mathrm{ct}(\Delta(d-a)/\Delta(A))$ *(in particular, $\Delta(d'-a)$ is not in $\Delta(A)$).*
- $d' < a \Longleftrightarrow d < a$.

*Proof.* The first condition is equivalent to $d'-a \in G(d/A) \smallsetminus H(d/A)$. The second condition ensures that $d'$ lies in the correct connected component. □

The two statements proposition 3.2.18 and proposition 3.2.19 yield a simple classification of the cuts over $A$. In order to consider non-forking extensions of types, we now have to deal with ways to refine those cuts to cuts over a larger parameter set $B$. As we classified the cuts over $A$ by looking at $\mathrm{Aut}(M/A)$-invariant convex subgroups, we have to compute what are the $\mathrm{Aut}(M/B)$-invariant convex subgroup corresponding to their refinements.

**Definition 3.2.21.** If $G$ is an $A$-type-definable convex subgroup, then we define the following $A$-$\vee$-definable convex subgroup:

$$\underline{G}_A = \bigcup_{a \in A \cap G} [-|a|, |a|]$$

If $H$ is an $A$-$\vee$-definable convex subgroup, then we define the following $A$-type-definable convex subgroup:

$$\overline{H}^A = \bigcap_{a \in A \smallsetminus H} ]-|a|, |a|[$$

*Remark* 3.2.22. Let $d \in M$. We have:

$$H(d/A) \leqslant \overline{H(d/A)}^M \leqslant \overline{H(d/A)}^B \leqslant \underline{G(d/A)}_B \leqslant \underline{G(d/A)}_M \leqslant G(d/A)$$

and :

$$\underline{G(d/A)}_A = H(d/A), \ \overline{H(d/A)}^A = G(d/A)$$

moreover, if $d \mathop{\downarrow}_A^{\mathbf{cut}} B$, then by proposition 3.2.9 we recall:

$$H(d/A) \leqslant H(d/B) \leqslant G(d/B) \leqslant G(d/A).$$

### 3.2.2 Weak orthogonality

Now we start to deal with global invariant extensions, in order to prove (C7).

**Lemma 3.2.23.** *Let $d \in M \setminus A$ be Archimedean over $A$, with $d \mathop{\downarrow}^{\mathbf{cut}}_{A} B$. Let $p$ be a global* $\mathrm{Aut}(M/A)$*-invariant extension of* $\mathrm{tp}(d/B)$*. Then $G(d/B) = G(d/A)$ (as $B$-type-definable sets), and $p(x) \models x - d \in G(d/A) \setminus \underline{G(d/A)}_M$.*

Note that the statement holds even if the element $d$ is ramified over $B$.

*Proof.* Suppose by contradiction $G(d/B) < G(d/A)$. Then there exists $b$ in $B \setminus \mathrm{Stab}(d/B)$ such that $b \in G(d/A)$. Let $b_1 \in B$ which lies strictly between $d$ and $d + b$. Like wise, as $-b$ is also in $G(d/A) \setminus \mathrm{Stab}(d/B)$, we find $b_2 \in B$ which lies strictly between $d$ and $d - b$. Then the closed interval with bounds $b_1, b_2$ contains $d$ and is strictly included in $d$ mod $G(d/A) = \mathrm{ct}(d/A)$, which contradicts cut-independence.

By proposition 3.2.18, as $p$ extends $\mathrm{tp}(d/A)$, we necessarily have:

$$p(x) \models x - d \in G(d/A)$$

Let $c \models p$ in some elementary extension. As $p$ is $\mathrm{Aut}(M/A)$-invariant, we have $c \mathop{\downarrow}^{\mathbf{cut}}_{A} M$, thus $G(c/M) = G(d/A)$ by the above paragraph ($c$ is of course Archimedean over $A$). As $c - d \in G(c/M)$, and $d \in M \not\ni c$, we have $d \in \mathrm{ram}(c/M)$, thus $c - d \notin H(c/M) = \underline{G(c/M)}_M = \underline{G(d/A)}_M$. We conclude that $p(x) \models x - d \notin \underline{G(d/A)}_M$. □

**Lemma 3.2.24.** *Let $d \in M$ be ramified over $A$, with $d \mathop{\downarrow}^{\mathbf{cut}}_{A} B$, and let $a \in \mathrm{ram}(d/A)$. Then $\Delta(d - a) \notin \Delta(B)$.*

*Moreover, if $p$ is a global* $\mathrm{Aut}(M/A)$*-invariant extension of* $\mathrm{tp}(d/B)$*, then, either ($p(x) \models x - a \in G(d/A) \setminus \underline{G(d/A)}_M$ and $G(d/B) = G(d/A)$) or ($p(x) \models x - a \in \overline{H(d/A)}^M \setminus H(d/A)$ and $H(d/B) = H(d/A)$).*

*Proof.* Suppose we have $\Delta(d - a) = \Delta(b)$ for some $b \in B$. Let $n > 0$ such that $\frac{1}{n}|b| \leqslant |d - a| \leqslant n|b|$. Then we have $d \in \left[a - n|b|, a - \frac{1}{n}|b|\right] \cup \left[a + \frac{1}{n}b, a + n|b|\right]$. As $b \in G(d/A) \setminus H(d/A)$, none of those two intervals has a point in $A$, which contradicts $d \mathop{\downarrow}^{\mathbf{cut}}_{A} B$.

As $d$ is ramified over $A$, and $p$ extends $\mathrm{tp}(d/A)$, it follows from proposition 3.2.19 that $p(x) \models x - a \in G(d/A) \setminus H(d/A)$.

Suppose by contradiction $p(x) \vDash x - a \in \underline{G(d/A)}_M \setminus \overline{H(d/A)}^M$. Then there exists $m_1, m_2 \in M$ such that $m_1 \notin H(d/A)$, $m_2 \in G(d/A)$, and:

$$p(x) \vDash |m_1| \leqslant |x - a| \leqslant |m_2|$$

i.e. $p(x) \vDash x \in [a - |m_2|, a - |m_1|] \cup [a + |m_1|, a + |m_2|]$. As $m_1 \notin H(d/A)$ and $m_2 \in G(d/A)$, each of those two intervals is included in each connected component of $a + (G(d/A) \setminus H(d/A))$. As a result, none of those intervals belongs to $p$ by corollary 1.2.12, a contradiction.

By lemma 3.2.23, and by remark 3.2.22, if we had $p(x) \vDash x - a \in \underline{G(d/A)}_B$ (resp. $p(x) \vDash x - a \notin \overline{H(d/A)}^B$), then we would have:

$$p(x) \vDash x - a \in \overline{H(d/A)}^M \setminus H(d/A)$$

and respectively: $p(x) \vDash x - a \in G(d/A) \setminus \underline{G(d/A)}_M$. To conclude, it now suffices to show the two following properties:

1. If $G(d/B) < G(d/A)$, then $\mathrm{tp}(d/B) \vDash x - a \in \underline{G(d/A)}_B$.
2. If $H(d/B) > H(d/A)$, then $\mathrm{tp}(d/B) \vDash x - a \notin \overline{H(d/A)}^B$.

Let us proceed:

1. By hypothesis, let $b \in B \setminus \mathrm{Stab}(d/B)$ such that $b \in G(d/A)$. Let $b' \in B$ which lies strictly between $d$ and $d + b$.

   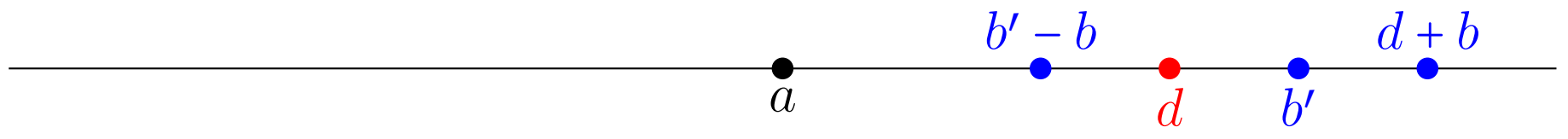
   

   then the interval $I$ with bounds $b' - a$, $b' - b - a$ is included in $G(d/A)$, and $\mathrm{tp}(d/B) \vDash x - a \in I$.

   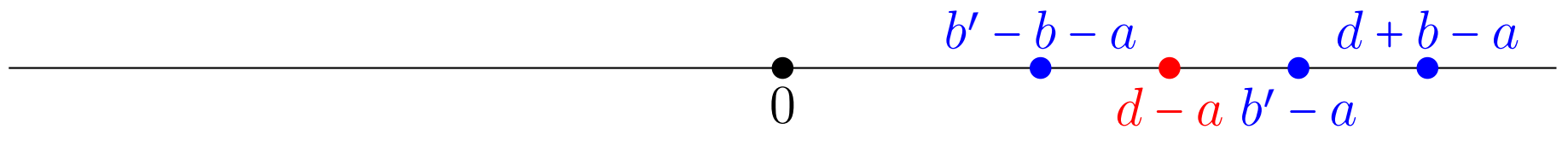
   

2. By hypothesis, let $b \in \mathrm{Stab}(d/B)$ such that $b \notin H(d/A)$. Then we have $\mathrm{ct}(d - b/B) = \mathrm{ct}(d/B) = \mathrm{ct}(d + b/B)$, thus $a \notin [d - |b|, d + |b|]$. As a result, $\mathrm{tp}(d/B) \vDash x - a \notin [-|b|, |b|]$. □

The two statements lemma 3.2.23, lemma 3.2.24 imply:

**Corollary 3.2.25.** *Let $d \in M$ such that $d \downarrow^{\mathbf{cut}}_A B$. Let $p \in S(M)$ be a global* $\mathrm{Aut}(M/A)$*-invariant extension of* $\mathrm{tp}(d/B)$*. Then there exists $f$ an $A$-definable map such that one of the three following conditions hold:*

- *$d \in A$, $f = id$.*
- *The cut over $M$ induced by $p$ is one of the two connected components of the translation by $f(d)$ of $G(d/A) \smallsetminus \underline{G(d/A)}_M$. In that case, $G(d/B) = G(d/A)$.*
- *The cut over $M$ induced by $p$ is one of the two connected components of the translation by $f(d)$ of $\overline{H(d/A)}^M \smallsetminus H(d/A)$. In that case, $H(d/B) = H(d/A)$.*

*Proof.* If $d \in A$, then the statement is trivial. If $d \notin A$ and $d$ is Archimedean over $A$, then the second condition holds with $f = id$, as stated in lemma 3.2.23. If $d$ is ramified over $A$, then either the second or the third condition holds with $f : x \longmapsto x - a$, by lemma 3.2.24. □

*Remark* 3.2.26*.* If $d \notin A$, while only one of the conditions from corollary 3.2.25 holds, we might still have $G(d/B) = G(d/A)$ and $H(d/B) = H(d/A)$ at the same time. For instance, this is the case if $\Delta(d) < \Delta(B \smallsetminus \{0\})$, and $B$ has no element that is infinitesimal compared to $A$.

In fact, one can see that $p$ is the type of an $M$-ramified point which adds an Archimedean class $\delta$ which is cut-independent over $\Delta(A)$: the two conditions of the statement expresses that $\delta$ either leans left or right with respect to $\Delta(A)$ and $\Delta(M)$.

*Remark* 3.2.27*.* For all $d \in M$, recall from remark 3.2.22 that $G(d/B) = \overline{H(d/B)}^B$. As a result, if $H(d/B) = H(d/A)$, then $G(d/B) = \overline{H(d/A)}^B$. In particular we see that, in either case of corollary 3.2.25, $p(x)$ implies $x - f(d) \in G(d/B)$. Likewise:

**Lemma 3.2.28.** *Let $d \in M$ such that $d \downarrow^{\mathbf{cut}}_A B$, and $d \notin A$. Let $p \in S(M)$ be a global* $\mathrm{Aut}(M/A)$*-invariant extension of* $\mathrm{tp}(d/B)$*. Then $p(x)$ implies $x \notin (M + H(d/B))$.*

*Proof.* As $d \notin A$, the first condition of corollary 3.2.25 fails. Let $f$ be the $A$-definable map witnessing corollary 3.2.25.

As $G(d/A) \geqslant G(d/B)$, we have:

$$\underline{G(d/A)}_M \geqslant \underline{G(d/B)}_M \geqslant H(d/B)$$

If the second condition holds, let $G = G(d/B), H = \underline{G(d/B)}_M$, else let $G = \overline{H(d/B)}^M, H = H(d/B)$. Either way, $G$ and $H$ have the same points in $M$, $p(x) \models x - f(d) \in G \setminus H$, and $H \geqslant H(d/B)$. Suppose by contradiction that there is $m \in M$ such that $p(x) \models x - m \in H(d/B)$. Then, as $p(x) \models x - f(d) \in G$, we have $m - f(d) \in G$. However, $m - f(d) \in M$, therefore $m - f(d) \in H$. This is a contradiction, as we cannot have at the same time $p(x) \models x - f(d) \notin H$, $p(x) \models x - m \in H$ and $f(d) - m \in H$. □

We established that realizations of unary global invariant types belong to some cosets of convex subgroups, and not others. As in lemma 3.1.15, we may relate those statements to orthogonality, in order to prove (C7).

**Lemma 3.2.29.** *Let $d_1, \dots, d_n$ be a family of tuples of $M$ such that, for all $i$,* $\mathrm{tp}(d_i/A)$ *and* $\mathrm{tp}(d_{<i}/A)$ *are weakly orthogonal. Then the family of types* $(\mathrm{tp}(d_1/A), \dots, \mathrm{tp}(d_n/A))$ *is weakly orthogonal with respect to definition 1.2.16.*

*Proof.* We show by induction on $i$ that $(\mathrm{tp}(d_1/A), \dots, \mathrm{tp}(d_i/A))$ is weakly orthogonal. The statement is trivial if $i = 1$. By induction hypothesis, we have $\bigcup_{j<i} \mathrm{tp}(d_j/A) \models \mathrm{tp}(d_{<i}/A)$. By hypothesis, $\mathrm{tp}(d_i/A) \cup \mathrm{tp}(d_{<i}/A) \models \mathrm{tp}(d_{\leqslant i}/A)$, and we conclude by induction. □

**Lemma 3.2.30.** *Let $d_1, d_2$ be two tuples from $M$. Suppose for all $f \in \mathrm{LC}(\mathbb{Q})$,* $\mathrm{tp}(f(d_1)/A)$ *is weakly orthogonal to* $\mathrm{tp}(d_2/A)$*. Then* $\mathrm{tp}(d_1/A)$ *is weakly orthogonal to* $\mathrm{tp}(d_2/A)$*.*

*Proof.* Let $d \equiv_A d_1$. It suffices to show that $d \equiv_{Ad_2} d_1$. By quantifier elimination, this holds if and only if, for every $f \in \mathrm{LC}(\mathbb{Q})$, we have $f(d) \equiv_{Ad_2} f(d_1)$, which is exactly the hypothesis of the statement. □

**Lemma 3.2.31.** *Let $d_1 \in M$ be a singleton, and let $d_2 \in M$ be a tuple. Let $V$ be the $\mathbb{Q}$-vector space* $\mathrm{dcl}(Ad_2)$*, and let $G$ be some* $\mathrm{Aut}(M/A)$*-invariant convex subgroup of $M$. Suppose we have $V \subseteq (A + G) \not\ni d_1$. Then* $\mathrm{tp}(d_1/A)$ *and* $\mathrm{tp}(d_2/A)$ *are weakly orthogonal.*

*Proof.* The proof is the same as that of lemma 3.1.15. $\square$

Let us now prove (C7):

**Proposition 3.2.32.** *Let $c = (c_{ij})_{ij}$ be a finite tuple from $M$ that is $\mathrm{val}^1_B$-separated, such that its $\mathrm{val}^1_B$-blocks are the $(c_i)_i = ((c_{ij})_j)_i$. Suppose we have $(p_i(x_i))_i$ a weak $\mathrm{val}^1_B$-block extension of the $c_{ij}$. Then $(p_i)_i$ is weakly orthogonal.*

Recall from definition 3.1.20 that types in a block extension are $\mathrm{Aut}(M/A)$-invariant. In particular, we must have $c_i \mathop{\downarrow}^{\mathbf{cut}}_A B$ for each $i$.

*Proof.* Enumerate the $c_i$ in increasing order of their values. Let $N$ be some $|M|^+$-saturated, strongly $|M|^+$-homogeneous elementary extension of $M$, and let $d_i \in N$ be a realization of $p_i$. By lemma 3.2.29, and lemma 3.2.30, it suffices to show that for all $i$, for all $f \in \mathrm{LC}(\mathbb{Q})$, $\mathrm{tp}(f(c_i)/M)$ is weakly orthogonal to $\mathrm{tp}(c_{<i}/M)$. This is trivial if $f = 0$.

Suppose $f \neq 0$. By $\mathrm{val}^1_B$-separatedness, we have $\mathrm{val}^1_B(f(c_i)) = \mathrm{val}^1_B(c_i)$. Let $H = H(f(c_i)/B)$. As $H$ is a $B$-$\vee$-definable group, $H(N)$ is $\mathrm{Aut}(M/N)$-invariant. As $\mathrm{tp}(f(d_i)/M)$ is a global $\mathrm{Aut}(M/A)$-invariant extension of $\mathrm{tp}(f(c_i)/M)$, we have $f(d_i) \notin (M + H(N))$ by lemma 3.2.28. Let $V = \mathrm{dcl}(Md_{<i})$. By lemma 3.2.31, it suffices to show that $V \subseteq (M + H(N))$. Let $v \in V$. Let $g \in \mathrm{LC}(\mathbb{Q})$ and $m \in M$ be such that $v = m + g(d_{<i})$. Let $h$ be the $A$-definable map witnessing corollary 3.2.25 applied to $g(d_{<i})$. As $\mathrm{val}^1_B$ is a $B$-valuation, and the values of the blocks before $c_i$ are strictly smaller than that of $c_i$, we have $G(g(c_{<i})/B) < G(c_i/B)$, thus $G(g(c_{<i})/B) < H$ by proposition 3.2.6. By remark 3.2.27, $g(d_{<i}) - h \circ g(c_{<i}) \in G(g(c_{<i})/B) < H$. It follows that $v = m + h \circ g(c_{<i}) + (g(d_{<i}) - h \circ g(c_{<i})) \in (M + H(N))$, concluding the proof. $\square$

Let us now build $\mathrm{val}^2_A$:

**Lemma 3.2.33.** *Let $d, d' \in M$ such that $G(d/A) = G(d'/A)$, $d$ is Archimedean over $A$ and $d'$ is ramified over $A$. Then $G(d + d'/A) = G(d/A)$ and $d + d'$ is Archimedean over $A$.*

*Proof.* Let $a \in \mathrm{ram}(d'/A)$. Then the (convex) coset $d + a \bmod G(d/A)$ has no point in $A$ and contains $d + d'$, in particular $d + d' \bmod G(d/A) \subseteq \mathrm{ct}(d + d'/A)$. By ultrametric inequality for the valuation $v^1_A$, $G(d + d'/A) \leqslant G(d/A)$. By lemma 3.2.17, $\mathrm{ct}(d + d'/A)$ is contained in $d + d' \bmod G(d + d'/A)$, therefore we have:

$$d{+}d' \bmod G(d/A) \subseteq \mathrm{ct}(d{+}d'/A) \subseteq d{+}d' \bmod G(d{+}d'/A) \subseteq d{+}d' \bmod G(d/A)$$

As a result, $G(d + d'/A) = G(d/A)$. This concludes the proof, as the coset $d + d' \bmod G(d + d'/A)$ has no point in $A$. □

**Corollary 3.2.34.** *The following satisfies the conditions of lemma 2.1.8:*

$$x \leqslant y$$
$$\Longleftrightarrow$$
*either* $\mathrm{val}^1_A(x) < \mathrm{val}^1_A(y)$*, or they have the same value and* $x$ *is ramified over* $A$*, or they have the same value and* $y$ *is Archimedean over* $A$*.*

**Definition 3.2.35.** We define $\mathrm{val}^2_A$ the valuation induced by the above preorder. It is clearly an $A$-valuation which refines $\mathrm{val}^1_A$.

This new valuation splits the $\mathrm{val}^1_A$-blocks in two, putting the $A$-ramified points below the $A$-Archimedean ones. For instance, in example 3.1.16, we have $\mathrm{val}^2_A(c_1) < \mathrm{val}^2_A(c_3) < \mathrm{val}^2_A(c_2)$.

Let us now prove (C6):

**Proposition 3.2.36.** *Let* $c = (c_{ij})_{ij}$ *be a finite tuple from* $M$ *that is* $\mathrm{val}^2_A$*-separated, such that its* $\mathrm{val}^2_A$*-blocks are the* $(c_i)_i = ((c_{ij})_j)_i$*. Then* $(\mathrm{tp}(c_i/A))_i$ *is weakly orthogonal.*

Note that the types which are weakly orthogonal are types over $A$ here, as opposed to the types in proposition 3.2.32 which were global types. The respective restrictions over $A$ of non-weakly orthogonal global types can very well be weakly orthogonal.

*Proof.* Just like in proposition 3.2.32, enumerate the blocks in increasing order of their value by $\mathrm{val}^2_A$. Choose $i$, and $f \in \mathrm{LC}(\mathbb{Q})$ with $f \neq 0$. As in the proof of proposition 3.2.32, it suffices to find $G$ an $\mathrm{Aut}(M/A)$-invariant convex subgroup of $M$ such that $f(c_i) \notin (A + G)$, and, for all $g \in \mathrm{LC}(\mathbb{Q})$, $g(c_{<i}) \in (A + G)$.

- Suppose the $(c_{ij})_j$ are ramified over $A$, then we choose $G = H(c_i/A)$. By separatedness, $G = H(f(c_i)/A)$. Now, for any $g \in \mathrm{LC}(\mathbb{Q})$, we have by proposition 3.2.6 $G(g(c_{<i})/A) < G$, therefore there exists $a$ in $A \smallsetminus \mathrm{Stab}(g(c_{<i})/A)$ such that $a \in G$. Let $a' \in A$ be such that $a' \in \,]g(c_{<i}) - |a|, g(c_{<i}) + |a|[$. Then $g(c_{<i}) \in \,]a' - |a|, a' + |a|[ \subseteq (A + G)$.

- Suppose the $(c_{ij})_j$ are Archimedean over $A$, and choose $g$. This time we let $G = G(c_i/A)$. By definition of being Archimedean, $f(c_i) \notin (A+G)$. Either $G(g(c_{<i})/A) < G$, in which case $g(c_{<i}) \in (A+H(c_i/A)) \subseteq (A+G)$ just as in the above item, or $G(g(c_{<i})/A) = G$, in which case $g(c_{<i})$ is ramified over $A$ as $\mathrm{val}^2_A(g(c_{<i})) < \mathrm{val}^2_A(c_i)$, therefore $g(c_{<i}) \in (A+G)$.

This concludes the proof. □

### 3.2.3 Tensor product

We saw that the blocks of a block extension for the coarsest valuation are weakly orthogonal, thus can be glued by just taking the union. It will no longer be the case for finer valuations, the tensor product of the blocks will no longer commute, and we will have to choose the order of the factors. For the third valuation, which is yet undefined, this choice will have to be taken with care, because some choices may give a global type which does not extend our type over $B$.

The following proposition shows us how to properly glue the elements of a $\mathrm{val}^2_B$-block extension.

**Proposition 3.2.37.** *Let $c$ be a finite tuple from $M$ such that $c \mathop{\downarrow}^{\mathbf{cut}}_A B$. Suppose $c = (c_i)_{i\in\{1,2\}}$ is a $\mathrm{val}^1_B$-block (i.e. a family having only one block) that is $\mathrm{val}^2_B$-separated, and $c_1$ (resp. $c_2$) is the block of $B$-ramified (resp. Archimedean) elements of $c$. Suppose $p_1, p_2$ is a weak $\mathrm{val}^2_B$-block extension of $c$. Then $c$ satisfies the restriction to $B$ of the tensor product of $p_1, p_2$ in any order.*

Note that, although the restrictions to $B$ of $p_1 \otimes p_2$ and $p_2 \otimes p_1$ coincide, these two global types are distinct, i.e. $(p_1, p_2)$ is not weakly orthogonal. This will be explained in the next section.

*Proof.* By proposition 3.2.36, $(\mathrm{tp}(c_i/B))_{i\in\{1,2\}}$ is weakly orthogonal, thus $\mathrm{tp}(c_1/B)$ is complete in $S(B+\mathbb{Q}c_2)$, and it coincides in particular with $(p_1)_{|B+\mathbb{Q}c_2}$. It follows that $c \vDash (p_1 \otimes p_2)_{|B}$, and the same reasoning can be applied to $p_2 \otimes p_1$. □

Let us now build $\mathrm{val}^3_B$:

**Lemma 3.2.38.** *Let $d_1, d_2 \in M$ be two points which are ramified over $A$, and have the same $\mathrm{val}^2_A$-value. Suppose $\delta(d_1/A) < \delta(d_2/A)$. Then $\mathrm{val}^2_A(d_1 + d_2) = \mathrm{val}^2_A(d_2)$, and $\delta(d_1 + d_2/A) = \delta(d_2/A)$ (as defined in definition 3.2.15).*

*Proof.* Let $a_i \in \mathrm{ram}(d_i/A)$. Then, by remark 2.1.4 we have $\Delta(d_1+d_2-a_1-a_2) = \delta(d_2/A)$, therefore $d_1 + d_2$ is ramified over $A$, and $\delta(d_1 + d_2/A) = \delta(d_2/A)$. By remark 3.2.16, we have:

$$\mathrm{ct}_>(\Delta(\mathrm{Stab}(d_2/A))/\Delta(A)) = \mathrm{ct}_>(\Delta(\mathrm{Stab}(d_1 + d_2/A))/\Delta(A))$$

However, those stabilizers are convex subgroups of $A$, thus the sets of their Archimedean classes are initial segments of $\Delta(A)$, so this equality of cuts implies $\Delta(\mathrm{Stab}(d_2/A)) = \Delta(\mathrm{Stab}(d_1 + d_2/A))$. It follows that $\mathrm{Stab}(d_2/A) = \mathrm{Stab}(d_1 + d_2/A)$, thus $H(d_2/A) = H(d_1 + d_2/A)$, thus $\mathrm{val}^1_A(d_2) = \mathrm{val}^1_A(d_1 + d_2)$. This concludes the proof, as $d_2$, $d_1 + d_2$ are both ramified over $A$. □

**Corollary 3.2.39.** *The following satisfies the conditions of lemma 2.1.8:*

$$x \leqslant y$$
$$\iff$$

*either $\mathrm{val}^2_A(x) < \mathrm{val}^2_A(y)$, or they have the same value, they are both ramified over $A$, and $\delta(x/A) \leqslant \delta(y/A)$.*

**Definition 3.2.40.** We define $\mathrm{val}^3_A$ the valuation induced by the above preorder. It is clearly an $A$-valuation which refines $\mathrm{val}^2_A$.

The valuation splits the ramified $\mathrm{val}^2_A$-blocks in a way which depends on the Archimedean classes $\delta$ added by their elements. Note that, while for $i \in \{1, 2\}$, $\mathrm{val}^i_A$ factors through the quotient by $A$-elementary equivalence (these are model-theoretic objects), $\mathrm{val}^3_A$ is an *algebraic* object that is no longer refined by types.

We need to define additional notions in order to better understand $\mathrm{val}^3_B$-block extensions. Contrary to $\mathrm{val}^2_B$ and $\mathrm{val}^1_B$, a weak $\mathrm{val}^3_B$-block extension of some tuple might not be strong. We need to find necessary and sufficient conditions for such a block extension to be strong.

**Proposition 3.2.41.** *Let $c = (c_{ij})_{ij}$ be a finite $\mathrm{val}^2_B$-block of $B$-ramified points from $M$ that is $\mathrm{val}^3_B$-separated, and whose $\mathrm{val}^3_B$-blocks are the $(c_i)_i = ((c_{ij})_j)_i$. Suppose the indices $i$ are ordered such that: $i < k \iff \delta(c_i/B) < \delta(c_k/B)$. Let $d = (d_{ij})_{ij}$ be another tuple from $M$. Then $c \equiv_B d$ if and only if the following conditions hold:*

1. *For each $i$, we have $d_i \equiv_B c_i$. In particular, each $d_i$ is a $\mathrm{val}^3_B$-block of $B$-ramified points that has the same $\mathrm{val}^2_B$-value as $c$.*

2. *For every $i$ and $k$, we have $i < k$ if and only if $\delta(d_i/B) < \delta(d_k/B)$.*

*Proof.* For each $i$ and $j$, let $b_{ij} \in \mathrm{ram}(c_{ij}/B)$.

Suppose $c \equiv_B d$ (then condition 1 holds). Let $i, k, j, l$ be indices. If $i < k$, then $\mathrm{tp}(c/B) \models |x_{ij} - b_{ij}| < \frac{1}{n}|x_{kl} - b_{kl}|$ for every $n$, therefore we have $\Delta(d_{ij} - b_{ij}) < \Delta(d_{kl} - b_{kl})$. Moreover, by corollary 3.2.20, we have $\delta(d_{ij}/B) = \Delta(d_{ij} - b_{ij})$ (and similarly for $kl$), and condition 2 holds.

Conversely, suppose that the conditions of the list hold. Let $f_i \in \mathrm{LC}^{|c_i|}(\mathbb{Q})$, and let $b' = \sum_i f_i(b_i)$. We have to prove that $c' = \sum_i f_i(c_i)$ and $d' = \sum_i f_i(d_i)$ have the same cut over $B$. This is trivial if each $f_i$ is zero, else let $k$ be the maximal index for which $f_k$ is non-zero. By $\mathrm{val}^3_B$-separatedness, $\mathrm{val}^3_B(f_k(c_k)) = \mathrm{val}^3_B(c_{kj})$ for any $j$. Moreover, $f_k((c_{kj} - b_{kj})_j)$ is a linear combination of elements of $G(c_k/B)$, hence it belongs to $G(c_k/B)$, and $f_k(b_k)$ is a ramifier of $f_k(c_k)$. Likewise, we have $f_i(c_i) - f_i(b_i) \in G(c_i/B)$ for any $i$. It follows by maximality of $k$ that:

$$\delta(c_k/B) = \Delta(f_k(c_k) - f_k(b_k)) > \Delta(f_i(c_i) - f_i(b_i)) \qquad (3)$$

for each $i \neq k$. By the conditions of the list, this also holds for $d$. Choose an arbitrary index $j$. There exists $n \in \omega_{>0}$ and $\varepsilon \in \{\pm 1\}$ such that $f_k(c_k) - f_k(b_k)$ lies between $\frac{\varepsilon}{n}(c_{kj} - b_{kj})$ and $\varepsilon n(c_{kj} - b_{kj})$. By the first condition, this must hold for $c$ replaced with $d$ (with the same $\varepsilon, n$). Moreover, as 3 holds for $d$, for all $i \neq k$ and $N > 0$, we have $|f_i(d_i) - f_i(b_i)| < \frac{1}{N}|d_{kj} - b_{kj}|$. By choosing $N$ large enough, and by summing everything, we get that $d' - b'$ lies between $\frac{\varepsilon}{N+1}(d_{kj} - b_{kj})$ and $\varepsilon(N+1)(d_{kj} - b_{kj})$, and we get the same property with $d$ replaced by $c$. Then we have:

$$\begin{aligned} \mathrm{ct}(\Delta(c' - b')/\Delta(B)) &= \mathrm{ct}(\Delta(c_{kj} - b_{kj})/\Delta(B)) \\ &= \mathrm{ct}(\Delta(d_{kj} - b_{kj})/\Delta(B)) \\ &= \mathrm{ct}(\Delta(d' - b')/\Delta(B)) \end{aligned}$$

The second equality is a consequence of corollary 3.2.20 applied to $c_{kj}$, $d_{kj}$ (by the hypothesis $c_k \equiv_B d_k$). The other equalities follow from the inequalities of the previous paragraph. We also have:

$$d' < b' \iff \varepsilon = -1 \iff c' < b'$$

as a result, by corollary 3.2.20, we have $c' \equiv_B d'$, concluding the proof. □

*Remark* 3.2.42. We recall some facts from previous subsections. Let $c = (c_{ij})_{ij}$ be a finite $\mathrm{val}^2_B$-block of $B$-ramified points from $M$ that is $\mathrm{val}^3_B$-separated, whose $\mathrm{val}^3_B$-blocks are the $(c_i)_i = ((c_{ij})_j)_i$. Let $(p_i)_i$ be a weak $\mathrm{val}^3_B$-block extension of $c$. By corollary 3.2.25, and by separatedness, there exists $G$ (resp. $H$) a unique $A$-type-definable(resp. $A$-$\vee$-definable) convex subgroup such that $H < G$, and exactly one of the following conditions holds for each $i$:

1. For every non-zero $f \in \mathrm{LC}^{|c_i|}(\mathbb{Q})$, there exists $b \in B$ such that we have $p_i(x_i) \models f(x_i) - b \in G \setminus \underline{G}_M$.
2. For every non-zero $f \in \mathrm{LC}^{|c_i|}(\mathbb{Q})$, there exists $b \in B$ such that we have $p_i(x_i) \models f(x_i) - b \in \overline{H}^M \setminus H$.

The following notion will mostly be used to understand strong $\mathrm{val}^3_B$-block extensions, and to prove (C5):

**Definition 3.2.43.** If the first condition of remark 3.2.42 holds, then we say that $p_i$ is *outer*, else $p_i$ is *inner*.

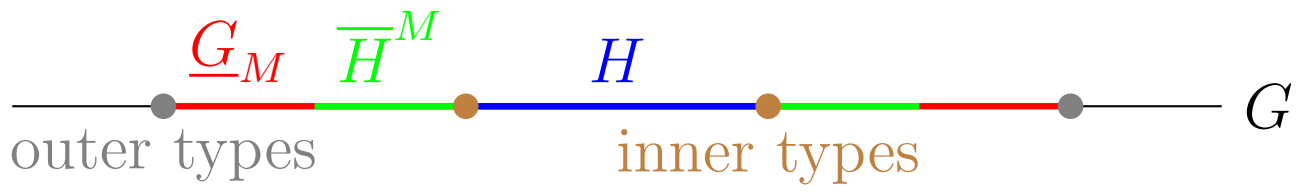

**Lemma 3.2.44.** *Let $c = (c_{ij})_{ij}$ be a finite $\mathrm{val}^2_B$-block of $B$-ramified points from $M$ that is $\mathrm{val}^3_B$-separated, such that its $\mathrm{val}^3_B$-blocks are the $(c_i)_i = ((c_{ij})_j)_i$. Suppose the indices $i$ are ordered such that: $i < k \iff \delta(c_i/B) < \delta(c_k/B)$. Let $(p_i)_i$ be a weak $\mathrm{val}^3_B$-block extension of $c$. Let $I$ (resp. $O$) be the set of indices $i$ such that $p_i$ is inner (resp. outer). Then the following are equivalent:*

- *$(p_i)_i$ is a strong $\mathrm{val}^3_B$-block extension of $c$.*
- *For every $i \in I, o \in O$, we have $i < o$.*

*Proof.* The realizations of $\bigcup_i p_i$ satisfy the first condition of proposition 3.2.41, so we are trying to characterize when some of them satisfy the second condition. Let $H = H(c/A) < G(c/A) = G$ be the convex subgroups witnessing remark 3.2.42.

Suppose we have $i \in I$ and $o \in O$ such that $o < i$. Let $d \models \bigcup_j p_j$. Then, as $p_i$ is inner, we have $\delta(d_i/M) \in \Delta(\overline{H}^M)$. By the same argument, we have $\delta(d_o/M) \notin \Delta(\underline{G}_M)$, therefore $\delta(d_i/M) < \delta(d_o/M)$, thus $d \not\equiv_B c$, and we proved the top-to-bottom implication.

Conversely, suppose we have $i < o$ for every $i \in I$, $o \in O$. In order to satisfy the second condition of proposition 3.2.41, it suffices to prove that the following partial type is consistent:

$$q(x) = \bigcup_i p_i(x_i) \cup \left\{ |x_{i,f(i)} - b_i| < \frac{1}{n}|x_{k,f(k)} - b_k| \ : \ i < k, n > 0 \right\}$$

with $f(i)$ an arbitrary choice of a coordinate of $(c_{ij})_j$ for each $i$, and with $b_i \in \mathrm{ram}(c_{i,f(i)}/B)$. Let us build a realization $\beta$ of $q$.

Let $(\alpha_i)_{i\in I} \models \bigcup_{i \in I} p_i$. We define by induction, for each $i \in I$, a tuple $\beta_i \models p_i$ such that, for $i < k$ in $I$, we have $\Delta(\beta_{if(i)} - b_i) < \Delta(\beta_{kf(k)} - b_k)$. Let $N$ be an $|M|^+$-saturated, strongly $|M|^+$-homogeneous elementary extension of $M$ containing $(\alpha_i)_{i\in I}$. Let $k \in I$, and suppose we defined $\beta_i \in N$ for every $i < k$. By saturation, let $\gamma \in N$ such that $\gamma \in \overline{H}^M \smallsetminus H$, and $\Delta(\gamma) > \Delta(\beta_{if(i)} - b_i)$ for all $i < k$ (in other words, $\gamma \in \overline{H}^M \smallsetminus \underline{(\overline{H}^M)}_{M + \sum_{i<k} \mathbb{Q}\beta_{i,f(i)}}$). Such a $\gamma$ exists, as we are only dealing with indices of inner types. Then, by corollary 3.2.20, there exists $\varepsilon \in \{\pm 1\}$ such that $b_k + \varepsilon\gamma \equiv_M \alpha_{k,f(k)}$. By strong homogeneity, there exists $\sigma \in \mathrm{Aut}(N/M)$ such that $\sigma(\alpha_{k,f(k)}) = b_k + \varepsilon\gamma$. We set $\beta_k = \sigma(\alpha_k)$.

We can similarly define $(\beta_o)_{o\in O}$. Then $(\beta_i)_{i\in I\cup O}$ satisfies $q$, concluding the proof. □

Now we can prove (C5):

**Proposition 3.2.45.** *Let $c = (c_{ij})_{ij}$ be a finite $\mathrm{val}^2_B$-block of $B$-ramified points from $M$ that is $\mathrm{val}^3_B$-separated, whose $\mathrm{val}^3_B$-blocks are the $(c_i)_i = ((c_{ij})_j)_i$. Let $(p_i)_i$ be a strong $\mathrm{val}^3_B$-block extension of $c$. Then some tensor product of the $p_i$ is consistent with $\mathrm{tp}(c/B)$.*

We will see much later, in lemma 3.3.25, that the tensor product in the following proof is in fact the only completion of $(p_i)_i$ (be it $\mathrm{Aut}(M/A)$-invariant or not) which is consistent with $\mathrm{tp}(c/B)$.

*Proof.* Suppose the indices $i$ are ordered such that $i < k$ if and only if we have $\delta(c_i/B) < \delta(c_k/B)$. Define $I$, $O$, $G$, $H$, just as in lemma 3.2.44. Let $\alpha_I \vDash p_I = \bigcup_{i \in I} p_i$, $\alpha_O \vDash p_O = \bigcup_{o \in O} p_o$. As the conditions of lemma 3.2.44 are satisfied, if we had $\alpha_I \equiv_B (c_i)_{i \in I}$ and $\alpha_O \equiv_B (c_o)_{o \in O}$, then by proposition 3.2.41 we would clearly have $\alpha_I \alpha_O \equiv_B c$. Moreover:

$$G(\alpha_I/M) = \overline{H}^M < \underline{G}_M = H(\alpha_O/M)$$

thus one clearly sees that $\alpha_I$ and $\alpha_O$ are $\mathrm{val}^2_M$-blocks of distinct values, thus their types over $M$ are weakly orthogonal by proposition 3.2.36. Therefore, if the respective types over $M$ of $\alpha_I$ and $\alpha_O$ are tensor products of the $p_i$, then so is the type of $\alpha_I \alpha_O$. As a result, we can assume that either $I$ or $O$ is empty.

Suppose $O = \varnothing$. For $k \in I$, and $c' = (c_i)_{i > k}$, let us show that $c_k \vDash p_{k|D}$, with $D$ the $\mathbb{Q}$-vector subspace generated by $Bc'$. This implies that $\mathrm{tp}(c/B)$ is consistent with the tensor product of the $p_i$ in increasing order. Let $f$ be non-zero in $\mathrm{LC}^{|c_k|}(\mathbb{Q})$. By $\mathrm{val}^3_B$-separatedness, we have $\delta(f(c_k)/B) = \delta(c_k/B)$ (call this Archimedean class $\delta$). Let $b \in \mathrm{ram}(f(c_k)/B)$. As $i \in I$, we have $p_k(x_k) \vDash f(x_k) - b \in \overline{H}^M \smallsetminus H$. It is enough to show that $f(c_k) - b \in \overline{H}^D$, for the cut over $M$ corresponding to the pushforward of $p_k$ by $f$ would be included in $\mathrm{ct}(f(c_k)/D)$, and we could conclude. For $i > k$, let $\delta_i = \delta(c_i/B)$. By $\mathrm{val}^3_B$-separatedness, we clearly have $\Delta(D) = \Delta(B) \cup \{\delta_i | i > k\}$. By definition of $I$, we have $\delta < \delta_i$ for every $i > k$. Moreover, we have by remark 3.2.16:

$$\mathrm{ct}(\delta_i/\Delta(B)) = \mathrm{ct}_>(\Delta(H) \cap \Delta(B)/\Delta(B))$$

as a result, $\Delta(H) \cap \Delta(D) = \Delta(H) \cap \Delta(B)$, and:

$$\mathrm{ct}(\delta/\Delta(D)) = \mathrm{ct}_>(\Delta(H) \cap \Delta(D)/\Delta(D))$$

thus $\delta \in \Delta(\overline{H}^D)$, concluding the proof.

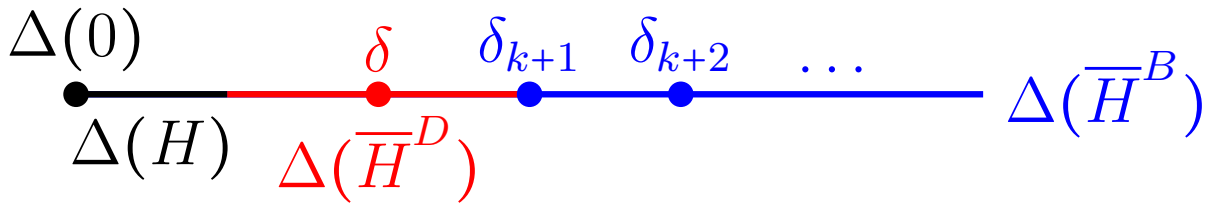

Now, if instead $I = \varnothing$, then a similar proof would show that the tensor product of the $p_o$ in decreasing order would be consistent with $\mathrm{tp}(c/B)$. □

**Corollary 3.2.46.** *The property (C5) holds.*

*Proof.* Proposition 3.2.45 deals with the particular case of a $\mathrm{val}^2_B$-block. Then, given the global invariant extensions of the types of each $\mathrm{val}^2_B$-block given by proposition 3.2.45, some tensor product of those global types witnesses the statement, by proposition 3.2.37 and proposition 3.2.32. □

## 3.3 How to build a block extension

In the previous section, we defined the key valuations $\mathrm{val}_A^i$, and we gave a proof of the properties (C5), (C6), (C7). In this section, we give a proof of (C4). The main technical goal is to show that $\downarrow_A^{\mathbf{cut}} B$ coincides with $\downarrow_A^{\mathbf{inv}} B$, when what we put on the left is a $\mathrm{val}_B^3$-separated $\mathrm{val}_B^3$-block. This will show the existence of a weak block extension, and then it will not be hard to show that such a block extension can be chosen strong. In fact, given some $\mathrm{val}_B^3$-separated tuple $c$, we describe explicitly the homeomorphism type of the space of global $\mathrm{Aut}(M/A)$-invariant extensions of $\mathrm{tp}(c/B)$. We adopt in this section assumptions 3.1.14

### 3.3.1 The Archimedean case

In this subsection, we deal with the case where $c$ is a $\mathrm{val}_B^3$-separated block of $B$-Archimedean points from $M$.

*Remark* 3.3.1. Suppose $c$ is a $\mathrm{val}_B^3$-separated block of $B$-Archimedean points from $M$. As $c \downarrow_A^{\mathbf{cut}} B$, the $c_i$ must be Archimedean over $A$. Indeed, by contraposition, if $c_i$ is ramified over $A$ for some $i$, then we cannot have $\delta(c_i/A) \in \Delta(B)$, otherwise there would exist $b \in B$, $a \in A$ such that $\delta(c_i/A) = \Delta(c_i - a) = \Delta(b)$, thus $c_i$ would lie between $a + \frac{1}{n}b$ and $a + nb$ for some $n \in \mathbb{Z}$, which would contradict $c \downarrow_A^{\mathbf{cut}} B$. Now, as $\Delta(c_i - a) \notin \Delta(B)$ we have a contradiction with the hypothesis that all the $c_i$ are Archimedean over $B$.

**Assumptions 3.3.2.** On top of assumptions 3.1.14, we assume that $c$ is a $\mathrm{val}_B^3$-separated block of $B$-Archimedean points from $M$. By remark 3.3.1 and lemma 3.2.23, let us also fix $G = G(c/A) = G(c/B)$.

**Lemma 3.3.3.** *Let $p_i$ be a global $\mathrm{Aut}(M/A)$-invariant extension of $\mathrm{tp}(c_i/B)$. Suppose $q$ is a complete global extension of $\mathrm{tp}(c/B)$ which extends $\bigcup_i p_i$. If we have $q(x) \models f(x) - f(c) \notin \underline{G}_M$ for all non-zero $f \in \mathrm{LC}^n(\mathbb{Q})$, then $q$ is $\mathrm{Aut}(M/A)$-invariant.*

*Proof.* Suppose by contradiction $q$ is not $\mathrm{Aut}(M/A)$-invariant. Then there exists a non-zero $f \in \mathrm{LC}^n(\mathbb{Q})$, and a cut $X$ over $M$, such that $q(x) \models f(x) \in X$, and the global 1-type induced by $X$ is not $\mathrm{Aut}(M/A)$-invariant. As $c \mathop{\downarrow}^{\mathbf{cut}}_{A} B$, we have $f(c) \mathop{\downarrow}^{\mathbf{inv}}_{A} B$. By lemma 3.2.23, we clearly see that the two global $\mathrm{Aut}(M/A)$-invariant extensions of $\mathrm{tp}(f(c)/B)$ correspond to the translates by $f(c)$ of the two connected components of $G \setminus \underline{G}_M$. We obtain a contradiction by showing that $X$ must coincide with one of those. It is enough to show that $q(x) \vDash f(x) - f(c) \in G \setminus \underline{G}_M$. By hypothesis, we just need to have $q(x) \vDash f(x) - f(c) \in G$, which follows from the fact that $c$ is a $\mathrm{val}^3_B$-separated block. Indeed, $\mathrm{val}^3_B(f(c)) = \mathrm{val}^3_B(c)$, therefore $f(c)$ is Archimedean over $B$, and $G(f(c)/B) = G$. This implies by proposition 3.2.18 that $\mathrm{tp}(f(c)/B) \vDash y - f(c) \in G$, concluding the proof. □

It turns out that the space of global $\mathrm{Aut}(M/A)$-invariant extensions of $\mathrm{tp}(c/B)$ in the Archimedean case has a very simple description:

**Proposition 3.3.4.** *Let $N$ be an $|M|^+$-saturated, strongly $|M|^+$-homogeneous elementary extension of $M$. Let $H = (\underline{G}_M)(N)$ be the convex subgroup of $N$ generated by $G(M)$, and let $N'$ be the ordered Abelian group ${}^{N}\!\diagup_{H}$, which is a model of* DOAG *by saturation. Let $F_1 \subseteq S^{|c|}(M)$ be the closed set of* $\mathrm{Aut}(M/A)$*-invariant extensions of* $\mathrm{tp}(c/B)$*, and $F_2 \subseteq S^{|c|}(\{0\})$ be the closed set of types of $\mathbb{Q}$-free families. Then the map:*

$$g : \ \mathrm{tp}^N(\alpha/M) \longmapsto \mathrm{tp}^{N'}((\alpha_i - c_i \bmod H)_i/\{0\})$$

*is a well-defined homeomorphism $F_1 \longrightarrow F_2$.*

*Proof.* Let $f \in \mathrm{LC}^n(\mathbb{Q})$ be non-zero, and let $\alpha \in N$ realize a global $\mathrm{Aut}(M/A)$-invariant extension of $\mathrm{tp}(c/B)$. Then, by separatedness, $f(c)$ is Archimedean over $A$, and $G(f(c)/A) = G$, thus $f(c) - f(\alpha) \notin \underline{G}_M$ by lemma 3.2.23, i.e. $f(c) - f(\alpha) \notin H$. No non-trivial linear combination of $(c_i - \alpha_i)_i$ belongs to $H$, hence $(c_i - \alpha_i \bmod H)_i$ is $\mathbb{Q}$-free.

Moreover, if $\beta \equiv_M \alpha$, then we have $f((\beta_i - c_i \bmod H)_i) > 0$ in $N'$ if and only if $f((\beta_i - c_i)_i) \in \bigcap_{m \in H} ]m, +\infty[$ in $N$, if and only if the same holds for $\alpha$, therefore $(\beta_i - c_i \bmod H)_i \equiv_{\{0\}} (\alpha_i - c_i \bmod H)_i$ in $N'$. It follows that $g$ is well-defined, and it is continuous by definition (and quantifier elimination).

Let $p_1, p_2 \in F_1$ be distinct. Then there exists two distinct cuts $X_1, X_2$ over $M$, and $f \in \mathrm{LC}^n(\mathbb{Q})$ non-zero, such that $p_i(x) \vDash f(x) \in X_i$. However, by lemma 3.2.23, $X_i$ must be the translate by $f(c)$ of a connected component of $G \setminus \underline{G}_M$. As a result, one of the $X_i$ (say, $X_1$) is the leftmost connected component, and the other is the rightmost one. We have:

$$p_1(x) \vDash f(x) - f(c) < \underline{G}_M$$

thus $g(p_1) \models f(x) < 0$. With the same reasoning, we get $g(p_2) \models f(x) > 0$, and we conclude that $g(p_1) \neq g(p_2)$, proving that $g$ is injective.

Now let us prove surjectivity. Let $\beta'$ be a $\mathbb{Q}$-free family from $N'$. Let $\beta$ be a preimage of $\beta'$ by $\pi: N \longrightarrow N'$. Let $d_1 \in N$ such that $\Delta(d_1) > \Delta(\beta_i)$ for every $i$. By saturation, let $d_2 \in N$ such that $\Delta(d_2) \in \Delta(G) \smallsetminus \Delta(\underline{G}_M)$. Then $|d_1|$ and $|d_2|$ have the same cut over $G(M)$: $+\infty$, the cut of positive elements larger than all the elements of $G(M)$. By strong homogeneity, there is an automorphism $\sigma \in \mathrm{Aut}(N/G(M))$ such that $\sigma(|d_1|) = |d_2|$. Since $\sigma$ fixes $G(M)$, we have $\sigma(H) = H$, thus the map $\sigma': \pi(x) \longmapsto \pi(\sigma(x))$ is a well-defined automorphism of the ordered group $N'$. As $\beta'$ is $\mathbb{Q}$-free, so is $\sigma'(\beta')$. As $\sigma(\beta)$ is a preimage of $\sigma'(\beta')$, none of its elements is in $H$. Moreover, $|\sigma(\beta_i)| \leqslant \sigma(|d_1|) = |d_2| \in G$ for every $i$, thus $\sigma(\beta_i) \in G$. As $\mathrm{tp}^{N'}(\beta'/\{0\}) = \mathrm{tp}^{N'}(\sigma'(\beta')/\{0\})$, and we are looking for a preimage of this type by $g$, we may assume $\sigma = id$, and thus, by $\mathbb{Q}$-freeness of $\beta'$, $f(\beta) \in G \smallsetminus \underline{G}_M$ for every non-zero $f \in \mathrm{LC}^n(\mathbb{Q})$. Let $p = \mathrm{tp}((c_i + \beta_i)_i/M)$. As $\beta_i \in G$ for all $i$, $p$ extends $\mathrm{tp}(c/B)$, and it is $\mathrm{Aut}(M/A)$-invariant by lemma 3.3.3. Thus $p \in F_1$, and it is a preimage of $\mathrm{tp}(\beta'/\{0\})$ by $g$.

We showed that $g$ is a continuous bijection, and we conclude by compactness and separation. □

Since $F_2$ is non-empty, we have in particular proved the following:

**Corollary 3.3.5.** *If $c$ is a $\mathrm{val}^3_B$-separated block of $B$-Archimedean points from $M$, then we have $c \mathop{\downarrow}^{\mathbf{inv}}_A B$.*

The next subsection gives tools to understand what $F_2$ looks like. For the most part, it is used to deal with the ramified blocks.

### 3.3.2 Archimedean groups

We need to define some algebraic notions about Archimedean groups in order to well-understand the global invariant extensions of the types of ramified blocks. Let us start with well-known facts:

**Proposition 3.3.6.** *Let $G$ be an ordered Abelian group. Then $G$ has at most two Archimedean classes (including $\Delta(0)$) if and only if it embeds into $(\mathbb{R}, +, <)$. In that case, for all $x \in G_{>0}$, for all $\mu \in \mathbb{R}_{>0}$, there exists a unique embedding of ordered groups $\sigma : G \longrightarrow \mathbb{R}$ sending $x$ to $\mu$.*

The idea of the proof is to send each $y \in G$ to $\mu \cdot (\sup\{\lambda \in \mathbb{Q} | \lambda y < x\})^{-1}$.

*Remark* 3.3.7. In the special case $G = \mathbb{R}$, we see that an embedding of ordered groups $\sigma : \mathbb{R} \longrightarrow \mathbb{R}$ is uniquely determined by the choice of $\sigma(1)$, hence $\sigma$ actually coincides with the automorphism $x \longmapsto \sigma(1) \cdot x$. This establishes an easy description of the group of ordered group automorphisms of $\mathbb{R}$, which is naturally isomorphic to the multiplicative group $\mathbb{R}_{>0}$.

**Definition 3.3.8.** We define $\mathbb{P}^+_n$, the space of half-lines in dimension $n$ (plus the origin) as the set of orbits of $\mathbb{R}^n$ under the automorphisms of the ordered group $\mathbb{R}$. We call $\mathbb{P}^+$ the canonical map $\mathbb{R}^n \longrightarrow \mathbb{P}^+_n$.

Suppose $G$ is an Archimedean ordered Abelian group, and let $u$ be a finite tuple from $G$. By the universal property of proposition 3.3.6, we see that $\mathbb{P}^+(\sigma(u))$ is always the same for any ordered group embedding $\sigma : G \longrightarrow \mathbb{R}$. As a result, we can define $\mathbb{P}^+(u) \in \mathbb{P}^+_{|u|}$ as this unique class.

*Remark* 3.3.9. If $f \in \mathrm{LC}^n(\mathbb{Q})$, then $f$ commutes with any automorphism, thus one can easily show that $f(\mathbb{P}^+(u)) = \mathbb{P}^+(f(u)) \in \mathbb{P}^+_1$ for any $u \in \mathbb{R}^n$. However, $\mathbb{P}^+_1$ is not very complicated to describe: $f(\mathbb{P}^+(u))$ is either $\mathbb{R}_{>0}$, $\mathbb{R}_{<0}$ or $\{0\}$.

**Definition 3.3.10.** Let $u \in \mathbb{R}^n$. We say that $\mathbb{P}^+(u)$ is $\mathbb{Q}$-*free* if $u$ is $\mathbb{Q}$-free. By remark 3.3.9, this definition does not depend on the choice of the representative $u$.

**Proposition 3.3.11.** *Let $G$ be an ordered Abelian group. Then for any subgroup $H$ of $G$, and $\Delta(0) \neq \delta \in \Delta(G)$, the quotient:*

$$H_\delta = {}^{H_{\leqslant\delta}} \diagup {}_{H_{<\delta}}$$

*as defined in definition 3.1.13, is Archimedean.*

In the literature, the $(H_\delta)_{\delta\in\Delta(H)}$ are called the “components” of $H$.

**Definition 3.3.12.** Let $\delta \in \Delta(M)$, and $u$ a tuple from $M_{=\delta}$. By proposition 3.3.11, we can define $\mathbb{P}^+(u) = \mathbb{P}^+(u \bmod M_{<\delta})$ with respect to the definition 3.3.8 in the Archimedean group ${}^{M_{\leqslant\delta}} \diagup {}_{M_{<\delta}}$.

Suppose $d = (d_1 \dots d_n)$ is a $\mathrm{val}^3_B$-block of $B$-ramified elements of $M$, and let $\delta = \delta(d_i/B) \in \Delta(M)$. Let $b_i \in \mathrm{ram}(d_i/B)$. Then we can define $\mathbb{P}^+(d/B) = \mathbb{P}^+((d_i - b_i)_i)$. As $H(d_i/B) \leqslant M_{<\delta}$, by remark 3.2.12, this definition does not depend on the choice of the $b_i$.

**Proposition 3.3.13.** *Let $d = (d_1 \dots d_n)$ be a finite $\mathrm{val}^3_B$-block of $B$-ramified points from $M$. Then $d$ is $\mathrm{val}^3_B$-separated if and only if $\mathbb{P}^+(d/B)$ is $\mathbb{Q}$-free.*

*Proof.* Let $b_i \in \mathrm{ram}(d_i/B)$, and let $\delta = \Delta(d/B)$. Then $d$ is not $\mathrm{val}^3_B$-separated if and only if $\mathrm{val}^3_B(f(d)) < \mathrm{val}^3_B(d)$ for some non-zero $f \in \mathrm{LC}^n(\mathbb{Q})$, if and only if there exists such an $f$ with $\Delta(f(d) - f(b)) < \delta$, if and only if we have $f$ such that $\{0\} = f(\mathbb{P}^+((d_i - b_i)_i)) = f(\mathbb{P}^+(d/B))$, if and only if $\mathbb{P}^+(d/B)$ is not $\mathbb{Q}$-free. □

### 3.3.3 The ramified case

**Assumptions 3.3.14.** In this subsection, on top of assumptions 3.1.14 suppose $c = (c_i)_i$ is a $\mathrm{val}^3_B$-separated $\mathrm{val}^3_B$-block of $B$-ramified points. Fix $G = G(c/B)$, $H = H(c/B)$, $\delta = \delta(c/B)$. For each $i$, let $b_i$ be in $\mathrm{ram}(c_i/B)$, and write $b = (b_i)_i$.

Recall that corollary 3.2.20 characterizes the type of a ramified singleton in terms of cuts of Archimedean classes. We can now extend this characterization to the type of our $\mathrm{val}^3_B$-separated ramified $\mathrm{val}^3_B$-block:

**Proposition 3.3.15.** *Let $d = (d_i)_i$ be a tuple from $M$. Then we have $d \equiv_B c$ if and only if the following conditions hold:*

- *$d_i \equiv_B c_i$ for every $i$.*
- *$d$ is a $\mathrm{val}^3_B$-block, i.e. $\delta(d_i/B) = \delta(d_j/B)$ for every $i$, $j$.*
- *$\mathbb{P}^+(d/B) = \mathbb{P}^+(c/B)$.*

*Proof.* Suppose the conditions of the list hold. Let $f \in \mathrm{LC}^n(\mathbb{Q})$ be non-zero. Then, as $c$ is $\mathrm{val}^3_B$-separated, $\mathbb{P}^+(c/B)$ is $\mathbb{Q}$-free, therefore $\Delta(f(c) - f(b)) = \delta$. As a result, $\mathrm{ct}(f(c)/B)$ is the translate by $f(b)$ of one of the connected components of $G \smallsetminus H$. The same holds for $d$, but we need to make sure that $f(c)$ and $f(d)$ lie in the same connected component. The thing is, $f(c)$ lies in the rightmost one if and only if $(f(c) - f(b) \bmod M_{<\delta}) > 0$ in the ordered group ${M_{\leq\delta}}/{M_{<\delta}}$, if and only if $\mathbb{P}^+(f(c) - f(b)) = \mathbb{R}_{>0}$, if and only if $\mathbb{R}_{>0} = f(\mathbb{P}^+((c_i - b_i)_i)) = f(\mathbb{P}^+(c/B))$. Now, by the hypothesis $\mathbb{P}^+(d/B) = \mathbb{P}^+(c/B)$, those conditions are equivalent to $f(d)$ lying in the rightmost connected component, and we get the bottom-to-top direction.

Conversely, suppose now $d \equiv_B c$, in particular the first condition holds. As for the second, $d$ must be a $\mathrm{val}^3_B$-block, otherwise it would be witnessed in its type by formulas $(|x_i - b_i| > n|x_j - b_j|)_n$ for some $i \neq j$, which would contradict the hypothesis on $c$. By strong homogeneity of $M$, let $\sigma \in \mathrm{Aut}(M/B)$ be such that $\sigma(c) = d$. The equivalence relation $\Delta(x) = \Delta(y)$ is an $\vee$-definable subset of $S^2(\varnothing)$, therefore it is invariant by $\sigma$. As a result, the map:

$$\sigma' : (x \bmod M_{<\delta}) \longmapsto (\sigma(x) \bmod M_{<\delta(d/B)})$$

is a well-defined automorphism of ordered groups from the quotient ${M_{\leqslant\delta}} / {M_{<\delta}}$ to ${M_{\leqslant\delta(d/B)}} / {M_{<\delta(d/B)}}$, which sends each $c_i - b_i$ to $d_i - b_i$. Let $\tau$ be an embedding of ordered groups ${M_{\leqslant\delta(d/B)}} / {M_{<\delta(d/B)}} \to \mathbb{R}$. Then $\mathbb{P}^+(c/B) \ni \tau \circ \sigma'((c_i - b_i)_i) = \tau((d_i - b_i)_i) \in \mathbb{P}^+(d/B)$, thus the last condition of the list holds, concluding the proof. □

**Lemma 3.3.16.** *Choose an arbitrary index $i$. Suppose that $p$ is a global $\mathrm{Aut}(M/A)$-invariant extension of $\mathrm{tp}(c_i/B)$. Then $q = p \cup \mathrm{tp}(c/B)$ is a complete global type in $S(M)$.*

Note that $q$ may not be $\mathrm{Aut}(M/A)$-invariant. We will see an example soon.

*Proof.* Proving that $q$ is consistent is easy: choose $d$ a realization of $p$ in some strongly $|B|^+$-homogeneous elementary extension $N$ of $M$, choose $\sigma$ in $\mathrm{Aut}(N/B)$ such that $\sigma(c_i) = d$, then $\sigma(c)$ is a realization of $q$.

Let $\alpha \vDash q$, and $\delta' = \Delta(\alpha_i - b_i)$. As $\alpha_i \vDash p$, we have:

$$\delta' \in (\Delta(G) \smallsetminus \Delta(\underline{G}_M)) \cup (\Delta(\overline{H}^M) \smallsetminus \Delta(H))$$

in particular it does not belong to $\Delta(M)$. For every $j \neq i$, not only do we have $\Delta(\alpha_j - b_j) \notin \Delta(B)$ as $\alpha_j \equiv_B c_j$, but in fact we have $\Delta(\alpha_j - b_j) = \delta'$ by the second condition of proposition 3.3.15. In particular $\alpha$ is a $\mathrm{val}^3_M$-block, and we have:

$$\mathbb{P}^+(\alpha/M) = \mathbb{P}^+((\alpha_j - b_j)_j) = \mathbb{P}^+(\alpha/B) = \mathbb{P}^+(c/B)$$

The last equality follows from proposition 3.3.15, and the others follow from the definitions. We established that, for all $\beta \vDash q$, $\beta$ is a $\mathrm{val}^3_M$-block, and $\mathbb{P}^+(\beta/M) = \mathbb{P}^+(\alpha/M)$. By proposition 3.3.15 (with $c, B$ replaced by $\alpha, M$), in order to check that $q$ is complete, it is enough to show that for every $j$, the pushforward $q_j$ of $q$ by the $j$-th projection, is a complete type in $S(M)$.

Let $f \in \mathrm{LC}^N(\mathbb{Q})$ be the projection on the $j$-th coordinate.

- On one hand, we have $q_j(x_j) \vDash x_j > b_j$ if and only if $f(\mathbb{P}^+(c/B)) > 0$.
- On the other hand, we have:

  $$q(x) \vDash \mathrm{ct}(\Delta(x_j - b_j)/\Delta(M)) = \mathrm{ct}(\Delta(x_i - b_i)/\Delta(M)) = \mathrm{ct}(\delta'/\Delta(M))$$

  The first equality follows from the fact that the realizations of $q$ are $\mathrm{val}^3_M$-blocks, and the last equality follows from corollary 3.2.20, as $p$ is complete. In particular, $q_j(x_j) \vDash \mathrm{ct}(\Delta(x_j - b_j)/\Delta(M)) = \mathrm{ct}(\delta'/\Delta(M))$.

Then, we can apply corollary 3.2.20 to show that $q_j$, and thus $q$, is complete. □

**Proposition 3.3.17.** *Let $i$ (if it exists) be such that $c_i$ is Archimedean over $A$. Then $\mathrm{tp}(c_i/B)$ has exactly one global $\mathrm{Aut}(M/A)$-invariant extension $p$, and $p \cup \mathrm{tp}(c/B)$ is a complete, $\mathrm{Aut}(M/A)$-invariant global type in $S(M)$.*

*Proof.* Suppose for simplicity $c_i > b_i$ (else apply the proposition with $c_i, b_i$ replaced by $-c_i, -b_i$). Recall that $c \mathop{\downarrow}^{\mathbf{cut}}_A B$, thus $\mathrm{tp}(c_i/B)$ must have at least one global $\mathrm{Aut}(M/A)$-invariant extension. By lemma 3.2.23, $G(c_i/A) = G$, and the types $p_>, p_<$ over $M$ corresponding respectively to the cuts:

$$\mathrm{ct}_>(\mathrm{ct}(c_i/A)/M) = \mathrm{ct}_>(b_i \bmod G(M)/M)$$

$$\mathrm{ct}_<(\mathrm{ct}(c_i/A)/M) = \mathrm{ct}_<(b_i \bmod G(M)/M)$$

are clearly global $\mathrm{Aut}(M/A)$-invariant extensions of $\mathrm{tp}(c_i/A)$. These two types are the only ones implying that $x - b_i \in G \smallsetminus \underline{G}_M$. As $c_i > b_i$, only $p_>$ extends $\mathrm{tp}(c_i/B)$. By lemma 3.2.23, there are no other global $\mathrm{Aut}(M/A)$-invariant extension of $\mathrm{tp}(c_i/B)$.

Let $q = p_> \cup \mathrm{tp}(c/B)$, which is complete by lemma 3.3.16. Let us show that $q$ is $\mathrm{Aut}(M/A)$-invariant. As $\mathbb{P}^+(c/B)$ is $\mathbb{Q}$-free, the realizations of $q$ are $\mathrm{val}^3_M$-separated by proposition 3.3.13. As a result, for every non-zero $f \in \mathrm{LC}^N(\mathbb{Q})$, $q(x) \vDash f(x) - f(b) \in G \smallsetminus \underline{G}_M$. As $G = G(c_i/A)$ and $\underline{G}_M$ are $\mathrm{Aut}(M/A)$-invariant, we only have to show that $X = f(b) \bmod \underline{G}_M$ is $\mathrm{Aut}(M/A)$-invariant as a subset of $M$, because:

$$q(x) \vDash \mathrm{ct}(f(x)/M) \in \{\mathrm{ct}_>(X/M), \mathrm{ct}_<(X/M)\}$$

However, as a subset of $M$, $X$ coincides with $f(b) \bmod G(M)$. If $f(c)$ is Archimedean over $A$, then $X = \mathrm{ct}(f(c)/A)(M)$, which is $\mathrm{Aut}(M/A)$-invariant. Else, by definition, there exists $a \in A$ such that $f(c) - a \in G(M)$.

However, we also have $f(c) - f(b) \in G(M)$, therefore $X = a \bmod G(M)$, which is $\mathrm{Aut}(M/A)$-invariant. This concludes the proof. $\square$

Note that $p_>$, and thus $q$, is outer (with respect to definition 3.2.43).

*Remark* 3.3.18. If $j \neq i$, and $c_j$ is ramified over $A$, then $\mathrm{tp}(c_j/B)$ might have two global $\mathrm{Aut}(M/A)$-invariant extensions, and the union of each of those extensions with $\mathrm{tp}(c/B)$ would be a consistent, complete type in $S(M)$. However, one of the two will not be $\mathrm{Aut}(M/A)$-invariant, as it will not be consistent with $p_>$. In the above proof, we would encounter a problem where we would want to show that $X$ is $\mathrm{Aut}(M/A)$-invariant, because we would have to deal with the case where $f(c)$ is Archimedean over $A$, and $X = f(b) \bmod H$.

This is exactly what happens in the following example:

*Example* 3.3.19. Let $A = \mathbb{Q}$, $B = \mathbb{Q} + \mathbb{Q}\sqrt{2}$, $\varepsilon$ a positive infinitesimal element, $c_1 = \sqrt{2} + \varepsilon$, and $c_2 = \varepsilon \cdot \sqrt{2}$. Then $c_1c_2 \mathop{\downarrow}^{\mathbf{cut}}_A B$, and $c_1c_2$ is a $\mathrm{val}^3_B$-separated block. Moreover, $\mathrm{tp}(c_1/B)$ has exactly one global $\mathrm{Aut}(M/A)$-invariant extension, and $\mathrm{tp}(c_2/B)$ has two, one of which is $p$, the type of a positive element that is infinitesimal with respect to $M$. Then $q = p \cup \mathrm{tp}(c_1c_2/B)$ is complete and consistent in $S(M)$, but it is not $\mathrm{Aut}(M/A)$-invariant. The reason is that we have $q(x,y) \vDash \sqrt{2} \leqslant x \leqslant \sqrt{2} + \varepsilon$.

**Proposition 3.3.20.** *Suppose $c_i$ is ramified over $A$ for every $i$. Choose an arbitrary index $i$. Let $H' = H(c/A), G' = G(c/A)$. If $H \neq H'$ or $G \neq G'$, then* $\mathrm{tp}(c_i/B)$ *has exactly one global* $\mathrm{Aut}(M/A)$*-invariant extension, else it has exactly two. Moreover, if $p$ is such an extension, then $p \cup \mathrm{tp}(c/B)$ is a complete* $\mathrm{Aut}(M/A)$*-invariant type in $S(M)$.*

*Proof.* We can assume $b_j \in A$ for every $j$. Suppose $c_i > b_i$. As $c \mathop{\downarrow}^{\mathbf{cut}}_A B$, recall that $H' \leqslant H \leqslant G \leqslant G'$. Let $p_<, p_>$ be the global 1-types corresponding to the respective cuts $\mathrm{ct}_<(\mathrm{ct}(c_i/A)/M) = \mathrm{ct}_>(b_i \bmod H'/M)$, $\mathrm{ct}_>(\mathrm{ct}(c_i/A)/M) = \mathrm{ct}_>(b_i \bmod \underline{G'}_M/M)$. One can see that these two types are exactly all the global $\mathrm{Aut}(M/A)$-invariant extensions of $\mathrm{tp}(c_i/A)$. Moreover $p_<$ (resp. $p_>$) is consistent with $\mathrm{tp}(c_i/B)$ if and only if $H = H'$ (resp. $G = G'$).

Now, let $p$ be a global $\mathrm{Aut}(M/A)$-invariant extension of $\mathrm{tp}(c_i/B)$, and $q = p \cup \mathrm{tp}(c/B)$. Then $q$ is complete by lemma 3.3.16. Let us show that $q$ is $\mathrm{Aut}(M/A)$-invariant.

Let $f \in \mathrm{LC}^n(\mathbb{Q})$ be non-zero. One just has to show that $f(b) \bmod H'(M)$ and $f(b) \bmod G'(M)$ are $\mathrm{Aut}(M/A)$-invariant. However, $H'(M)$ and $G'(M)$ are both $\mathrm{Aut}(M/A)$-invariant, and the $b_i$ were now chosen in $A$. This concludes the proof. □

Note that $p_<$ (thus the extension $p_< \cup \mathrm{tp}(c/B)$) is inner, while $p_>$ (thus $p_> \cup \mathrm{tp}(c/B)$) is outer.

*Remark* 3.3.21. The groups $H'$ and $G'$ have the same points in $A$. By quantifier elimination, it is not hard to show that $G'$ does not admit an $A$-($\vee$/type)-definable proper convex subgroup that strictly contains $H'$. Recall that $H' \leqslant H \leqslant G \leqslant G'$ by cut-independence. Then, if $G \neq G'$ then $G$ cannot be $A$-type-definable. Similarly, if $H \neq H'$, then $H$ is not $A$-$\vee$-definable.

**Corollary 3.3.22.** *If one of the $c_i$ is Archimedean over A, or H is not $A$-$\vee$-definable, or G is not A-type-definable, then* $\mathrm{tp}(c/B)$ *has exactly one global* $\mathrm{Aut}(M/A)$*-invariant extension. Else, it has exactly two. More precisely, if G is A-type-definable, then* $\mathrm{tp}(c/B)$ *admits a global* $\mathrm{Aut}(M/A)$*-invariant extension that is outer ; while if H is A-$\vee$-definable and none of the $c_i$ is Archimedean over A, then* $\mathrm{tp}(c/B)$ *has a global* $\mathrm{Aut}(M/A)$*-invariant extension that is inner.*

### 3.3.4 Gluing everything together

We recall that we adopt assumptions 3.1.14. We now have enough tools to prove (C4). We can actually prove the following more precise statement:

**Proposition 3.3.23.** *Suppose $c = (c_{ik})_{ik}$ is a finite $\mathrm{val}^2_B$-block of B-ramified points from M that is $\mathrm{val}^3_B$-separated, such that its $\mathrm{val}^3_B$-blocks are the $(c_i)_{i \in E} = ((c_{ik})_k)_{i \in E}$. Let $G = G(c/B)$, $H = H(c/B)$. Suppose the indices $i \in E$ are ordered such that $i < l$ if and only if $\delta(c_i/B) < \delta(c_l/B)$. Define $I, O \subseteq E$ as follows:*

- *If H is not A-$\vee$-definable, then $I = \varnothing, O = E$.*
- *If G is not A-type-definable, then $I = E, O = \varnothing$.*
- *Else, $I = \varnothing$, and O is the least final segment of E containing the set of indices $o \in E$ such that one of the $(c_{ok})_k$ is Archimedean over A.*

*Either way, let $J = E \setminus (I \cup O)$. Then the following conditions hold:*

1. *For all $i \in I$,* $\mathrm{tp}(c_i/B)$ *has exactly one global* $\mathrm{Aut}(M/A)$*-invariant extension $p_i$, and it is inner.*

2. *For all* $o \in O$, $\mathrm{tp}(c_o/B)$ *has at least one global* $\mathrm{Aut}(M/A)$*-invariant extension, of which exactly one is outer, denote it by* $q_o$.

3. *For all* $j \in J$, $\mathrm{tp}(c_j/B)$ *has exactly two global* $\mathrm{Aut}(M/A)$*-invariant extensions:* $r_j$*, which is inner, and* $s_j$*, which is outer.*

4. *The map* $K \longmapsto (p_i)_{i\in I}, (r_j)_{j\in K}, (s_j)_{j\in J\smallsetminus K}, (q_o)_{o\in O}$ *is a bijection from the set of initial segments of* $J$ *to the set of strong* $\mathrm{val}^3_B$*-block extensions of* $c$*. In particular,* $c$ *admits exactly* $|J|+1$ *strong* $\mathrm{val}^3_B$*-block extensions.*

*Proof.* Note that $I, O, J$ are pairwise-disjoint convex subsets of $E$ which cover $E$. In fact, $I$ is an initial segment of $E$, $O$ is a final segment, and $J$ is in-between. The first two conditions are easy consequences of corollary 3.3.22. The third condition also easily follows from the fact that for all $j \in J$, none of the $(c_{jk})_k$ is Archimedean over $A$. There remains to prove the last condition.

Suppose we have $o \in O$ such that $\mathrm{tp}(c_o/B)$ has a global $\mathrm{Aut}(M/A)$-invariant extension $q'$ which is inner. Let us show that $q'$ cannot be extended to a strong $\mathrm{val}^3_B$-block extension of $c$. Let $(q'_e)_{e\in E}$ be some strong block extension of $c$. By definition of $O$, there exists $o'$, $k$ such that $o' \leqslant o$ and $c_{o'k}$ is Archimedean over $A$. As a result, $q_{o'}$ is the unique global $\mathrm{Aut}(M/A)$-invariant extension of $\mathrm{tp}(c_{o'}/B)$ by corollary 3.3.22, therefore $q'_{o'} = q_{o'}$. As $o' \leqslant o$, by lemma 3.2.44, $q'_o$ must be outer, thus $q'_o \neq q'$ and we are done.

As a result, all the strong block extensions of $c$ must extend $(p_i)_{i\in I}, (q_o)_{o\in O}$. Then, each such extension $(q'_e)_e$ identifies with a subset of $J$: the set of all $j$ such that $q'_j = r_j$. This concludes the proof, as lemma 3.2.44 tells us that the valid choices are exactly the initial segments of $J$. □

**Corollary 3.3.24.** *If* $c$ *is a finite* $\mathrm{val}^3_B$*-separated family from* $M$ *such that* $c \mathop{\downarrow}^{\mathbf{cut}}_A B$*, then* $c$ *admits a strong* $\mathrm{val}^3_B$*-block extension.*

*Proof.* By proposition 3.2.36, it is enough to show that each $\mathrm{val}^2_B$-block of $c$ has a $\mathrm{val}^3_B$-block extension. The Archimedean $\mathrm{val}^2_B$-blocks of $c$ coincide with its Archimedean $\mathrm{val}^3_B$-blocks, by definition of $\mathrm{val}^3_B$. We find a $\mathrm{val}^3_B$-block extension of those blocks by corollary 3.3.5 (the valid choices are built explicitly in proposition 3.3.4), and we find a block extension of the $B$-ramified $\mathrm{val}^2_B$-blocks by proposition 3.3.23. □

We just proved the property (C4), but we are actually very close to get a full classification of the space of global $\mathrm{Aut}(M/A)$-invariant extensions of $\text{tp}(c/B)$. We achieve this in the remainder of this section, but these results are not necessary for our main results about forking.

**Lemma 3.3.25.** *Suppose $c = (c_{ij})_{ij}$ is a finite $\text{val}^2_B$-block of $B$-ramified points from $M$ that is $\text{val}^3_B$-separated, whose $\text{val}^3_B$-blocks are the $(c_i)_i = ((c_{ij})_j)_i$. Let $(p_i)_i$ be a strong $\text{val}^3_B$-block extension of $c$. Then $q = \bigcup_i p_i \cup \text{tp}(c/B)$ is a complete type in $S(M)$.*

*Proof.* Suppose the indices $i$ are ordered such that $i < k$ if and only if we have $\delta(c_i/B) < \delta(c_k/B)$. Let $b_{ij} \in \text{ram}(c_{ij}/B)$. We know from proposition 3.2.45 that $q$ is consistent. Let $\alpha, \beta \vDash q$. Then each $\alpha_i, \beta_i$ is a $\text{val}^3_M$-block of $M$-ramified points, and:

$$\Delta(\alpha_{ij} - b_{ij}) = \delta(\alpha_i/M), \Delta(\beta_{ij} - b_{ij}) = \delta(\beta_i/M)$$

for all $i$ and $j$. Moreover, for $i < k$, and arbitrary indices $jl$, we have $\text{tp}(c/B) \vDash |x_{kj} - b_{kj}| > n|x_{il} - b_{kl}|$ for all $n > 0$, therefore $\delta(\alpha_i/M) < \delta(\alpha_k/M)$ and $\delta(\beta_i/M) < \delta(\beta_k/M)$. We can apply proposition 3.2.41 to get $\alpha \equiv_M \beta$, which concludes the proof. □

*Remark* 3.3.26. In the setting of lemma 3.3.25, by proposition 3.3.23, the Stone space of all the global $\text{Aut}(M/A)$-invariant extensions of $\text{tp}(c/B)$ is finite Hausdorff, hence it is discrete.

*Remark* 3.3.27. Let $X$ be a topological space. Suppose $(O_i)_i$ is an open partition of $X$. Then the open subsets of $X$ are exactly the subsets that can be written $\bigcup_i O'_i$, with $O'_i$ an open subset of $O_i$. As a result, $X = \coprod_i O_i$.

**Corollary 3.3.28.** *Let $X$ be a topological space, and $Y$ a discrete topological space such that we have a continuous map $X \longrightarrow Y$. Then $X$ is the coproduct of the fibers.*

**Lemma 3.3.29.** *Suppose $c = (c_{ij})_{ij}(c'_k)_k$ is a finite, $\text{val}^3_B$-separated, $\text{val}^1_B$-block such that $c'$ is the family of its $B$-Archimedean points, and its $B$-ramified $\text{val}^3_B$-blocks are the $(c_i)_i = ((c_{ij})_j)_i$. Let $(p_i)_i$ be a strong $\text{val}^3_B$-block extension of $(c_i)_i$, and suppose all the $p_i$ are outer. Let $G = G(c/B) = G(c'/A)$. For each $i$ and $j$, choose $b_{ij} \in \text{ram}(c_{ij}/B)$. Let $d = (b_{ij})_{ij}(c'_k)_k$. Let $q_k$ be a global $\text{Aut}(M/A)$-invariant extension of $\text{tp}(c'_k/B)$, and let $q$ be a complete global extension of $\text{tp}(c/B)$ which extends $\bigcup_i p_i$ and $\bigcup_k q_k$. If for all non-zero $f \in \text{LC}^n(\mathbb{Q})$ we have $q(x) \vDash f(x) - f(d) \notin \underline{G}_M$, then $q$ is $\text{Aut}(M/A)$-invariant.*

*Proof.* Let $f \in \mathrm{LC}^n(\mathbb{Q})$ be non-zero. As is the proof of lemma 3.3.3, it suffices to show that $q(x) \models f(x) \in X$, with $X$ some $\mathrm{Aut}(M/A)$-invariant cut over $M$. By $\mathrm{val}^3_B$-separatedness, we have $G(f(c)/B) = G$. Write $f = g + h$, with $g \in \mathrm{LC}^{|(c_{ij})_{ij}|}(\mathbb{Q})$, $h \in \mathrm{LC}^{|(c'_k)_k|}(\mathbb{Q})$. By lemma 3.2.44, $(p_i)_i$ is a strong block extension of $(c_{ij})_{ij}$. By lemma 3.3.25, the pushforward of $q$ by the projection over the coordinates $(x_{ij})_{ij}$ is $\mathrm{Aut}(M/A)$-invariant, thus there is $X$ an $\mathrm{Aut}(M/A)$-invariant cut over $M$ such that $q(x) \models g((x_{ij})_{ij}) \in X$. If $h = 0$, then $q(x) \models f(x) \in X$. Now suppose $h \neq 0$. If $g = 0$, then we are done as the proof is the same as in lemma 3.3.3, so suppose $g \neq 0$. By $\mathrm{val}^3_B$-separatedness, $f(c)$ is Archimedean over $B$, and $G(f(c)/B) = G$. By lemma 3.2.23, the two global $\mathrm{Aut}(M/A)$-invariant extensions of $\mathrm{tp}(h((c'_k)_k)/B)$ are the translates by $h((c'_k)_k)$ of the connected components of $G \smallsetminus \underline{G}_M$. By lemma 3.3.3, the pushforward of $q$ by the projection on the coordinates $(x'_k)_k$ is $\mathrm{Aut}(M/A)$-invariant, therefore there is $Y$ one of those two translates such that $q(x) \models h((x'_k)_k) \in Y$. As the $p_i$ are outer, and by $\mathrm{val}^3_B$-separatedness, $X$ is the translate by $g((b_{ij})_{ij})$ of some connected component of $G \smallsetminus \underline{G}_M$. Then, we have $q(x) \models f(x) - f(d) \in G \smallsetminus \underline{G}_M$, i.e. $q(x) \models f(x) \in L \cup R$, with $L = \mathrm{ct}_<\big(\big[g((b_{ij})_{ij}) + h((c'_k)_k) \bmod G(M)\big]/M\big)$, and $R = \mathrm{ct}_>\big(\big[g((b_{ij})_{ij}) + h((c'_k)_k) \bmod G(M)\big]/M\big)$. By $A$-invariance of $X$ and $Y$, both $L$ and $R$ are $\mathrm{Aut}(M/A)$-invariant, and concluding the proof. □

**Lemma 3.3.30.** *Suppose $c = (c_{ij})_{ij}(c'_k)_k$ is a finite, $\mathrm{val}^3_B$-separated, $\mathrm{val}^1_B$-block such that $c'$ is the family of its $B$-Archimedean points, and its $B$-ramified $\mathrm{val}^3_B$-blocks are the $(c_i)_i = ((c_{ij})_j)_i$. Let $(p_i)_i$ be a strong $\mathrm{val}^3_B$-block extension of $(c_i)_i$. Let $q$ be the global $\mathrm{Aut}(M/A)$-invariant extension of $\mathrm{tp}((c_{ij})_{ij}/B)$ which extends $\bigcup_i p_i$ ($q$ exists and is unique, by lemma 3.3.25).*

*Then the following description yields a full classification of the Stone space $F_1$ of global $\mathrm{Aut}(M/A)$-invariant types extending $q \cup \mathrm{tp}(c/B)$:*

*Suppose $c'$ is not empty (or else $q \cup \mathrm{tp}(c/B)$ is already complete). Let $O$ be the set of indices $o$ such that $p_o$ is outer. Let $G = G(c/B) = G(c'/A)$. Define $N$, $H$, $N'$ as in the statement of proposition 3.3.4. For each $i$ and $j$, choose $b_{ij} \in \mathrm{ram}(c_{ij}/B)$, and choose $\alpha$ a realization of the projection of $q$ onto the coordinates $(x_o)_{o\in O}$. Let $F_2$ be the closed subspace of $S_{(x_o)_{o\in O},x'}(\{0\})$ corresponding to the types of $\mathbb{Q}-free$ families satisfying $\mathrm{tp}^{N'}((\alpha_{oj} - b_{oj} \bmod H)_{oj}/\{0\})$. Then the map:*

$$\mathrm{tp}^N((\gamma_{ij})_{i\notin O,j}\alpha(\beta_k)_k/M) \longmapsto \mathrm{tp}^{N'}((\alpha_{oj}-b_{oj} \bmod H)_{o\in O,j}(\beta_k - c'_k \bmod H)_k/\varnothing)$$

*is a well-defined homeomorphism $F_1 \longrightarrow F_2$.*

*Proof.* Let $d = ((b_{oj})_{o\in O,j}(c'_k)_k)$, and let $\pi$ be the projection:

$$(x_i)_i, x' \longmapsto (x_o)_{o\in O}, x'$$

Let us show first that the map from the statement factors through a homeomorphism $\pi(F_1) \longrightarrow F_2$. The proof is very similar to that of proposition 3.3.4.

- Let $\beta \in N$ be such that $\mathrm{tp}(\alpha\beta/M) \in \pi(F_1)$. Let $f_1$ in $\mathrm{LC}^{|\alpha|}(\mathbb{Q})$, $f_2$ in $\mathrm{LC}^{|\beta|}(\mathbb{Q})$ be such that $f_1 + f_2 \neq 0$. Just as in proposition 3.3.4, in order to show that the map $g : \pi(F_1) \longrightarrow F_2$ given by the statement is well-defined, it suffices to show that $f_1(\alpha) + f_2(\beta) - (f_1 + f_2)(d) \notin \underline{G}_M$. If $f_2 \neq 0$, then $(f_1 + f_2)(d)$ is Archimedean over $A$, which concludes the proof by lemma 3.2.23 just as in proposition 3.3.4. As all the $p_o$ are outer, and $\alpha$ is $\mathrm{val}^3_B$-separated (as $(c_{oj})_{oj} \equiv_B \alpha$ is), we have $f_1(\alpha) - (f_1 + 0)(d) \notin \underline{G}_M$, proving that $g$ is well-defined.

- Let us show that $g$ is injective. Let $p_1$, $p_2 \in \pi(F_1)$ be distinct, and $f_1$, $f_2 \in \mathrm{LC}(\mathbb{Q})$ be as in the above paragraph. We use the same case disjunction ($f_2 \neq 0$ or $f_2 = 0$) as in the above paragraph to show that the pushforward of $p_i$ by $f_1 + f_2$ is (as a cut over $M$) a translate by $(f_1 + f_2)(d)$ of one of the two connected components of $G \setminus \underline{G}_M$. Then, it follows from the same reasoning as in proposition 3.3.4 that $g$ is injective.

- Let us now prove surjectivity. Let $\beta'$ be a tuple from $N'$ such that $r' = \mathrm{tp}((\alpha_{oj} - b_{oj} \bmod H)_{oj}\beta'/\{0\}) \in F_2$. Let $\beta$ be a preimage of $\beta'$ by the projection $N \longrightarrow N'$. Let $d_1, d_2 \in N$ be such that we have $\Delta(d_1) > \Delta(\beta_i)$, $\Delta(d_k) > \Delta(\alpha_{oj} - b_{oj})$, and $d_2 \notin \underline{G}_M$ for all $i, j, k, o$. Let $\sigma \in \mathrm{Aut}\left(N/H + \sum_{oj} \mathbb{Q}\cdot(\alpha_{oj} - b_{oj})\right)$ be such that $\sigma(|d_1|) = |d_2|$, and let $r = \mathrm{tp}(\alpha(c'_k + \sigma(\beta_k))_k/M)$. Then $r$ is $\mathrm{Aut}(M/A)$-invariant by lemma 3.3.29, and, just as in proposition 3.3.4, $r$ is a preimage of $r'$ by $g$.

- We conclude that $g$ is a homeomorphism $\pi(F_1) \longrightarrow F_2$.

It remains to show that $F_1 \longrightarrow \pi(F_1)$ is homeo. By lemma 3.3.25, recall that $p = \mathrm{tp}((c_i)_{i\notin O}/B) \cup \bigcup_{i\notin O} p_i$ is a complete global type. Let $q' \in \pi(F_1)$. Then, if $\gamma_1 \vDash p, \gamma_2 \vDash q'$, one can note that $G(\gamma_2/M) = G$, while $G(\gamma_1/M) =$ $\overline{H(c/B)}^M < G$ (this is the difference between being inner and outer). As a result, the $\mathrm{val}^1_M$(and thus $\mathrm{val}^2_M$)-values of the elements of $\gamma_1$ (there is only one such value, but we do not care) are different from the values of the elements of $\gamma_2$. By proposition 3.2.36, $p$ and $q$ are weakly orthogonal, which shows that $F_1 \longrightarrow \pi(F_1)$ is homeo. □

When put together with remark 3.3.26 and corollary 3.3.28, we get that the space of global $\mathrm{Aut}(M/A)$-invariant extensions of the type of any finite, $\mathrm{val}^3_B$-separated $\mathrm{val}^1_B$-block $(c_i)_i, c'$ that is cut-independent is the finite coproduct of the topological spaces described in lemma 3.3.30. The continuous map we consider for corollary 3.3.28 is the projection $\pi: \ ((x_i)_i, x') \longmapsto (x_i)_i$. Note that, by 3.3.30, any global $\mathrm{Aut}(M/A)$-invariant extension of the type of $(c_i)_i$ over $B$ extends to one of the type of $((c_i)_i, c')$ (as one can always extend the type of some $\mathbb{Q}$-free family to an arbitrary type of a larger $\mathbb{Q}$-free family), so $\pi$ is onto the finite space described in remark 3.3.26. Now, given a finite, $\mathrm{val}^3_B$-separated tuple $c$ that is cut-independent, we just described $(S_i)_i$, the Stone spaces of the global $\mathrm{Aut}(M/A)$-invariant extensions of the type of each $\mathrm{val}^1_B$-block of $c$. By (C7) and lemma 3.1.19, the Stone space of every global $\mathrm{Aut}(M/A)$-invariant extension of the type of the full tuple $c$ is exactly the finite direct product of the $S_i$. Let us rephrase all of these conclusions (with their hypothesis) in a single statement for good measure:

*Remark* 3.3.31. In DOAG, let $M$ be a sufficiently saturated and strongly homogeneous model, and $C \geqslant A \leqslant B$ be $\mathbb{Q}$-vector subspaces, such that $C \mathop{\downarrow}^{\mathbf{cut}}_A B$. Let $c = (c_{ij})_{ij}$ be some $\mathrm{val}^3_B$-separated family of $C$, with $(c_i)_i = ((c_{ij})_j)_i$ its $\mathrm{val}^1_B$-blocks. For each $i$, let $d_i$ be the tuple of $B$-ramified points of $c_i$. Let $F_i$ be the space of all global $\mathrm{Aut}(M/A)$-invariant extensions of $\mathrm{tp}(d_i/B)$, which is finite, and fully described by proposition 3.3.23 and lemma 3.3.25. For each $i$, and for each $q \in F_i$, let $S_q$ be the space (described in lemma 3.3.30) of global $\mathrm{Aut}(M/A)$-invariant extensions of $\mathrm{tp}(c_i/B)$ which extend $q$. We recall that $S_i$ is homeomorphic to some closed set of parameter-free DOAG-types. Then the space $S_i$ of all global $\mathrm{Aut}(M/A)$-invariant extensions of $\mathrm{tp}(c_i/B)$ is of course the finite disjoint union $\bigcup_{q \in F_i} S_q$, which is naturally homeomorphic to the coproduct $\coprod_{q \in F_i} S_q$. Lastly, the whole space of global $\mathrm{Aut}(M/A)$-invariant extensions of $\mathrm{tp}(c/B)$ is naturally homeomorphic to $\prod_i S_i$.

In conclusion, we established a natural and explicit homeomorphism be-tween the space of global $\mathrm{Aut}(M/A)$-invariant extensions of $\mathrm{tp}(c/B)$, and a topological space of the form $\prod_i \coprod_j F_{ij}$, with $(F_{ij})_{ij}$ some closed sets of parameter-free DOAG-types.

Moreover, as definable bijections are homeomorphisms in $S(M)$, this homeomorphism type is kept when replacing $c$ by an $A$-interdefinable tuple. Note that we show in the next section that every cut-independent tuple is $A$-interdefinable with a $\mathrm{val}^3_B$-separated, cut-independent family. As a result, this will yield a complete description of the Stone space of global $\mathrm{Aut}(M/A)$-invariant extensions of any non-forking type.

## 3.4 Normal forms

We built in the last two sections a nice framework to get a fine understanding of a particular class of tuples, the $\mathrm{val}^3_B$-separated families. What remains for us to do in order to prove theorem 3.1.5 is to prove (C3), which allows us to always reduce to the case where we deal with such a nice family.

We adopt assumptions 3.1.22 in this section.

We have to show that $c$ is $A$-interdefinable with a $\mathrm{val}^3_B$-separated family $d$, which is cut-independent by remark 3.1.21. In fact, $d$ will also be $\mathrm{val}^3_A$-separated. This could be useful, as we saw in proposition 3.3.23 that we have to keep track of which points from $d$ are Archimedean or ramified over $A$ in order to have a good understanding of the global $\mathrm{Aut}(M/A)$-invariant extensions of its type. By remark 3.1.23, $d$ is the image of $c$ by a composition of $A$-translations, and a map from $\mathrm{GL}_n(\mathbb{Q})$. In order to build $d$, we start with a particular basis of $\mathrm{dcl}(Ac)$ (remember that the definable closure of a set is the $\mathbb{Q}$-vector subspace generated by said set in DOAG), then we go through a process of several steps, where we apply $A$-translations and maps from $\mathrm{GL}_n(\mathbb{Q})$ to this basis. At the end of each step, our current tuple satisfies an additional property from a list of conditions, the conjunction of which is a sufficient condition to be a $\mathrm{val}^3_B$ and $\mathrm{val}^3_A$-separated family.

*Remark* 3.4.1. Recall that, if $d \in M$ is a $B$-Archimedean singleton such that $d \mathop{\downarrow}^{\mathbf{cut}}_A B$, then $d$ must be Archimedean over $A$ (see remark 3.3.1). In particular, $G(d/B) = G(d/A)$.

The following definition allows us to enumerate our families in a way that helps us to simultaneously consider their $\mathrm{val}^3_B$-blocks and their $\mathrm{val}^3_A$-blocks.

The elements of our families are indexed by the subgroups $G(-/B)$, and the eventual Archimedean classes $\delta(-/B)$ that they bring. We use the single quote $'$ to refer to elements that are Archimedean over $A$ and ramified over $B$. We use the tilde $\tilde{}$ to refer to elements that are ramified over $A$ (and thus over $B$ if the family is cut-independent from $B$ over $A$). By elimination, the other elements are those that are Archimedean over $B$ and $A$.

**Definition 3.4.2.** We say that the family $(d_{Gi})_{Gi}(d'_{G'\delta'i'})_{G'\delta'i'}(\tilde{d}_{\tilde{G}\tilde{\delta}\tilde{i}})_{\tilde{G}\tilde{\delta}\tilde{i}}$ is *under normal enumeration* (with respect to $A$ and $B$) if the following conditions hold:

- All the $d_{Gi}$ are Archimedean over $B$.
- All the $d'_{G'\delta'i'}$ are Archimedean over $A$ and ramified over $B$.
- All the $\tilde{d}_{\tilde{G}\tilde{\delta}\tilde{i}}$ are ramified over $A$.
- For all $G$, $i$, we have $G(d_{Gi}/B) = G$.
- For all $G'$, $\delta'$, $i'$, we have $G(d'_{G'\delta'i'}/B) = G'$ and $\delta(d'_{G'\delta'i'}/B) = \delta'$.
- For all $\tilde{G}$, $\tilde{\delta}$, $\tilde{i}$, we have $G(\tilde{d}_{\tilde{G}\tilde{\delta}\tilde{i}}/B) = \tilde{G}$, and $\delta(\tilde{d}_{\tilde{G}\tilde{\delta}\tilde{i}}/A) = \tilde{\delta}$.

Given such a normal enumeration, the ramified $\mathrm{val}^3_B$-blocks of the family are the $((d'_{G'\delta'i'})_{i'}(\tilde{d}_{G'\delta'\tilde{i}})_{\tilde{i}})_{G'\delta'}$, its Archimedean $\mathrm{val}^3_B$-blocks are the $((d_{Gi})_i)_G$, its ramified $\mathrm{val}^3_A$-blocks are the $((\tilde{d}_{\tilde{G}\tilde{\delta}\tilde{i}})_{\tilde{i}})_{\tilde{G}\tilde{\delta}}$, and its Archimedean $\mathrm{val}^3_A$-blocks are the $((d_{Gi})_i(d'_{G\delta'i'})_{\delta'i'})_G$. One can note that the presence of the elements of $d'$ is the reason why none of the valuations $\mathrm{val}^3_A$, $\mathrm{val}^3_B$ refines the other.

Let us now give the list of properties that we want to satisfy:

**Definition 3.4.3.** Let $dd'\tilde{d}$ be a family under normal enumeration with respect to $A$ and $B$. We say that $dd'\tilde{d}$ is *under normal form* with respect to $A$ and $B$ if it is a lift of a $\mathbb{Q}$-free family from ${}^{M}\!/\!{}_{B}$, and the following conditions hold:

(P1) For all $\tilde{G}$ and $\tilde{\delta}$, $(\tilde{d}_{\tilde{G}\tilde{\delta}\tilde{i}})_{\tilde{i}}$ is a lift of a $\mathbb{Q}$-free family from ${}^{M_{\leqslant\tilde{\delta}}}\!/\!{}_{M_{<\tilde{\delta}}}$.

(P2) Any non-trivial linear combination of $dd'$ is Archimedean over $A$.

(P3) For all $G$, any non-trivial linear combination $e$ of the $(d_{Gi})_i(d'_{G\delta'i'})_{\delta'i'}$ satisfies $G(e/A) = G$.

(P4) For all $G$, any non-trivial linear combination $c$ of the tuple $(d_{Gi})_i$ is Archimedean over $B$ (which implies $G(c/B) = G$ if the family is cut-independent, by (P3) and remark 3.4.1).

(P5) For all $G'$ and $\delta'$, any non-trivial linear combination $e'$ of the tuple $(d'_{G'\delta' i'})_{i'}(\tilde{d}_{G'\delta'\tilde{i}})_{\tilde{i}}$ is ramified over $B$, and satisfies $\delta(e'/B) = \delta'$.

*Remark* 3.4.4. (P1) implies that each ramified $\mathrm{val}^3_A$-block is $\mathrm{val}^3_A$-separated. The conjunction of (P2) and (P3) is equivalent to each Archimedean $\mathrm{val}^3_A$-block being $\mathrm{val}^3_A$-separated. Likewise, (P4) (resp. (P5)) is equivalent to each Archimedean (resp. ramified) $\mathrm{val}^3_B$-block being $\mathrm{val}^3_B$-separated. As a result, if we manage to prove that $c$ is $A$-interdefinable with a family under normal form, then this family would be simultaneously $\mathrm{val}^3_A$-separated and $\mathrm{val}^3_B$-separated, and (C3) would follow.

*Remark* 3.4.5. From the definition of the normal form, if $dd'\tilde{d}$ is under normal form with respect to $A$ and $B$, then one can note that any subfamily of $dd'\tilde{d}$ is also under normal form with respect to $A$ and $B$.

**Definition 3.4.6.** Define $C$ to be the $\mathbb{Q}$-vector space generated by $Ac$. We define a *basis* (resp. *free family*) *of* $C$ *over* $A$ as a lift of a basis (resp. free family) of ${C}\diagup{A}$.

*Remark* 3.4.7. Remember that $dim\left({C}\diagup{A}\right)$ is finite. It is easy to see that any concatenation of lifts of bases of ${C_{\leqslant\tilde{\delta}}}\diagup{C_{<\tilde{\delta}}}$, for each $\tilde{\delta} \in \Delta(C) \smallsetminus \Delta(A)$, is a free family over $A$, hence a finite family. In particular, $\Delta(C) \smallsetminus \Delta(A)$ is finite.

**Lemma 3.4.8.** *For each $\tilde{\delta} \in \Delta(C) \smallsetminus \Delta(A)$, let $\tilde{d}_{\tilde{\delta}}$ be a lift of a basis of ${C_{\leqslant\tilde{\delta}}}\diagup{C_{<\tilde{\delta}}}$, and let $\tilde{d}$ be the concatenation of all the $\tilde{d}_{\tilde{\delta}}$. Then, for all $e$ in $C \smallsetminus \mathrm{dcl}(A\tilde{d})$, there exists $\tilde{e} \in \mathrm{dcl}(A\tilde{d})$ such that $e + \tilde{e}$ is Archimedean over $A$.*

*Moreover, if $e$ is ramified over $A$, then we can choose $\tilde{e}$ such that there exists $\tilde{\delta} \in \Delta(C) \smallsetminus \Delta(A)$ for which $\Delta(e + \tilde{e}) < \tilde{\delta} \leqslant \Delta(e)$.*

*Proof.* Let us build $\tilde{e}$ by induction. Let $\tilde{e}_0 = 0 \in \mathrm{dcl}(A\tilde{d})$. If $e + \tilde{e}_n$ is not Archimedean over $A$, then let $a \in \mathrm{ram}(e + \tilde{e}_n/A)$, and $\tilde{\delta}_n = \delta(e + \tilde{e}_n/A)$. By definition of $\tilde{d}$, there exists $v$ a linear combination of $\tilde{d}_{\tilde{\delta}_n}$ for which we have $\Delta(e + \tilde{e}_n - a - v) < \tilde{\delta}_n$, so let us set $\tilde{e}_{n+1} = \tilde{e}_n - a - v$.

By finiteness of $\Delta(C) \smallsetminus \Delta(A)$, the sequence of the $\tilde{\delta}_n$, and hence that of the $\tilde{e}_n$, must eventually stop. The ultimate value $\tilde{e}$ of the $\tilde{e}_n$ witnesses the first part of the statement. If $e$ is ramified over $A$, then $\tilde{e} \neq 0 = \tilde{e}_0$, so $\tilde{\delta}_0$ exists and witnesses the second part of the statement. □

**Lemma 3.4.9.** *For each $\tilde{\delta} \in \Delta(C) \smallsetminus \Delta(A)$, let $\tilde{d}_{\tilde{\delta}}$ be a lift of a basis of $C_{\leqslant\tilde{\delta}} / C_{<\tilde{\delta}}$, and let $\tilde{d}$ be the concatenation of all the $\tilde{d}_{\tilde{\delta}}$. Let $(v_n)_n$ be a sequence such that $(\tilde{d}, v_n)$ is a basis of $C$ over $A$ for each $n$. Suppose that, for all $n$, there exists $v'_1$ a term from $v_n$, $v'_2$ a term from $v_{n+1}$, and $\delta \in \Delta(C) \smallsetminus \Delta(A)$, such that $\Delta(v'_2) < \delta \leqslant \Delta(v'_1)$, and $v_{n+1}$ is obtained from $v_n$ by replacing $v'_1$ by $v'_2$.*

*Then the sequence $(v_n)_n$ must eventually stop.*

*Proof.* Let $N = |\Delta(C) \smallsetminus \Delta(A)| < \omega$, $(\delta_n)_{n<N}$ the strictly increasing enumeration of $\Delta(C) \smallsetminus \Delta(A)$, and $\mathcal{I}$ be the set of intervals:

$$\{]-\infty, \delta_0[\,, [\delta_0, \delta_1[\,, [\delta_1, \delta_2[\,, \dots [\delta_{N-2}, \delta_{N-1}[\,, [\delta_{N-1}, +\infty[\}$$

Define on $\mathcal{I}$ the only total order that makes the canonical projection $f : \Delta(C) \longrightarrow \mathcal{I}$ an order-preserving map. With respect to that order, let $g$ be the unique order-isomorphism $\mathcal{I} \longrightarrow \{0, 1 \dots N\}$. For each $i$ between 0 and $N$, let $F(n,i)$ be the number of terms $v$ from $v_n$ such that $g(f(\Delta(v))) = i$. Finally, let $\rho$ be the following function:

$$n \longmapsto F(n,N) \times \omega^N + F(n, N-1) \times \omega^{N-1} + \dots + F(n,1) \times \omega + F(n,0)$$

then $\rho$ has its values in the well-ordered set $\omega^{N+1}$. By hypothesis, for each $n$, there exists $m < m'$ such that $F(n,m') = F(n+1,m') + 1$, $F(n,m) + 1 = F(n+1,m)$, and $F(n,i) = F(n+1,i)$ for every $i \neq m$, $m'$. As a result, $\rho$ is strictly decreasing, so the well-ordering forces the sequence $(v_n)_n$ to stop eventually. □

**Lemma 3.4.10.** *For each $\tilde{\delta} \in \Delta(C) \smallsetminus \Delta(A)$, let $u_{\tilde{\delta}}$ be a lift of a basis of $C_{\leqslant\tilde{\delta}} / C_{<\tilde{\delta}}$, and let $u$ be the concatenation of all the $u_{\tilde{\delta}}$. Then there exists $dd'\tilde{d}$ a basis of $C$ over $A$, under normal enumeration with respect to $A$ and $B$, for which (P2) holds, and $\tilde{d} = u$ (so (P1) holds as well).*

*Proof.* Let $w_0$ be some basis of $C$ over $A$ containing $u$. Note that $w_n = (u, v_n)$ does not witness the statement if and only if there exists a non-trivial linear combination $v'$ of $v_n$ that is ramified over $A$. In that case, let $v$ be some term of $v_n$ that appears in $v'$, chosen such that $\Delta(v)$ is maximal. By valuation inequality, we must have $\Delta(v) \geqslant \Delta(v')$. Apply lemma 3.4.8 to find $u' \in \mathrm{dcl}(Au)$ and $\delta \in \Delta(C) \setminus \Delta(A)$ for which $\Delta(v' + u') < \delta \leqslant \Delta(v') \leqslant \Delta(v)$. Replace $v$ by $v' + u'$ in $w_n$, let $w_{n+1}$ be the new family obtained. Now, $w_{n+1}$ is also a basis of $C$ over $A$ containing $u$.

The sequence $(v_n)_n$ witnesses the hypothesis of lemma 3.4.9, so it must eventually stop, and its last element witnesses the statement. □

*Remark* 3.4.11. If $w = (u, v)$ witnesses lemma 3.4.10, and $M$ is an invertible matrix of size $|v|$ with coefficients in $\mathbb{Q}$, then $(u, (Mv))$ also witnesses lemma 3.4.10.

**Lemma 3.4.12.** *Let $dd'\tilde{d} = (v, \tilde{d})$ be a free family over $A$ that witnesses (P1) and (P2), and $N = |v|$. Then there exists $M \in GL_N(\mathbb{Q})$ such that $((Mv), \tilde{d})$ witnesses (P3) (and thus (P1), (P2) by remark 3.4.11).*

*Proof.* This is by induction over $N$. (P3) trivially holds if $N = 0, 1$.

For every $G$, let $v_G$ be the $\mathrm{val}^1_A$-block of $v$ of value $G$. Let $G$ be maximal such that $v_G$ is not $\mathrm{val}^3_A$-separated. As long as $v_G$ is not $\mathrm{val}^3_A$-separated, there must exist $e$ a non-trivial linear combination of $v_G$ such that $G(e/A) < G$, then replace some term of $v_G$ that appears in $e$ by $e$ (making that replacement corresponds to replacing $v$ by $Mv$ for some $M \in \mathrm{GL}_N(\mathbb{Q})$). As we have $G(e/A) < G$, these replacements always decrease the size of $v_G$, and leave $v_{>G}$ untouched as the value of $e$ is smaller, so the $\mathrm{val}^1_A$-blocks of value larger than $G$ remain separated, thus we keep the maximality of $G$. As a result, these replacements eventually stop, and they stop before we removed all the terms from $v_G$, because if there is only one term left in $v_G$, then $v_G$ is clearly separated.

By induction hypothesis, we can replace $v_{<G}$ by $Mv_{<G}$ such that the family $((Mv_{<G}), \tilde{d})$ witnesses (P3). It is now clear that $((Mv_{<G}), v_{\geqslant G}, \tilde{d})$ witnesses (P3). The only changes we made to the original family $v$ are operations from $\mathrm{GL}_N(\mathbb{Q})$, concluding the proof. □

*Remark* 3.4.13. Let $w = dd'\tilde{d} = (v, \tilde{d})$ be a free family over $A$ that witnesses (P1), (P2), (P3). Let $v_G$ be the $\mathrm{val}^1_A$-block of $v$ of value $G$. If we replace $v_G$ by $Mv_G$ in $w$ for some invertible matrix $M$, then the new family $w'$ still witnesses (P1), (P2) and (P3).

**Lemma 3.4.14.** *There exists a basis of $C$ over $A$ that witnesses (P1), (P2), (P3), (P4).*

*Proof.* Let $w = dd'\tilde{d}$ be a basis of $C$ over $A$ that witnesses (P1), (P2) and (P3). For each $G$, as long as there is a non-trivial linear combination $e$ of the $(d_{Gi})_i$ that is not Archimedean over $B$, then replace some term of $(d_{Gi})_i$ that appears in $e$ by $e$. Each of these replacements removes a term from $(d_{Gi})_i$ and adds a new term to $(d'_{G\delta' i'})_{\delta' i'}$. This eventually stops as the number of terms of $(d_{Gi})_i$ is finite and strictly decreases at each step.

By remark 3.4.13, the new family $w'$ that we obtain after all these replacements still witnesses (P1), (P2), (P3). We perform these replacements for each $G$, and the new family obtained clearly witnesses (P4). □

*Remark* 3.4.15. Let $w = dd'\tilde{d}$ be a free family over $A$ that witnesses (P1), (P2), (P3), (P4), and $G'$ some index. Then, if we replace $(d'_{G'\delta' i'})_{\delta' i'}$ by $M(d'_{G'\delta' i'})_{\delta' i'}$ in $w$ for some invertible matrix $M$, it is clear that the new family still witnesses (P1), (P2), (P3), (P4).

**Lemma 3.4.16.** *Let $w = dd'\tilde{d}$ be a basis of $C$ over $A$ that witnesses (P1), (P2), (P3), (P4). Let $d'_{G'\delta' k'}$ be a term from $d'$, and $\tilde{e}$ be a linear combination of the $(\tilde{d}_{G'\tilde{\delta}\tilde{i}})_{\tilde{\delta}\tilde{i}}$. If we replace $d'_{G'\delta' k'}$ by $e = d'_{G'\delta' k'} + \tilde{e}$ in $w$, then the new family obtained is still a basis of $C$ over $A$ that witnesses (P1), (P2), (P3), (P4).*

*Proof.* First of all, we have $\mathrm{val}^2_A(d'_{G'\delta' k'}) > \mathrm{val}^2_A(\tilde{e})$, hence $e$ is Archimedean over $A$, and $G(e/A) = G'$. As $\mathrm{val}^2_B(e) \leqslant \max(\mathrm{val}^2_B(d'_{G'\delta' k'}), \mathrm{val}^2_B(\tilde{e}))$, $e$ cannot be Archimedean over $B$. This implies that the $A$-ramified and $B$-Archimedean blocks of our new family are left untouched, hence we keep the properties (P1), (P4). Let $u'$ be the new block of $A$-Archimedean and $B$-ramified points in the normal enumeration. Then $e$ lies in the block $u'_{G'}$. In order to show that (P2) and (P3) hold in the new family $du'\tilde{d}$, we just have to show that the $\mathrm{val}^3_A$-block $d_{G'}u'_{G'}$ is $\mathrm{val}^3_A$-separated. Indeed, it would follow that the whole $du'$ is $\mathrm{val}^3_A$-separated, which is clearly equivalent to (P2) and (P3) being true in $du'\tilde{d}$.

Any non-trivial linear combination of $d_{G'}u'_{G'}$ can be written $u + \lambda \cdot \tilde{e}$, with $u$ a non-trivial linear combination of $d_{G'}d'_{G'}$, and $\lambda \in \mathbb{Q}$. As (P2) and (P3) hold in $dd'\tilde{d}$, we have $G(u/A) = G'$, and $u$ is Archimedean over $A$. Then, we have $\mathrm{val}^3_A(u) > \mathrm{val}^3_A(\tilde{e})$, hence $\mathrm{val}^3_A(u + \lambda \cdot \tilde{e}) = \mathrm{val}^3_A(u)$. This value equals that of the corresponding block, $\mathrm{val}^3_A(d_{G'}d'_{G'})$, as (P2) and (P3) hold in $dd'\tilde{d}$, thus it equals $\mathrm{val}^3_A(d_{G'}u'_{G'})$, concluding the proof. □

**Lemma 3.4.17.** *Let $w = dd'\tilde{d}$ be a basis of $C$ over $A$ that witnesses (P1), (P2), (P3), (P4), and not (P5). Then there exists $G'$, $\delta'$, and two $\mathbb{Q}$-linear maps $f, g \in \mathrm{LC}(\mathbb{Q})$, such that $v = f((d'_{G'\delta' i'})_{i'}) + g((\tilde{d}_{G'\delta'\tilde{i}})_{\tilde{i}})$ is ramified over $B$, and $\delta(v/B) < \delta'$.*

*Moreover, if $d'_{G'\delta' k'}$ appears in $v$, then replacing it by $v$ in $w$ gives a new family $w'$ which is still a basis of $C$ over $A$ that witnesses (P1), (P2), (P3), (P4).*

*Proof.* By definition of (P5), there exists $G'$, $\delta'$, and $\mathbb{Q}$-linear maps $f, g$ in $\mathrm{LC}(\mathbb{Q})$ such that, either $v = f((d'_{G'\delta' i'})_{i'}) + g((\tilde{d}_{G'\delta'\tilde{i}})_{\tilde{i}})$ is not ramified over $B$, or $\delta(v/B) \neq \delta'$.

Note that $f$ must be non-zero, otherwise the above hypothesis fails by (P1). By (P1), (P2), (P3), $dd'\tilde{d}$ is $\mathrm{val}^3_A$-separated. As $f$ is non-zero, $\mathrm{val}^3_A(v) = \mathrm{val}^3_A(d'_{G'})$, thus $v$ is Archimedean over $A$, and $G(v/A) = G' = G(v/B)$. For all $i'$, $\tilde{i}$, we have $d'_{G'\delta' i'}, \tilde{d}_{G'\delta'\tilde{i}} \in B + G'$ (by the definition of being ramified over $B$), so $v \in B + G' = B + G(v/B)$, therefore $v$ must be ramified over $B$, thus by hypothesis we have $\delta(v/B) \neq \delta'$. By valuation inequality (for $\mathrm{val}^3_B$), we must have $\delta(v/B) < \delta'$, and we get the first part of the statement.

Suppose $d'_{G'\delta' k'}$ appears in $v$. Then replacing $d'_{G'\delta' k'}$ by $v' = f((d'_{G'\delta' i'})_{i'})$ is just a replacement of $(d'_{G'\gamma' i'})_{\gamma' i'}$ by $M(d'_{G'\gamma' i'})_{\gamma' i'}$ for some invertible matrix $M$. By remark 3.4.15, the new family is still a basis of $C$ over $A$ that witnesses (P1), (P2), (P3), (P4). By (P2) and (P3), $v'$ is Archimedean over $A$, and $G(v'/A) = G'$. Moreover, as all the $(d'_{G'\delta' i'})_{i'}$ are in $B + G'$, $v'$ must be ramified over $B$. Then, $w'$ is obtained by replacing $v'$ by $v$ in (a normal enumeration of) $w''$. This is exactly the operation described in lemma 3.4.16, concluding the proof. □

**Theorem 3.4.18.** *There exists $w = dd'\tilde{d}$ a basis of $C$ over $A$ that is under normal form with respect to $A$ and $B$.*

*Proof.* Let $w_0 = dd'^0\tilde{d}$ be a basis of $C$ over $A$ that witnesses (P1), (P2), (P3), (P4).

Suppose $w_n = dd'^n\tilde{d}$ does not witness (P5). Then the hypothesis of lemma 3.4.17 hold in $w_n$, so we do the replacement from this lemma. Let $w_{n+1}$ be the new family obtained.

To prove the theorem, we have to prove that the sequence $(w_n)_n$ must stop.

To do this, we once again build a strictly decreasing map with values in $\omega^\omega$. Let $N = |\Delta(C + B) \setminus \Delta(B)|$, and $(\delta_i)_{i<N}$ be the strictly increasing

enumeration of $\Delta(CB) \setminus \Delta(B)$. Define $F(n,i)$ as the number of terms $c$ from $d'^n$ for which $\delta(c/B) = \delta_i$. We define $g$ as the map:

$$n \longmapsto F(n,N) \times \omega^N + F(n,N-1) \times \omega^{N-1} + \ldots + F(n,1) \times \omega + F(n,0)$$

For each $n$, the family $w_{n+1}$ is obtained from $w_n$ by replacing some term $v$ from $d'^n$ by a vector $v'$ from $d'^{n+1}$ for which $\delta(v'/B) < \delta(v/B)$, so $g$ is strictly decreasing. This concludes the proof. $\square$

With that, we proved theorem 3.1.5, the main theorem of this chapter stating that $\downfree^{\mathbf{cut}} = \downfree^{\mathbf{inv}}$ in DOAG. The reasoning is explained in the end of section 3.1, and the different steps of the proof are carried out in proposition 3.2.32, proposition 3.2.36, proposition 3.2.45, corollary 3.3.24 and theorem 3.4.18.

# Chapter 4

# Regular ordered Abelian groups

This chapter corresponds to sections 5 and 6 of our preprint ([Hos23b]).

## 4.1 Quantifier-free types in the Presburger language

**Definition 4.1.1.** The *Presburger language* is defined as :

$$\mathfrak{L}_P = \{+, -, <, 0, \mathbb{1}, (\mathfrak{d}_{l^n})_{l \text{ prime}, n>0}\}$$

Given an ordered Abelian group $G$, we see $G$ as an $\mathcal{L}_P$-structure by setting $\mathfrak{d}_{l^n}(x)$ as the predicate $\exists y,\ l^n y = x$, and setting $\mathbb{1} = \min(]0, +\infty[)$ if $G$ is discrete, else we interpret $\mathbb{1}$ as the only $\varnothing$-definable choice, $\mathbb{1} = 0$ (this is the standard interpretation of the language). A *special* subgroup of $G$ is a pure subgroup (that is a relatively divisible subgroup) $A \leqslant G$ containing $\mathbb{1}$. Special subgroups will correspond to definably closed sets in the groups we are interested in (see corollary 4.2.6).

Define $S^m_{\mathrm{qf}}(A)$ to be the Stone space of quantifier-free types over $A$ in $m$ variables. Given an $m$-tuple $c$, we write $\mathrm{tp}_l(c/A)$ as the partial type generated by the set of quantifier-free formulas satisfied by $c$, with parameters in $A$, which only involve the predicates $(\mathfrak{d}_{l^n})_{n>0}$ (equality does not appear in these formulas either). We write $S^m_l(A)$ for the Stone space of all such partial types, and we call its elements the *l-types*. We define similarly $\mathrm{tp}_<(c/A)$, $S^m_<(A)$ for the quantifier-free formulas which only involve the predicates $=$, $<$.

For ease of notation, define the set of indices $J$ as the union of the set of primes with $\{<\}$. We have natural restriction maps $S^m_{\mathrm{qf}}(A) \longrightarrow S^m_j(A)$ for every $j \in J$, and they are all surjective. We write $(\pi_j)_j$ for those maps.

**Assumptions 4.1.2.** We fix an ordered Abelian group $M$, and a special subgroup $A \leqslant M$, such that $M$ is $|A|^+$-saturated and strongly $|A|^+$-homogeneous.

**Lemma 4.1.3** (Standard variant of the Chinese remainder theorem)**.** *Let $L$ be a finite set of primes, and $N > 0$. Then for all $(a_l)_l \in M^L$, there exists $b \in M$ such that $M \vDash \mathfrak{D}_{l^N}(a_l - b)$ for all $l \in L$.*

**Corollary 4.1.4.** *The natural map $S^m_{\mathrm{qf}}(A) \longrightarrow \prod_{l \text{ prime}} S^m_l(A)$ is surjective.*

*Proof.* For each $l$, let $a_l = (a_{i,l})_i \in M^m$. By lemma 4.1.3 and compactness, for each $i$, there must exist $b_i \in M$ such that $\mathfrak{D}_{l^N}(b_i - a_{i,l})$ for all $l$, $N$. In particular, $\mathrm{tp}_l(b/A) = \mathrm{tp}_l(a_l/A)$ for all $l$, which concludes the proof. □

*Remark* 4.1.5. The map $\pi : S^m_{\mathrm{qf}}(A) \longrightarrow \prod_{j \in J} S^m_j(A)$ is obviously injective. Note that $S^m_<(A)$ is not a factor of the product of lemma 4.1.3, thus $\pi$ may not be surjective.

## 4.2 Basic properties of regular groups

The class of regular ordered Abelian groups (ROAG) has several equivalent definitions. The ones that will matter for our proofs are quantifier elimination in $\mathcal{L}_P$, and some compatibility conditions between the $<$-types and the $l$-types. The definitions involving Archimedean groups and definable convex subgroups will be relevant, because they give a motivation as to why we are interested in ROAG.

The equivalence between these different definitions is folklore. In this section, we prove the easy implications between these equivalent definitions, and we give references for the harder implications.

**Definition 4.2.1.** Let $G$ be an ordered Abelian group. We say that $G$ is *regular* if, for every positive integer $n$, every interval of $G$ that contains at least $n$ elements intersects $nG$.

For any ordered Abelian group $G$, we recall that $\operatorname{div}(G)$ refers to the divisible closure of $G$, which we see as an ordered group to which the order on $G$ naturally extends.

**Lemma 4.2.2.** *Let $G \vDash \mathrm{ROAG}$, and $A$ a special subgroup of $G$. Let $F$ be the closed subspace of $S^m_{\mathrm{qf}}(A)$ of types of $m$-tuples that are $\mathbb{Q}$-free over $A$ (i.e. lifts in $G^m$ of $\mathbb{Q}$-free families from $\operatorname{div}(G) \diagup \operatorname{div}(A)$). Then the map $F \longrightarrow \pi_<(F) \times \prod_{l \text{ prime}} S^m_l(A)$ is a homeomorphism.*

In particular, $\pi_l(F) = S^m_l(A)$ for every prime $l$.

*Proof.* Two elements having the same type must have the same $j$-type for all $j \in J$, therefore the map is injective, and it is clearly continuous. Let us prove that it is surjective.

Let $p = (p_j)_j$ be in the direct product. Let $c = (c_i)_i$ be a realization of $p_<$ which is $\mathbb{Q}$-free over $A$. By corollary 4.1.4, let $c'$ be some tuple which realizes simultaneously all the $(p_l)$ for $l$ prime ($c'$ might not be $\mathbb{Q}$-free over $A$, but it does not matter). Let us show that there exist $d = (d_i)_{i<m}$ such that $d_i - c'_i \in \bigcap\limits_{N>0} NG$ for all $i$, and $d \vDash p_<$ (in particular, $d$ will be $\mathbb{Q}$-free over $A$). The existence of $d$ will be established by a disjunction of two cases.

Suppose $G$ is discrete. Let $N > 0$. Then, by regularity of $G$ applied to the interval $[c_i - c'_i, c_i - c'_i + N \cdot \mathbb{1}]$, there must exist an integer $k_i$ between 0 and $N$ such that $c_i - c'_i + k_i \cdot \mathbb{1} \in NG$. Let us show that $d' = (d'_i)_i = (c_i + k_i \cdot \mathbb{1})_i \vDash p_<$. If not, then there must exist some atomic formula with predicate = or < which is satisfied by $c$ and not $d'$, or vice-versa. Without loss, there must exist $a \in A$, and some non-zero $f \in \mathrm{LC}^m(\mathbb{Z})$ such that $f(c) < a \leqslant f((c_i + k_i \cdot \mathbb{1})_i)$. By subtracting $f((k_i \cdot \mathbb{1})_i)$ to the second inequality, $f(c)$ belongs to the interval $[a - f((k_i \cdot \mathbb{1})_i), a]$, which is included in $A$ as $A$ is special, contradicting $\mathbb{Q}$-freeness of $c$ over $A$. It follows that $d' \vDash p_<$, and we conclude by compactness that $d$ exists.

Now, suppose $G$ is dense. Let $D$ be the special subgroup generated by $A$, and all the $(c_i)_i$. For every positive element $a$ of $D$, the open interval $]c_i - c'_i, c_i - c'_i + a[$ is infinite, thus contains an element of $NG$ by regularity, for any $N > 0$. By compactness, there exists, for each $i$, an element $d'_i$ such that $d'_i \in ]c_i - c'_i, c_i - c'_i + a[$ for every positive $a \in D$, and $d'_i \in \bigcap\limits_{N>0} NG$. Let $d_i = d'_i + c'_i$. Notice that $d_i - c_i$ is a positive element which is infinitesimal with respect to $D$. It remains to show that $d \vDash p_<$. Let $f \in \mathrm{LC}^m(\mathbb{Z})$ be non-zero, and $a \in A$, such that $f(c) > a$. As each $d_i - c_i$ is infinitesimal with respect to $D$, this is also the case for $f(d-c) = f(d) - f(c)$, thus $|f(d) - f(c)| < f(c) - a$, thus $f(c) - f(d) < f(c) - f(a)$, and we conclude that $f(d) > f(a)$, proving that $d \vDash p_<$.

In both cases, we have that $d \vDash p_<$. As $d_i - c_i' \in l^N G$ for every $l$ prime and $N > 0$, we have $d \vDash p_l$ for every $l$, which concludes the proof. □

*Remark* 4.2.3. Note that, in the above proof, $d_i$ can always be chosen such that $d_i > c_i$ for each $i$ (in the discrete case, replace $k_i$ by $k_i + N > 0$). In particular, if we have $c = c'$ in the construction of the proof, then the tuple $d$ that we build is distinct from $c$, and re-iterating this operation generates infinitely many pairwise-distinct tuples having the same quantifier-free type as $d$. It follows that any consistent quantifier-free type over $A$ whose realizations are $\mathbb{Q}$-free over $A$ has infinitely many realizations.

**Definition 4.2.4.** Let $G$ be an ordered Abelian group, $n < \omega$ and $g \in G$. Define $H_n(g)$ to be the largest convex subgroup of $G$ for which we have $(g + H_n(g)) \cap nG = \varnothing$ (by convention, if $g \in nG$, then $H_n(g) = \{0\}$). This convex subgroup is definable:

$$H_n(g) = \{0\} \text{ or } \{x \in G | \forall y \in G,\ (|y| \leqslant n|x| \Longrightarrow g + y \notin nG)\}$$

for if $x \in G$, $z \in G$ satisfy $\Delta(x) = \Delta(z)$ ($\Delta$ is defined in definition 3.1.11) and $z + g \in nG$, then $x \notin H_n(g)$ is witnessed by choosing $y = z \pm mn|x|$, with $m$ the least natural integer for which $|z| - mn|x| \leqslant n|x|$.

In the literature, the set of the $H_n(g)$ for all $g \in G$ is called the *n-spine* of $G$. This is a definable family of definable convex subgroups. One of the most general "complexity class" of ordered Abelian groups that is still considered rather "nice" is the class of ordered Abelian groups with finite spines. This class contains in particular the dp-finite ordered Abelian groups. For reference, section 2 of [Far17] gives nice characterizations and a quantifier elimination result for this class.

**Proposition 4.2.5.** *Let $G$ be an ordered Abelian group. Then the following conditions are equivalent:*

1. *$G$ does not have any proper non-trivial definable convex subgroup.*

2. *For all $g \in G$, for all $n < \omega$ $H_n(g) = \{0\}$.*

3. *$G$ is regular.*

4. *The theory of $G$ eliminates quantifiers in the Presburger language.*

5. *There exists an Archimedean ordered Abelian group elementarily equivalent to $G$.*

*Proof.* The directions that we show in order to establish the equivalence are $1 \implies 2 \implies 3$, $3 \implies 4 \implies 1$ and $3 \implies 5 \implies 1$.

The implications $1 \implies 2$ and $5 \implies 1$ are trivial, $3 \implies 4$ is due to Weispfenning ([Wei81], Theorems 2.3 and 2.6), and Presburger in the discrete case, and $3 \implies 5$ follows from ([Zak61], Theorem 2.5) and ([RZ60], Theorem 4.7).

Let us show $2 \implies 3$ by contraposition. We have an interval $I$ of $G$ and an integer $n$ such that $|I| \geqslant n$ and $I \cap nG = \varnothing$. If $|I| < 2n+1$, then $G$ is discrete, and the points from $I + n \cdot \mathbb{1}$ have the same cosets $\bmod nG$ as those of $I$, so we can assume $|I| \geqslant 2n+1$ by replacing $I$ by $I \cup (I + n \cdot \mathbb{1}) \cup (I + 2n \cdot \mathbb{1})$. Let $f$ be a strictly increasing map $\{0, \dots 2n\} \longrightarrow I$, and $h = \min\{f(i+1) - f(i) | 0 \leqslant i \leqslant 2n-1\}$. As $f$ is strictly increasing, we have $h \neq 0$. As $h$ is minimal, we can sum $n$ many inequalities to get $n \cdot h \leqslant \sum_{i=0}^{n-1} (f(i+1) - f(i)) = f(n) - f(0)$. With the same argument, we have $n \cdot h \leqslant f(2n) - f(n)$, thus $[f(n) - n \cdot h, f(n) + n \cdot h] \subseteq I$, which is disjoint from $nG$. It follows that $h \in H_n(f(n)) \smallsetminus \{0\}$, so 2 fails.

Let us prove $4 \implies 1$ by contraposition. Let $G$ be an ordered Abelian group with a proper non-trivial definable convex subgroup $H$. Let $A$ be a definably closed subset of $G$ such that $H$ is $A$-definable. Let $A_1 = A_{\geqslant 0} \cap H$, $A_2 = (A_{>0} \smallsetminus H) \cup \{+\infty\}$. Let $X_< = \left( \bigcap_{a \in A_1, b \in A_2} ]a, b[ \right)$, an $A$-type definable set which is non-empty by compactness. For each prime $l$, let $X_l = \bigcap_N l^N G$, which corresponds to $\mathrm{tp}_l(0/A)$. Let us show that the partial type defined as $q(x) : \ x \in X_< \cap \bigcap_l X_l$ is consistent with $H$ and $\neg H$. Let $a \in A_1$, $b \in A_2$, $N > 0$. By compactness, it suffices to show that $Y = ]a, b[ \cap NG$ intersects $H$ and $\neg H$, which is witnessed by $N \cdot a \in Y \cap H$, and $N \cdot (b - a) \in Y \smallsetminus H$. Now let $h \in H, g \notin H$ be two realizations of $q$. Then $\mathrm{qftp}(h/A) = q = \mathrm{qftp}(g/A)$, but $\mathrm{tp}(h/A) \neq \mathrm{tp}(g/A)$, and the theory of $G$ does not eliminate quantifiers in $\mathcal{L}_P$. □

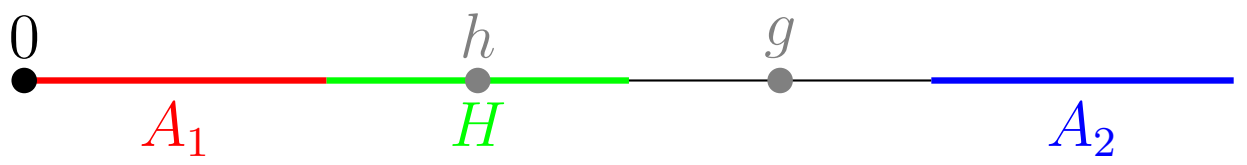

Note that in particular, the domain of the homeomorphism given by lemma 4.2.2 is in fact a space of *complete types*.

**Corollary 4.2.6.** *The definable closure of a parameter set $A$ coincides with the special subgroup it generates.*

*Proof.* The special subgroup generated by $A$ is clearly included in $\mathrm{dcl}(A)$. Conversely, if $A$ is a special subgroup, and $c \notin A$, then the 1-tuple $c$ is $\mathbb{Q}$-free over $A$, thus $c \notin \mathrm{acl}(A)$ by remark 4.2.3. $\square$

There is another important corollary which allows us to better understand independence:

**Corollary 4.2.7.** *Let $M \vDash$ROAG, $A \leqslant B$ special subgroups of $M$ such that $M$ is $|B|^+$ saturated and strongly $|B|^+$-homogeneous, and $c = c_1 \dots c_n \in M$. Suppose $c_1 \dots c_n$ is $\mathbb{Q}$-free over $B$, and let $F$ be the closed subspace of $S^n(M)$ of tuples that are $\mathbb{Q}$-free over $M$. For each $j \in J$, consider the action of $\mathrm{Aut}(M/A)$ on $\pi_j(F)$. Then the following conditions are equivalent:*

1. *$c \mathop{\downarrow}^{\mathbf{inv}}_{A} B$ (resp. $c \mathop{\downarrow}^{\mathbf{bo}}_{A} B$).*
2. *For each $j \in J$, $\mathrm{tp}_j(c/B)$, which is a partial type over $B$, can be extended to an $\mathrm{Aut}(M/A)$-invariant (resp. of bounded orbit) element of $\pi_j(F)$.*

*Proof.* First of all, note that the global extensions of $\mathrm{tp}(c/B)$ realized by tuples that are not $\mathbb{Q}$-free over $M$ divide over $A$, hence have an unbounded orbit by fact 1.1.19, so we do have to restrict ourselves to $F$.

For each $j \in J$, $F \longrightarrow \pi_j(F)$ is an equivariant surjection, thus the stabilizer of a point of $F$ is a subgroup of that of its image, so the cardinal of its orbit is larger than the supremum of the cardinals of the orbits of each of its images, and we get the top-to-bottom direction by contraposition.

Suppose for each $j$, we have $p_j \in \pi_j(F)$ a witness of the second condition: it extends $\mathrm{tp}_j(c/B)$, and its orbit is a singleton (resp. bounded). Note that the orbit of $(p_j)_j$ under $\mathrm{Aut}(M/A)$ is contained in the Cartesian product of the orbits of the $(p_j)_j$, therefore it is also a singleton (resp. bounded). Now, by lemma 4.2.2 and quantifier elimination, the map $F \longrightarrow \prod_j \pi_j(F)$ is an equivariant bijection, so the consistent global type corresponding to $(p_j)_j$ witnesses the first condition, and we get the bottom-to-top direction. $\square$

Note that, given special subgroups $A \leqslant B$, and a tuple $c$, $c$ is always $A$-interdefinable with a subtuple $d$ which is $\mathbb{Q}$-free over $A$ (choose any of those subtuples, of maximal length), thus there is a natural equivariant homeomorphism between the space of global extensions of $\mathrm{tp}(c/B)$ and that of $\mathrm{tp}(d/B)$. Corollary 4.2.7 completes the picture, and gives us easy conditions on $d$ to check whether $c \downarrow^{\mathbf{inv}}_A B$ or $c \downarrow^{\mathbf{bo}}_A B$.

*Remark* 4.2.8. Corollary 4.2.7 raises one interesting question: does the orbit of an element of $\prod_j \pi_j(F)$ coincide with the Cartesian product of the orbits of its components ? In other words, for $p_j \in \pi_j(F)$, is the natural (injective) map:

$$\mathrm{Aut}(M/A) \big/ \bigcap_j Stab(p_j) \longrightarrow \prod_j \mathrm{Aut}(M/A) \big/ Stab(p_j)$$

a bijection ?

## 4.3 Invariance and boundedness of the global extensions of the partial types

**Assumptions 4.3.1.** Let $M \vDash \mathrm{ROAG}$, let $A \leqslant B$ be special subgroups of $M$, let $\lambda = \max(|B|, 2^{\aleph_0})^+$, let $c = c_1 \dots c_n$ be a tuple from $M$, and suppose that $M$ is $\lambda$-saturated and strongly $\lambda$-homogeneous. Suppose $c$ is $\mathbb{Q}$-free over $B$, and let $F$ be the closed space of global types whose realizations are $\mathbb{Q}$-free over $M$. For each $j \in J$, let $p_j = \mathrm{tp}_j(c/B)$. A set of cardinality $\kappa$ is called *small* if $M$ is $\kappa^+$-saturated and strongly $\kappa^+$-homogeneous, else it is *large*.

By corollary 4.2.7, in order to understand the global extensions of $\mathrm{tp}(c/B)$ which are invariant or have a bounded orbit under the action of $\mathrm{Aut}(M/A)$, one has to understand, for each $j \in J$ the global extensions in $\pi_j(F)$ of $p_j$ which are invariant or have a bounded orbit.

### 4.3.1 Partial types using equality and order

**Assumptions 4.3.2.** In this subsection, on top of assumptions 4.3.1 we fix $\overline{M} \vDash \mathrm{DOAG}$ some $|M|^+$-saturated, strongly $|M|^+$-homogeneous elementary extension of $\mathrm{div}(M)$.

*Remark* 4.3.3. Let $d$ be some tuple from $M$. Then $d \vDash p_<$ if and only if $c$ and $d$ have the same type over $B$ in $\overline{M}$.

**Lemma 4.3.4.** *Let $D$ be a small special subgroup of $M$. Let $\alpha = \alpha_0 \dots \alpha_{m-1}$ be a tuple from $\overline{M}$ which is $\mathbb{Q}$-free over $D$. Suppose for all non-zero $f \in \mathrm{LC}^m(\mathbb{Q})$, and for all $d \in D$, we have $f(\alpha) \notin ]d, d+\mathbb{1}[$ (this interval of $\overline{M}$ being by convention empty when $\mathbb{1} = 0$). Then there exists $\sigma \in \mathrm{Aut}(\overline{M}/D)$ such that $\sigma(\alpha)$ is a tuple from $M$.*

*Proof.* Suppose by induction we have $\sigma_i \in \mathrm{Aut}(\overline{M}/D)$ sending $\alpha_{<i}$ to a tuple from $M$ for some $i < m$ ($\sigma_0 = id$ will do for $i = 0$). Let $D_i$ be the special subgroup of $M$ generated by $D\sigma_i(\alpha_{<i})$, which is exactly the relative divisible closure in $M$ of $D + \sum_{k<i} \mathbb{Z} \cdot \sigma_i(\alpha_k)$. Then, for all non-zero $f \in \mathrm{LC}^{n-i}(\mathbb{Q})$, for all $d \in D_i$, we have by hypothesis $f(\alpha_{\geqslant i}) \notin ]d, d+\mathbb{1}[$. Let us find $\tau \in \mathrm{Aut}(\overline{M}/D_i)$ such that $\tau(\sigma_i(\alpha_i)) \in M$. Then we could set $\sigma_{i+1} = \tau \circ \sigma_i$, and conclude by induction.

By strong homogeneity of $\overline{M}$, it is enough to show that any interval of $\overline{M}$ with bounds in $\mathrm{div}(D_i) \cup \{\pm\infty\}$ containing $\sigma_i(\alpha_i)$ has a point $\beta$ in $M$. Let $I$ be such an interval. If either the lower or the upper bound of $I$ is in $\{\pm\infty\}$, then $\beta$ clearly exists, one may choose some multiple of $d$ with large enough absolute value if $d \neq 0$, else choose any non-trivial element with correct sign. Now, suppose $I = \left]\frac{d_1}{N_1}, \frac{d_2}{N_2}\right[$, with $d_k \in D_i$, $N_k > 0$ (by the $\mathbb{Q}$-freeness assumption on $\alpha$, $\sigma_i(\alpha_i) \notin \mathrm{div}(D_i)$, thus it does not matter whether the bounds of $I$ belong to $I$). Then $I$ has a point in $M$ if and only if $]N_2 d_1, N_1 d_2[$ has a point in $N_1 N_2 M$. If $M$ is dense, then $I$ has infinitely many points in $M$, and we can use the axioms of ROAG to conclude.

There remains to deal with the case where $M$ is discrete. We established earlier that the $\mathbb{Q}$-linear combinations of $\sigma_i(\alpha_{\geqslant i})$ do not belong to $]d, d+\mathbb{1}[$ for any $d \in D_i$. In particular, for every $N \in \mathbb{Z}$:

$$N_1 N_2 \sigma_i(\alpha_i) \notin \bigcup_{N \in \mathbb{Z}} ]N_2 d_1 + N \cdot \mathbb{1}, N_2 d_1 + (N+1) \cdot \mathbb{1}[$$

and by $\mathbb{Q}$-freeness we also have $N_1 N_2 \sigma_i(\alpha_i) \notin \{N_2 d_1 + N \cdot \mathbb{1} | N \in \mathbb{Z}\}$. As a result, $N_1 N_2 \sigma_i(\alpha_i) \in ]N_2 d_1, N_1 d_2[ \smallsetminus ]N_2 d_1, N_2 d_1 + (N_1 N_2 + 2) \cdot \mathbb{1}]$, which must imply that $N_2 d_1 + (N_1 N_2 + 1) \cdot \mathbb{1} < N_1 d_2$, so $I$ has at least $N_1 N_2 + 1$ (in fact, infinitely many) points in $M$, and we can also conclude with the axioms of ROAG. □

**Lemma 4.3.5.** *Let $\alpha$ be a tuple from $M^n$ which is $\mathbb{Q}$-free over $B$, and such that $p = \mathrm{tp}^{\mathrm{DOAG}}(\alpha/B)$ does not fork over $A$. Let $D$ be a small special subgroup of $M$ containing $B$, and let $q \in S^n_{\mathrm{DOAG}}(\mathrm{div}(D))$ be some extension of $p$ which does not fork over $A$. Then $q$ has a realization in $M^n$.*

*Proof.* First of all, the realizations of $q$ must be $\mathbb{Q}$-free over $D$. Let $d \in D$, and let $f$ be a non-zero element of $\mathrm{LC}^n(\mathbb{Q})$. By lemma 4.3.4, one just has to prove that $q(x) \vDash f(x) \notin ]d, d+\mathbb{1}[$. If not, then by theorem 3.1.5, the interval $[d, d+\mathbb{1}]$ must have a point in $\mathrm{div}(A)$. By multiplying everything by a sufficiently large $N > 0$, the interval $]Nd, Nd + N\cdot\mathbb{1}[$ has a point in $A$, and $Nd \in D$, thus $Nd \in A$. In particular, $q_{|A}(x) \vDash Nf(x) \in ]Nd, Nd + N\cdot\mathbb{1}[$, therefore $p(x) \vDash Nf(x) \in ]Nd, Nd+N\cdot\mathbb{1}[$, thus $Nf(\alpha) \in ]Nd, Nd+N\cdot\mathbb{1}[$. As $Nd$ and $\alpha$ are in $A$, it follows that $Nf(\alpha) \in A$, contradicting $\mathbb{Q}$-freeness of $\alpha$. □

**Proposition 4.3.6.** *By remark 4.3.3, let $h$ be the natural injection:*

$$\pi_<(F) \longrightarrow S^n_{\mathrm{DOAG}}(\mathrm{div}(M))$$

*Consider the action of* $\mathrm{Aut}(M/A)$ *on* $Q = \pi_<(F)$. *Then the following conditions are equivalent:*

1. *Some element of $Q$ is invariant and extends $p_<$.*
2. *Some element of $Q$ has a bounded orbit and extends $p_<$.*
3. *The partial type $p_<$ does not fork over $A$.*
4. *The partial type $p_<$ does not divide over $A$.*
5. *Every closed bounded interval with bounds in $B$ containing a $\mathbb{Z}$-linear combination of $c$ also has a point in* $\mathrm{div}(A)$.
6. *For some $p \in Q$ extending $p_<$, $h(p)$ extends to some* $\mathrm{Aut}(\overline{M}/A)$*-invariant type over $\overline{M}$.*

*Proof.* The directions $1 \Longrightarrow 2$ and $3 \Longrightarrow 4$ are immediate.

Let us prove $2 \Longrightarrow 3$. Let $q_< \in Q$ witness 2. For each prime $l$, let $q_l(x)$ in $S^n_l(M)$ be the partial type $\{\mathfrak{d}_{l^N}(f(x)) | n > 0, f \in \mathrm{LC}(\mathbb{Z}), f \neq 0\} = \mathrm{tp}_l(0/M)$. Note that $q_l$ is clearly invariant. Then the complete global type in $F$ corresponding to $(q_<, (q_l)_l)$ (remember that $F \longrightarrow \pi_<(F) \times \prod_{l \text{ prime}} S^n_l(M)$ is a homeomorphism by lemma 4.2.2, thus this type is consistent) has a bounded orbit, and thus does not fork over $A$, and we get 3.

Let us prove $4 \implies 5$. Suppose we have $b_k \in B$, $f \in \mathrm{LC}(\mathbb{Z})$ such that the formula $f(c) \in [b_1, b_2]$ witnesses the failure of 5. Let us show that the formula $f(x) \in [b_1, b_2]$ divides over $A$, witnessing the failure of 4. As in the proof of corollary 1.2.12, it is enough to find $d \equiv_A b_1$ such that $d > b_2$. As $[b_1, b_2]$ has no point in $\mathrm{div}(A)$, we have $b_k \notin A$, thus the singletons $b_1$ and $b_2$ are both $\mathbb{Q}$-free over $A$. Now, by lemma 4.2.2, one just has to show that $\mathrm{tp}_<(b_1/A)$ is consistent with $]b_2, +\infty[$. If not, then all the elements in $M$ of $]b_2, +\infty[$ have a type over $A$ in $\overline{M}$ which is distinct from that of $b_1$. As a result, by the characterization of 1-types in DOAG, there must exist in $\overline{M}$ a point in $\mathrm{div}(A) \cap ]b_1, b_2]$, a contradiction.

Let us show $5 \implies 6$. Suppose 6 fails. Then, by theorem 3.1.5, $p_<$ is inconsistent with the partial type $\{f(x) \notin I | I \in \mathcal{I}, f \in \mathrm{LC}(\mathbb{Z})\}$, with $\mathcal{I}$ the set of closed bounded intervals with bounds in $M$ that have no points in $\mathrm{div}(A)$. By compactness, there exist finite subsets $\mathcal{I}' \subseteq \mathcal{I}$, $G \subseteq \mathrm{LC}^n(\mathbb{Z})$ such that $p_<(x) \models \bigvee_{I \in \mathcal{I}', g \in G} g(x) \in I$. Let $D$ be the special subgroup of $M$ generated by $B$ and the bounds of the elements of $\mathcal{I}'$. Then, for all $p \in Q$ extending $p_<$, the restriction of $p$ to $D$ corresponds to a type in $S^n_{\mathrm{DOAG}}(D)$ which forks over $A$ in DOAG.

Suppose now, by contradiction, that 5 holds. Let $q$ be the type of $S^n_{\mathrm{DOAG}}(B)$ corresponding to $p_<$. Then $q$ does not fork over $A$. By extension, let $q' \in S^n_{\mathrm{DOAG}}(D)$ be an extension of $q$ which does not fork over $A$. By lemma 4.3.5, $q'$ must admit a realization $\beta$ in $M^n$ (in particular $\beta \models p_<$). As $q'$ does not fork over $A$, we have $g(\beta) \notin I$ for all $g \in G$, $I \in \mathcal{I}'$. As a result, $q_< = \mathrm{tp}_<(\beta/D)$ cannot be extended to any element of $F$. This means that in any elementary extension of $M$, no realization of $q_<$ is $\mathbb{Q}$-free over $M$. By compactness, this means that there exist finite subsets $P \subseteq M$, $G' \subseteq \mathrm{LC}^n(\mathbb{Q}) \setminus \{0\}$, such that $q_<(x) \models \bigvee_{m \in P, g \in G'} g(x) = m$.

Then we can finally reach a contradiction by extending $q'$ again: let $\tilde{D}$ be the special subgroup of $M$ generated by $D \cup P$, and let $\tilde{q} \in S^n_{\mathrm{DOAG}}(\tilde{D})$ be some extension of $q'$ which does not fork over $A$. Then with the same argument, $\tilde{q}$ has a realization $\gamma \in M^n$, and by hypothesis $\gamma$ is not $\mathbb{Q}$-free over $\tilde{D}$, a contradiction with the fact that $\tilde{q}$ does not fork over $A$.

Suppose $p \in Q$ witnesses 6, and let us show that $p$ is $\mathrm{Aut}(M/A)$-invariant, which would allow us to conclude the whole proof with the direction $6 \implies 1$.

Let $q \in S^n_{\mathrm{DOAG}}(\overline{M})$ be some global $\mathrm{Aut}(\overline{M}/A)$-invariant extension of $h(p)$, and $\sigma \in \mathrm{Aut}(M/A)$. Then $\sigma$ extends uniquely to $\sigma'$, an automorphism of the ordered group $\mathrm{div}(M)$, pointwise-fixing $A$. We clearly have $h(\sigma(p)) = \sigma'(h(p))$ (look at the atomic formulas with predicate $=, <$ that belong to $X_<$, check that their image by $\sigma$ is satisfied by the realizations of $\sigma'(h(p))$). By quantifier elimination in DOAG, $\sigma'$ is a partial elementary map in $\overline{M}$. By strong homogeneity of $\overline{M}$, $\sigma'$ extends to some $\tilde{\sigma} \in \mathrm{Aut}(\overline{M})$. As $\sigma'$ pointwise-fixes $A$, so does $\tilde{\sigma}$, thus $\tilde{\sigma}(q) = q$. Now, we have $\sigma(p) = h^{-1}(\sigma'(h(p)) = h^{-1}(\sigma'(q_{|\mathrm{div}(M)})) = h^{-1}(\tilde{\sigma}(q_{|\mathrm{div}(M)})) = h^{-1}(\tilde{\sigma}(q)_{|\tilde{\sigma}(\mathrm{div}(M))}) = h^{-1}(q_{|\mathrm{div}(M)}) = h^{-1}(h(p)) = p$, concluding the proof. □

The fifth condition of proposition 4.3.6 is a very simple geometric condition, the kind of statement that would be very satisfactory for a characterization of forking. In the subsections 4.3.3 and 4.3.4, we look for similar conditions for the $(p_l)_l$.

### 4.3.2 Mapping monsters to monsters

*Remark* 4.3.7. Up to isomorphism, the unique discrete subgroup of $\mathbb{R}$ is $\mathbb{Z}$, and all its elements are $\emptyset$-definable in $\mathcal{L}_P$. As a result, any discrete model of ROAG is an elementary extension of $\mathbb{Z}$, which we identify with $\mathrm{dcl}(\emptyset)$.

If $d$ is a tuple in some elementary extension of $M$, then proposition 4.3.6 shows that $p_<$ does not fork/divide over $A$ if and only if the corresponding partial type (via $h$) in $S_{\mathrm{DOAG}}(\mathrm{div}(M))$ does not fork over $A$. Moreover, one can easily see that $h$ is continuous, and the restriction of $h$ to the closed subspace of $\pi_<(F)$ of extensions of $p_<$ which do not fork over $A$ is, by lemma 4.3.5, onto the closed subspace of $S(\mathrm{div}(M))$ of extensions of $h(p_<)$ which do not fork over $A$. It follows that they are homeomorphic.

One may think that we just proved that this topological space is homeomorphic to the one described in remark 3.3.31, but there is actually a subtle obstruction, coming from the fact that $\mathrm{div}(M)$ is not necessarily $|B|^+$-saturated. For the remainder of this subsection, we show that the space of extensions of $p_<$ which do not fork over $A$ is indeed homeomorphic to the space described remark 3.3.31. For this, we have to do a case disjunction depending on whether $M$ is dense or not.

If $M$ is dense, we show that there exists a sufficiently saturated elementary extension $N$ of $M$ such that $\mathrm{div}(N)$ is also sufficiently saturated.

If $M$ is discrete, then $\operatorname{div}(M)$ is never $\aleph_1$-saturated, for it does not admit non-zero elements having smaller Archimedean value than $\mathbb{1}$. The homeomorphism $h$ that we build does not come from the divisible closure, but from the quotient $M/\mathbb{Z}$, which is a model of DOAG. We show that this structure is sufficiently saturated.

Note that we do not care about strong homogeneity here. Indeed, one may easily see that a partial global $A$-invariant type on an $|A|^+$-saturated first-order structure extends uniquely to an $A$-invariant type over any $|A|^+$-saturated elementary extension, regardless of strong homogeneity.

**Definition 4.3.8.** Suppose $M$ is discrete. We define $\widetilde{M}$ as the quotient $M/\mathbb{Z}$, and $\pi$ the canonical surjection. As $\mathbb{Z}$ is a convex subgroup, $\widetilde{M}$ is naturally an ordered Abelian group by definition 1.2.13, and $\pi$ is an order-preserving map.

**Proposition 4.3.9.** *If $M$ is discrete, then $\widetilde{M} \vDash$ DOAG.*

*Proof.* By saturation, $M \neq \mathbb{Z}$, therefore $\widetilde{M}$ is non-trivial. Let $\pi(a) \in \widetilde{M}$, and $m > 0$. The quotient and remainder for the Euclidean division by $m$ are $\varnothing$-definable functions in $M$, therefore there exists $b \in M$ such that $a = m\cdot b + k\cdot\mathbb{1}$ for some $k \in \mathbb{Z}$ ($k < m$), in particular $\pi(a) = m\cdot\pi(b)$, concluding the proof. □

**Lemma 4.3.10.** *The map $h' : \operatorname{tp}_<(c/M) \longrightarrow \operatorname{tp}\left(\pi(c)/\widetilde{M}\right)$ is a well-defined continuous injection $\pi_<(F) \longrightarrow S^n_{\mathrm{DOAG}}\left(\widetilde{M}\right)$.*

*Proof.* Let $f \in \mathrm{LC}^n(\mathbb{Q})$, $b \in M$, and $c, d \in M^n$ such that $c$ and $d$ are $\mathbb{Q}$-free over $\operatorname{dcl}(b)$. Let $m > 0$ be such that $m \cdot f \in \mathrm{LC}^n(\mathbb{Z})$.

Let us show that $h'$ is well-defined. Suppose $\operatorname{tp}_<(c/\operatorname{dcl}(b)) = \operatorname{tp}_<(d/\operatorname{dcl}(b))$. Then we have:

$$\begin{aligned} f(\pi(c)) > \pi(b) &\iff m\cdot f(\pi(c)) > m\cdot\pi(b) \\ &\iff \forall k \in \mathbb{Z}\; m\cdot f(c) + k\cdot\mathbb{1} > b \\ &\iff \forall k \in \mathbb{Z}\; m\cdot f(d) + k\cdot\mathbb{1} > b \\ &\iff f(\pi(d)) > \pi(b) \end{aligned}$$

$$\begin{aligned} f(\pi(c)) = \pi(b) &\iff m\cdot f(\pi(c)) = m\cdot\pi(b) \\ &\iff m\cdot f(c) \in (m\cdot b + \mathbb{Z}) \\ &\implies m\cdot f(c) \in \operatorname{dcl}(b) \\ &\implies \begin{cases} f = 0 \\ b \in \mathbb{Z} \end{cases} \text{by } \mathbb{Q}\text{-freeness} \\ &\implies f(\pi(d)) = \pi(b) \end{aligned}$$

it follows that $\text{tp}(\pi(c)/\pi(\text{dcl}(b))) = \text{tp}(\pi(d)/\pi(\text{dcl}(b)))$, therefore $h'$ is well-defined. The above computations also show us that $h'$ is continuous.

Now, if $f \in \text{LC}^n(\mathbb{Z})$, then, by $\mathbb{Q}$-freeness, we have:

$$f(c) > b \iff \forall k \in \mathbb{Z}\ f(c) > b + k \cdot \mathbb{1} \iff f(\pi(c)) > \pi(b)$$

and:

$$f(c) = b \implies f = 0 \implies f(\pi(c)) = 0$$

it follows that $h'$ is an injection. □

**Lemma 4.3.11.** *Suppose $M$ is discrete. Let $h'$ be the map from lemma 4.3.10, and let $q \in \pi_<(F)$. Then the following are equivalent:*

1. *For all $f \in \text{LC}^n(\mathbb{Z})$, for all closed bounded intervals $I$ with bounds in $M$ having no point in $\text{div}(A)$, $q(x) \models f(x) \notin I$.*

2. *For all $f \in \text{LC}^n(\mathbb{Q})$, for all closed bounded intervals $I$ with bounds in $\pi(M)$ having no point in $\pi(A)$, $h'(q)(x) \models f(x) \notin I$.*

*Proof.* Suppose 1 fails, and let $f$, $I$ be a witness. As $A$ is a special subgroup, we have $\mathbb{Z} \subseteq A$, therefore the interval $[\min(I) - m \cdot \mathbb{1}, \max(I) + m \cdot \mathbb{1}]$ has no point in $A$ for every $m < \omega$. It follows that $I' = [\pi(\min(I)), \pi(\max(I))]$ has no point in $\pi(A)$. As a result, $f$, $I'$ is a witness of the failure of 2.

Conversely, suppose 2 fails, and let $f, I$ be a witness. We must have $f \neq 0$, otherwise $I$ would have a point in $\pi(A)$. Let $b_1$, $b_2$ be respective preimages of $\min(I)$, $\max(I)$ by $\pi$. Let $m > 0$ be such that $m \cdot f \in \text{LC}^n(\mathbb{Z})$. As $q \in \pi_<(F)$, and $f \neq 0$, $q(x) \models m \cdot f(x) \notin M$. This would not be the case if we had $m \cdot \pi(b_1) = m \cdot \pi(b_2)$, therefore we have $m \cdot \pi(b_1) < m \cdot \pi(b_2)$, thus $m \cdot b_1 < m \cdot b_2$, and $q(x) \models m \cdot f(x) \in I' = [m \cdot b_1, m \cdot b_2]$. In order to conclude the proof, it suffices to show that $I'$ has no point in $\text{div}(A)$, which would establish the failure of 1, with the witness $m \cdot f$, $I'$. Let $a \in A$, and $k > 0$. If we had $\frac{a}{k} \in I'$, then we would have $\frac{1}{mk}\pi(a) \in I$. However, as $A$ is definably closed, and the quotient for the Euclidean division by $mk$ is $\varnothing$-definable, we would have $\frac{1}{mk}\pi(a) \in \pi(A)$, a contradiction. □

**Corollary 4.3.12.** *Suppose $M$ is discrete, and let $q \in \pi_<(F)$. Then the following are equivalent:*

1. *$q$ does not fork over $A$.*

2. *$h'(q)$ does not fork over $\pi(A)$.*

*Proof.* By condition 5 of proposition 4.3.6, $q$ does not fork over $A$ if and only if condition 1 of lemma 4.3.11 holds. By definition 3.1.3, the second condition of lemma 4.3.11 exactly states that $h'(q)$ avoids all the definable sets witnessing cut-dependence (i.e. forking-dependence by theorem 3.1.5) from $\pi(M)$ over $\pi(A)$, therefore $h'(q)$ does not fork over $\pi(A)$ if and only if condition 2 of lemma 4.3.11 holds. □

**Lemma 4.3.13.** *Suppose $M$ is discrete. Let $q \in S^n_{\mathrm{DOAG}}(\widetilde{M})$ be the type of a family which is $\mathbb{Q}$-free over $\widetilde{M}$. Then there exists $p \in \pi_<(F)$ such that $h'(p) = q$.*

*Proof.* Let $I$ be a finite set. For each $i \in I$, let $b_i, b_i' \in M$, and $f_i \in \mathrm{LC}(\mathbb{Q})$, such that $q(x) \vDash \pi(b_i) < f_i(x) < \pi(b_i')$. By compactness, it suffices to find a tuple $d$ in $M$ such that $\pi(b_i) < f_i(\pi(d)) < \pi(b_i')$ for every $i$. By consistency of $q$, there exists a tuple $\delta$ in $\widetilde{M}$ such that $\pi(b_i) < f_i(\delta) < \pi(b_i')$ for every $i$. We simply choose $d$ an arbitrary lift of $\delta$ by $\pi$. □

It follows that every type over $\widetilde{M}$ which extends $h'(p_<)$ and which does not fork over $\pi(A)$ has a lift in $\pi_<(F)$ by $h'$ which does not fork over $A$. As a result, $h'$ restricts to a homeomorphism between the space of types in $\pi_<(F)$ which extend $p_<$ and which do not fork over $A$, and the space of types in $S_{\mathrm{DOAG}}(\widetilde{M})$ which extend $h'(p_<)$ and which do not fork over $\pi(A)$.

Now we show that $\mathrm{div}(M)$ or $\widetilde{M}$ is sufficiently saturated:

**Proposition 4.3.14.** *Let $\kappa$ be an infinite cardinal such that $M$ is $\kappa$-saturated. If $M$ is dense, then $\mathrm{div}(M)$ is $\kappa$-saturated.*

*Proof.* Suppose $M$ is dense. Let $(I_i)_i$ be a consistent family of strictly less than $\kappa$ many intervals from $\mathrm{div}(M)$. It suffices to show that $\bigcap_i I_i$ has a point in $\mathrm{div}(M)$ (as saturation for unary types implies saturation). We may assume every intersection of finitely many intervals from the family is also in the family. If some $I_i$ has empty interior, then it is a singleton, and the proof is trivial, therefore we may assume that every $I_i$ is open of non-empty interior.

As $M$ is $\kappa$-saturated and dense, there exists a non-zero $\varepsilon \in M$ such that $\Delta(\varepsilon) < \Delta(\sup(I_i) - \inf(I_i))$ ($\Delta$ is defined in definition 3.1.11) for all $i$. In particular:

$$\inf(I_i) < \inf(I_i) + |\varepsilon| < \sup(I_i) - |\varepsilon| < \sup(I_i)$$

For each $i$, let $m_i < \omega$ be sufficiently large, such that $m_i \cdot \inf(I_i), m_i \cdot (\inf(I_i) + |\varepsilon|), m_i \cdot (\sup(I_i) - |\varepsilon|), m_i \cdot \sup(I_i)$ are all in $M$. By regularity, as both the intervals $]m_i \cdot \inf(I_i), m_i \cdot (\inf(I_i) + |\varepsilon|)[$, $]m_i \cdot (\sup(I_i) - |\varepsilon|), m_i \cdot \sup(I_i)[$ are infinite, they have a point in $m_i \cdot M$, therefore there exists $a_i, b_i \in M$ such that:

$$\inf(I_i) < a_i < \inf(I_i) + |\varepsilon| < \sup(I_i) - |\varepsilon| < b_i < \sup(I_i)$$

Let $I'_i = ]a_i, b_i[$, which is an interval of $M$. We have $I'_i \subseteq I_i$ for every $i$. It suffices to show that the family $(I'_i)_i$ has a point in $M$. By saturation and density, it suffices to show that $a_i < b_j$ for every $i, j$. It suffices to show that $\inf(I_i) + |\varepsilon| < \sup(I_j) - |\varepsilon|$. As $(I_i)_i$ is closed under finite intersection, and by definition of $\varepsilon$, we have:

$$\begin{aligned} 2 \cdot |\varepsilon| &< \sup(I_i \cap I_j) - \inf(I_i \cap I_j) \\ &= \min(\sup(I_i), \sup(I_j)) - \max(\inf(I_i), \inf(I_j)) \\ &\leqslant \sup(I_j) - \inf(I_i) \end{aligned}$$

hence the result is proved. □

**Proposition 4.3.15.** *Let $\kappa$ be an infinite cardinal such that $M$ is $\kappa$-saturated. If $M$ is discrete, then $\widetilde{M}$ is $\kappa$-saturated.*

*Proof.* Suppose $M$ is discrete. Let $(I_i)_i$ be a consistent family of strictly less than $\kappa$ many open intervals of non-empty interior of $\widetilde{M}$. Let $a_i, b_i \in M$ be such that $I_i = ]\pi(a_i), \pi(b_i)[$. Then the family of every interval of the form $]a_i + k \cdot \mathbb{1}, b_i + k' \cdot \mathbb{1}[$, with $k, k' \in \mathbb{Z}$ is consistent in $M$. It can be defined with as many parameters as the initial family, therefore its intersection has a point $d \in M$. We clearly have $\pi(d) \in \bigcap_i I_i$, concluding the proof. □

**Proposition 4.3.16.** *The space of types in $\pi_<(F)$ which extend $p_<$ and do not fork over $A$ is homeomorphic to the space described in remark 3.3.31.*

*Proof.* The homeomorphism is given by $h$ if $M$ is dense, and $h'$ if it is discrete.
□

### 4.3.3 Partial types using a prime of finite index

For the rest of this section, we fix a prime $l$.

**Proposition 4.3.17.** *Let $G$ be a torsion-free Abelian group. Then, for all $N > 0$, we have in $G^{eq}$ a $\varnothing$-definable group isomorphism between the two definable groups $G/_{lG}$ and $l^N G/_{l^{N+1}G}$.*

*Proof.* Let $x, y \in G$, and $N > 0$. Suppose $l^N x - l^N y \in l^{N+1}G$. Then, as $G$ is torsion-free, we have $x - y \in lG$, thus the map:

$$f : l^N x \bmod l^{N+1}G \longmapsto x \bmod lG$$

is a well-defined group homomorphism, which is clearly surjective. Now, if $x \notin lG$, then $l^N x \notin l^{N+1}G$, so the map is injective. □

**Corollary 4.3.18.** *If $lG$ has a finite index $d$, then $\bigcap_N l^N G$ has index at most equal to $2^{\aleph_0}$.*

*Proof.* Define a tree structure on $\bigcup_{N \geqslant 0} G/_{l^N G}$ of root $G$, such that the parent of some node $x \bmod l^{N+1}G$ is $x \bmod l^N G$. Then every node has $d$ children and the tree has $\aleph_0$ levels, thus it admits at most $2^{\aleph_0}$ branches. To conclude, the map $x \bmod \bigcap_N l^N G \longmapsto (x \bmod l^N G)_N$ is a natural bijection between $G/_{\bigcap_N l^N G}$ and the set of the branches. □

**Lemma 4.3.19.** *Let $D$ be some small special subgroup of $M$. Suppose we have $e_1, e_2 \in M$, such that $e_1, e_2 \notin D + l^N M$, and $e_1 - e_2 \in l^{N-1}M$, for some $N > 0$. Then we have $\mathrm{tp}_l(e_1/D) = \mathrm{tp}_l(e_2/D)$.*

*Proof.* Any atomic formula with parameters in $D$ and predicate in $\{\mathfrak{d}_{l^m} | m > 0\}$ can clearly be written $\mathfrak{d}_{l^m}(\lambda x - d)$, with $\lambda \in \mathbb{Z}, m > 0, d \in D$. We have several cases:

- $\lambda \in l\mathbb{Z}$, $d \notin lM$, in which case the formula always fails.
- $\lambda \in l\mathbb{Z}$, $d \in lM$, $m = 1$, in which case the formula always holds.
- $\lambda \notin l\mathbb{Z}$, in which case $\lambda$ and $l$ are coprime. Therefore, by Bézout's identities, the formula is equivalent to $\mathfrak{d}_{l^m}(x - d')$ for some $d' \in D$. By hypothesis on the $e_i$, this formula is satisfied by $e_1$ if and only if it is satisfied by $e_2$.

- $\lambda \in l\mathbb{Z}$, $d \in lM$, $m > 1$, in which case the formula is equivalent to $\mathfrak{d}_{l^{m-1}}\left(\frac{\lambda}{l}x - \frac{d}{l}\right)$, and we reduce by induction to one of the above three cases.

either way, we clearly see that those formulas are satisfied by $e_1$ if and only if they are satisfied by $e_2$, which concludes the proof. □

**Corollary 4.3.20.** *Let $d$ be a tuple from $M$, which is $\mathbb{Q}$-free over $B$. Then $d \vDash p_l$ if and only if, for each $f \in \mathrm{LC}^n(\mathbb{Z})$, for each $N > 0$, either we have $f(c) - f(d) \in l^N M$, or neither $f(c)$ nor $f(d)$ are in $B + l^N M$.*

*Proof.* Any atomic formula $\varphi(x_1 \dots x_n)$ with predicate in $\{\mathfrak{d}_{l^N} | N > 0\}$ can be written $\psi(f(x))$, with $f \in \mathrm{LC}(\mathbb{Z})$, and $\psi$ an atomic formula on the same predicate with one variable. As a result, $d \vDash p_l$ if and only if $\mathrm{tp}_l(f(d)/B) = \mathrm{tp}_l(f(c)/B)$ for all $f \in \mathrm{LC}(\mathbb{Z})$. Fix $f \in \mathrm{LC}^n(\mathbb{Z})$.

Suppose that we have $f(c) \in B + l^N M$ for all $N > 0$. For each $N$, let $b_N \in B$ such that $f(c) - b_N \in l^N M$. Then $\mathrm{tp}_l(f(c)/B)$ contains $\{\mathfrak{d}_{l^N}(x - b_N) | N > 0\}$. This inclusion is an equality, for if $\mathfrak{d}_{l^N}(f(d) - b_N)$ for all $N$, then we have $f(d) - f(c) \in l^N M$ for all $N$, which clearly implies $\mathrm{tp}_l(f(c)/B) = \mathrm{tp}_l(f(d)/B)$. In particular, $\mathrm{tp}_l(f(d)/B) = \mathrm{tp}_l(f(c)/B)$ if and only if $f(d) - f(c) \in \bigcap_N l^N M$.

Suppose now that there exists $N > 0$ such that $f(c) \notin B + l^N M$. The set of all such $N$ is a final segment of $\omega$, choose $N$ its least element. Let $b \in B$ such that $f(c) - b \in l^{N-1} M$ (if $N = 1$, then $b = 0$ will do). Then, by lemma 4.3.19, $\mathrm{tp}_l(f(c)/B)$ is generated by the partial type:

$$\{\mathfrak{d}_{l^{N-1}}(x - b)\} \cup \{\neg \mathfrak{d}_{l^N}(x - b') | b' \in B\}$$

this concludes the proof. □

*Remark* 4.3.21. If $[M : lM]$ is finite, then by corollary 4.3.18 the $\varnothing$-type-definable equivalence relation $\{\mathfrak{d}_{l^N}(x - y) | N > 0\}$ is bounded (has a small number of classes) and finer than the equivalence relation of having the same $l$-type over any parameter set, thus $S_l^n(M)$ is small, and each orbit of $S_l^n(M)$ under $\mathrm{Aut}(M/A)$ is obviously bounded.

Moreover, the classes of this $\varnothing$-type-definable equivalence relation are all $\mathrm{Aut}(M^{eq}/\mathrm{acl}^{eq}(\varnothing))$-invariant, thus $\mathrm{Aut}(M^{eq}/\mathrm{acl}^{eq}(\varnothing))$ acts trivially over $S_l^n(M)$.

Let us recall that we adopted assumptions 4.3.1 in the section.

**Proposition 4.3.22.** *Suppose $[M : lM]$ is finite. Consider the action of* $\mathrm{Aut}(M/A)$ *on* $Q = S_l^n(M)$*. Then* $p_l$ *has an invariant extension in* $Q$ *if and only if every* $\mathbb{Z}$*-linear combination of* $c$ *belongs to* $\bigcap_N (A + l^N M)$*.*

*Proof.* Suppose every $\mathbb{Z}$-linear combination of $c$ belongs to $\bigcap_N (A+l^N M)$. Then $p_l(M)$ can be written as the intersection of $A$-definable sets of the form:

$$\{x | \mathfrak{d}_{l^N} (f(x) - a)\},$$

with $N > 0, a \in A$, $f \in \mathrm{LC}(\mathbb{Z})$. As a result, $p_l(M)$ is $A$-type definable (thus $\mathrm{Aut}(M/A)$-invariant), and all its realizations $d \in M^n$ satisfy the conditions $f(c) - f(d) \in \bigcap_N l^N M$ for all $f \in \mathrm{LC}^n(\mathbb{Z})$. By corollary 4.3.20, $p_l$ is complete as an element of $Q$, so we get the right-to-left direction.

Suppose now that we have $f \in \mathrm{LC}^n(\mathbb{Z})$ and $k > 0$ such that $f(c)$ does not belong to $A + l^k M$. Let $p \in Q$ be an extension of $p_l$. By compactness and finiteness of $[M : lM]$, let $\alpha \in M$ be such that $p(x) \models \mathfrak{d}_{l^N} (f(x) - \alpha)$ for all $N > 0$. As $[M : A + l^k M] > 1$, we must have by proposition 4.3.17 the inequalities $1 < [M : l^k M] = [M : lM]^k$, thus $1 < [M : lM] = [l^k M : l^{k+1} M]$, so let $m \in l^k M \smallsetminus l^{k+1} M$. Then, by lemma 4.3.19, $\mathrm{tp}_l(\alpha/A) = \mathrm{tp}_l(\alpha+m/A)$. Let $D$ be a special subgroup of $M$ of size $\leqslant \max(2^{\aleph_0}, |A|)$ containing $A$, and a system of representatives of the cosets of $\bigcap_N l^N M$. By saturation of $M$, we can find $\beta \in M$ such that $\mathrm{tp}_l(\beta/D) = \mathrm{tp}_l(\alpha - m/D)$ (which implies $\beta - \alpha - m \in \bigcap_N l^N M$), and $\mathrm{tp}_j(\beta/D) = \mathrm{tp}_j(\alpha/D)$ for all $j \neq l$ (which implies $\mathrm{tp}_j(\beta/A) = \mathrm{tp}_j(\alpha/A)$). Then, by strong homogeneity of $M$, there exists $\sigma \in \mathrm{Aut}(M/A)$ such that $\sigma(\alpha) = \beta$. Now, we have $\sigma(p)(x) \models \mathfrak{d}_{l^{k+1}} (f(x) - \alpha - m)$, a formula that is inconsistent with $\mathfrak{d}_{l^{k+1}} (f(x) - \alpha)$. As a result, $\sigma(p) \neq p$, concluding the proof. □

*Remark* 4.3.23. Just like in subsection 4.3.2, we would like to describe the space of global invariant/bounded extensions of $p_l$. In case the condition of proposition 4.3.22 holds, we saw that $p_l$ is $\mathrm{Aut}(M/A)$-invariant and complete in $S_l^n(M)$, thus this space is a single point. In case the condition does not hold, the space of invariant extensions of $p_l$ is obviously empty, but we would like to describe the topological space of all the extensions that have a bounded orbit. This is what we do in the remainder of this subsection, we exhibit a homeomorphism between this orbit and some closed subspace of the Cantor space.

**Definition 4.3.24.** Let $N > 0$. We define $\mathcal{T}_N$ the tree $\{0, \dots N-1\}^{<\omega}$, and $\mathcal{B}_N = \{0, \dots N-1\}^{\omega}$ the set of its branches. The $m$-th-level of $\mathcal{T}_N$ is the set $\{0, \dots N-1\}^m$. We see $\mathcal{B}_N$ as a topological space, with the elementary open sets being the $O_w = \{b \in \mathcal{B}_N | w \text{ is a prefix of } b\}$ for every $w \in \mathcal{T}_N$. For every $m$, the opens $(O_w)_{w \in \{0,\dots N-1\}^m}$ form a finite disjoint open cover of $\mathcal{B}_N$, thus each $O_w$ is in fact clopen.

**Lemma 4.3.25.** *For all $N, d > 0$, $\mathcal{B}_N$ is naturally homeomorphic to $\mathcal{B}_{N^d}$.*

*Proof.* The natural bijection $\{0, \dots N^d - 1\} \longrightarrow \{0, \dots N-1\}^d$ extends to a bijection $h\colon \mathcal{B}_{N^d} \longrightarrow \mathcal{B}_N$. For each $m < \omega$, we have:

$$\{h(O_w) | w \text{ in the } m\text{-th level of } \mathcal{T}_{N^d}\} = \{O_w | w \text{ in the } md\text{-th level of } \mathcal{T}_N\},$$

thus $h^{-1}$ is continuous. To conclude, for each $w \in \mathcal{T}_N$, $O_w$ can be written as a finite union of elementary opens $O_v$, with $v$ in the $d\lceil |w|/d \rceil$-th level of $\mathcal{T}_N$, thus $h^{-1}(O_w)$ can be written as a finite union of some $O_v$, with $v$ in the $\lceil |w|/d \rceil$-th level of $\mathcal{T}_{N^d}$, therefore $h$ is continuous. □

**Lemma 4.3.26.** *The two spaces $\mathcal{B}_l$ and $\mathbb{Z}_l$ (the topological space of $l$-adic integers) are naturally homeomorphic.*

Recall that with our conventions, the $l$-adic valuation of $l^m$ is $-m$.

*Proof.* The tree of non-empty balls of $\mathbb{Z}_l$ (including the whole $\mathbb{Z}_l$) ordered by inclusion is isomorphic to $\mathcal{T}_l$ via the isomorphism mapping each $w \in \mathcal{T}_l$ to the ball of radius $-|w|$ around $\sum_m w_m l^m$. The corresponding bijection $\mathcal{B}_l \longrightarrow \mathbb{Z}_l$ is clearly $h\colon w \longmapsto \sum_{m<\omega} w_m l^m$, which concludes the proof. □

**Lemma 4.3.27.** *Suppose $[M : lM] = l^d$ is finite. Then $S_l^n(M)$ is naturally homeomorphic to the direct product $\mathcal{B}_{l^d}^n$ (with $n = |c|$).*

*Proof.* Note that each point of $S_l^n(M)$ is a partial type that can be written $\{x \in Y_m | m < \omega\}$, with $Y_m \in \left(M \diagup l^m M\right)^n$. In particular, $S_l^n(M)$ is in natural bijection with $\left(M \diagup \bigcap_{m<\omega} l^m M\right)^n$. For any $\mathcal{L}_P(M)$-term $t(x)$, the definable subset of $M$ defined by $\mathfrak{d}_{l^m}(t(x))$ is a union of cosets mod $l^m M$, thus the elementary clopen subsets of $\left(M \diagup \bigcap_{m<\omega} l^m M\right)^n$ induced by $S_l^n(M)$ are the Cartesian products of cosets of $l^m M$ for $m < \omega$. We build by induction a sequence of bijections $(f_m)_m$ between the $m$-th level of $\mathcal{T}_{l^d}$ and ${}^{M}\!/\!_{l^m M}$. Suppose by induction hypothesis that we have $f_m$ ($f_0$ is the bijection between the point $\left\{{}^{M}\!/\!_{M}\right\}$ and the point $\{\varnothing\}$). For each $w$ in the $m$-th level, $w$ has exactly $l^d$ children in the $m+1$-th level, and $f_m(w)$ is the union of exactly $l^d$ cosets $\bmod\ l^{m+1}M$, so there is a bijection between them for each $w$, and we set $f_{m+1}$ as the union of those bijections. Then the $(f_m)_m$ define an homeomorphism between $\mathcal{B}_{l^d}$ and ${}^{M}\!/\!_{\bigcap_m l^m M}$. This concludes the proof. □

**Corollary 4.3.28.** *The closed space of every global extension of $p_l$ of bounded orbit under the action of* $\mathrm{Aut}(M/A)$ *is naturally homeomorphic to some closed subspace of $\mathbb{Z}_l^n$.*

### 4.3.4 Partial types using a prime of infinite index

We still fix a prime $l$, and we assume here that $[M : lM]$ is infinite.

**Lemma 4.3.29.** *Let $D$ be a small special subgroup of $M$, and $\alpha = \alpha_1 \dots \alpha_m$ a finite tuple from $M$ which is $\mathbb{Q}$-free over $D$. Let $\beta = \beta_1 \dots \beta_m$ be another tuple from $M$ such that* $\mathrm{tp}_l(\beta/D) = \mathrm{tp}_l(\alpha/D)$. *Then there exists $\gamma = \gamma_1 \dots \gamma_m \equiv_D \alpha$ such that $\gamma_i - \beta_i \in \bigcap_N l^N M$ for all $i$.*

*Proof.* As $\bigcap_N l^N M$ is large, choose by induction $\varepsilon_i \in \bigcap_N l^N M$ such that $\varepsilon_i$ is not in the special subgroup generated by $D\alpha\beta\varepsilon_{<i}$, and let $\beta'_i = \beta_i + \varepsilon_i$. Then we have $\beta_i - \beta'_i \in \bigcap_N l^N M$, and $\varepsilon$ is $\mathbb{Q}$-free over the special subgroup generated by $D\alpha\beta$. As a result, $\alpha$ (and hence $(\alpha_i - \beta'_i)_i$) is $\mathbb{Q}$-free over $D'$, the special subgroup generated by $D\beta'$. By lemma 4.2.2, there exists in $M$ a tuple $e = e_1 \dots e_m$ such that $e_i \in \bigcap_N l^N M$ for all $i$, and $\mathrm{tp}_j(e/D') = \mathrm{tp}_j((\alpha_i - \beta'_i)_i/D')$ for all $j \neq l$. Let $\gamma_i = \beta'_i + e_i$. Then, for all $j \neq l$, we have $\mathrm{tp}_j(\gamma/D') = \mathrm{tp}_j(\alpha/D')$. Moreover, we have $\gamma_i - \beta'_i \in \bigcap_N l^N M$ for all $i$, thus $\mathrm{tp}_l(\gamma/D) = \mathrm{tp}_l(\beta'/D) = \mathrm{tp}_l(\beta/D) = \mathrm{tp}_l(\alpha/D)$, so we conclude that $\gamma \equiv_D \alpha$. □

*Remark* 4.3.30. Recall that any discrete model of ROAG is an elementary extension of $\mathbb{Z}$. In particular, if $[M : lM]$ is infinite for some $l$, then $M$ must be dense. In that case, the special subgroups of $M$ (which are exactly the definably closed sets) are its pure subgroups.

**Lemma 4.3.31.** *Let $(\alpha_i)_{i\in I}$ be a finite tuple from $M$ which is $\mathbb{Q}$-free over $B$. Suppose that for each $i \in I$, there exists $N > 0$ such that $\alpha_i \notin B + l^N M$, and choose $N_i$ the least of those integers for each $i \in I$. Then, there exists $(\beta_i)_i \equiv_B (\alpha_i)_i$ such that $\beta_i - \alpha_k \notin l^{N_i} M$ for each $i, k \in I$.*

*Proof.* Note that by remark 4.3.30, the special subgroups of $M$ are exactly its pure subgroups.

Let $b_i \in B$ such that $b_i - \alpha_i \in l^{N_i - 1} M$. Suppose we can find a witness $(\gamma_i)_i$ of the statement with $(\alpha_i)_i$ replaced by $\left(\frac{\alpha_i - b_i}{l^{N_i-1}}\right)_i$ ($N_i$ will be replaced by 1). Then we can choose $\beta_i = l^{N_i-1}\gamma_i + b_i$. For the remainder of the proof, we can suppose without loss that $N_i = 1$ for every $i$ (in particular, $b_i = 0$).

Let $V$ be the large $\mathbb{F}_l$-vector space ${M}\diagup{lM}$, and $U \leqslant V$ the $\mathbb{F}_l$-vector subspace ${(B + lM)}\diagup{lM}$. Note that $U$ is small, as it is the image by $B$ of the map $M \longrightarrow {M}\diagup{lM}$. Note also that we have $\alpha_i \in V \setminus U$ for every $i \in I$. Now, let $I_0 \subseteq I$ such that $(\alpha_i)_{i\in I_0}$ is a lift of an $\mathbb{F}_l$-basis of the image in ${V}\diagup{U}$ of $U + \sum_{i\in I} \mathbb{F}_l \cdot (\alpha_i \bmod lM)$. Let $D$ be the special subgroup of $M$ generated by $B(\alpha_i)_{i\in I}$, and let $U' = {(D + lM)}\diagup{lM}$, a small vector subspace of $V$ containing $U$. As ${V}\diagup{U'}$ is large, one can choose an arbitrary family $(\beta_i)_{i\in I_0} \in M^{I_0}$ whose image in ${V}\diagup{U'}$ is $\mathbb{F}_l$-free. Now, for each non-zero $f \in \mathrm{LC}^{|I_0|}(\mathbb{Z})$, if $m$ is the largest integer such that all the coefficients of $f$ lie in $l^m\mathbb{Z}$, then we have:

$$f((\beta_i)_{i\in I_0}) \in l^m M \ni f((\alpha_i)_{i\in I_0})$$

$$f((\beta_i)_{i\in I_0}) \notin B + l^{m+1} M \not\ni f((\alpha_i)_{i\in I_0})$$

thus we have $\mathrm{tp}_l((\alpha_i)_{i\in I_0}/B) = \mathrm{tp}_l((\beta_i)_{i\in I_0}/B)$ by the same reasoning as in corollary 4.3.20. For each $i \in I_0$, as $\beta_i \bmod lM \notin U'$, we have $\beta_i \notin D + lM$, thus $\beta_i - \alpha_k \notin lM$ for each $k \in I$. By lemma 4.3.29, we can suppose without loss that $(\beta_i)_{i\in I_0} \equiv_B (\alpha_i)_{i\in I_0}$.

Choose $\sigma \in \mathrm{Aut}(M/B)$ such that $\sigma(\alpha_i) = \beta_i$ for every $i \in I_0$. Let $\beta_i = \sigma(\alpha_i)$ for each $i \in I \setminus I_0$. Then we have $(\beta_i)_{i\in I} \equiv_B (\alpha_i)_{i\in I}$. In order to conclude, we need to show that $\beta_i - \alpha_k \notin lM$ for each $i \in I \setminus I_0$, $k \in I$. Choose $i \in I \setminus I_0$. There must exist $f_i \in \mathrm{LC}(\mathbb{Z})$ such that $\alpha_i - f_i((\alpha_k)_{k\in I_0}) \in B + lM$. Choose $e_i \in B$ such that $\alpha_i - f_i((\alpha_k)_{k\in I_0}) - e_i \in lM$. Then $\beta_i - f_i((\beta_k)_{k\in I_0}) - e_i \in lM$, so we just have to show that $f_i((\beta_k)_{k\in I_0}) \notin D + lM$. As $\alpha_i \notin B + lM$, some coefficient of $f_i$ must be coprime with $l$, which concludes the proof as $(\beta_k \bmod lM)_{k\in I_0}$ is $\mathbb{F}_l$-free over $U'$. ☐

**Lemma 4.3.32.** *Let $(\alpha_i)_i$ be some tuple from $M$. Then there exists in $M$ a tuple $(\beta_i)_i$ which is $\mathbb{Q}$-free over $B$, and for which $\alpha_i - \beta_i \in \bigcap_N l^N M$ for every $i$.*

*Proof.* One merely has to define by induction $\varepsilon_i \in \bigcap_N l^N M$ which does not lie in the small special subgroup generated by $B(\alpha_k)_k \varepsilon_{<i}$, and choose $\beta_i = \alpha_i + \varepsilon_i$. □

**Lemma 4.3.33.** *Let $N > 0$. Suppose we have $f_1, \dots f_N \in \mathrm{LC}^n(\mathbb{Z})$, as well as $Y_1, \dots Y_N \in \bigcup_m M \diagup l^m M$, such that every tuple of $M$ which realizes $p_l$ satisfies the formula $\bigvee_i f_i(x) \in Y_i$. Then at least one of the $Y_i$ must have a point in $B$.*

*Proof.* Let $\alpha_i \in Y_i$ for each $i$. By lemma 4.3.32, we can assume without loss that $(\alpha_i)_i$ is $\mathbb{Q}$-free over $B$. Suppose by contradiction that none of the $Y_i$ has a point in $B$. Then we can apply lemma 4.3.31 to $(\alpha_i)_i$, we find $(\beta_i)_i \equiv_B (\alpha_i)_i$ such that, for all $m > 0$ and for all $i$, $k$, if $\alpha_i \notin B + l^m M$, then $\beta_i - \alpha_k \notin l^m M$. In particular, as $Y_k$ does not intersect $B$, we have $\beta_i \notin Y_k$ for all $i$, $k$. Let $\sigma \in \mathrm{Aut}(M/B)$ such that $\sigma(\alpha_i) = \beta_i$ for each $i$. Then we must have $\left(\bigcup_i Y_i\right) \cap \left(\bigcup_i \sigma(Y_i)\right) = \varnothing$. However, as $p_l$ is $\mathrm{Aut}(M/B)$-invariant, we must have $p_l(x) \vDash \bigvee_i f_i(x) \in \sigma(Y_i)$, a contradiction. □

Let us recall that we adopted assumptions 4.3.1 in the section.

**Proposition 4.3.34.** *Consider the action of $\mathrm{Aut}(M/A)$ on $Q = S^n_l(M)$. Let $C$ be the group of all $\mathbb{Z}$-linear combinations of $c$. Then the following conditions are equivalent:*

1. *Some invariant element of $Q$ extends $p_l$.*
2. *Some element of $Q$ of bounded orbit extends $p_l$.*
3. *$p_l$ does not fork over $A$.*
4. *$p_l$ does not divide over $A$.*
5. *For each $N > 0$, we have $(C + l^N M) \cap (B + l^N M) = A + l^N M$.*

*Proof.* The directions $1 \Longrightarrow 2$ and $3 \Longrightarrow 4$ are trivial.

Let us prove $2 \Longrightarrow 3$. For each prime $l' \neq l$, we can define $q_{l'}$ just as in proposition 4.3.6. There remains to show that there exists $q_< \in \pi_<(F)$ whose orbit is $\mathrm{Aut}(M/A)$-invariant. In some elementary extension of $M$, define by induction $(a_i)_i, (D_i)_i$ such that $D_i$ is the special subgroup generated by $Ma_{<i}$, and $a_i \in \bigcap_{b \in D_i} ]b, +\infty[$. Define $q_< = \mathrm{tp}_<(a/M)$. Then, for all $(\lambda_1 \dots \lambda_n) \in \mathbb{Z}^n$, and $b \in M$, whether $\sum_i \lambda_i a_i > b$ depends uniquely on the sign of $(\lambda_i)_i$ in the anti-lexicographic sum $\mathbb{Z}^n$, and it does not depend on $b$. Moreover, $(a_i)_i$ is clearly $\mathbb{Q}$-free over $M$. As a result, $q_<$ is an $\mathrm{Aut}(M)$-invariant (and hence $\mathrm{Aut}(M/A)$-invariant) element of $\pi_<(F)$, and we conclude.

Let us prove $4 \implies 5$. Let $f \in \mathrm{LC}(\mathbb{Z})$ and $N > 0$ be such that $f(c)$ is in $B + l^N M \smallsetminus A + l^N M$. Let $b \in B$ such that $f(c) - b \in l^N M$ (we know that $b \notin A + l^N M$ by hypothesis, in particular it is $\mathbb{Q}$-free over $A$ as a singleton). Let us show that the formula $\varphi(x, b) \coloneqq \mathfrak{d}_{l^N}(f(x) - b)$ divides over $A$, which implies that 4 fails. We can repeatedly apply lemma 4.3.31 to find $(b_i)_{i<\omega}$ such that $b_0 = b$, $b_i \equiv_B b_{i+1}$, and $b_i - b_j \notin l^N M$ for all $j \neq i$. Then the set of formulas $\{\varphi(x, b_i) | i < \omega\}$ is clearly 2-inconsistent, so we can conclude.

Let us prove $5 \implies 1$. Suppose 5 holds, and let us build an explicit invariant element of $Q$ which extends $p_l$, witnessing 1. Define:

$$F = \left\{(f, N) \in \mathrm{LC}^n(\mathbb{Z}) \times \omega | f(c) \in A + l^N M\right\}$$

For each $(f, N) \in F$, let $a_{f,N} \in A$ such that $f(c) - a_{f,N} \in l^N M$. Define the following partial type:

$$p(x) = \{\mathfrak{d}_{l^N}(f(x) - a_{f,N}) | (f, N) \in F\} \cup \{\neg \mathfrak{d}_{l^N}(f(x) - e) | (f, N) \notin F, e \in M\}$$

Suppose by contradiction that $p$ is not consistent with $p_l$. Then, by compactness, there exists a finite tuple $(f_i, N_i, e_i)_i$ such that, for every $i$, $(f_i, N_i) \notin F$, and $p_l(x) \vDash \bigvee_i \mathfrak{d}_{l^{N_i}}(f_i(x) - e_i)$. By lemma 4.3.33 applied to $Y_i = e_i \bmod l^{N_i} M$, there must exist $i$ such that $e_i \in B + l^{N_i} M$. As $(f_i, N_i) \notin F$, we must have $e_i \notin A + l^{N_i} M$, thus $e_i \in B + l^{N_i} M \smallsetminus A + l^{N_i} M$. By hypothesis 5, we have $p_l(x) \vDash \neg \mathfrak{d}_{l^{N_i}}(f_i(x) - e_i)$, therefore $p_l(x) \vDash \bigvee_{k \neq i} \mathfrak{d}_{l^{N_k}}(f_k(x) - e_k)$. We keep decreasing by induction the size of this disjunction, until we inevitably reach a contradiction. As a result, $p$ must be consistent with $p_l$, and it is clearly $\mathrm{Aut}(M/A)$-invariant. This partial type is of course complete as an element of $Q$, because every atomic formula with predicate in $\{\mathfrak{d}_{l^N} | N > 0\}$ and parameters in $M$ either lies in $p$, or its negation lies in $p$. □

*Remark* 4.3.35. When the conditions of proposition 4.3.34 hold, the explicit type in the proof of $5 \implies 1$ is actually the unique global extension of $p_l$ which has a bounded orbit. Indeed, if $q$ is another global extension, and $d$ is a realization of $q$ in some elementary extension, then we must have $f(d) - \alpha \in l^N M$ for some $(f, N) \notin F$, and $m \in M$. In that case, $m \notin A + l^N M$, thus the fifth condition from proposition 4.3.34 fails with $c, B$ replaced by $d, M$. In particular, $q$ divides over $A$, thus its orbit is unbounded.

## 4.4 Computation of forking

**Theorem 4.4.1.** *Let $M \vDash$ROAG, let $A$, $B$ be parameter subsets of $M$, let $\kappa = \max(|AB|, 2^{\aleph_0})^+$, and let $c = c_1 \dots c_n \in M^n$. Suppose $M$ is $\kappa$-saturated and strongly $\kappa$-homogeneous. Let $C$ be the subgroup of $M$ generated by $c$, and $A', B'$ the special subgroups generated by $A, AB$. Then the following conditions are equivalent:*

- $c \mathop{\downarrow}^{\mathbf{f}}_{A} B$
- $c \mathop{\downarrow}^{\mathbf{d}}_{A} B$
- $c \mathop{\downarrow}^{\mathbf{bo}}_{A} B$
- $\mathrm{tp}(c/AB)$ *has a global* $\mathrm{Aut}(M^{eq}/A \cup \mathrm{acl}^{eq}(\varnothing))$*-invariant extension.*
- *the following conditions hold:*
    1. *Every closed bounded interval of $B'$ that has a point in $C$ already has a point in* $\mathrm{div}(A')$.
    2. *For all prime $l$, if $[M : lM]$ is infinite, then we have, for all $N > 0$, $(C + l^N M) \cap (B' + l^N M) = A' + l^N M$.*

*Moreover, $c \mathop{\downarrow}^{\mathbf{inv}}_{A} B$ if and only if the above conditions hold, and, additionally, for every prime $l$ for which $[M : lM]$ is finite, we have $C \subseteq A' + \bigcap_{N>0} l^N M$.*

Note that, if $c \mathop{\not\downarrow}^{\mathbf{d}}_{A} B$, then $\mathrm{tp}(c/AB)$ divides over $\mathrm{acl}^{eq}(A)$, thus it does not admit any global $\mathrm{Aut}(M^{eq}/A \cup \mathrm{acl}^{eq}(\varnothing))$-invariant extension.

*Proof.* Firstly, we have $A' = \mathrm{dcl}(A)$ and $B' = \mathrm{dcl}(AB)$ by corollary 4.2.6, thus we have $c \mathop{\downfree}_A B$ if and only if $c \mathop{\downfree}_{A'} B'$ for every $\downfree \in \{\downfree^{\mathbf{f}}, \downfree^{\mathbf{d}}, \downfree^{\mathbf{bo}}, \downfree^{\mathbf{inv}}\}$. Secondly, given $c' = c'_1 \dots c'_m$ a maximal subtuple of $c$ which is $\mathbb{Q}$-free over $A'$, $c$ and $c'$ are $A$-interdefinable, thus one can show that for all of those independence notions we have $c \mathop{\downfree}_{A'} B'$ if and only if $c' \mathop{\downfree}_{A'} B'$. Likewise, $\mathrm{tp}(c/AB)$ has an $\mathrm{Aut}(M^{eq}/A \cup \mathrm{acl}^{eq}(\varnothing))$-invariant extension if and only if $\mathrm{tp}(c'/B')$ has an $\mathrm{Aut}(M^{eq}/A' \cup \mathrm{acl}^{eq}(\varnothing))$-invariant extension.

If $c'$ was not $\mathbb{Q}$-free over $B'$, then one could show that $c' \mathop{\not\downfree}^{\mathbf{d}}_{A'} B'$ (hence $c' \mathop{\not\downfree}^{\mathbf{f}}_{A'} B'$, $c' \mathop{\not\downfree}^{\mathbf{bo}}_{A'} B'$, and $\mathrm{tp}(c/AB)$ does not admit any global extension which is $\mathrm{Aut}(M^{eq}/A \cup \mathrm{acl}^{eq}(\varnothing))$-invariant). Moreover, condition 1 would fail on some singleton of $B \smallsetminus A$ witnessing the fact that $c'$ is not $\mathbb{Q}$-free over $B$.

Suppose $c'$ is $\mathbb{Q}$-free over $B'$. For each $j \in J$, let $p_j = \mathrm{tp}_j(c'/B')$, and let $F_j$ be the image by $\pi_j$ of the space of complete global types with realizations that are $\mathbb{Q}$-free over $M$.

If condition 1 fails, then $p_<$ divides over $A'$ by proposition 4.3.6, and if condition 2 fails for some $l$, then $p_l$ divides over $A'$ by proposition 4.3.34. Either way, we have $c' \mathop{\not\downfree}_{A'} B'$ for every $\downfree \in \{\downfree^{\mathbf{f}}, \downfree^{\mathbf{d}}, \downfree^{\mathbf{bo}}, \downfree^{\mathbf{inv}}\}$.

Suppose conditions 1 and 2 hold, then $p_<$ extends to an $\mathrm{Aut}(M/A')$-invariant (thus $\mathrm{Aut}(M^{eq}/A' \cup \mathrm{acl}^{eq}(\varnothing))$-invariant) element of $F_<$ by proposition 4.3.6, and, for every prime $l$ for which $[M : lM]$ is infinite, $p_l$ extends to some $\mathrm{Aut}(M/A')$-invariant extension of $F_l$ by 4.3.34. By remark 4.3.21, for each prime $l$ of finite index, any global extension of $p_l$ is automatically $\mathrm{Aut}(M^{eq}/A' \cup \mathrm{acl}^{eq}(\varnothing))$-invariant, thus it follows that $\mathrm{tp}(c'/A'B')$ has a global $\mathrm{Aut}(M^{eq}/A' \cup \mathrm{acl}^{eq}(\varnothing))$-invariant extension. Moreover, by remark 4.3.21, we know that $p_l$ extends to some element of $F_l$ of bounded orbit under $\mathrm{Aut}(M/A')$ for every prime $l$ for which $[M : lM]$ is finite. By corollary 4.2.7, $c' \mathop{\downfree}^{\mathbf{bo}}_{A'} B'$, thus $c' \mathop{\downfree}^{\mathbf{f}}_{A'} B', c' \mathop{\downfree}^{\mathbf{d}}_{A'} B'$. Moreover, by corollary 4.2.7, we have $c' \mathop{\downfree}^{\mathbf{inv}}_{A'} B'$ if and only if, for every prime $l$ for which $[M : lM]$ is finite, $p_l$ extends to some $\mathrm{Aut}(M/A')$-invariant element of $F_l$. Then we can conclude using proposition 4.3.22. □

**Corollary 4.4.2.** *We have* $\downfree^{\mathbf{f}} = \downfree^{\mathbf{Sh}}$ *in* ROAG.

Recall from definition 1.3.6 that $c \mathop{\downfree}^{\mathbf{Sh}}_A B$ is equivalent to $\mathrm{tp}(c/AB)$ having a global $\mathrm{Aut}(M^{eq}/\mathrm{acl}^{eq}(A))$-invariant extension, which is a weaker condition than being $\mathrm{Aut}(M^{eq}/A \cup \mathrm{acl}^{eq}(\varnothing))$-invariant.

We know (see remark 3.1.1) that in any NIP theory, $\mathop{\downarrow}^{\mathbf{f}}$ coincides with $\mathop{\downarrow}^{\mathbf{inv}}$ over models. We can now refine that result for ROAG:

**Corollary 4.4.3.** *Let $A'$ be the special subgroup generated by $A$. Then we have $\mathop{\downarrow}^{\mathbf{f}}_A = \mathop{\downarrow}^{\mathbf{inv}}_A$ if and only if $A' + \bigcap_N l^N M = M$ for all primes $l$ of finite index.*

*Moreover, $\mathop{\downarrow}^{\mathbf{f}} = \mathop{\downarrow}^{\mathbf{inv}}$ if and only if, for each prime $l$, either $M$ is $l$-divisible, or the index $[M : lM]$ is infinite.*

*Remark* 4.4.4. Suppose that $M$ is dense. Then, by proposition 4.3.16, remark 4.3.35, corollary 4.3.28, and remark 3.3.31, we get the three following properties :

- If $c \mathop{\downarrow}^{\mathbf{f}}_A B$, then every global extension of $p_<$ which does not fork over $A$ is $\mathrm{Aut}(M/A)$-invariant. Moreover, the Stone space of those extensions can be written as $S_1$, a finite product of finite coproducts of closed subspaces of $S_{\mathrm{DOAG}}(\{0\})$.
- Let $l$ be prime such that $[M : lM]$ is finite. If $C \subseteq A' + \bigcap_N l^N M$, then $p_l$ has only one global extension, and it is $\mathrm{Aut}(M/A)$-invariant. Else, none of the global extensions of $p_l$ is $\mathrm{Aut}(M/A)$-invariant, but they are all non-forking over $A$, and the space of those extensions can be written $S_{2,l}$, a closed subspace of $\mathbb{Z}_l^{|c'|}$.
- Let $l$ be prime such that $[M : lM]$ is infinite. If $c \mathop{\downarrow}^{\mathbf{f}}_A B$, then $p_l$ has only one global extension which does not fork over $A$, and this extension is $\mathrm{Aut}(M/A)$-invariant.

As a result, if we define:

$$L = \left\{ l \text{ prime } \middle| [M : lM] < \omega, C \nsubseteq A' + \bigcap_N l^N M \right\},$$

then the space $S$ of global extensions of $\mathrm{tp}(c/AB)$ which do not fork over $A$ is either empty, or homeomorphic to $S_1 \times \prod_{l \in L} S_{2,l}$.

Moreover, if $S \neq \varnothing$ and $L = \varnothing$, then all the elements of $S$ are $\mathrm{Aut}(M/A)$-invariant, thus $S$ is homeomorphic to $S_1$, and it coincides with the space of global $\mathrm{Aut}(M/A)$-invariant extensions of $\mathrm{tp}(c/AB)$.

## 4.5 Examples

### 4.5.1 Non-forking versus finite satisfiability

Let us start with the very basic remark that the cut independence notion does not coincide with finite satisfiability in DOAG. In an elementary extension of the ordered group $\mathbb{R}$, let $\varepsilon > 0$, $c > 0$ such that $\Delta(c) < \Delta(\varepsilon) < \Delta(1)$. Let $A = \mathbb{Q}$, $B = \mathbb{Q} + \mathbb{Q}\varepsilon$. Then we clearly have $c \mathop{\downfree}^{\mathbf{cut}}_A B$, but the type of $c$ over $B$ is not finitely satisfiable over $A$: the $B$-definable set $]0, \varepsilon]$ contains $c$ and has no point in $A$.

### 4.5.2 Non-forking versus sequential independence

We would like to compare non-forking independence with a stronger independence notion:

**Definition 4.5.1.** Let $M$ be any first-order structure, and $A$, $B$, $C \subseteq M$ parameter sets. We write $C \mathop{\downfree}^{\mathbf{seq}}_A B$ when, for any finite tuple $c$ from $C$, there exists a finite tuple $d = (d_0, \dots d_n)$ from $\operatorname{acl}(AC)$ such that $\operatorname{acl}(Ac) \subseteq \operatorname{acl}(Ad)$ and $d_i \mathop{\downfree}^{\mathbf{f}}_{Ad_{\leqslant i}} B$ for every $i \leqslant n$.

This relation, which we call *sequential independence*, is stronger than $\downfree^{\mathbf{f}}$ in any theory (an easy application of left-transitivity of non-forking). One can show that it coincides with $\downfree^{\mathbf{f}}$ in simple theories, and this is also the case in some non-simple theories, like DLO (a parameter set coincides with its algebraic closure, and it is easy to guess the global invariant extensions of the type of some finite tuple). In a theory where $\downfree^{\mathbf{seq}} = \downfree^{\mathbf{f}}$, it is really easy to understand forking, because one can reduce the problem to forking in dimension one. We showed in the previous section that forking in ROAG is essentially a dimension one phenomenon, so it makes sense to ask whether $\downfree^{\mathbf{seq}} = \downfree^{\mathbf{f}}$ in these structures. We show that the answer is no with a counterexample in DOAG. The only other explicit example we know of $\downfree^{\mathbf{f}} \neq \downfree^{\mathbf{seq}}$ is in ([HHM05], example 13.6), in the theory of algebraically closed valued fields. This example does not occur exclusively in the value group, one needs field-theoretic tools to make it work, which is why our example is arguably less complicated.

In the reduct to the ordered group structure of an elementary extension of the real-closed field $\mathbb{R}$, let $\varepsilon > 0$ such that $\Delta(\varepsilon) < \Delta(1)$. Let $A = \mathbb{Q}$, $B = \mathbb{Q} + \mathbb{Q}\sqrt{2} + \mathbb{Q}\sqrt{3}$, and $c = (c_1, c_2) = (\sqrt{2} + \varepsilon, \sqrt{3} + \varepsilon\sqrt{2})$.

Let us show that we have $c \mathop{\downarrow}^{\mathbf{f}}_{A} B$. We know that any two distinct real numbers have a different cut over $\mathbb{Q} = A$, so any bounded $B$-definable closed interval that does not have any point in $A$ must be a singleton. As we have $\mathrm{dcl}(Ac) \cap \mathrm{dcl}(B) = A$, we have in fact $c \mathop{\downarrow}^{\mathbf{cut}}_{A} B$. In particular, we have $c \mathop{\downarrow}^{\mathbf{f}}_{A} B$.

Now let us show $c \not\mathop{\downarrow}^{\mathbf{seq}}_{A} B$. Let $d = (d_1, \dots d_n)$ be a finite tuple from $\mathrm{acl}(Ac)$ such that $\mathrm{acl}(Ac) \subseteq \mathrm{acl}(Ad)$. We can assume that the family $d$ is $\mathbb{Q}$-free over $A$, for if $d_i \in \mathrm{dcl}(Ad_{<i})$, then we trivially have $d_i \mathop{\downarrow}^{\mathbf{f}}_{Ad_{<i}} B$. We can also assume $d_i \in \mathbb{Q}c_1 + \mathbb{Q}c_2$, for if $a \in A$, then $\mathrm{acl}(Ad_{\leqslant i}) = \mathrm{acl}(Ad_{<i}(d_i + a))$, so we can freely translate any $d_i$ by an element of $A$. By basic linear algebra, we have $n = 2$, and $(d_1, d_2)$ is a $\mathbb{Q}$-basis of $\mathbb{Q}c_1 + \mathbb{Q}c_2$. Let $f_i$ be the unique map of $\mathrm{LC}^2(\mathbb{Q})$ such that $f_i(c_1, c_2) = d_i$. As $(c_1, c_2)$ is also $\mathbb{Q}$-free, $(f_1, f_2)$ must be $\mathbb{Q}$-free. For any $g \in \mathrm{LC}^2(\mathbb{Q})$, we have:

$$d_2 - g(1, d_1) = f_2(\sqrt{2}, \sqrt{3}) - g(1, f_1(\sqrt{2}, \sqrt{3})) + \varepsilon\left[f_2(1, \sqrt{2}) - g(0, f_1(1, \sqrt{2}))\right]$$

As $(f_1, f_2)$ and $(1, \sqrt{2}, \sqrt{3})$ are both $\mathbb{Q}$-free, we must have:

$$f_2(\sqrt{2}, \sqrt{3}) - g(1, f_1(\sqrt{2}, \sqrt{3})) \neq 0$$

so $\Delta(d_2 - g(1, d_1)) = \Delta(1)$. We established that, for all $x \in \mathrm{dcl}(Ad_1)$, we have $\Delta(d_2 - x) > \Delta(\varepsilon)$, which implies $\varepsilon \in H(d_2/\mathrm{dcl}(Ad_1))$. Now, if we had $d_2 \mathop{\downarrow}^{\mathbf{f}}_{Ad_1} B$, then we would necessarily have $\varepsilon \in H(d_2/\mathrm{dcl}(Bd_1))$ by proposition 3.2.9. Let us show that this fails. It is enough to show that $\Delta(d_2 - b) = \Delta(\varepsilon)$ for some $b \in \mathrm{dcl}(Bd_1)$. Choose $b = d_1 - f_1(\sqrt{2}, \sqrt{3}) + f_2(\sqrt{2}, \sqrt{3})$. Then, as $d_i = f_i(c_1, c_2) = f_i(\sqrt{2}, \sqrt{3}) + f_i(\varepsilon, \varepsilon\sqrt{2})$, we have $d_2 - b = f_2(\varepsilon, \varepsilon\sqrt{2}) - f_1(\varepsilon, \varepsilon\sqrt{2})$. As $(f_1, f_2)$ is $\mathbb{Q}$-free, $f_2 - f_1 \neq 0$, therefore $d_2 - b$ is a non-trivial linear combination of $(\varepsilon, \varepsilon\sqrt{2})$. This latter family is clearly $\Delta$-separated, thus $\Delta(d_2 - b) = \Delta(\varepsilon)$, and we can conclude.

This holds for any choice for the family $d$, so we have $c \not\mathop{\downarrow}^{\mathbf{seq}}_{A} B$, which concludes the example.

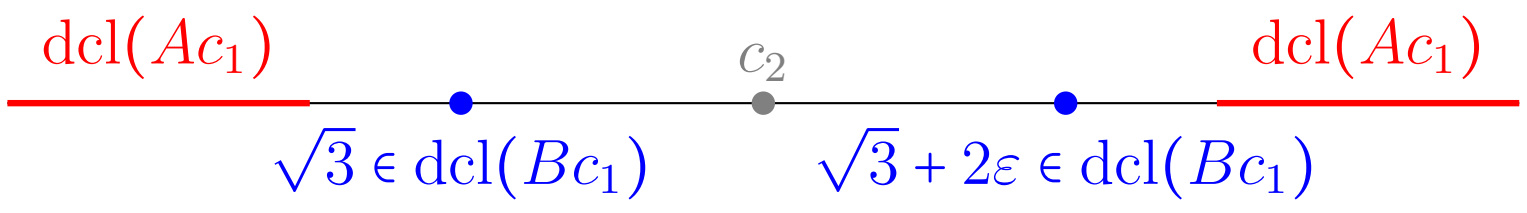

### 4.5.3 Weak versus strong versions of sequential independence

In the definition of sequential independence, we require that there exists a good enumeration $d$ that witnesses independence. Does it change anything to ask this for all the enumerations? We show that the answer is yes with a counterexample in DOAG.

**Definition 4.5.2.** Define $C \mathop{\downfree}^{\mathbf{seq}*}_A B$ when, for every finite tuple $c = (c_0..c_n)$ from $\mathrm{acl}(AC)$, we have $c_i \mathop{\downfree}^{\mathbf{f}}_{Ac_{<i}} B$ for all $i \leqslant n$.

This notion is at least as strong as $\downfree^{\mathbf{seq}}$, which is strictly stronger than $\downfree^{\mathbf{f}}$ in DOAG by the above example. Let us show that $\downfree^{\mathbf{seq}*}$ is strictly stronger than $\downfree^{\mathbf{seq}}$ in DOAG. Let $\varepsilon > 0$, $c_2 > 0$ such that $\Delta(\varepsilon) < \Delta(c_2) < \Delta(1)$. Let $A = \mathbb{Q}$, $B = \mathbb{Q} + \mathbb{Q}\sqrt{2}$, and $c_1 = \sqrt{2} + \varepsilon$.

With the same reasoning as in the previous example, we have $c_1 c_2 \mathop{\downfree}^{\mathbf{cut}}_A B$. In particular, we have $c_1 \mathop{\downfree}^{\mathbf{f}}_A B$ and $c_2 \mathop{\downfree}^{\mathbf{f}}_A B$.

As $\Delta(c_2) \notin \Delta(\mathrm{dcl}(Ac_1))$, $c_2$ is ramified over $\mathrm{dcl}(Ac_1)$. The only points of $\mathrm{dcl}(Bc_1)$ that have the same cut over $\mathrm{dcl}(Ac_1)$ as that of $c_2$ are the $\lambda\varepsilon$, with $\lambda \in \mathbb{Q}_{>0}$. They are all smaller than $c_2$, so $c_2$ leans right with respect to $\mathrm{dcl}(Ac_1)$ and $\mathrm{dcl}(Bc_1)$, and $c_2 \mathop{\downfree}^{\mathbf{f}}_{Ac_1} B$, which implies $c_1 c_2 \mathop{\downfree}^{\mathbf{seq}}_A B$.

Let us show $c_1 c_2 \mathop{\not\downfree}^{\mathbf{seq}*}_A B$. The interval $[\sqrt{2}, \sqrt{2} + c_2]$ contains $c_1$ and has no point in $\mathrm{dcl}(Ac_2)$, which implies $c_1 \mathop{\not\downfree}^{\mathbf{f}}_{Ac_2} B$.

### 4.5.4 Forking versus geometric cut-dependence

Let us look at a different definition for cut independence:

**Definition 4.5.3.** Let $M$ be the expansion of some totally ordered set, and $A$, $B$, $C$ parameter sets. We write $C \mathop{\not\downfree}^{\mathbf{cut}*}_A B$ when there exists $c \in \mathrm{dcl}(AC)$ and $b_1, b_2 \in \mathrm{dcl}(AB)$ such that $b_1 \leqslant c \leqslant b_2$, $b_1 \equiv_A b_2$ and $b_i \notin \mathrm{dcl}(A)$.

We clearly have the implications $\downfree^{\mathbf{d}} \subseteq \downfree^{\mathbf{cut}*}$ and:

$$C \mathop{\downfree}^{\mathbf{cut}*}_A B \Longleftarrow \mathrm{div}(\mathrm{dcl}(AC)) \mathop{\downfree}^{\mathbf{cut}}_{\mathrm{div}(\mathrm{dcl}(A))} \mathrm{div}(\mathrm{dcl}(AB))$$

in any ordered Abelian group. It is not hard to prove that these two implications are equivalences in DOAG. In a dense model of ROAG $M$, if $B$ is an $|A|^+$-saturated model containing $A$, then one may show that the second implication is an equivalence, i.e. $C \mathop{\downfree}^{\mathbf{cut}*}_A B$ holds in $M$ if and only if $\operatorname{div}(C) \mathop{\downfree}^{\mathbf{cut}}_{\operatorname{div}(A)} \operatorname{div}(B)$ holds in $\operatorname{div}(M)$.

Whether the first inclusion is an equality is a little more complicated. Let $G$ be an ordered Abelian group, $A$ a parameter set, and $B$ an $|A|^+$-saturated model containing $A$. In case $G$ is regular and $|G/lG|$ is infinite for some prime $l$ (this is possible, take for instance $G = \sum_{n<\omega} \mathbb{Z}r_n$, with $(r_n)_n$ a $\mathbb{Q}$-free subfamily of $\mathbb{R}$), the characterization of forking and dividing that we have in theorem 4.4.1 implies that $\mathop{\downfree}^{\mathbf{d}}_A B$ does not coincide with $\mathop{\downfree}^{\mathbf{cut}*}_A B$, because we have partial types $\bmod \equiv_l$ that divide over $A$. However, if $G$ is dp-minimal, we have indeed $\mathop{\downfree}^{\mathbf{f}}_A B = \mathop{\downfree}^{\mathbf{cut}*}_A B$. In fact, Simon showed in ([Sim11], Proposition 2.5) that, for $c$ a singleton, we have $c \mathop{\downfree}^{\mathbf{f}}_A B$ if and only if $c \mathop{\downfree}^{\mathbf{cut}*}_A B$ if $G$ is dp-minimal, even when $G$ is not regular. However, we will see an example where $G$ is a dp-minimal, poly-regular, non-regular dense ordered Abelian group where $c \mathop{\downfree}^{\mathbf{cut}*}_A B$ (and thus $c \mathop{\downfree}^{\mathbf{f}}_A B$) holds in $G$, but $c \mathop{\downfree}^{\mathbf{cut}}_{\operatorname{div}(A)} \operatorname{div}(B)$ does not hold in $\operatorname{div}(G)$.

Let $G$ be any ordered group elementarily equivalent to the lexicographical product $\mathbb{Z}\left[\frac{1}{2}\right] \times \mathbb{Z}\left[\frac{1}{2}\right]$. Then $G$ is clearly dp-minimal, poly-regular and dense. In $\mathbb{Z}\left[\frac{1}{2}\right] \times \mathbb{Z}\left[\frac{1}{2}\right]$, we have $H_3((1,0)) = \{0\} \times \mathbb{Z}\left[\frac{1}{2}\right] \neq \{0\}$, so $G$ is not regular, and it has exactly one proper non-trivial definable convex subgroup, let us call it $H(G)$. Let $A$ be a definably closed parameter subset of $G$, $B$ an $|A|^+$-saturated elementary extension of $G$, and $M$ a $|B|^+$-saturated elementary extension of $B$. Let $c \in H(M)$ such that $c > H(B)$. Then we clearly have $c \mathop{\downfree}^{\mathbf{cut}*}_A B$, hence $c \mathop{\downfree}^{\mathbf{f}}_A B$ by [Sim11]. Let us show $c \mathop{\not\downfree}^{\mathbf{cut}}_{\operatorname{div}(A)} \operatorname{div}(B)$ in $\operatorname{div}(M)$ (in fact, we even show $c \mathop{\not\downfree}^{\mathbf{cut}}_{\operatorname{div}(G)} \operatorname{div}(B)$). Let $b_1 \in H(B)$ such that $b_1 > H(G)$. Let $b_2 \in B$ such that $b_2 \notin H(B)$, but $b_2 < g$ for every $g \in G_{>0} \setminus H(G)$. Then the interval $I = [b_1, b_2]$ contains $c$.

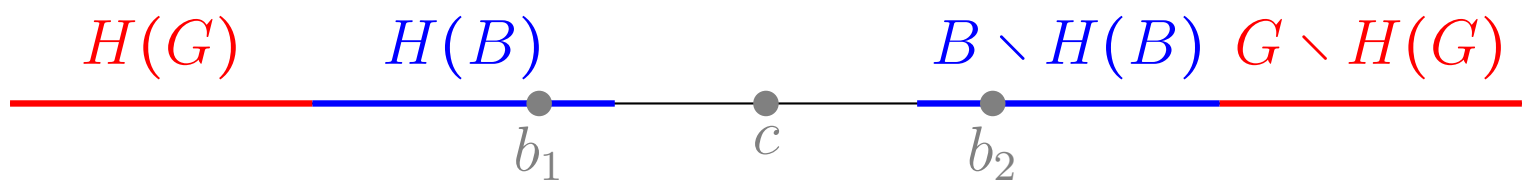

Let us show that $I$ does not have a point in $\mathrm{div}(G)$. Let $g \in \mathrm{div}(G)$ such that $0 < g < b_2$, and suppose by contradiction $g > b_1$. Let $0 < N < \omega$ such that $Ng \in G$. We have $Ng \geqslant g > b_1$, so $Ng \notin H(G)$. The ordered group $G \diagup HG$ is dense (as $\mathbb{Z}\left[\frac{1}{2}\right]$ is), so there must exist $g' \in G$ such that $g' \notin H(G)$ and $Ng' < (Ng \bmod H(G))$ (choose $N$ distinct cosets between $H(G)$ and $(Ng \bmod H(G))$, choose $g'$ a representative of the coset that corresponds to the minimal distance between them in $G \diagup H(G)$). As $g' \notin H(G)$ and $g' > 0$, we have $b_2 < g'$. However, we have $Ng' < Ng < Nb_2$, a contradiction. This concludes the example.

For future work, we should try to generalize our results to dp-minimal ordered groups. In such a group, one can hope that a global type over a monster model $B$ does not fork over $A$ if and only if its realizations $C$ satisfy $C \downarrow^{\mathbf{cut}*}_A B$. The dimension one case has already been proved by Simon.

### 4.5.5 Forking versus Presburger-quantifier-free dependence

For this last example, we study an ordered Abelian group $G$ with definably closed parameter sets $A \subseteq B$, and a singleton $c \in G$ for which $c \downarrow^{\mathbf{cut}}_{\mathrm{div}(A)} \mathrm{div}(B)$ in $\mathrm{div}(G)$ (hence $c \downarrow^{\mathbf{cut}*}_A B$), and:

$$\left(\mathrm{dcl}(Ac) + l^N G\right) \cap \left(B + l^N G\right) = A + l^N G$$

for every prime $l$, $N > 0$, but $c \not\downarrow^{\mathbf{d}}_A B$. The dividing definable set will not be built from the cuts or the cosets mod $l$, but from the chain of definable convex subgroups of $G$. Our example needs to have an infinite spine.

Let $\Gamma$ be the lexicographical product $\mathbb{Q}_{\leqslant 0} \times \mathbb{Z}_{\leqslant 0}$. Let $G$ be the ordered group of all Hahn series with coefficients in $\mathbb{Z}$ and powers in $\Gamma$ (the valuation coincides with the Archimedean valuation). Then $G$ is a dense ordered Abelian group. For each $g \in G$, let $S(g) \subseteq \Gamma$ be the support of $g$. Each $g \in G$ can be written $\sum_{\delta \in S(g)} \lambda_\delta t^\delta$, with $\lambda_\delta \in \mathbb{Z}$, and $t$ an indeterminate. If $\sigma$ is an automorphism of the ordered set $\mathbb{Q}_{\leqslant 0}$, then let $\bar{\sigma} \in \mathrm{Aut}(G)$ be defined as $\sum_{(\mu,n)\in\Gamma} \lambda_{(\mu,n)} t^{(\mu,n)} \longmapsto \sum_{(\mu,n)\in\Gamma} \lambda_{(\mu,n)} t^{(\sigma(\mu),n)}$. Let $A = \{a \in G | S(a) \subseteq (\{0\} \times \mathbb{Z}_{\leqslant 0})\}$ and $B = \mathrm{dcl}\left(A \cup \{t^{(-1,0)}\}\right)$. The action on $\mathbb{Q}_{\leqslant 0}$ of its automorphism group has only two orbits: $\{0\}$ and $\mathbb{Q}_{<0}$. As a result, we must have $\mathrm{dcl}(\varnothing) \subseteq A$. By ([CH11], Corollary 1.10), $\mathrm{dcl}(A)$ is exactly the relative divisible closure in $G$ of $A$, in other words $A$ is definably closed. With the same argument, we have $B = A + \mathbb{Z}t^{(-1,0)}$. Let $c = t^{(-1,-1)}$, so $\mathrm{dcl}(Ac) = A + \mathbb{Z}c$. We have $\Delta(c) < \Delta(B_{\neq 0})$, so $c$ clearly leans left with respect to $\mathrm{div}(A)$ and $\mathrm{div}(B)$ in $\mathrm{div}(G)$, which implies $c \mathop{\downarrow}^{\mathbf{cut}}_{A} B$ in $\mathrm{div}(G)$. Moreover, for every prime $l$, and $N > 0$, the elements of $l^N G$ are exactly the series $\sum_{\delta} l^N \lambda_\delta t^\delta$ with $\lambda_\delta \in \mathbb{Z}$, so we clearly have $(\mathrm{dcl}(Ac) + l^N G) \cap (B + l^N G) = A + l^N G$.

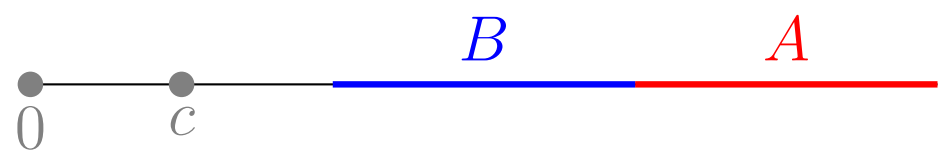

Now let us prove $c \not\mathop{\downarrow}^{\mathbf{d}}_{A} B$. For each $n < \omega$, let $\sigma_n \in \mathrm{Aut}(\mathbb{Q}_{\leqslant 0})$ such that $\sigma_n(1) = -n$. For each prime $l$, for each $g = \sum_{\delta} \lambda_\delta t^\delta \in G$, if $g \notin l^N G$, then one can show (see for instance [CH11], example 4.2) that $H_l^N(g)$ is the group of elements of $G$ with an Archimedean class smaller or equal to $\Delta(t^\delta)$, where $\delta$ is the largest element of $S(g)$ for which $\lambda_\delta \notin l^N \mathbb{Z}$. As a result, there does not exist $g \in G$ for which $H_2(c) < H_2(g) < H_2(t^{(-1,0)})$. The following $B$-definable set which contains $c$:

$$X = \left\{ g \in G | g \notin 2G \text{ and } \left( \forall h \in G,\ \left( H_2(h) \leqslant H_2(g) \text{ or } H_2(t^{(-1,0)}) \leqslant H_2(h) \right) \right) \right\}$$

divides over $A$, because the $(\bar{\sigma_n}(X))_n$ are pairwise-disjoint. This concludes the example.

# Chapter 5

# Short exact sequences and expansions

## 5.1 Preliminaries and resplendent quantifier elimination

In this section, we define what are pure short exact sequences of Abelian structures, and we recall results from previous literature.

**Notations** The notations that we use in this chapter differ from the rest of this report. Just as in [ACGZ22], we denote by $A$, $B$, $C$ the terms of the short exact sequences at hand, therefore they will refer to multi-sorted structures that are reducts of the short exact sequence. Parameter sets will be denoted by $\Gamma$, $\Delta$, $U$. Our tuples will always be tuples of reals, and they will be denoted by $u$, $v$, $w$... The real sorts will be in the disjoint union of $A^*$, $B$, $C^*$, with $A^*$ an expansion of $A$, and $C^*$ an expansion of $C$. As a result, a tuple $u$ will be the concatenation of three disjoint subtuples $u_A$, $u_B$, $u_C$, with of course $u_A \in (A^*)^{<\omega}$, $u_B \in B^{<\omega}$, $u_C \in (C^*)^{<\omega}$. We shall always use those subscripts in that fashion, and we will not redefine what they refer to.

**Definition 5.1.1.** Fix $\mathcal{L}$ a first-order language. We define the set of *p.p.-formulas* (for Primitive Positive) as the closure of the set of parameter-free atomic formulas under finite conjunction and existential quantification.

Let $A$ be a multi-sorted first-order structure with language $\mathcal{L}$. A *fundamental system* for $A$ is a set $\Phi$ of p.p.-formulas such that, modulo the theory of $A$, every p.p.-formula is equivalent to a finite conjunction of formulas of the form $\phi(t(x))$, with $\phi \in \Phi$, and $t$ a tuple of parameter-free terms.

We define $A$ to be an *Abelian structure* when the following conditions hold:

- On each sort $s$, the structure $A$ expands a structure of Abelian group on $s(A)$ .
- Any parameter-free term $t(x_{s_1} \dots x_{s_n})$ (seen as a $\varnothing$-definable function $s_1(A) \times \dots \times s_n(A) \longrightarrow s(A)$) is a group homomorphism.
- For every p.p.-formula $\phi(x_{s_1} \dots x_{s_n})$, the set $\phi(A)$ is a subgroup of $s_1(A) \times \dots \times s_n(A)$.

Given Abelian structures $A$, $B$ on the language $\mathcal{L}$, an embedding (for the language $\mathcal{L}$) $\iota : A \longrightarrow B$ is said to be *pure* when we have $\iota^{-1}(\phi(B)) = \phi(A)$ for every p.p.-formula $\phi$.

Let $\Phi$ be a set of p.p.-formulas. A *pure short exact sequence of Abelian structures expanded by* $\Phi$ ($\Phi$-PSES) is a first-order structure

$$\begin{array}{ccccc} A & \xrightarrow{\quad\iota\quad} & B & \xrightarrow{\quad\nu\quad} & C \\ & {\scriptstyle (\pi_\phi)_{\phi\in\Phi}} \searrow & & \swarrow {\scriptstyle (\rho_\phi)_{\phi\in\Phi}} & \\ & & \left({}^{A}\!\big/_{\phi(A)}\right)_{\phi\in\Phi} & & \end{array}$$

consisting of:

- Three Abelian structures $A$, $B$, $C$ on copies of $\mathcal{L}$ such that the respective sets of sorts of $A$, $B$, $C$ are pairwise-disjoint, together with the quotient groups ${}^{A}\!\big/_{\phi(A)}$ for each $\phi \in \Phi$.
- A pure embedding $\iota : A \longrightarrow B$, and a surjective morphism $\nu : B \longrightarrow C$ such that, on each sort of $\mathcal{L}$, we have $\mathrm{Im}(\iota) = \mathrm{Ker}(\nu)$.
- For each $\phi \in \Phi$, the quotient map $\pi_\phi : A \longrightarrow {}^{A}\!\big/_{\phi(A)}$.
- For each $\phi \in \Phi$, the map $\rho_\phi : B \longrightarrow {}^{A}\!\big/_{\phi(A)}$, which is zero outside of $\phi(B) + \iota(A)$, and extends the group homomorphism:

$$\phi(B) + \iota(A) \longrightarrow {}^{(\phi(B) + \iota(A))}\!\big/_{\phi(B)} \simeq {}^{A}\!\big/_{\phi(A)}$$

Note that it follows from the axioms of PSES that we have $\nu^{-1}(\phi(C)) = \phi(B) + \iota(A)$ for every p.p.-formula $\phi$. Equivalently, $\nu$ is a *pure projection*, that is for every p.p.-formula $\phi$, we have $\phi(C) = \nu(\phi(B))$.

For the remainder of this chapter, let $M^-$ be a $\Phi$-PSES:

$$0 \longrightarrow A \underset{\iota}{\longrightarrow} B \underset{\nu}{\longrightarrow} C \longrightarrow 0$$

with $\Phi$ a fundamental system for $B$. Let $\mathcal{L}_B$ be the language of $B$, and let $\mathcal{L}_A$, $\mathcal{L}_C$ be arbitrary enrichments of the respective languages of $A$ and $C$, such that the sort $A/_{\phi(A)}$ and the map $\pi_\phi$ are in $\mathcal{L}_A$ for every $\phi \in \Phi$. Let $A^*$ be an expansion of $A$ with language $\mathcal{L}_A$, and let $C^*$ an expansion of $C$ with language $\mathcal{L}_C$. Note that $A^*$, $C^*$ are completely arbitrary expansions of $A$ and $C$, they are in general not interpretable. Let $M$ be the expansion

$$0 \longrightarrow A^* \underset{\iota}{\longrightarrow} B \underset{\nu}{\longrightarrow} C^* \longrightarrow 0$$

of $M^-$ obtained by expanding $A$, $C$ to $A^*$, $C^*$. Our goal is to obtain a characterization for dividing and forking in $M$ relative to dividing and forking in $A^*(M)$ and $C^*(M)$.

A specific case of this very general context is important to understand (expansions of) valued fields. One can often reduce a problem taking place in a valued field $K$ to a problem on the short exact sequence:

$$0 \longrightarrow k \longrightarrow \mathrm{RV}^* \longrightarrow \Gamma \longrightarrow 0$$

with $k$ the multiplicative group of the residue field of $K$ expanded to the field structure $k$, $\mathrm{RV}^*$ the group of non-zero leading terms of $K$ (as in definition 2.2.4), and $\Gamma$ the value group of $K$ expanded to its ordered group structure and $-\infty$. In that case, the Abelian structures at hand are merely Abelian groups, and the fundamental system that we choose for $B$ is usually the set of formulas of the form $\exists y\ n\cdot y = x$, with $n$ a natural integer (including 0). The notion of pure embedding that we defined corresponds in this context to the usual notion of pure embeddings of Abelian groups. Reductions to RV also exist in common expansions of valued fields, such as differential valued fields, ordered valued fields or difference valued fields, and they may correspond to richer expansions of $k$ and $\Gamma$.

Now we describe known quantifier elimination results from previous literature.

**Definition 5.1.2.** Let $\Gamma \subseteq M$ be a substructure. We see $A^*$, $B$, $C^*$ as (tuples of) sorts, and we recall that $A^*(\Gamma)$, $B(\Gamma)$, $C^*(\Gamma)$ refer to the sets of elements of $\Gamma$ that are in $A^*$, $B$, $C^*$.

A $C$-*formula* with parameters in $\Gamma$ is a formula of the form:

$$\psi_C(\nu(x_B), x_C, \gamma)$$

with $\psi_C \in \mathcal{L}_C$ parameter-free, and $\gamma \in C^*(\Gamma)$. If $u$ is a tuple from $M$, then we write $\mathrm{tp}_C(u/\Gamma)$ for the $C$-type of $u$ over $\Gamma$, that is the partial type of every $C$-formula with parameters in $\Gamma$ satisfied by $u$.

An $A$-*formula* with parameters in $\Gamma$ is a formula of the form:

$$\psi_A(\alpha, x_A, (\rho_{\phi_i}(t_i(x_B) - \beta_i))_i)$$

with $\psi_A \in \mathcal{L}_A$ parameter-free, $(\phi_i)_i$ a finite family of elements of $\Phi$, $(t_i)_i$ a finite family of parameter-free $\mathcal{L}_B$-terms, $\alpha \in A^*(\Gamma)$, and $\beta_i \in B(\Gamma)$. We define $\mathrm{tp}_A(u/\Gamma)$ similarly.

Note that, given $t(x_1 \dots x_n, y_1 \dots y_m)$ a parameter-free $\mathcal{L}_B$-term, and given $\beta_1 \dots \beta_m \in B(\Gamma)$, $t(x_1 \dots x_n, \beta_1 \dots \beta_m)$ is equivalent to the term:

$$t(x_1 \dots x_n, 0 \dots 0) - t(0 \dots 0, -\beta_1 \dots - \beta_m)$$

thus it may be written just as in the definition of $A$-formulas, as a term of the form $s(x_1 \dots x_n) - \beta$, with $s$ parameter-free, and $\beta \in B(\Gamma)$.

We will be using the following quantifier elimination result:

**Theorem 5.1.3** ([ACGZ22], Corollary 4.20)**.** *If $\Gamma$ is generated by reals (i.e. it is a substructure generated by actual points from $A^*(M)$, $B(M)$ and $C^*(M)$), then we have:*

$$\mathrm{tp}(u/\Gamma) = \mathrm{tp}(v/\Gamma) \iff \begin{cases} \mathrm{tp}_A(u/\Gamma) = \mathrm{tp}_A(v/\Gamma) \\ \mathrm{tp}_C(u/\Gamma) = \mathrm{tp}_C(v/\Gamma) \end{cases}$$

We specify that $\Gamma$ is generated by reals because, in the next sections, we work with imaginary parameter sets.

In fact, they prove this result under weaker assumptions: they allow the expansion to contain functions and relations that involve both $A$ and $C$. However, under these assumptions, $A^*$ and $C^*$ may not be orthogonal, making it unreasonable (and arguably impossible at this level of generality) to try to obtain Ax-Kochen-Ershov principles. For the same reason, we do not work in the setting of weakly pure short exact sequences of Abelian structures (see their subsubsection 4.5.2), as $\nu \circ \iota$ would be a non-trivial map from $A^*$ to $C^*$. In our setting, $A^*$ and $C^*$ are orthogonal, and they are stably embedded in a strong sense:

**Corollary 5.1.4.** *If* $\Gamma$ *is generated by reals, then we have* $u_A u_C \equiv_\Gamma v_A v_C$ *in* $M$ *if and only if* $u_A \equiv_{A^*(\Gamma)} v_A$ *in the reduct* $A^*(M)$*, and* $u_C \equiv_{C^*(\Gamma)} v_C$ *in* $C^*(M)$*. In particular,* $u_A \equiv_\Gamma v_A$ *in* $M$ *if and only if* $u_A \equiv_{A^*(\Gamma)} v_A$ *in* $A^*(M)$*, and likewise for* $C$*.*

*Remark* 5.1.5. We should point out that $A$-types and $C$-types are not as independent from each other as one would expect, due to the technical definition of $A$-formulas. For instance, the $A$-formula $\rho_\phi(x) \neq 0$ must imply the $C$-formula $\phi(\nu(x))$. Because of that, one has to be careful when using orthogonality to study forking in this setting. We build in section 5.5 a type $\mathrm{tp}(u/\Delta)$ over a larger substructure $\Delta \supseteq \Gamma$ which divides over $\Gamma$, while $\mathrm{tp}_A(u/\Delta)$ and $\mathrm{tp}_C(u/\Delta)$ do not fork over $\Gamma$.

However, even though the Ax-Kochen-Ershov principle fails for forking and dividing in the full Stone space of types over $\Delta$, it holds in a subspace which only depends on the base parameter set $\Gamma$. More precisely, given $u$ and $\Delta$, we will define $p_0$ a partial type over $\Gamma$ realized by $u$ such that, in the Stone space of types over $\Delta$ containing $p_0$, $A$-types and $C$-types *are* independent from each other. We do manage to classify, in an Ax-Kochen-Ershov fashion, the types in this Stone space which do not fork/divide over $\Gamma$, in terms of the induced $A$-types and $C$-types. This yields in particular a classification of the space of non-forking extensions of $\mathrm{tp}(u/\Gamma)$.

# 5.2 Orthogonality

Just as in the previous section, we have $M$ an expansion of a $\Phi$-PSES:

$$0 \longrightarrow A^* \underset{\iota}{\longrightarrow} B \underset{\nu}{\longrightarrow} C^* \longrightarrow 0$$

with $\Phi$ a fundamental system for $B$, we have substructures $\Gamma \subseteq \Delta \subseteq M^{eq}$ (we allow imaginaries this time), and $u = u_A u_B u_C \in M$ a tuple. We defined $A$-types and $C$-types in the last section, and we noted that they are not completely independent from each other. We get around that problem in this section.

As we want to describe forking, we may assume that $M$ is $|\Delta|^+$-saturated and strongly $|\Delta|^+$-homogeneous. Likewise, we may freely assume that $\Gamma = \mathrm{dcl}^{eq}(\Gamma)$, $\Delta = \mathrm{dcl}^{eq}(\Delta)$. For technical reasons, we also make the following two mild hypothesis:

(H1) We assume that $\Gamma$ has enough algebraic imaginaries: we assume that for every $\phi \in \Phi$, each coset from ${}^{C}\!/\!_{\phi(C)}$ which belongs to $\mathrm{acl}^{eq}(\Gamma)$ belongs in fact to $\Gamma$.

(H2) We assume that $\Delta$ is generated by reals:

$$\Delta = \mathrm{dcl}^{eq}(A^*(\Delta) \cup B(\Delta) \cup C^*(\Delta))$$

**Notations** Our goal is to describe the Stone space of extensions in $S(\Delta)$ of $\mathrm{tp}(u/\Gamma)$ which do not fork/divide over $\Gamma$. To achieve this, we define various partial types, Stone spaces and definable functions. The partial types that we define typically state that $\nu(x_B)$ does or does not belong to various cosets of ${}^{C}\!/\!_{\phi(C)}$. One partial type, $p_0$, depends solely on $\mathrm{tp}(u/\Gamma)$, and plays an important role in our main results. We write with the subscript $_0$ the various objects (such as Stone spaces) which depend on $p_0$ in some way. Other objects will depend on the parameter set $\Delta$, we will typically write them with the subscript/superscript $_\Delta$.

**Definition 5.2.1.** Let $\mathcal{T}$ be the set of parameter-free $\mathcal{L}_B$-terms. We define the following sets:

$$E_0^+ = \left\{(\phi, t) \in \Phi \times \mathcal{T} | \nu(t(u_B)) + \phi(C) \in {}^{C}\!/\!_{\phi(C)}(\Gamma)\right\}$$

$$E_0^- = (\Phi \times \mathcal{T}) \setminus E_0^+$$

$$E_\Delta^+ = \{(\phi, t) \in E_0^+ | \exists \beta \in B(\Delta)\ \nu(t(u_B) - \beta) \in \phi(C)\}$$

For each $(\phi, t) \in E_\Delta^+$, choose $\beta_{(\phi,t)} \in B(\Delta)$ a witness of the last equation.

By definition of $E_0^+$ and $E_0^-$, the following partial type:

$$p_0(x_A, x_B, x_C) = \{\nu(t(x_B - u_B)) \in \phi(C)) | (\phi, t) \in E_0^+\}$$

$$\cup \left\{\nu(t(x_B)) + \phi(C) \neq X | (\phi, t) \in E_0^-, X \in {}^{C}\!/\!_{\phi(C)}(\Gamma)\right\}$$

can be written with formulas over $\Gamma$. We also define the following partial type over $\Delta$:

$$p_\Delta(x_A,x_B,x_C) = p_0(x_A,x_B,x_C) \cup \{\neg\phi(\nu(t(x_B)-\beta)) | (\phi,t) \in E_0^-,\ \beta \in B(\Delta)\}$$

Finally, we define the $\varnothing$-definable function:

$$f_C:\ (x_A,x_B,x_C) \longmapsto (\nu(x_B),x_C)$$

and we define the (infinite) tuple of $\Delta$-definable functions:

$$f_A^\Delta:\ (x_A,x_B,x_C) \longmapsto (x_A,(\rho_\phi(t(x_B)-\beta_{(\phi,t)}))_{(\phi,t)\in E_\Delta^+})$$

The graph of $f_A^\Delta$ is $\Delta$-$\star$-definable.

Our goal in this chapter is to classify the types over $\Delta$ which contain $p_0$ and do not fork/divide over $\Gamma$. We see $f_C$ and $f_A^\Delta$ as projections to $A^*$ and $C^*$. With respect to those projections, the sets of realisations of $p_0$ and $p_\Delta$ are preimages by $f_C$, i.e. if $u \vDash p_0$ (resp. $u \vDash p_\Delta$), and $f_C(u) = f_C(v)$, then it follows from the definitions that $v \vDash p_0$ (resp. $v \vDash p_\Delta$).

We call $f_A^\Delta$ and $f_C$ projections, but $f_A^\Delta \times f_C$ may not be onto the full direct product, i.e. its direct image may not be a rectangle, as remark 5.1.5 points out. The following statements shows that the direct images of $p_0$ and $p_\Delta$ by $f_A^\Delta \times f_C$ are, in fact, rectangles:

**Proposition 5.2.2.** *Let $v,v' \in M$ be realizations of $p_0$. Then there exists $w \in M$ such that $f_A^\Delta(v) = f_A^\Delta(w)$, and $f_C(v') = f_C(w)$ (in particular, $w \vDash p_0$).*

*Proof.* Let $(\phi_i,t_i)_i$ be a finite family from $E_\Delta^+$. We find $w \in M$ such that $f_C(w) = f_C(v')$, and $\rho_{\phi_i}(t_i(w_B)-\beta_{(\phi_i,t_i)}) = \rho_{\phi_i}(t_i(v_B)-\beta_{(\phi_i,t_i)})$ for every $i$, and the statement will follow by saturation.

Let $\psi(x)$ be the following p.p.-formula:

$$\bigwedge_i \phi_i(t_i(x))$$

then the map $x_B \longmapsto (\rho_{\phi_i}(t_i(x_B)-\beta_{(\phi_i,t_i)}))_i$ factors through the quotient $B/\psi(B)$, thus it suffices to find $w_B$ such that $\nu(w_B) = \nu(v'_B)$ (in other words, $w_B - v'_B \in \iota(A)$) and $w_B - v_B \in \psi(B)$, and set $w = (v_A, w_B, v'_C)$. Such a $w_B$ exists if and only if $v_B - v'_B \in \psi(B) + \iota(A)$, if and only if $\nu(v_B - v'_B) \in \psi(C)$, if and only if $\nu(t_i(v_B - v'_B)) \in \phi_i(C)$ for each $i$. Now, $v$ and $v'$ both realize $p_0$, thus we have $\nu(t_i(u_B - v_B)) \in \phi_i(C) \ni \nu(t_i(u_B - v'_B))$, concluding the proof. □

As the set of realizations of $p_\Delta$ is contained in that of $p_0$, and is also a preimage by $f_C$, it immediately follows that:

**Corollary 5.2.3.** *Proposition 5.2.2 also holds if we replace $p_0$ by $p_\Delta$.*

Note that proposition 5.2.2 is not stated in its most optimal form: it holds for the following weaker partial type:

$$\{\nu(t(x_B - u_B)) \in \phi(C)) | (\phi, t) \in E_0^+\}$$

*Remark* 5.2.4. Let $S_0^\Delta$ (resp. $S^\Delta$) be the Stone space of every type from $S^u(\Delta)$ which extends $p_0$ (resp. $p_\Delta$). Then $(f_A^\Delta, f_C)$ yield continuous maps $S_0^\Delta \longrightarrow f_A^\Delta(S_0^\Delta) \times f_C(S_0^\Delta)$, and $S^\Delta \longrightarrow f_A^\Delta(S^\Delta) \times f_C(S^\Delta)$. By corollary 5.1.4 and (H2), the factors of those direct products can be identified with closed subspaces of $S(A^*(\Delta))$ and $S(C^*(\Delta))$ in the first-order structures $A^*(M)$ and $C^*(M)$. Then, proposition 5.2.2 and corollary 5.2.3 state that those two continuous maps are *surjective*.

**Proposition 5.2.5.** *The map $S^\Delta \longrightarrow f_A^\Delta(S^\Delta) \times f_C(S^\Delta)$ is a homeomorphism.*

*Proof.* By compactness of the domain, separation of the image, and surjectivity, we just need to show injectivity. Let $v$, $w$ be realizations of $p_\Delta$ such that $f_A^\Delta(v) \equiv_{A^*(\Delta)} f_A^\Delta(w)$ and $f_C(v) \equiv_{C^*(\Delta)} f_C(w)$. We need to show that $v \equiv_\Delta w$. By corollary 5.1.4 and (H2), we have $f_A^\Delta(v) f_C(v) \equiv_\Delta f_A^\Delta(w) f_C(w)$. As a result, by strong homogeneity, we can freely assume that $f_A^\Delta(v) = f_A^\Delta(w)$ and $f_C(v) = f_C(w)$.

Now let us show $v \equiv_\Delta w$. Let $\Delta_r = A^*(\Delta) \cup B(\Delta) \cup C^*(\Delta)$. By (H2), it suffices to show that $v \equiv_{\Delta_r} w$. By theorem 5.1.3, it suffices to show that $\mathrm{tp}_A(v/\Delta_r) = \mathrm{tp}_A(w/\Delta_r)$ and $\mathrm{tp}_C(v/\Delta_r) = \mathrm{tp}_C(w/\Delta_r)$.

The fact that $\mathrm{tp}_C(v/\Delta_r) = \mathrm{tp}_C(w/\Delta_r)$ immediately follows from the fact that $f_C(v) = f_C(w)$.

Let us show $\mathrm{tp}_A(v/\Delta_r) = \mathrm{tp}_A(w/\Delta_r)$. Let $t \in \mathcal{T}$, $\phi \in \Phi$, and $\beta \in B(\Delta_r) = B(\Delta)$. It suffices to show that $\rho_\phi(t(v_B) - \beta) = \rho_\phi(t(w_B) - \beta)$.

- Suppose $(\phi, t) \in E_\Delta^+$, and $\beta - \beta_{(\phi,t)} \in \phi(B) + \iota(A)$. Then:

$$\begin{array}{rcl} \rho_\phi(t(v_B) - \beta) & = & \rho_\phi(t(v_B) - \beta_{(\phi,t)}) + \rho_\phi(\beta_{(\phi,t)} - \beta) \\ & = & \rho_\phi(t(w_B) - \beta_{(\phi,t)}) + \rho_\phi(\beta_{(\phi,t)} - \beta) \\ & = & \rho_\phi(t(w_B) - \beta) \end{array}$$

  The second equality follows from the fact that $f_A^\Delta(v) = f_A^\Delta(w)$, the other equalities follow from the definition of $\rho_\phi$.

- Suppose $(\phi, t) \in E^+_\Delta$, and $\beta - \beta_{(\phi,t)} \notin \phi(B) + \iota(A)$. Then:

$$t(v_B) - \beta \notin \phi(B) + \iota(A) \not\ni t(w_B) - \beta$$

  thus $\rho_\phi(t(v_B) - \beta) = 0 = \rho_\phi(t(w_B) - \beta)$.

- Suppose $(\phi, t) \in E^+_0 \setminus E^+_\Delta$. By definition of $p_0$, we have:

$$t(v_B - u_B) \in \phi(B) + \iota(A) \ni t(w_B - u_B)$$

  however, by definition of $E^+_\Delta$, we have $t(u_B) - \beta \notin \phi(B) + \iota(A)$, thus we conclude just as in the above item that $\rho_\phi(t(v_B) - \beta) = 0 = \rho_\phi(t(w_B) - \beta)$.

- Suppose $(\phi, t) \in E^-_0$. Then, by definition of $p_\Delta$, we have:

$$t(v_B) - \beta \notin \phi(B) + \iota(A) \not\ni t(w_B) - \beta$$

  just as in the above item.

This conclude the proof, since all cases were considered. □

*Remark* 5.2.6. Let $v$, $w$ be realizations of $p_\Delta$. By definition, $\mathrm{tp}_C(v/\Delta) = \mathrm{tp}_C(w/\Delta)$ if and only if $f_C(v) \equiv_\Delta f_C(w)$, and the above proof shows us that the (non-trivial) same equivalence holds for $A$: $\mathrm{tp}_A(v/\Delta) = \mathrm{tp}_A(w/\Delta)$ if and only if $f^\Delta_A(v) \equiv_\Delta f^\Delta_A(w)$. So we have a homeomorphism between the space of $A$-types of realizations of $p_\Delta$, and $f^\Delta_A(S^\Delta)$. In particular, $S^\Delta$ is homeomorphic to the direct product of the respective spaces of $A$-types and $C$-types of realizations of $p_\Delta$.

## 5.3 Dividing

We defined in the previous section a space of types $S^\Delta_0$ which is closed in $S(\Delta)$. We classify here the subspace of $S^\Delta_0$ of types which do not divide over $\Gamma$ in terms of the projection maps $f^\Delta_A$, $f_C$. With respect to those maps, we saw that $S^\Delta_0$ maps onto a rectangle, while $S^\Delta$ *is* a rectangle. Thus the elements of $S^\Delta$ are in some sense easier to classify than in $S^\Delta_0$. Thankfully, the space that we are trying to describe is a subspace of $S^\Delta$:

**Proposition 5.3.1.** *Every type in $S^\Delta_0 \setminus S^\Delta$ divides over $\Gamma$.*

*Proof.* Let $p$ be such a type. By definition, there must exist $\beta \in B(\Delta)$, and $(\phi, t) \in E_0^-$ such that $p(x) \models \phi(\nu(t(x_B) - \beta))$. Let us show that this formula divides over $\Gamma$. Let $Y$ be the coset $\nu(\beta) + \phi(C)$. Then we have $p(x) \models \nu((t(x_B))) + \phi(C) = Y$. As $p$ extends $p_0$, we have $Y \notin \Gamma$. By (H1), we have $Y \notin \operatorname{acl}^{eq}(\Gamma)$. Let $(Y_n)_{n<\omega}$ be a sequence of pairwise-distinct $\Gamma$-conjugates of $Y$. Then the definable sets $(\nu(t(x_B)) + \phi(C) = Y_n)_n$ are $\Gamma$-conjugates and pairwise-disjoint. This implies that they divide over $\Gamma$, which concludes the proof. □

**Corollary 5.3.2.** *The space of $A$-types (resp. $C$-types) over $\Delta$ of realizations $u$ of $p_0$ such that $u \mathop{\downarrow}^{\mathbf{d}}_{\Gamma} \Delta$ is naturally homeomorphic to the direct image by $f_A^{\Delta}$ (resp. $f_C$) of the closed subspace of $S_0^{\Delta}$ of extensions of $p_0$ to $\Delta$ which do not divide over $\Gamma$. The same holds if we replace dividing by forking.*

*Proof.* As the spaces of types at hand are included in $S^{\Delta}$, we apply remark 5.2.6. □

Now, one wonders what does the space of types which do not divide over $\Gamma$ look like. It turns out that it is a *subrectangle* of $S^{\Delta}$:

**Theorem 5.3.3.** *Let $v, v' \in M$ be realizations of $p_0$ such that $v \mathop{\downarrow}^{\mathbf{d}}_{\Gamma} \Delta$ and $v' \mathop{\downarrow}^{\mathbf{d}}_{\Gamma} \Delta$. Let $w \in M$ be such that $f_A^{\Delta}(w) = f_A^{\Delta}(v)$, $f_C(w) = f_C(v')$. Then we have $w \mathop{\downarrow}^{\mathbf{d}}_{\Gamma} \Delta$.*

*Proof.* Let $\phi(x)$ be a formula with parameters in $\Delta$ which is realized by $w$. By proposition 5.2.5, corollary 5.1.4, and the fact that $p_0$ and $p_{\Delta}$ are preimages by $f_C$, the partial type:

$$p_0 \cup \left\{ f_A^{\Delta}(x) \models \operatorname{tp}(f_A^{\Delta}(w)/A^*(\Delta)) \right\} \cup \{ f_C(x) \models \operatorname{tp}(f_C(w)/C^*(\Delta)) \}$$

generates the complete type $\operatorname{tp}(w/\Delta)$. By compactness, there must exist:

- A formula $\psi(x)$ with parameters in $\Gamma$ satisfied by all the realizations of $p_0$.
- A formula $\phi_A(y)$ in $\mathcal{L}_A$, and with parameters in $A^*(\Delta)$, satisfied by $f_A^{\Delta}(w)$.
- A formula $\phi_C(z)$ in $\mathcal{L}_C$, and with parameters in $C^*(\Delta)$, satisfied by $f_C(w)$.

such that:

$$M^{eq} \vDash \forall x \ \left[\left(\psi(x) \wedge \phi_A\left(f_A^{\Delta}(x)\right) \wedge \phi_C(f_C(x))\right) \Longrightarrow \phi(x)\right]$$

Note that only finitely many of the components of $f_A^{\Delta}$ appear in the formula $\phi_A(f_A^{\Delta}(x))$. Let $\delta$ be a finite tuple from $\Delta$ such that the formula $\phi_A(f_A^{\Delta}(x)) \wedge \phi_C(f_C(x))$ is $\delta$-definable, say, $\phi_A(f_A^{\Delta}(x)) = \theta_A(x,\delta)$ and $\phi_C(f_C(x)) = \theta_C(x,\delta)$, with $\theta_A(x,y)$ and $\theta_C(x,z)$ parameter-free formulas. It suffices to show that the formula $\psi(x) \wedge \theta_A(x,\delta) \wedge \theta_C(x,\delta)$ does not divide over $\Gamma$.

Let $I = (\delta_n)_{n<\omega}$ be a $\Gamma$-indiscernible sequence containing $\delta$. As $v \mathop{\downarrow}^{\mathbf{d}}_{\Gamma} \Delta$, there must exist by fact 1.1.5 some $\overline{v} \equiv_{\Delta} v$ such that $I$ is $\Gamma\overline{v}$-indiscernible. Likewise, let $\overline{v'} \equiv_{\Delta} v'$ be such that $I$ is $\Gamma\overline{v'}$-indiscernible. Then $\overline{v}$ and $\overline{v'}$ are still realizations of $p_0$, thus there exists, by proposition 5.2.2 a tuple $\overline{w} \in M$ such that $f_A^{\Delta}(\overline{w}) = f_A^{\Delta}(\overline{v})$ and $f_C(\overline{w}) = f_C(\overline{v'})$ (in particular, $\overline{w} \vDash p_0$, thus $M^{eq} \vDash \psi(\overline{w})$). Let us show that $M \vDash \theta_A(\overline{w},\delta_n) \wedge \theta_C(\overline{w},\delta_n)$ for all $n < \omega$, which will conclude the proof. Choose $N < \omega$.

As $\overline{v'} \equiv_{\Delta} v'$, we have $M \vDash \theta_C(\overline{v'},\delta)$. As $(\delta_n)_n$ is a $\Gamma\overline{v'}$-indiscernible sequence containing $\delta$ and $\delta_N$, we have $M \vDash \theta_C(\overline{v'},\delta_N)$. As $f_C(\overline{w}) = f_C(\overline{v'})$, we have $M \vDash \theta_C(\overline{w},\delta_N)$.

Likewise, we have $M \vDash \theta_A(\overline{v},\delta_N)$, but there is a subtlety to get around if we want to prove that $M \vDash \theta_A(\overline{w},\delta_N)$: since $f_A^{\Delta}$ depends on $\Delta$, when replacing $\delta$ with its $\Gamma$-conjugate $\delta_N$, we also replace $f_A^{\Delta}$ with a $\Gamma$-conjugate... By strong homogeneity, let $\sigma \in \mathrm{Aut}(M^{eq}/\Gamma)$ be such that $\sigma(\delta) = \delta_N$. We have to show that $\sigma\left(f_A^{\Delta}\right)(\overline{v}) = \sigma\left(f_A^{\Delta}\right)(\overline{w})$.

Since the definition of $f_A^{\Delta}$ is rather technical, we recall the various definitions:

- $E_0^+ = \left\{(\phi,t) \in \Phi \times \mathcal{T} | \nu(t(u_B)) + \phi(C) \in {}^{C}\!/\!_{\phi(C)}(\Gamma)\right\}$
- $E_{\Delta}^+ = \{(\phi,t) \in E_0^+ | \exists \beta \in B(\Delta) \ \nu(t(u_B) - \beta) \in \phi(C)\}$
- For each $(\phi,t) \in E_{\Delta}^+$, $\beta_{(\phi,t)} \in B(\Delta)$ is such that $\nu(t(u_B) - \beta_{(\phi,t)}) \in \phi(C)$.
- $f_A^{\Delta}: \ (x_A,x_B,x_C) = (x_A, (\rho_\phi(t(x_B) - \beta_{(\phi,t)}))_{(\phi,t)\in E_{\Delta}^+})$.

Choose $(\theta,t) \in E_{\Delta}^+$, and let us show that:

$$\rho_\theta(t(v_B) - \sigma(\beta_{(\theta,t)})) = \rho_\theta(t(w_B) - \sigma(\beta_{(\theta,t)}))$$

By definition of $E_\Delta^+$ and $E_0^+$, the coset $\nu(\beta_{(\theta,t)}) + \theta(C)$ belongs to $\Gamma$, thus $\nu(\sigma(\beta_{(\theta,t)}))$ belongs to that coset as well. It follows that $\beta_{(\theta,t)} - \sigma(\beta_{(\theta,t)})$ is in $\theta(B) + \iota(A)$. As $\overline{v}$ and $\overline{w}$ both realize $p_0$, their image by $\nu$ belongs to that coset as well, and we conclude that:

$$\begin{aligned}\rho_\theta(\overline{v}_B - \sigma(\beta_{(\theta,t)})) &= \rho_\theta(\overline{v}_B - \beta_{(\theta,t)}) + \rho_\theta(\beta_{(\theta,t)} - \sigma(\beta_{(\theta,t)})) \\ &= \rho_\theta(\overline{w}_B - \beta_{(\theta,t)}) + \rho_\theta(\beta_{(\theta,t)} - \sigma(\beta_{(\theta,t)})) \\ &= \rho_\theta(\overline{w}_B - \sigma(\beta_{(\theta,t)}))\end{aligned}$$

The second equality follows from the fact that $f_A^\Delta(\overline{v}) = f_A^\Delta(\overline{w})$, and the other equalities follow from the fact that $\rho_\theta$ is always applied here to elements of $\theta(B) + \iota(A)$, over which it restricts to a group homomorphism. □

**Corollary 5.3.4.** *Let $S_{\mathbf{d}}^\Delta$ be the space of every type in $S_0^\Delta$ which do not divide over $\Gamma$. Then the map $S_{\mathbf{d}}^\Delta \longrightarrow f_A^\Delta(S_{\mathbf{d}}^\Delta) \times f_C(S_{\mathbf{d}}^\Delta)$ is a homeomorphism.*

*Proof.* This map is clearly continuous. It is surjective by theorem 5.3.3. By proposition 5.3.1, we have $S_{\mathbf{d}}^\Delta \subseteq S^\Delta$, therefore this map is injective as a restriction of an injective map, by proposition 5.2.5. We conclude by compactness and separation. □

**Corollary 5.3.5.** *Let $\mathcal{A}_{\mathbf{d}}^\Delta$ (resp. $\mathcal{C}_{\mathbf{d}}^\Delta$) be the space of $A$-types (resp. $C$-types) of realizations $u$ of $p_0$ such that $u \mathop{\downarrow}\limits_{\Gamma}^{\mathbf{d}} \Delta$. Then $S_{\mathbf{d}}^\Delta$ is naturally homeomorphic to $\mathcal{A}_{\mathbf{d}}^\Delta \times \mathcal{C}_{\mathbf{d}}^\Delta$.*

*Proof.* This follows immediately from corollary 5.3.2. □

It follows that $u \mathop{\downarrow}\limits_{\Gamma}^{\mathbf{d}} \Delta$ if and only if each of the partial types $\mathrm{tp}_C(u/\Delta)$ and $p_0 \cup \mathrm{tp}_A(u/\Delta)$ have a realization $v$ such that $v \mathop{\downarrow}\limits_{\Gamma}^{\mathbf{d}} \Delta$. Note that this statement is in general not equivalent to: $\mathrm{tp}_C(u/\Delta)$ and $p_0 \cup \mathrm{tp}_A(u/\Delta)$ do not divide over $\Gamma$. This is related to the issue of forking being potentially stronger than dividing, or equivalently, that it is possible for a formula not to divide while containing only types that divide. More precisely, we should remind the reader that:

*Remark* 5.3.6. In any first-order theory, given a set of parameters $\Gamma$, the following are equivalent:

- For every set of parameters $\Delta$, for every partial type $p$ over $\Gamma\Delta$, $p$ does not divide over $\Gamma$ if and only if $p$ admits a realization $v$ such that $v \mathop{\downarrow}\limits_{\Gamma}^{\mathbf{d}} \Gamma\Delta$.

- For every set of parameters $\Delta$, for every $\Gamma\Delta$-definable set $X$, $X$ does not divide over $\Delta$ if and only if $X$ admits a realization $v$ such that $v \underset{\Gamma}{\overset{\mathbf{d}}{\perp\!\!\!\perp}} \Gamma\Delta$.

- Forking coincides with dividing over $\Gamma$, i.e. every definable set which forks over $\Gamma$ divides over $\Gamma$.

A common example of a definable set which does not divide while all its realizations divide is in the well-known cyclical order $(\mathbb{Q}, \mathrm{cyc})$, with $\Gamma = \varnothing$, and $\Delta = \{*\}$ a singleton. In this setting, the empty partial type does not divide over $\Gamma$ (because it never does), but every type over $\Delta$ divides over $\Gamma$, thus all the realizations $v$ of the empty partial type satisfy $v \underset{\Gamma}{\overset{\mathbf{d}}{\not\perp\!\!\!\perp}} \Delta$.

Nonetheless, it is true in any theory that a partial type $p$ over $\Gamma\Delta$ does not fork over $\Gamma$ if and only if there exists $v$ a realization of $p$ such that $v \underset{\Gamma}{\overset{\mathbf{f}}{\perp\!\!\!\perp}} \Delta$. As a consequence, the non-trivial result that we can get from this section is:

**Corollary 5.3.7.** *If both partial types $p_0 \cup \mathrm{tp}_A(u/\Delta)$ and $\mathrm{tp}_C(u/\Delta)$ do not fork over $\Gamma$, then we have $u \underset{\Gamma}{\overset{\mathbf{d}}{\perp\!\!\!\perp}} \Delta$.*

And the converse is trivial whenever forking coincides with dividing over $\Gamma$. In fact, one can see that the converse holds whenever $- \underset{\Gamma}{\overset{\mathbf{f}}{\perp\!\!\!\perp}} \Delta = - \underset{\Gamma}{\overset{\mathbf{d}}{\perp\!\!\!\perp}} \Delta$. In particular, it holds when $\Delta$ is a $|\Gamma|^+$-saturated model. Let us state it formally, as it will come in handy in the next section:

**Corollary 5.3.8.** *Suppose $\Delta$ is a $|\Gamma|^+$-saturated model. Then we have $u \underset{\Gamma}{\overset{\mathbf{f}}{\perp\!\!\!\perp}} \Delta$ if and only if both partial types $p_0 \cup \mathrm{tp}_A(u/\Delta)$ and $\mathrm{tp}_C(u/\Delta)$ do not fork over $\Gamma$.*

Now, corollary 5.3.7 also suggests that corollary 5.3.7 holds for any $\Delta$. This is what we prove in the next section.

## 5.4 Forking

Recall that a partial type does not fork over $\Gamma$ if and only if it admits a global completion which does not divide over $\Gamma$. Note that the monster model $M^{eq}$ satisfies the hypothesis (H2), thus all the results that we obtained for $\Delta$ also hold for $M^{eq}$, by switching to some big elementary extension of $M$. By corollary 1.1.17, forking for types over $\Delta$ involves dividing for types over $M^{eq}$, thus the results of the previous sections will help us to establish here a classification of the space of every type in $S_0^{\Delta}$ which do not fork over $\Gamma$.

**Definition 5.4.1.** Let $\mathcal{B}_0^{\Delta}$ be the Boolean algebra associated to $S_0^{\Delta}$. Let $e^{\Delta}$ be the embedding $\mathcal{B}_0^{\Delta} \longrightarrow \mathcal{B}_0^{M^{eq}}$. Define $F_A^{\Delta}$ (resp. $F_C^{\Delta}$) as the set of every element of $\mathcal{B}_0^{\Delta}$ which is equivalent to some formula of the form $\phi(f_A^{\Delta}(x))$ (resp. $\phi(f_C(x))$), with $\phi \in \mathcal{L}_A$ (resp. $\mathcal{L}_C$) with parameters in $\Delta$.

There are several ways to formalize the Stone duality. We identify by convention an element of $\mathcal{B}_0^{\Delta}$ with a (clopen) subset of $S_0^{\Delta}$. With this formalism, given $I$ an ideal of $\mathcal{B}_0^{\Delta}$ (for instance: the ideal of definable sets $X$ such that $p_0 \cup \{X\}$ forks over $\Gamma$), a type $p \in S_0^{\Delta}$ is inconsistent with $I$ (i.e. $p$ does not fork over $\Gamma$) if and only if $p \notin \cup I$.

Recall that in a Boolean algebra $\mathcal{B}$, a subset is an ideal (as in ring theory) if and only if it is non-empty, downward-closed and closed under finite join. In other words, the ideal generated by some set $P$ coincides with:

$$\{X \in \mathcal{B} | \exists n < \omega \ \exists Y_1 \dots Y_n \in P \ X \leqslant (Y_1 \vee \dots \vee Y_n)\}$$

moreover, if $P$ is already downward-closed, then this set coincides with:

$$\{X \in \mathcal{B} | \exists n < \omega \ \exists Y_1 \dots Y_n \in P \ X = (Y_1 \vee \dots \vee Y_n)\}$$

in particular, if $I$, $J$ are ideals of $\mathcal{B}$, then:

$$I + J = \{X \vee Y | X \in I, Y \in J\}$$

therefore in our setting:

$$\cup(I + J) = (\cup I) \cup (\cup J)$$

**Lemma 5.4.2.** *Let $I_A$ (resp. $I_C$) be an ideal of $\mathcal{B}_0^{M^{eq}}$ generated by some subset of $F_A^{M^{eq}}$ (resp. $F_C^{M^{eq}}$). Let $X_A \in F_A^{\Delta}$, and $X_C \in F_C^{\Delta}$. Then we have $e(X_A \wedge X_C) \in I_A + I_C$ if and only if $e(X_A) \in I_A$ or $e(X_C) \in I_C$.*

*Proof.* The right-to-left direction is trivial.

Let us show the other direction by contraposition. Suppose $e(X_A) \notin I_A$ and $e(X_C) \notin I_C$. Then there exists $p \in e(X_A)$, $q \in e(X_C)$ such that $p \notin \cup I_A$, $q \notin \cup I_C$. By proposition 5.2.2, let $r \in S_0^{M^{eq}}$ be such that $f_A^{M^{eq}}(r) = f_A^{M^{eq}}(p)$, $f_C(r) = f_C(q)$. As $f_A^{\Delta}$ is a restriction of $f_A^{M^{eq}}$, we have $f_A^{\Delta}(r) = f_A^{\Delta}(p)$. As $p \in e(X_A)$, and $X_A \in F_A^{\Delta}$, we have $r \in e(X_A)$. Likewise, we have $r \in e(X_C)$, thus $r \in e(X_A \wedge X_C)$.

Now, it suffices to show that $r \notin \cup(I_A + I_C)$. Note that $\cup(I_A + I_C) = (\cup I_A) \cup (\cup I_C)$, thus it suffices to show that $r \notin \cup I_A$ and $r \notin \cup I_C$.

By hypothesis on $I_A$, we have $\cup I_A = \cup(I_A \cap F_A^{M^{eq}})$. As $f_A^{M^{eq}}(r) = f_A^{M^{eq}}(p)$, for every $X \in F_A^{M^{eq}}$, we have $r \in X$ if and only if $p \in X$. As $p \notin \cup I_A$, we have $r \notin \cup(I_A \cap F_A^{M^{eq}})$, thus $r \notin \cup I_A$.

Likewise, $r \notin \cup I_C$, concluding the proof. □

Note that this kind of reasoning is very general, it may be used in more abstract settings which are not limited to expansions of PSES. The key property that we use to make this lemma work is that the Stone space at hand maps onto a rectangle by our two projection maps $f_A$, $f_C$. However, this lemma does not apply to every definable set, it only applies to *actual* rectangles, i.e. definable sets of the form $X_A \wedge X_C$ as in the statement. This is where the rectangle $S^{\Delta}$ comes handy, because every element of the Boolean algebra corresponding to this Stone space may be written as a finite join of rectangles. Then, with proposition 5.3.1, and the fact that $S^{\Delta}$ may be written as the intersection of a subset of $F_C^{\Delta}$, we have all the tools that we need to obtain our main result about forking. Let us formalize all this:

**Theorem 5.4.3.** *Let $p_A$ (resp. $p_C$) be the partial type generated by all formulas in $F_A^{\Delta}$ (resp. $F_C^{\Delta}$) satisfied by $u$. Then we have $u \underset{\Gamma}{\overset{\mathbf{f}}{\downarrow}} \Delta$ if and only if both partial types $p_0 \cup p_A$ and $p_C$ do not fork over $\Gamma$.*

*Proof.* The left-to-right direction is trivial. Let us show the other direction, so suppose that $p_0 \cup p_A$ and $p_C$ do not fork over $\Gamma$.

Let $I_A$ (resp. $I_C$) be the ideal of $\mathcal{B}_0^{M^{eq}}$ generated by formulas $\phi$ of $F_A^{M^{eq}}$ (resp. $F_C^{M^{eq}}$) such that $p_0 \cup \{\phi\}$ (resp. $\phi$) forks over $\Gamma$. By corollary 5.3.8, the (open) space of every type in $S_0^{M^{eq}}$ which forks over $\Gamma$ is exactly $\cup(I_A + I_C)$. It follows that the space of every type in $S_0^{\Delta}$ which forks over $\Gamma$ is $\cup e^{-1}(I_A + I_C)$.

If we had $u \nvDash p_{\Delta}$, then no realization of $p_C$ would realize $p_{\Delta}$, therefore, by 5.3.1, we would have $v \underset{\Gamma}{\overset{\mathbf{d}}{\not\downarrow}} \Delta$ for every $v \vDash p_C$, contradicting our hypothesis that $p_C$ does not fork over $\Gamma$. So $u$ realizes $p_{\Delta}$.

Let $X$ be a $\Delta$-definable set containing $u$. As $u$ realizes $p_{\Delta}$, by proposition 5.2.5 and compactness, there exists $X_A \in F_A^{\Delta}$, $X_C \in F_C^{\Delta}$, and $Y$ a $\Delta$-definable set containing all the realizations of $p_{\Delta}$, such that, in $\mathcal{B}_0^{\Delta}$, we have $(Y \wedge X_A \wedge X_C) \leqslant X$. As $p_\Delta$ is a preimage by $f_C$, there must exist $Z \in F_C^\Delta$, containing $u$, and such that $Z \leqslant Y$. We have $(X_C \wedge Z) \in F_C^\Delta$, thus:

$$p_C \vDash x \in X_C \wedge x \in Z$$

As $p_C$ does not fork over $\Gamma$, we have $(X_C \wedge Z) \notin e^{-1}(I_A + I_C)$, in particular $e(X_C \wedge Z) \notin I_C$. Likewise, we have $e(X_A) \notin I_A$, therefore by lemma 5.4.2 we have $e(X_A \wedge X_C \wedge Z) \notin I_A + I_C$. In particular, this $\Delta$-definable set does not fork over $\Gamma$, and neither does $X$. This holds for all $X \in \mathcal{B}_0^\Delta$ containing $u$, so we conclude that $u \mathop{\downarrow}_\Gamma^{\mathbf{f}} \Delta$. □

**Corollary 5.4.4.** *Let $S_{\mathbf{f}}^\Delta$ be the space of all types in $S_0^\Delta$ which do not fork over $\Gamma$. Then the map $S_{\mathbf{f}}^\Delta \longrightarrow f_A^\Delta(S_{\mathbf{f}}^\Delta) \times f_C(S_{\mathbf{f}}^\Delta)$ is a homeomorphism.*

*Proof.* Just as in corollary 5.3.4, this map is a continuous injection between compact separated spaces. There remains to prove surjectivity.

Let $(q_A(y), q_C(y)) \in f_A^\Delta(S_{\mathbf{f}}^\Delta) \times f_C(S_{\mathbf{f}}^\Delta)$. By proposition 5.2.2, let $q \in S_0^\Delta$ be such that $q(x) \vDash q_A(f_A^\Delta(x))$ and $q(x) \vDash q_C(f_C(x))$, and let $v$ realize $q$. As $(q_A(y), q_C(y)) \in f_A^\Delta(S_{\mathbf{f}}^\Delta) \times f_C(S_{\mathbf{f}}^\Delta)$, each of the partial types $p_0 \cup q_A(f_A^\Delta(x))$ and $q_C(f_C(x))$ admits a realization $w$ such that $w \mathop{\downarrow}_\Gamma^{\mathbf{f}} \Delta$. In particular, those partial types do not fork over $\Gamma$. By theorem 5.4.3, we have $v \mathop{\downarrow}_\Gamma^{\mathbf{f}} \Delta$, therefore $q$ is a preimage of $(q_A(y), q_C(y))$, concluding the proof. □

**Corollary 5.4.5.** *Let $\mathcal{A}_{\mathbf{f}}^\Delta$ (resp. $\mathcal{C}_{\mathbf{f}}^\Delta$) be the space of $A$-types (resp. $C$-types) of realizations $u$ of $p_0$ such that $u \mathop{\downarrow}_\Gamma^{\mathbf{f}} \Delta$. Then $S_{\mathbf{f}}^\Delta$ is naturally homeomorphic to $\mathcal{A}_{\mathbf{f}}^\Delta \times \mathcal{C}_{\mathbf{f}}^\Delta$.*

*Proof.* This follows immediately from corollary 5.3.2. □

**Corollary 5.4.6.** *We have $u \mathop{\downarrow}_\Gamma^{\mathbf{f}} \Delta$ if and only if both partial types $\mathrm{tp}_C(u/\Delta)$ and $p_0 \cup \mathrm{tp}_A(u/\Delta)$ do not fork over $\Gamma$.*

## 5.5 A silly example

In this section, we show a very basic example of a PSES of Abelian groups $M$ (without any expansion), with substructures $\Gamma \leqslant \Delta$ satisfying (H1) and (H2), and $u \in M$ such that $u \not\mathop{\downarrow}_\Gamma^{\mathbf{d}} \Delta$, but both partial types $\mathrm{tp}_A(u/\Delta)$ and $\mathrm{tp}_C(u/\Delta)$ do not fork over $\Gamma$. This illustrates why the partial type $p_0$ is absolutely necessary in our main result theorem 5.4.3.

Let $M$ be the following $\Phi$-PSES of Abelian groups:

$$\mathbb{Q} \times \{0\} \underset{\iota}{\longrightarrow} \mathbb{Q} \times \mathbb{Q} \underset{\nu}{\longrightarrow} \{0\} \times \mathbb{Q}$$

with $\Phi = \{x = 0\}$, which is indeed a fundamental system in torsion-free divisible Abelian groups. In fact, we should remind the reader that in the theory of such groups, any formula with one variable $x$ and parameters in some set $P$ is equivalent to a Boolean combination of formulas of the form $x - e = 0$, with $e$ lying in the $\mathbb{Q}$-vector space generated by $P$.

Let $\Gamma = \mathrm{dcl}^{eq}(\{\nu(0,1)\})$, $\Delta = \mathrm{dcl}^{eq}((1,1)) \supseteq \Gamma$, and let $u = (1,1)$. Let $N$ be some $\aleph_1$-saturated, strongly $\aleph_1$-homogeneous elementary extension of $M$.

By quantifier elimination for torsion-free divisible Abelian groups, we have $C(\mathrm{acl}^{eq}(\Gamma)) = \mathbb{Q} \cdot \nu(0,1) = C(\Gamma)$, thus (H1) holds, and (H2) holds by definition.

One can check that $u \notin \mathrm{acl}^{eq}(\Gamma)$, thus $u \not\downarrow^{\mathbf{d}}_{\Gamma} \Delta$.

The $C$-type of $u$ over $\Delta$ is isolated by a formula with parameters in $\Gamma$: $\nu(x) = \nu(0,1)$. Let $v \in B(N)$ be such that $v \in \iota(A)$, but $v \neq w$ for every $w \in M$. By theorem 5.1.3, the type of $v + (0,1)$ over $M$ is generated by the following partial type:

$$\begin{gathered}\{\nu(x) = \nu(0,1)\} \\ \cup \{\rho_{x=0}(x - \beta) \neq 0 | \beta \in B(M), \nu(\beta) = \nu(0,1)\} \\ \cup \{\rho_{x=0}(x - \beta) = 0 | \beta \in B(M), \nu(\beta) \neq \nu(0,1)\}\end{gathered}$$

As this partial type is $\mathrm{Aut}(M^{eq}/\Gamma)$-invariant, $\mathrm{tp}(v + (0,1)/M)$ is also $\mathrm{Aut}(M^{eq}/\Gamma)$-invariant. In particular, this type does not fork over $\Gamma$. As $v + (0,1)$ realizes $\mathrm{tp}_C(u/\Delta)$, we can deduce that $\mathrm{tp}_C(u/\Delta)$ does not fork over $\Gamma$.

Using theorem 5.1.3, one may show that for all $\lambda, \mu \in \mathbb{Q}$, we have $(\lambda, \mu) \equiv_\Gamma (2\lambda, 2\mu)$, thus $B(\Gamma) = \{\iota(0,0)\}$. In particular, as $\rho_{x=0}(\lambda \cdot u - \iota(0,0)) = \rho_{x=0}(\lambda, \lambda) = (0,0)$ for every $\lambda \in \mathbb{Q}$, we have $\mathrm{tp}_A(u/\Delta) = \mathrm{tp}_A(\iota(0,0)/\Delta)$. As $\iota(0,0) \in \Gamma$, we have of course $\iota(0,0) \downarrow^{\mathbf{f}}_{\Gamma} \Delta$, thus $\mathrm{tp}_A(u/\Delta)$ does not fork over $\Gamma$, which concludes the example.

## 5.6 Models

In theorem 5.4.3, we reduce the problem of characterizing forking in general to the problem of characterizing forking for certain partial types which encode data related to $A$ and $C$. However, the partial types at hand are still types in the full PSES, it would be better if we could reduce the problem to checking forking in the reducts $A^*$ and $C^*$. While this is probably not possible in full generality, we isolate conditions on $\Gamma$ which make it possible. Those conditions hold in particular when $\Gamma$ is a model.

**Definition 5.6.1.** In this section, we assume that for every $\phi \in \Phi$, for every coset $X \in {}^{C}\!/_{\phi(C)}(\Gamma)$, there exists $\beta \in B(\Gamma)$ such that $\nu(\beta) \in X$. We also assume that $\Gamma$ is generated by reals, i.e. $\Gamma = \mathrm{dcl}^{eq}(A^*(\Gamma) \cup B(\Gamma) \cup C^*(\Gamma))$.

*Remark* 5.6.2. With that assumption, we have $E^+_\Delta = E^+_0$, and the $\beta_{(\phi,t)}$ can be (and are, from now on) chosen in $B(\Gamma)$ for any $(\phi, t) \in E^+_\Delta$. In particular, $f^\Delta_A$ is $\Gamma$-$*$-definable and does not depend on $\Delta$, so we might as well call it $f_A$ in the rest of this section.

**Lemma 5.6.3.** *We have* $A^*(\mathrm{dcl}^{eq}(\Gamma u)) = A^*(\mathrm{dcl}^{eq}(A^*(\Gamma)f_A(u)))$, *as well as* $C^*(\mathrm{dcl}^{eq}(\Gamma u)) = C^*(\mathrm{dcl}^{eq}(C^*(\Gamma)f_C(u)))$.

*Proof.* The part of the statement about $C$ follows immediately from corollary 5.1.4, and the fact that $\Gamma$ is generated by reals.

As for the part about $A$, by theorem 5.1.3, it suffices to show that $\rho_\phi(t(u_B) - \beta) \in A^*(\mathrm{dcl}^{eq}(A^*(\Gamma)f_A(u)))$ for every term $t$, $\phi \in \Phi$, $\beta \in B(\Gamma)$. We follow the same case disjunction as in the proof of proposition 5.2.5:

- If $(\phi, t) \in E^+_0$, and $\beta - \beta_{(\phi,t)} \in \phi(B) + \iota(A)$, then:

$$\rho_\phi(t(u_B) - \beta) = \rho_\phi(t(u_B) - \beta_{(\phi,t)}) + \rho_\phi(\beta_{(\phi,t)} - \beta)$$

  We conclude as $\rho_\phi(t(u_B) - \beta_{(\phi,t)})$ is a projection of $f_A(u)$, and as $\rho_\phi(\beta_{(\phi,t)} - \beta)$ is in $A^*(\Gamma)$.

- If $(\phi, t) \in E^+_0$, and $\beta - \beta_{(\phi,t)} \notin \phi(B) + \iota(A)$, then $\rho_\phi(t(u_B) - \beta) = 0$.

- If $(\phi, t) \in E^-_0$, then $\rho_\phi(t(u_B) - \beta) = 0$.

this concludes the proof. ∎

**Theorem 5.6.4.** *The following are equivalent:*

- $u \underset{\Gamma}{\downfree}^{\mathbf{d}} \Delta$.
- $f_A(u) \underset{A^*(\Gamma)}{\downfree}^{\mathbf{d}} A^*(\Delta)$ *in the reduct* $A^*(M)$*, and* $f_C(u) \underset{C^*(\Gamma)}{\downfree}^{\mathbf{d}} C^*(\Delta)$ *in the reduct* $C^*(M)$*.*

*Proof.* The proof may be seen as an easier version than that of theorem 5.3.3.

Suppose $f_A(u) \underset{A^*(\Gamma)}{\not\downfree}^{\mathbf{d}} A^*(\Delta)$ in $A^*(M)$. Let $\phi_A(y,z) \in \mathcal{L}_A$, and $\delta \in A^*(\Delta)$ be such that $A^*(M) \vDash \phi_A(f_A(u),\delta)$, and $\phi_A(y,\delta)$ divides over $A^*(\Gamma)$. Let $(\delta_i)_i$ be a $A^*(\Gamma)$-indiscernible sequence witnessing division. By corollary 5.1.4, and as $\Gamma$ is generated by reals, $(\delta_i)_i$ is actually $\Gamma$-indiscernible. It follows that $\phi_A(y,\delta)$ divides over $\Gamma$, therefore $\phi_A(f_A(x),\delta)$ divides over $\Gamma$ by proposition 1.1.7, and $u \underset{\Gamma}{\not\downfree}^{\mathbf{d}} \Delta$.

Likewise, if $f_C(u) \underset{C^*(\Gamma)}{\not\downfree}^{\mathbf{d}} C^*(\Delta)$ in $C^*(M)$, then $u \underset{\Gamma}{\not\downfree}^{\mathbf{d}} \Delta$, therefore the first condition of the statement implies the second one.

Conversely, suppose by contradiction that the first condition fails, and the second one holds. As $f_C(u) \underset{C^*(\Gamma)}{\downfree}^{\mathbf{d}} C^*(\Delta)$, $u$ realizes $p_\Delta$ by the proof of proposition 5.3.1. By proposition 5.2.5, there must exist a $\Delta$-definable set $X$ implied by $p_\Delta$, $\alpha$ a tuple from $A^*(\Delta)$, $\gamma$ a tuple from $C^*(\Delta)$ and formulas $\phi_A(y,\alpha) \in \mathcal{L}_A, \phi_C(z,\gamma) \in \mathcal{L}_C$ such that $M \vDash \phi_A(f_A(u),\alpha) \wedge \phi_C(f_C(u),\gamma)$, and the formula $x \in X \wedge \phi_A(f_A(x),\alpha) \wedge \phi_C(f_C(x),\gamma)$ divides over $\Gamma$. As $\mathrm{tp}_C(u/\Delta)$ implies $p_\Delta$, we may assume $X$ is implied by $\phi_C(f_C(x),\gamma)$, therefore the formula $\phi_A(f_A(x),\alpha) \wedge \phi_C(f_C(x),\gamma)$ divides over $\Gamma$. Let $(\alpha_i\gamma_i)_i$ be a $\Gamma$-indiscernible sequence containing $\alpha\gamma$, and let us show that it does not witness division for this definable set, yielding a contradiction. By the second condition, and fact 1.1.5, let $(\alpha'_i)_i \equiv_{A^*(\Delta)} (\alpha_i)_i$, $(\gamma'_i)_i \equiv_{C^*(\Delta)} (\gamma_i)_i$ be such that $(\alpha'_i)_i$ is $A^*(\Gamma)f_A(u)$-indiscernible, and $(\gamma'_i)_i$ is $C^*(\Gamma)f_C(u)$-indiscernible. By lemma 5.6.3, and corollary 5.1.4 applied to $\mathrm{dcl}^{eq}(\Gamma u)$, the sequence $(\alpha'_i\gamma'_i)_i$ is in fact $\Gamma u$-indiscernible. However, as there exists $j$ such that $\alpha_j\gamma_j = \alpha\gamma$, we also have $\alpha'_j\gamma'_j = \alpha\gamma$. In particular, we have $M \vDash \phi_A(f_A(u),\alpha'_i) \wedge \phi_C(f_C(u),\gamma'_i)$ for every $i$ by indiscernibility. It follows that $\{\phi_A(f_A(x),\alpha'_i) \wedge \phi_C(f_C(x),\gamma'_i) | i\}$ is consistent, hence the same can be said of $\{\phi_A(f_A(x),\alpha_i) \wedge \phi_C(f_C(x),\gamma_i) | i\}$, concluding the proof. □

**Corollary 5.6.5.** *Let $U$ be the substructure of $M$ generated by $u$ and by $A^*(\Gamma) \cup B(\Gamma) \cup C^*(\Gamma)$. Then we have $U \underset{\Gamma}{\downarrow^{\mathbf{d}}} \Delta$ if and only if we have $A^*(U) \underset{A^*(\Gamma)}{\downarrow^{\mathbf{d}}} A^*(\Delta)$ and $C^*(U) \underset{C^*(\Gamma)}{\downarrow^{\mathbf{d}}} C^*(\Delta)$ in the respective reducts $A^*(M)$ and $C^*(M)$.*

*Proof.* This follows from the fact that $f_A$ and $f_C$ are defined using terms from the language of $M$ with parameters in $\Gamma$. □

**Theorem 5.6.6.** *The following are equivalent:*

- $u \underset{\Gamma}{\downarrow^{\mathbf{f}}} \Delta$.
- $f_A(u) \underset{A^*(\Gamma)}{\downarrow^{\mathbf{f}}} A^*(\Delta)$ *in the reduct* $A^*(M)$*, and* $f_C(u) \underset{C^*(\Gamma)}{\downarrow^{\mathbf{f}}} C^*(\Delta)$ *in the reduct* $C^*(M)$.

*Proof.* Let us use the same notations as in section 5.4.

By theorem 5.6.4, the ideal of $\mathcal{B}_0^{M^{eq}}$ generated by definable sets which fork over $\Gamma$ coincides with $I_A + I_C$, where:

- $I_A$ is the ideal generated by formulas of the form $\phi_A(f_A(x), \alpha) \in F_A^{M^{eq}}$ such that $\phi_A(y, \alpha)$ divides over $A^*(\Gamma)$ in the reduct $A^*(M)$.
- Likewise for $I_C$, by replacing $A^*$, $f_A$ and $F_A^{M^{eq}}$ by $C^*$, $f_C$ and $F_C^{M^{eq}}$.

Note that, without the hypothesis of this section, those ideals are in general smaller than the ones defined in the proof of theorem 5.4.3.

Suppose $f_A(u) \underset{A^*(\Gamma)}{\not\downarrow^{\mathbf{f}}} A^*(\Delta)$ in $A^*(M)$. Then there are finitely many formulas $(\phi_i(f_A(x), \alpha))_i$ in $I_A$ such that $\mathrm{tp}(f_A(u)/A^*(\Delta)) \vDash \bigvee_i \phi_i(y, \alpha)$. It follows that $\mathrm{tp}(u/\Delta) \in \cup e^{-1}(I_A)$, therefore we have $u \underset{\Gamma}{\not\downarrow^{\mathbf{f}}} \Delta$.

Likewise, if $f_C(u) \underset{C^*(\Gamma)}{\not\downarrow^{\mathbf{f}}} C^*(\Delta)$, then $u \underset{\Gamma}{\not\downarrow^{\mathbf{f}}} \Delta$, which proves the top-to-bottom direction of the statement.

Conversely, suppose $u \underset{\Gamma}{\not\downarrow^{\mathbf{f}}} \Delta$. We may assume $u$ realizes $p_\Delta$, for otherwise $f_C(u) \underset{C^*(\Gamma)}{\not\downarrow^{\mathbf{d}}} C^*(\Delta)$. As in theorem 5.6.4, by proposition 5.2.5, and the fact that $\mathrm{tp}_C(u/\Delta)$ implies $p_\Delta$, there exists $\Delta$-definable sets $X_A \in F_A^\Delta$, $X_C \in F_C^\Delta$ such that $X_A \wedge X_C$ forks over $\Gamma$, and $u \in X_A \wedge X_C$. By lemma 5.4.2, either $X_A \in I_A$, in which case $f_A(u) \mathop{\not\downarrow}^{\mathbf{f}}_{A^*(\Gamma)} A^*(\Delta)$ in $A^*(M)$, or $X_C \in I_C$, in which case $f_C(u) \mathop{\not\downarrow}^{\mathbf{f}}_{C^*(\Gamma)} C^*(\Delta)$ in $C^*(M)$. This concludes the proof. □

**Corollary 5.6.7.** *Let $U$ be the substructure of $M$ generated by $u$ and by $A^*(\Gamma) \cup B(\Gamma) \cup C^*(\Gamma)$. Then we have $U \mathop{\downarrow}^{\mathbf{f}}_{\Gamma} \Delta$ if and only if we have $A^*(U) \mathop{\downarrow}^{\mathbf{f}}_{A^*(\Gamma)} A^*(\Delta)$ and $C^*(U) \mathop{\downarrow}^{\mathbf{f}}_{C^*(\Gamma)} C^*(\Delta)$ in the respective reducts $A^*(M)$ and $C^*(M)$.*

# Chapter 6

# Extension bases

This chapter corresponds to our preprint ([Hos23a]).

We saw a question come up in the last chapter, the question of whether forking equals dividing over some parameter set. The reason why this question came up, is the same reason as to why the literature is more interested in forking than dividing: forking satisfies better behavioral properties, notably extension, and the fact that the set of formulas which fork over some parameter set is an ideal, as explained in the end of subsection 1.1.1. There is however an issue in the other direction which further motivates the question of whether forking equals dividing, one very basic property expected of forking which may sometimes fail, while it always holds for dividing: a tuple $c$ might not be independent from $A$ over $A$.

**Definition 6.0.1.** Let $A$ be a subset of some first-order structure. Then $A$ is an *extension base* if, for every parameter set $C$, we have $C \mathop{\smile\!\!\!\!\!\mid}^{\mathbf{f}}_{A} A$.

Note that if forking equals dividing over $A$, then $A$ is an extension base. In particular, every set is an extension base in a simple theory. The fact that $A$ is an extension base may be considered as a weak condition, the least that we expect when we compute forking, because it becomes very difficult to understand forking over sets that are not extension bases, and it may not contribute much to the initial problem of understanding elementary properties of first-order structures. Even if this property is weak and necessary, it happens to be a sufficient condition for various deep results about NIP and $NTP_2$ theories, such as in ([HP07], Corollary 2.10, Proposition 2.13, Proposition 3.2, Proposition 4.7. . . ), ([BYC14]). . . We should insist in particular on the main theorem from ([CK12]), which states that forking equals dividing over extension bases in any $\mathrm{NTP}_2$ theory. This goes to show that extension bases play an important role in the general study of unstable theories, and it suggests that they should play an important role in valued fields.

While it is already known in the literature that forking equals dividing in some classical valued fields, such as in ACVF or RCVF (it basically follows from Corollary 2.14 of [HP07] with some extra arguments: corollary 1.1.25, fact 1.1.19, and Theorem 1.2 from [CK12]), this fact remained unknown in $\mathrm{PL}_0$ (which we recall is the theory of the non-principal ultraproducts of the $p$-adic fields) up to the current work, which is $\mathrm{NTP}_2$. We find in this chapter various weak sufficient conditions for parameter sets in a Henselian valued field to be extension bases, and we show in particular that forking equals dividing over any parameter set (potentially with imaginaries) in $\mathrm{PL}_0$. Our results also unify the results from ACVF and RCVF from previous literature.

*Remark* 6.0.2. In the example involving cyclical orders ([TZ12], Exercise 7.1.6), one can see that the empty set is not an extension base.

Contrary to this classical non-example, note that models are extension bases: any type over a model is finitely satisfiable in that model, and hence does not fork over that model. In particular, from the fact that in any elementary extension of the ordered group $\mathbb{Z}$, any definably closed parameter set is a model, we know that any set in the theory of $(\mathbb{Z}, +, <)$ is an extension base.

Lastly, we observe that if $M$ is a field, and $\mathop{\smile\!\!\!\!\!\mid}^{\mathbf{f}} = \mathop{\smile\!\!\!\!\!\mid}^{\mathbf{alg}}$, then every set is an extension base.

Let $M$ be some first-order structure which is sufficiently saturated and strongly homogeneous. The existence of extension bases in $M$ can be proved using unary forking as follows:

**Lemma 6.0.3.** *Let $\mathfrak{A}$ be a large subset of $M$, and let $\mathfrak{C}$ be a class of small subsets of $M$. Suppose that, for every singleton $c \in \mathfrak{A}$, for all $B \in \mathfrak{C}$, we have $Bc \in \mathfrak{C}$. Then the following conditions are equivalent:*

- *For every $B \in \mathfrak{C}$, for every small subset $C$ of $\mathfrak{A}$, we have $C \mathop{\smile\!\!\!\!\!\mid}^{\mathbf{f}}_B B$*
- *For every $B \in \mathfrak{C}$, for every singleton $c \in \mathfrak{A}$, we have $c \mathop{\smile\!\!\!\!\!\mid}^{\mathbf{f}}_B B$.*

*Proof.* The first condition clearly implies the second one.

Suppose the first condition fails. Let $C$ be a small subset of $\mathfrak{A}$ such that $C \not\downarrow^{\mathbf{f}}_{B} B$. By left finite character, there exist singletons $c_1 \dots c_n \in C$ such that $c_1 \dots c_n \not\downarrow^{\mathbf{f}}_{B} B$. By (the contraposite of) proposition 1.1.18, there must exist $i$ such that $c_{i+1} \not\downarrow^{\mathbf{f}}_{Bc_1 \dots c_i} B$, so $c_{i+1} \not\downarrow^{\mathbf{f}}_{Bc_1 \dots c_i} Bc_1 \dots c_i$ by definition. However, by hypothesis, we have $Bc_1 \dots c_i \in \mathfrak{C}$, so the second condition fails, and we get the equivalence. □

**Lemma 6.0.4.** *Suppose $M$ is a Henselian valued field of residue characteristic zero. Let $A$ be a small subset of $M^{eq}$. Suppose $C \downarrow^{\mathbf{f}}_{A} A$ for every small subset $C$ of $K(M)$. Then $A$ is an extension base.*

*Proof.* Suppose $M$ is a Henselian valued field of residue characteristic zero. Let $C$ be a small subset of $M^{eq}$. For each $c \in C$, the sort of $M^{eq}$ that contains $c$ is the image of $K$ under some surjective $\varnothing$-definable function $f_c$ (in fact $M^{eq} = \mathrm{dcl}^{eq}(K(M))$), so we can find $d_c \in K(M)$ such that $f_c(d_c) = c$. Let $D$ be the union of all the $d_c$. Then $D$ is a small subset of $K(M)$ for which $C \subseteq \mathrm{dcl}(D)$, so $\mathrm{acl}(AC) \subseteq \mathrm{acl}(AD)$. By hypothesis, we have $D \downarrow^{\mathbf{f}}_{A} A$, so we get $C \downarrow^{\mathbf{f}}_{A} A$ by using corollary 1.1.22. This holds for every $C$, so $A$ is an extension base. □

Note that lemma 6.0.4 applies more generally to every first-order structure, where the sort $K$ would be replaced by the union of the home sorts of the structure at hand.

**Corollary 6.0.5.** *Let $\mathfrak{C}$ be a class of small subsets of $M$ for which, for every $A \in \mathfrak{C}$, $c \in K(M)$, we have $Ac \in \mathfrak{C}$. Then the following conditions are equivalent:*

- *Every element of $\mathfrak{C}$ is an extension base.*
- *For every $A \in \mathfrak{C}$, for every singleton $c \in K(M)$, we have $c \downarrow^{\mathbf{f}}_{A} A$.*

*Proof.* This is an immediate consequence of lemma 6.0.4 and lemma 6.0.3 (with $\mathfrak{A} = K(M)$). □

From the rest of this chapter, we assume $M$ is a Henselian valued field of residue characteristic zero.

Our work in section 5.6 should convince the reader that forking and dividing are easier to understand when some imaginaries can be lifted to reals in the parameter sets at hand. In section 6.2 and section 6.3, we look for optimal sufficient conditions for various imaginaries (balls, values in the value group, leading terms in RV...) to be lifted to a point in $K$.

## 6.1 More technical preliminaries

The contents of this section are well-known technical facts about the model theory of Henselian valued fields.

**Proposition 6.1.1** ([Fle11], Propositions 3.6, 4.3 and 5.1)**.** *Let $M \vDash \mathrm{HVF}_{0,0}$, and $A \subseteq M$ such that $A = K(A)^{\mathbf{alg}} \cap K(M)$. Let $n, m < \omega$, and $X$ an $A$-definable subset of $K^n \times \mathrm{RV}^m$. Then $X$ is definable by an $\mathcal{L}_{\mathrm{RV}}$-formula without quantifiers over the sort $K$. Moreover, if $n = 1$ and $m = 0$, then we can assume that the $K$-variable $x$ in the formula only occurs in terms of the form $\mathrm{rv}(x - a)$, with $a \in K(A)$.*

*Proof.* Proposition 4.3 of [Fle11] is exactly the first part of the statement without our hypothesis $A = A^{\mathbf{alg}} \cap K(M)$.

In the proof of Proposition 5.1 of the same paper, the $\alpha_i$ are given by applying Proposition 3.6 to polynomials with coefficient in $A$. The statement of this Proposition 3.6 shows that the $\alpha_i$ can be chosen as roots in $K(M)$ of some iterated derivatives of these polynomials, so they belong to $A^{\mathbf{alg}} \cap K(M)$, which concludes the proof. □

**Corollary 6.1.2.** *Let $M \vDash \mathrm{HVF}_{0,0}$, and $A \subseteq K(M)$. Then $K(\mathrm{acl}(A)) = A^{\mathbf{alg}} \cap K(M)$.*

We refer to this statement as the fact that models of $\mathrm{HVF}_{0,0}$ are algebraically bounded.

*Proof.* Let $X$ be a finite $A$-definable subset of $K$. Proposition 6.1.1 gives us a formula defining $X$ with parameters in $A^{\mathbf{alg}} \cap K(M)$, which satisfies the second part of proposition 6.1.1. Let $(a_i)_i$ be the finite family of points from $A^{\mathbf{alg}} \cap K(M)$ that appear in this formula. Suppose towards contradiction that $X$ has a point $a$ distinct from all the $a_i$. As the valuation is non-trivial, let $\gamma \in \Gamma^*(M)$ be strictly smaller than all the $\mathrm{val}(a - a_i)$. The closed ball $Y$ of radius $\gamma$ around $a$ is infinite, as the residue field is of characteristic zero. However, all its elements belong to $X$, for if $b \in Y$, then we have $\mathrm{rv}(b - a_i) = \mathrm{rv}(a - a_i)$ for all $i$. As a result, $X$ is infinite, a contradiction. □

*Remark* 6.1.3. If $M \vDash \mathrm{HVF}_{0,0}$, then corollary 6.1.2 implies in particular that $K(\mathrm{acl}(\varnothing)) \subseteq \mathbb{Q}^{\mathbf{alg}}$, which is trivially valued as the residue characteristic is zero.

**Lemma 6.1.4.** *Let* $M \vDash \mathrm{HVF}_{0,0}$*, and* $a_1 \dots a_n, b_1 \dots b_n \in K(M)$*. Suppose* $(\mathrm{val}(a_1), \dots \mathrm{val}(a_n))$ *is a* $\mathbb{Q}$*-free family in the divisible closure* $\mathrm{div}(\Gamma^*(M))$ *of* $\Gamma^*(M)$*, and* $\mathrm{rv}(a_i) = \mathrm{rv}(b_i)$ *for each* $i$*. Then* $a_1 \dots a_n \equiv_\varnothing b_1 \dots b_n$*.*

*Proof.* By 6.1.1, we just have to prove $\mathrm{rv}(P(a_1 \dots a_n)) = \mathrm{rv}(P(b_1 \dots b_n))$ for each polynomial $P$ with $\varnothing$-definable coefficients. Let $P$ be such a polynomial. By remark 6.1.3, the non-zero coefficients of $P$ all have value 0. As the $\mathrm{val}(a_i)$ form a $\mathbb{Q}$-free family, the values of every monomial involved in $P(a_1 \dots a_n)$ are pairwise-distinct. Let $Q$ be the monomial in $P$ such that $Q(a_1 \dots a_n)$ has largest value. Then we have $\mathrm{rv}(P(a_1 \dots a_n)) = \mathrm{rv}(Q(a_1 \dots a_n))$. As $\mathrm{rv}(a_i) = \mathrm{rv}(b_i)$, we have $\mathrm{val}(a_i) = \mathrm{val}(b_i)$, so $Q(b_1 \dots b_n)$ is also the monomial of $P(b_1 \dots b_n)$ of largest value, and $\mathrm{rv}(P(b_1 \dots b_n)) = \mathrm{rv}(Q(b_1 \dots b_n))$. Now, $\mathrm{rv}(Q(a_1 \dots a_n))$ is just a product of powers of the $\mathrm{rv}(a_i)$ multiplied by a constant. As $\mathrm{rv}(a_i) = \mathrm{rv}(b_i)$, we have $\mathrm{rv}(Q(a_1 \dots a_n)) = \mathrm{rv}(Q(b_1 \dots b_n))$. That concludes the proof. □

**Fact 6.1.5** (see for instance [vdD14], Corollary 5.25)**.** *Let* $M \vDash \mathrm{HVF}_{0,0}$*, and* $A \subseteq M$ *such that* $A = \mathrm{acl}(K(A))$*. Then any* $A$*-definable subset of* $\Gamma(M)$ *is actually* $\mathrm{val}(A)$*-definable in the first-order structure* $(\Gamma(M), +, <)$*. Moreover, any* $A$*-definable subset* $X$ *of* $k(M)$ *is actually* $k(M)$*-definable in the ring* $k(M)$*, and, if* ac *is an angular component on* $M$*, then* $X$ *is* $\mathrm{ac}(A)$*-definable in the ring* $k(M)$*.*

This shows that $k$ and $\Gamma$ are both stably embedded reducts, but the way that $\Gamma$ is stably embedded leaves for control on the parameters, which is crucial when one works with forking. For instance, one property that we get for $\Gamma$ and not for $k$ is:

**Corollary 6.1.6.** *Let* $M \vDash \mathrm{HVF}_{0,0}$*,* $A \leqslant M$ *a substructure, and let* $\gamma$ *and* $\delta$ *be two tuples from* $\Gamma(M)$*, such that* $A = \mathrm{acl}(K(A))$*. If* $\gamma \equiv_{\mathrm{val}(A)} \delta$ *in* $(\Gamma(M), +, <)$*, then* $\gamma \equiv_A \delta$ *in the valued field* $M$*.*

### 6.1.1 Lifts of the residue field

**Definition 6.1.7.** In a valued field $M$, a *lift of the residue field* is a ring morphism $k(M) \longrightarrow \mathcal{O}(M)$ that is a section of res.

In the Hahn field $k'((t^G))$, the standard lift maps an element of $k'$ to the corresponding constant polynomial.

**Proposition 6.1.8** ([vdD14], subsection 2.4)**.** *Let $M \vDash \mathrm{HVF}_{0,0}$. Then there is a canonical 1-1 correspondence between the lifts of the residue field in M and the maximal proper subfields of $\mathcal{O}(M)$.*

*Moreover, if $N \vDash \mathrm{HVF}_{0,0}$, $f : N \longrightarrow M$ is an embedding, and $l$ is a lift of the residue field in $N$, then there exists at least one lift of the residue field $l'$ in $M$ such that $l'$ extends $flf^{-1} : f(k(N)) \longrightarrow \mathcal{O}(M)$.*

**Corollary 6.1.9.** *In any valued field $M \vDash \mathrm{HVF}_{0,0}$, the residue field can be lifted.*

*Proof.* Apply proposition 6.1.8 to the (canonical) embedding into $M$ whose domain is the trivially valued field $\mathbb{Q}$ (which is of course Henselian). □

*Remark* 6.1.10. If $M \vDash \mathrm{HVF}_{0,0}$, $a \in K(M)$, and $\alpha \in k(M)$, then the following conditions are equivalent:

- There exists a lift of the residue field sending $\alpha$ to $a$.
- $\alpha$ and $a$ have the same ideal (i.e. the same algebraic type) over $\mathbb{Q}$.

**Lemma 6.1.11.** *Let $M \vDash \mathrm{HVF}_{0,0}$. Suppose we have a cross-section $s$, and let* ac *be the corresponding angular component. Let $l$ be a lift of the residue field. Let $a \in K(M)^*$ such that* $\mathrm{val}(a) \neq 0$. *Then $a \equiv_\varnothing l(\mathrm{ac}(a))s(\mathrm{val}(a))$ in $M$.*

This elementary equivalence is in $M$, without the cross-section and the angular component.

*Proof.* Let $b = l(\mathrm{ac}(a))s(\mathrm{val}(a))$. We have $\mathrm{ac}(a) = \mathrm{ac}(b)$ and $\mathrm{val}(a) = \mathrm{val}(b)$, so $\mathrm{rv}(a) = \mathrm{rv}(b)$ (as ac and $s$ exist, we have a group isomorphism between $\mathrm{RV}(M)$ and $k(M)^* \times \Gamma(M)$). Moreover, as $\mathrm{val}(a) \neq 0$, we can apply 6.1.4 to get $a \equiv_\varnothing b$. □

## 6.2 Generics of chains of balls

In this section we give definitions and show properties which formalize in the model theory of valued field the intuition that we had in remark 2.1.7.

**Definition 6.2.1.** Let $X$ be a ball in $M$. If $A \subseteq M^{eq}$, then $X$ is an *$A$-ball* if the imaginary canonical parameter of $X$ is an element of $\mathrm{dcl}^{eq}(A)$. Moreover, we say that an $A$-ball $X$ is *pointed* if it has a point in $K(A)$.

We define two notions: genericity and weak genericity. Let us first give an illustration before the formal details. Suppose $A = \mathrm{acl}^{eq}(A) \subseteq M^{eq}$, and let $c \in K(M)$. The set $\mathfrak{B}$ of every $A$-ball that contains $c$ turns out to be a chain with respect to inclusion, this is a well-known fact in ultrametric geometry. Now, not every $A$-ball is necessarily pointed, so let $\mathfrak{B}'$ be the chain of pointed $A$-balls containing $c$. The chain $\mathfrak{B}'$ is a final segment of $\mathfrak{B}$ which might be strict. Then, the type-definable set $\cap\mathfrak{B}'$ is an approximation of $c$ that is weaker than $\cap\mathfrak{B}$. Even if this approximation is weaker, it is witnessed by points from $K(A)$ (because the balls of $\mathfrak{B}'$ are pointed), hence it is often easier to manipulate than $\cap\mathfrak{B}$, especially when $A$ is generated by field parameters.

Let us now define formally our notions of genericity. Let $\mathfrak{B}$ be a chain (with respect to inclusion) of $A$-balls, with $A$ an arbitrary parameter set not necessarily algebraically closed, and let $c \in K$ be in an elementary extension of $M$. By convention, if $\mathfrak{B} = \varnothing$, then $\cap\mathfrak{B}$, the intersection of every ball in $\mathfrak{B}$, is the definable set $K$. We say that $c$ is *$A$-generic* of $\mathfrak{B}$ if the following conditions hold:

- $c$ is in the $A$-type-definable set $\cap\mathfrak{B}$.
- For every $A$-ball $X$ strictly contained in $\cap\mathfrak{B}$, we have $c \notin X$.

The above conditions clearly yield a partial type over $A$. We can use basic saturation arguments to show that such a type is consistent. In case $\mathfrak{B}$ has a least element, then we have to use the fact that the residue field is infinite (and maybe use a density argument in case the least ball is not closed). Else, consistency is simply established by compactness.

Likewise, if $\mathfrak{B}$ is a chain of **pointed** $A$-balls, then $c$ is said to be *weakly $A$-generic* of $\mathfrak{B}$ when the following conditions hold:

- $c \in \cap\mathfrak{B}$.

- For every **pointed** $A$-ball $X$ strictly contained in $\cap\mathfrak{B}$, we have $c \notin X$.

Note that our notion of weak genericity only makes sense when $\mathfrak{B}$ is a chain of pointed $A$-balls, otherwise it is undefined. If $\mathfrak{B}$ has a non-pointed element, then an $A$-generic point of $\mathfrak{B}$ is weakly $A$-generic of some chain of $A$-balls whose intersection is strictly coarser than $\cap\mathfrak{B}$, and it is not weakly $A$-generic of $\mathfrak{B}$. However, if every ball of $\mathfrak{B}$ is pointed, then $A$-genericity also implies weak $A$-genericity.

Geometrically speaking, the partial type corresponding to the $A$-generics of $\mathfrak{B}$ can be written as an intersection of some decreasing sequence of $A$-definable crowns (a crown being a ball from which we removed some subball that can be empty). In case of weak genericity, less crowns are involved, so the set of realizations is larger.

If $\mathfrak{B}$ has a least element which is a closed ball, then we say that $\mathfrak{B}$ is *residual*. If not, then $\mathfrak{B}$ is called *ramified* (with respect to $A$) if $\cap\mathfrak{B}$ strictly contains an $A$-ball or a point in $K(A)$, and $\mathfrak{B}$ is called *immediate* (with respect to $A$) otherwise. We refer the reader to fact 2.2.13 for further explanation.

If $X$ is a ball, $a \in X$ and $b \notin X$, then we have $\mathrm{rv}(b-a) = \mathrm{rv}(b-a')$ for all $a' \in X$. As a result, we can define $\mathrm{rv}(b-X)$ to be $\mathrm{rv}(b-a)$ for any $a \in X$.

Recall that if the value group is discrete, then every ball is both open and closed. If it is dense however, then there is no ball that is both open and closed.

*Remark* 6.2.2. One can note that the sets of $A$-generic points of each chain of $A$-balls are the classes of some equivalence relation (to belong to the same $A$-balls), which is even $A$-type-definable. As a result, two chains of $A$-balls $\mathfrak{B}$ and $\mathfrak{B}'$ have the same $A$-generic points if and only if they have at least one $A$-generic point in common.

Using the definition of genericity, one can easily show that $\mathfrak{B}$ and $\mathfrak{B}'$ have the same $A$-generic points if and only if the $A$-type-definable sets $\cap\mathfrak{B}$ and $\cap\mathfrak{B}'$ coincide.

Moreover, if the balls of $\mathfrak{B}$ and $\mathfrak{B}'$ are pointed, then the same two characterizations hold when we look at their $A$-weakly generic points.

When the parameter set is algebraically closed, and generated by field elements, we have a good understanding of when the notions of genericity and weak genericity do not coincide:

**Proposition 6.2.3.** *Let $A \subseteq M^{eq}$ such that $A = \mathrm{acl}^{eq}(K(A))$, and let $\mathfrak{B}$ be a chain of pointed $A$-balls. Suppose there is $c \in K(M)$, $A$-weakly generic of $\mathfrak{B}$, but not $A$-generic of $\mathfrak{B}$.*

*Then $\mathfrak{B}$ is residual. Moreover, if $Y$ is the largest open ball that contains $c$ which is strictly contained in $\cap\mathfrak{B}$, then $Y$ is an $A$-ball (which must not be pointed by weak genericity of $c$), and $c$ is $A$-generic of $\{Y\}$.*

*Proof.* By hypothesis, there exists $Y$ an $A$-ball containing $c$ that is not pointed. In particular, $Y$ is strictly contained in $\cap\mathfrak{B}$. By proposition 6.1.1, there exist a finite family $(a_i)_{i<N}$ of elements of $K(A)$, and $Z$ an $A$-definable subset of $\mathrm{RV}^N$, such that:

$$Y = \{x \in K | (\mathrm{rv}(x - a_i))_i \in Z\}$$

Clearly, as $Y$ is not pointed, none of the $a_i$ belongs to $Y$. Let $\gamma \in \Gamma(A)$ be the largest of the $\mathrm{val}(a_i - Y)$. There does not exist any $\delta \in \Gamma(M)$ such that $\gamma > \delta > \mathrm{rad}(Y)$, for if not, then choose $b \in K(M)$ for which $\mathrm{val}(b-c) = \delta$. For every $i$, we have $\mathrm{rv}(a_i - b) = \mathrm{rv}(a_i - c)$, so $b \in Y$, which contradicts $\mathrm{val}(b-c) > \mathrm{rad}(Y)$.

If $\Gamma(M)$ is dense, this must imply that $Y$ is open and $\mathrm{rad}(Y) = \gamma$. If $\Gamma(M)$ is discrete, then of course $Y$ is open and closed, and the radius of $Y$ as a closed ball is the successor of $\gamma$.

As $Y$ is open, let $X$ be the least closed ball that strictly contains $Y$. As $X$ is definable over the canonical parameter of $Y$, $X$ is an $A$-ball. The radius of $X$ is $\gamma$, so at least one of the $a_i$ must be in $X$, so $X$ is pointed. There is no ($M^{eq}$-)ball between $X$ and $Y$, and $c \in Y$, so $c$ must be $A$-weakly generic of $\{X\}$. Now, by remark 6.2.2 applied to $\{X\}$ and $\mathfrak{B}$, $X$ must be a lower bound of $\mathfrak{B}$. Then, $X$ must be the least element of $\mathfrak{B}$, otherwise a standard compactness argument could be used to show that the inclusion of type-definable sets $X \subseteq \cap\mathfrak{B}$ would be strict, which would contradict remark 6.2.2. As a result, we notice that $X$ is canonical, its construction does not depend on the choice of $Y$, but $Y$ does depend on $X$: it is the largest open ball containing $c$ and strictly contained in $X$. Therefore, the choice of $Y$ must be unique: $Y$ is the only $A$-ball containing $c$ that is not pointed. There is no $A$-ball smaller than $Y$ containing $c$, so $c$ is $A$-generic of $\{Y\}$, this concludes the proof. □

**Corollary 6.2.4.** *Let $A \subseteq M^{eq}$ be such that $A = \mathrm{acl}^{eq}(K(A))$. If $\Gamma(M)$ is dense, then any closed $A$-ball is pointed.*

*Proof.* Let $X$ be a closed $A$-ball, and $c \in K(M)$ an $A$-generic point of $\{X\}$. Let $\mathfrak{B}$ be the set of every pointed $A$-ball containing $c$. Suppose by contradiction that $X$ is not pointed. Then $c$ is not $A$-generic of $\mathfrak{B}$, so we can apply proposition 6.2.3 to show that $c$ is $A$-generic of $\{Y\}$ for some open ball $Y$. However, $X$ is not open, therefore $\cap\{X\} = X \neq Y = \cap\{Y\}$. By remark 6.2.2, $c$ cannot be $A$-generic of $\{X\}$ and $\{Y\}$ at the same time, which concludes our proof by contradiction. ☐

**Corollary 6.2.5.** *Let $A \subseteq M^{eq}$ be such that $A = \mathrm{acl}^{eq}(K(A))$. Let $\mathfrak{B}$ be a chain of $A$-balls. If $\mathfrak{B}$ is ramified, then $\cap\mathfrak{B}$ actually has a point in $K(A)$.*

*Proof.* Let $Y$ be an $A$-ball strictly contained in $\cap\mathfrak{B}$. If $Y$ is closed and $\Gamma(M)$ is dense, then we are done by corollary 6.2.4.

Else, let $X$ be the least closed ball strictly containing $Y$ (which is in fact the least ball strictly containing $Y$). The ball $X$ is clearly an $A$-ball, and we have $X \subseteq \cap\mathfrak{B}$ by minimality of $X$. This inclusion is strict as $\mathfrak{B}$ is not residual. Let $c$ be an $A$-generic point of $\{Y\}$. We are done if $Y$ is pointed, so suppose it is not. Let $\mathfrak{B}'$ be the chain of pointed $A$-balls containing $c$. By compactness, $Y$ is strictly included in $\mathfrak{B}'$, therefore $c$ is not $A$-generic of $\mathfrak{B}'$, despite being weakly $A$-generic. By proposition 6.2.3, $X$ is the least element of $\mathfrak{B}'$, therefore $X$ is pointed, and its points belong $\cap\mathfrak{B}$, concluding the proof. ☐

**Proposition 6.2.6.** *Let $A \subseteq M^{eq}$ be such that $A = \mathrm{acl}^{eq}(K(A))$. Let $c, d$ be singletons in $K(M)$, and let $\mathfrak{B}$ be the set of pointed $A$-balls that contain $c$. Suppose $d \in \cap\mathfrak{B}$. Then:*

- *If $\cap\mathfrak{B}$ has no point in $K(A)$, then $c \equiv_A d$.*
- *If $a \in K(A)$ is in $\cap\mathfrak{B}$, and $\mathrm{rv}(c-a) = \mathrm{rv}(d-a)$, then $c \equiv_A d$.*

*Proof.* As $A = \mathrm{acl}^{eq}(K(A))$, it is sufficient to show $c \equiv_{K(A)} d$. By proposition 6.1.1, it is sufficient to show that $\mathrm{rv}(c-a') = \mathrm{rv}(d-a')$ for all $a' \in K(A)$. Let $a' \in K(A)$.

- Suppose $a' \notin \cap\mathfrak{B}$. Let $X \in \mathfrak{B}$ such that $a' \notin X$. Then $\mathrm{rv}(a'-c) = \mathrm{rv}(a'-X) = \mathrm{rv}(a'-d)$.
- Suppose $a' \in \cap\mathfrak{B}$ (so the point $a$ from the statement exists). Let $X$ be the least closed ball that contains $a$ and $a'$. If $X$ is a singleton, then $a = a'$, therefore $\mathrm{rv}(c - a') = \mathrm{rv}(c - a) = \mathrm{rv}(d - a) = \mathrm{rv}(d - a')$. Suppose $X$ is not a singleton.

Let $Y$ (resp. $Y'$) be the largest open ball strictly contained in $X$ such that $a \in Y$ (resp. $a' \in Y'$). By minimality of $X$, we have $a' \notin Y$. As $a' \in \cap\mathfrak{B}$, and $a' \notin Y$, by definition of $\mathfrak{B}$, we have $c \notin Y$. Likewise, we have $c \notin Y'$.

By definition 2.2.4, in order to show that $\mathrm{rv}(c - a') = \mathrm{rv}(d - a')$, it suffices to show that $\mathrm{val}(c - d) < \mathrm{val}(c - a')$. By hypothesis, we have $\mathrm{val}(c-d) < \mathrm{val}(c-a)$. As $c \notin Y$ and $c \notin Y'$, we have $\mathrm{val}(c-a) = \mathrm{val}(c-a')$, and we are done.

either way, we have $\mathrm{rv}(c - a') = \mathrm{rv}(d - a')$, concluding the proof. □

An extensive study of unary types over imaginary parameters has been carried out in [HMRC18], in the setting of pseudo-local fields. In particular, there is an analogue of proposition 6.2.6 with imaginary parameters, when the value group is a $\mathbb{Z}$-group:

**Proposition 6.2.7** ([HMRC18], Lemma 3.7)**.** *Suppose that $\Gamma(M) \equiv \mathbb{Z}$. Let $A \subseteq M^{eq}$ be such that $A = \mathrm{acl}^{eq}(A)$. Let $c, d \in K(M)$ be singletons, and let $\mathfrak{B}$ be the set of $A$-balls that contain $c$. Suppose $d \in \cap\mathfrak{B}$, and $\mathfrak{B}$ is $A$-immediate. Then $c \equiv_A d$.*

**Proposition 6.2.8.** *Let $A \subseteq B \subseteq M^{eq}$ be such that $A = \mathrm{acl}^{eq}(A)$, let $\mathfrak{B}$ be a chain of $A$-balls, and let $c \in K(M)$ be an $A$-generic point of $\mathfrak{B}$. Then there exists $c' \equiv_A c$ (in an elementary extension of $M$) that is $B$-generic of $\mathfrak{B}$.*

*Proof.* Suppose not. Then, by compactness, there exists a finite family $(X_i)_i$ of $B$-balls strictly contained in $\cap\mathfrak{B}$ such that, for all $c' \equiv_A c$, $c'$ is in $\bigcup_i X_i$. We may assume that each $X_i$ contains at least one $A$-conjugate of $c$. Replace $M$ by a $|B|^+$-saturated, strongly $|B|^+$-homogeneous elementary extension.

- Suppose $\mathfrak{B}$ is residual. Let $Y_i$ be the largest open ball containing $X_i$, and strictly contained in $\min(\mathfrak{B})$ (the least element of $\mathfrak{B}$). Then the $Y_i$ are exactly the maximal proper subballs of $\min(\mathfrak{B})$ containing some $A$-conjugate of $c$, therefore the definable set $\bigcup_i Y_i$ is $\mathrm{Aut}(M/A)$-invariant, so each $Y_i$ is an $\mathrm{acl}^{eq}(A) = A$-ball, so $c \notin Y_i$ by $A$-genericity. However, we have $c \in \bigcup_i Y_i$, a contradiction.

- Suppose $\mathfrak{B}$ is not residual. Let $X$ be the least closed ball that contains all the $X_i$. Then $X$ is the least closed ball that contains every $A$-conjugate of $c$, so $X$ is $\mathrm{Aut}(M/A)$-invariant, therefore $X$ is an $A$-ball. As $c \in X$, $X \in \mathfrak{B}$ (we can replace $\mathfrak{B}$ by the set of every $A$-ball containing $\cap\mathfrak{B}$, without loss of generality). By hypothesis, $X$ is not the least element of $\mathfrak{B}$, so the largest open ball containing $c$ and strictly contained in $X$ (let us call it $Y$) is an $A$-ball. This contradicts the minimality of $X$, as each $X_i$ is strictly contained in $Y$, so the least closed ball that contains all of them must be (weakly) contained in $Y$. □

Before the next section, we prove a technical lemma giving sufficient conditions to lift the radius of $A$-balls:

**Lemma 6.2.9.** *Let $M \models HF_{0,0}$, and $A = \mathrm{acl}^{eq}(K(A)) \subseteq M^{eq}$. Suppose the following conditions hold:*

- *$M$ admits a cross-section $s$. Denote by* ac *the corresponding angular component.*
- *For every $\gamma \in \Gamma(\mathrm{dcl}(\varnothing))$, there exists an $a \in K(A)$ such that* $\mathrm{val}(a) = \gamma$.
- *For every $N > 0$, $\alpha \in \mathrm{ac}(K(A))$, there exists $a \in \mathcal{O}^* \cap A$ such that* $\mathrm{res}(a) \bmod (k^*)^N = \alpha \bmod (k^*)^N$.

*then, for every $\delta \in \Gamma(A)$, there exists some $\alpha \in K(A)$ such that* $\mathrm{val}(\alpha) = \delta$.

In particular, the radius of each $A$-ball is equal to $\mathrm{val}(a)$ for some $a \in K(A)$.

*Proof.* Let $\delta \in \Gamma(A)$. Then $\delta$ is $K(A)$-definable, so by fact 6.1.5 it is actually $\mathrm{val}(K(A))$-definable. By ([CH11], Corollary 1.10), the definable closure of $\mathrm{val}(K(A))$ in the ordered group $\Gamma(M)$ is exactly the relative divisible closure of the subgroup of $\Gamma(M)$ generated by $\mathrm{val}(K(A))$ and $\Gamma(\mathrm{dcl}(\varnothing))$. In other words, there exist $\gamma \in \Gamma(\mathrm{dcl}(\varnothing))$, $a \in K(A)$, and $N > 0$ such that $\delta = \frac{\mathrm{val}(a)+\gamma}{N}$. By hypothesis, there exists $a' \in K(A)$ such that $\mathrm{val}(a') = \gamma$.

Let $\bar{a} \in \mathcal{O}_A^*$ such that $\mathrm{res}(\bar{a}) \bmod (k^*)^N = \mathrm{ac}(aa') \bmod (k^*)^N$. Let $a'' = aa'\bar{a}^{-1}$. We have $a'' \in K(A)$ and $\mathrm{val}(a'') = \mathrm{val}(a) + \gamma$. We conclude the proof by showing that $a''$ has an $N$-th root in $K(M)$, this root will have value $\delta$, and it will also belong to $K(A)$, as $A = \mathrm{acl}(A)$.

We have $\mathrm{ac}(a'') = \frac{\mathrm{ac}(aa')}{\mathrm{res}(\bar{a})} \in (k^*)^N$. Let $r' \in k^*$ such that $r'^N = \mathrm{ac}(a'')$, and let $l$ be a lift of the residue field. We have $(l(r'))^N = l(\mathrm{ac}(a''))$, so $l(\mathrm{ac}(a''))$ has an $N$-th root in $K(M)$. As a result, $l(\mathrm{ac}(a''))s(\mathrm{val}(a''))$ also has an $N$-th root in $K(M)$, namely $l(r')s(\delta)$. By lemma 6.1.11, $l(\mathrm{ac}(a''))s(\mathrm{val}(a'')) \equiv_\varnothing a''$, so $a''$ also has an $N$-th root in $K(M)$, concluding the proof. □

## 6.3 The strength of the residue field's stable embeddedness

In fact 6.1.5, one cannot help but notice that the result on definable subsets of $k$ is weaker than that on $\Gamma$. In particular, we don't have an analogue of corollary 6.1.6 in $k$. This property will sometimes be needed for our results on unary forking, so we would like to understand when it holds. In this section, we show a counterexample in which this analogue does not hold, and we give sufficient conditions to obtain the same nice behavior as what we have in $\Gamma$. The strategy for the proofs is to work in an elementary extension that is a Hahn field, where the specific structure gives us more tools to build elementary maps.

We invite the reader to recall what we said in remark 2.3.2 to not get confused in the notations. For instance, in definition 6.3.6, we have an explicit field $k'$, but we use the notation $k$ to define some definable subset $X$ of $k'$, because $X$ has points in elementary extensions of $k'$ that are not in $k'$.

### 6.3.1 Algebraic technicalities

In order to build our example and prove our version of corollary 6.1.6 in $k$, we first need to state simple field-theoretic facts.

**Definition 6.3.1.** Let $G$ be an ordered Abelian group, and $k'$ a field of characteristic zero. We define two (injective) group homomorphisms:

$$\begin{array}{rccc} \mathrm{aut}^1_{k',G} : & \mathrm{Aut}(k') & \longrightarrow & \mathrm{Aut}(k'((t^G))) \\ \mathrm{aut}^2_{k',G} : & \mathrm{Hom}(G, k'^*) & \longrightarrow & \mathrm{Aut}(k'((t^G))) \end{array}$$

The automorphisms of $k'((t^G))$ that we consider are of course valued field automorphisms.

For every $\sigma \in \mathrm{Aut}(k')$, $\mathrm{aut}^1_{k',G}(\sigma)$ is defined as $\sum_\alpha a_\alpha t^\alpha \longmapsto \sum_\alpha \sigma(a_\alpha) t^\alpha$.

For every $\sigma \in \mathrm{Hom}(G, k'^*)$, $\mathrm{aut}^2_{k',G}(\sigma)$ is defined as $\sum_\alpha a_\alpha t^\alpha \longmapsto \sum_\alpha a_\alpha \sigma(\alpha) t^\alpha$.

The map $\mathrm{aut}^1_{k',G}$ is classical, and the map $\mathrm{aut}^2_{k',G}$ is denoted by $P$ in the Definition 3.3.5 of [KS22], to which we refer the interested reader for a more elaborate study of automorphisms of Hahn fields (and some of their subfields).

**Lemma 6.3.2.** *Let $G$ be a divisible group, and $g \in G$. Then there exists $\sigma \in \mathrm{Hom}(\mathbb{Q}, G)$ such that $\sigma(1) = g$.*

*Proof.* Let $g_1 = g$. By induction on $n > 0$, choose $g_{n+2} \in G$ an $(n+2)$-th root of $g_{n+1}$ (i.e. $g_{n+2}^{n+2} = g_{n+1}$). In particular, $g_{n+1}$ is always an $((n+1)!)$-th root of $g$. For each $(n,m) \in \mathbb{Z} \times \mathbb{Z}_{>0}$, define $f(n,m) = (g_m)^{n(m-1)!}$.

Let us show that $f$ factors through the canonical surjection $\mathbb{Z}\times\mathbb{Z}_{>0} \longrightarrow \mathbb{Q}$. Let $n, n' \in \mathbb{Z}$, $m, m' \in \mathbb{Z}_{>0}$ be such that $nm' = n'm$. We have to show that $(g_m)^{n(m-1)!} = (g_{m'})^{n'(m'-1)!}$. This is trivial if $m = m'$, so we may assume $m > m'$. Let $N = \dfrac{m!}{m'!} = \prod_{i=m'+1}^{m} i$. Then we have $g_{m-1}^{mN} = g_{m'}$. As a result:

$$\begin{aligned}
(g_{m'})^{n'(m'-1)!} &= g_{m-1}^{mNn'(m'-1)!} \\
&= g_{m-1}^{m'nN(m'-1)!} \\
&= g_{m-1}^{nN(m'!)} \\
&= g_{m-1}^{n(m!)} \\
&= g_m^{n(m-1)!}
\end{aligned}$$

It follows that $f$ factors through the canonical surjection $\mathbb{Z} \times \mathbb{Z}_{>0} \longrightarrow \mathbb{Q}$. The induced map $\overline{f} : \mathbb{Q} \longrightarrow G$ clearly sends 1 to $g$. It remains to show that $\overline{f}$ is a group homomorphism. This follows from the computation:

$$
\begin{aligned}
\overline{f}\left(\frac{n}{m}\right)\overline{f}\left(\frac{n'}{m'}\right) &= (g_m)^{n(m-1)!}(g_{m'})^{n'(m'-1)!} \\
&= g_{mm'}^{[(mm')!/m!]n(m-1)!} g_{mm'}^{[(mm')!/m'!]n'(m'-1)!} \\
&= g_{mm'}^{(mm')!([n/m]+[n'/m'])} \\
&= g_{mm'}^{(mm')![(nm'+n'm)/mm']} \\
&= g_{mm'}^{(nm'+n'm)(mm'-1)!} \\
&= \overline{f}\left(\frac{n}{m}+\frac{n'}{m'}\right)
\end{aligned}
$$

this concludes the proof. □

**Lemma 6.3.3.** *Let $k'$ be a field of characteristic zero, and $n > 0$. Then, in the field $k'(u)$ ($u$ being any transcendental element), $\frac{u}{u+n}$ is not a square.*

*Proof.* If $\frac{u}{u+n}$ were a square, then there would exist $P$, $Q \in k'[u]$ such that $Q \neq 0$ and $Q^2 \cdot u = P^2 \cdot (u+n)$. Now, $k'[u]$ is a UFD where $u$ and $u+n$ are coprime and irreducible. As a result, $u$ appears an even amount of times in the factorization of $P^2 \cdot (u+n)$ in $k'[u]$, and it appears an odd amount of times in that of $Q^2 \cdot u$, so we cannot have $Q^2 \cdot u = P^2 \cdot (u+n)$. □

**Lemma 6.3.4.** *Let $k'$ be an algebraically closed field of characteristic zero, $1 \neq \lambda \in \mathbb{Q}_{>0}$, and $P$, $Q$ coprime non-zero polynomials in $k'[u]$. If $\frac{P(u)}{Q(u)} = \frac{P(\lambda u)}{Q(\lambda u)}$ in $k'(u)$, then both $P$ and $Q$ are constant.*

*Proof.* Suppose, say, $P$ is non-constant (the same proof works if $Q$ is non-constant), and admits a non-zero root $\alpha$. As $\lambda \in \mathbb{Q}^*$, and $k'$ has characteristic zero, the orbit of $\alpha$ under the bijection $h : x \longmapsto \lambda^{-1}x$ is infinite. Therefore, there exists $n \in \mathbb{Z}$ such that $h^n(\alpha)$ is a root of $P(u)$, but $h^{n+1}(\alpha)$ is not. As $h^n(\alpha)$ is a root of $P(u)$, $h^{n+1}(\alpha)$ is a root of $P(\lambda u)$. As $P(u)Q(\lambda u) = P(\lambda u)Q(u)$, $h^{n+1}(\alpha)$ is a root of $P(u)Q(\lambda u)$, so it is a root of $Q(\lambda u)$. This contradicts the coprimality hypothesis, as $h^n(\alpha)$ is a common root of $P(u)$ and $Q(u)$.

Now, $P$ and $Q$ are coprime and do not admit any non-zero root. In particular, as $k'$ is algebraically closed, $P$ and $Q$ are either both constant, or one of them is constant and the other can be written $\beta X^N$ for some $\beta \in k'^*$, $N > 0$. Then, if $P$ and $Q$ are not both constant, it is easy to see that $\frac{P(u)}{Q(u)} \neq \frac{P(\lambda u)}{Q(\lambda u)}$, which concludes the proof. □

**Corollary 6.3.5.** *Let $k'$ be a field of characteristic zero, $0 < N \neq M > 0$, and $F \in k'(u) \setminus k'$. Then $F(Mu) \neq F(Nu)$.*

*Proof.* We can assume $k'$ is algebraically closed. We can write $F(Nu) = \frac{P(u)}{Q(u)}$, with $P$, $Q \in k'[u]$ coprime and non-zero (because $F \neq 0$). If we had $F(Mu) = F(Nu)$, then we would have $\frac{P(u)}{Q(u)} = \frac{P(\frac{M}{N}u)}{Q(\frac{M}{N}u)}$, so both $P$ and $Q$ would be constant by lemma 6.3.4, so $F \in k'$, contradicting the hypothesis. □

### 6.3.2 An example of weak stable embeddedness

Let us build an example of a Henselian valued field $M$ of residue characteristic zero, a subfield $A < K(M)$, and an $A$-definable set $X \subseteq k(M)$ such that $X$ is not $k(\mathrm{acl}(A))$-definable (and hence in particular not $\mathrm{res}(K(\mathrm{acl}(A)))$-definable). We will also have $\Gamma(\mathrm{acl}(A)) > \Gamma(\mathrm{dcl}(\varnothing))$, such that the example is not too degenerate.

**Definition 6.3.6.** Let $k' = \mathbb{Q}^{\mathbf{alg}}(u)$, $G = \mathbb{Q}$, and $M = k'((t^G))$. We write ac for the natural angular component map of $M$, and we identify $k'$ with a subfield of $K(M)$ via the natural lift of the residue field. Let $A \subseteq K(M)$ be the field generated by $u \cdot t^{-1}$, and $X = u \bmod (k^*)^2 = \{x \in k | x \cdot u^{-1} \in (k^*)^2\}$.

**Proposition 6.3.7.** *The set $X$ is $A$-definable.*

Recall that with our conventions, the elements of the valuation ring have negative value, in particular $t^{-1}$ is in the valuation ring.

*Proof.* Consider the $A$-definable set:

$$Y = \left\{x \in k | \exists y \in K \ \mathrm{val}(y) = \mathrm{val}(u \cdot t^{-1}) \wedge y \in (K^*)^2 \wedge x = \mathrm{res}(u \cdot t^{-1} \cdot y^{-1})\right\}$$

Let us show that $X = Y$.

If $x \in X(M)$, then we have $x = uv^2$ for some $v \in k'^*$. We have $x \in Y(M)$ by choosing $y = (v^{-1} \cdot t^{\frac{-1}{2}})^2$.

Let $x \in Y(M)$, and $y$ be as in the definition of $Y$. Then $\mathrm{ac}(y)$ must be in $(k'^*)^2$, and we have $x = u \cdot \mathrm{ac}(y)^{-1}$, so $x \in X(M)$.

We have $X(M) = Y(M)$, so $X = Y$. □

**Proposition 6.3.8.** *We have $k(\mathrm{acl}(A)) = \mathbb{Q}^{\mathbf{alg}}$.*

*Proof.* Let $N > 0$. Let $\sigma_N \in \mathrm{Aut}(k')$ be defined as $\sigma_N(F(u)) = F(Nu)$. Apply lemma 6.3.2 to the divisible group $(\mathbb{Q}^{\mathbf{alg}})^*$ to find $\sigma'_N \in \mathrm{Hom}(\mathbb{Q}, k'^*)$ such that $\sigma'_N(1) = N$. Let $\tau_N = \mathrm{aut}^2_{k',G}(\sigma'_N) \circ \mathrm{aut}^1_{k',G}(\sigma_N)$. Then $\tau_N(ut^{-1}) = (\mathrm{aut}^2_{k',G}(\sigma'_N))(Nut^{-1}) = N^{-1}Nut^{-1} = ut^{-1}$, therefore $\tau_N \in \mathrm{Aut}(M/A)$. By corollary 6.3.5 (with $k' = \mathbb{Q}^{\mathbf{alg}}$), for each $F \in k' \setminus \mathbb{Q}^{\mathbf{alg}}$, $(\tau_N(F))_N = (F(N \cdot u))_N$ is an infinite family of pairwise-distinct $A$-conjugates of $F$, so $F \notin \mathrm{acl}(A)$. This concludes the proof. □

**Proposition 6.3.9.** *The definable set $X$ has infinitely many $\mathbb{Q}^{\mathbf{alg}}$-conjugates.*

*Proof.* For each $N < \omega$, let $\sigma_N \in \mathrm{Aut}(k'/\mathbb{Q}^{\mathbf{alg}})$ be defined as $\sigma_N(F(u)) = F(u + N)$, and $\tau_N = \mathrm{aut}^1_{k',G}(\sigma_N)$.[1] By lemma 6.3.3, the $(\tau_N(X))_N$ are all pairwise-distinct. □

As a result, the $A$-definable set $X$ is not $k(\mathrm{acl}(A))$-definable, which concludes our example of a Henselian valued field of residue characteristic zero where the properties on $A$-definable subsets of $\Gamma$ stated in fact 6.1.5 fail for $k$. In particular, the analogue for $k$ of corollary 6.1.6 fails, as we have $u \equiv_{k(\mathrm{acl}(A))} u + 1$, but $u \not\equiv_A u + 1$.

### 6.3.3 Elementary embeddings into Hahn fields

The content of this subsection is mostly an application of some of the material covered in the notes of van den Dries [vdD14].

**Fact 6.3.10** ([vdD14], Corollary 4.29)**.** *Let $M$ be any valued field of residue characteristic zero. Let $\sigma_1$ and $\sigma_2$ be two valued field embeddings with domain $M$, and $M_i$ a maximal immediate extension of $\sigma_i(M)$. Then there exists $\tau$ a valued field isomorphism $M_1 \longrightarrow M_2$ such that $\sigma_2 = \tau \circ \sigma_1$.*

**Proposition 6.3.11.** *Let $M$ be a Henselian valued field of residue characteristic zero with a cross-section, and $l$ a lift of the residue field. Then there exists a valued field embedding $\sigma : M \longrightarrow k(M)((t^{\Gamma(M)}))$ such that the following conditions hold:*

- *$\sigma$ is the identity on $k(M) \cup \Gamma(M)$.*
- *$\sigma$ is a $\mathcal{L}_{cs}$-embedding, i.e. for every $\alpha \in k(M)$ and for every $\gamma \in \Gamma(M)$, $\sigma(l(\alpha))$ is the constant series $\alpha$ and $\sigma(s(\gamma)) = t^\gamma$.*

[1]The map $\sigma_N$ is $\mathbb{Q}^{\mathbf{alg}}$-invariant, but it is not $A$-invariant anymore.

- $\sigma$ *is an elementary* $\mathcal{L}_{Pas}$*-embedding.*

*Proof.* For the construction of $\sigma$ and the proof that it satisfies the two first conditions from the list, we refer the reader to the discussion from the notes of van den Dries between Lemmas 4.30 and 4.31.

For the proof that $\sigma$ is $\mathcal{L}_{Pas}$-elementary, we refer the reader to section 5.5 of the notes, especially the definitions at the beginning, and theorem 5.21. It is clear that $\sigma$ is a "good map" with respect to these definitions. ☐

We actually need a construction that is a little more specific to prove lemma 6.3.13, which is a technical lemma about unary definable subsets of RV:

**Proposition 6.3.12.** *Let $G$ be an ordered Abelian group, $k'$ a field of characteristic zero, and $M$ a $\mathcal{L}_{cs}$-structure which extends $k'((t^G))$, and such that the inclusion $k'((t^G)) \longrightarrow M$ is $\mathcal{L}_{Pas}$-elementary. Then there exists $\sigma$ a $\mathcal{L}_{cs}$-embedding of Pas fields $M \longrightarrow k(M)((t^{\Gamma(M)}))$ which is $\mathcal{L}_{Pas}$-elementary, such that the following conditions hold:*

- *$\sigma$ is the identity on $\Gamma(M) \cup k(M)$.*
- *The restriction of $\sigma$ to $k'((t^G))$ is the inclusion map.*

*Proof.* The proof is written in the same spirit as in the notes of van den Dries.

Let $l$ be a lift of the residue field in $M$ extending that of $k'((t^G))$ (by proposition 6.1.8). Let $K'$ be some maximal immediate extension of $M$, and $L$ be the subfield of $M$ generated by $k'((t^G)) \cup s(\Gamma(M)) \cup l(k(M))$. We have a valued field embedding $\tau : L \longrightarrow k(M)((t^{\Gamma(M)}))$ sending each $s(\gamma)$ to $t^\gamma$, sending each $l(\alpha)$ to the constant polynomial $\alpha$, and sending each element of $k'((t^G))$ to itself. Now, the extension $L \leqslant M$ is clearly immediate, so $K'$ and $k(M)((t^{\Gamma(M)}))$ are both maximal immediate extensions of $L$. By fact 6.3.10, we have $\sigma$ a valued field isomorphism $K' \longrightarrow k(M)((t^{\Gamma(M)}))$ that extends $\tau$. By construction of $\tau$, $\sigma$ is also a Pas isomorphism. Notice now that the inclusion $\sigma(M) \longrightarrow k(M)((t^{\Gamma(M)}))$ is a "good map", hence an elementary embedding of Pas fields. As a result, $M \longrightarrow \sigma(M) \longrightarrow k(M)((t^{\Gamma(M)}))$ witnesses the statement. ☐

**Lemma 6.3.13.** *Let $M \vDash HF_{0,0}$, and $X$ a definable subset of* $\mathrm{RV}^*$. *Then there exists a finite family $\gamma_1 \dots \gamma_n \in \Gamma^*$, definable subsets $X_i \subseteq \mathrm{RV}^*$ of* $\mathrm{val}^{-1}(\gamma_i)$*, and a bound* $N < \omega$*, such that* $X_i \subseteq X$ *for all* $i$*, and* $X \setminus \left(\bigcup_i X_i\right)$ *is a union of cosets of the* $\varnothing$*-definable subgroup of* $\mathrm{RV}^*$ *corresponding to the embedding:*

$$(k^*)^N \longrightarrow k^* = \mathcal{O}^* \diagup 1 + \mathfrak{M} \longrightarrow K^* \diagup 1 + \mathfrak{M} = \mathrm{RV}^*$$

This lemma is used for technical computations in the next section. There are hints of the proof in ([HMRC18], proofs of Proposition 5.9 and Lemma 5.10, and Remark 5.14). It should also follow from the material on short exact sequences in [ACGZ22]. The techniques used in the following proof are quite similar to those used in subsection 6.3.4. We give less formal details here than in the next subsection.

*Proof.* Let $H$ be the $\varnothing$-type-definable group $\bigcap_{N<\omega} (k^*)^N$. For each $N < \omega$, let $X_N$ be the definable set:

$$\left\{x \in \mathrm{RV}^* | x \bmod (k^*)^N \subseteq X\right\}$$

Let $X_\omega = \{x \in \mathrm{RV}^* | x \bmod H \subseteq X\}$. We just have to prove that the definable set $\mathrm{val}(X \setminus X_N)$ is finite for some $N < \omega$. If not, then by compactness $\mathrm{val}(X \setminus X_\omega)$ must be infinite. Let us show that this is impossible.

We can apply proposition 6.3.11 to some $\mathcal{L}_{Pas}$-field that is an elementary extension of $M$ (which exists by fact 2.3.5) to find $M_1$ a Hahn field that is an elementary extension of $M$. Let $\mathcal{N}$ be a $|M_1|^+$-saturated $\mathcal{L}_{cs}$-elementary extension of $M_1$. By proposition 6.3.12, let $M_2$ be a Hahn field that contains $M_1$ such that there exists some elementary embedding (of valued fields) $\sigma : \mathcal{N} \longrightarrow M_2$ which restricts to the identity on $M_1$. In particular, we have:

$$M \preceq M_1 \preceq \sigma(\mathcal{N}) \preceq M_2$$

with $M_i$ a Hahn field, and $\mathcal{N}$ a sufficiently saturated model, such that the map $M_1 \longrightarrow M_2$ is compatible with the power series. By compactness and saturation, if $\mathrm{val}(X \setminus X_\omega)$ was infinite, then there would exist $x \in \mathrm{RV}^*(\sigma(\mathcal{N}))$ such that $x \in X \setminus X_\omega$ and $\mathrm{val}(x) \notin \Gamma(M_1)$. By definition of $X_\omega$, let $y$ be in $H(\sigma(\mathcal{N}))$ such that $xy \notin X$. By lemma 6.3.2 applied to the divisible group $H$, let $\tau \in \mathrm{Hom}(\mathbb{Q} \cdot \mathrm{val}(x), H)$ such that $\tau(\mathrm{val}(x)) = y$. Let $A \leqslant \mathrm{div}(\Gamma(M_2))$ such that (as groups) $\mathrm{div}(\Gamma(M_2)) = \mathrm{div}(\Gamma(M_1)) \oplus A \oplus \mathbb{Q} \cdot \mathrm{val}(x)$. Let $\tau'$ be the group homomorphism $\mathrm{div}(\Gamma(M_2)) \longrightarrow k(M_2)^*$ that extends $\tau$ which is trivial on $\mathrm{div}(\Gamma(M_1)) \oplus A$. Then $\mathrm{aut}^2_{k(M_2),\Gamma(M_2)}\left(\tau'_{|\Gamma(M_2)}\right)$ is a valued field automorphism of $M_2$ leaving $M_1$ pointwise-invariant, and sending $x$ to $xy$. As $X$ is $M_1$-definable, we must have $xy \in X$, a contradiction. $\square$

### 6.3.4 A sufficient condition for strong stable embeddedness

**Lemma 6.3.14.** *Suppose $K$ is a valued field with an angular component. Let $A \leqslant B$ be subfields of $K$ such that* $\mathrm{ac}(A) = \mathrm{res}(A)$*, and* $\mathrm{val}(A) = \mathrm{val}(B)$*. Then* $\mathrm{ac}(B) = \mathrm{res}(B)$*.*

*Proof.* Let $b \in B^*$. Then there exists $a \in A$ such that $\mathrm{val}(a) = \mathrm{val}(b)$. Let $a' \in A \cap \mathcal{O}$ be such that $\mathrm{ac}(a) = \mathrm{res}(a')$. Then $\mathrm{ac}(b) = \mathrm{res}\left(\frac{b}{a} \cdot a'\right) \in \mathrm{res}(B)$. □

**Lemma 6.3.15.** *Let $M \vDash HF_{0,0}$ be $\aleph_1$-saturated, and $A \subseteq M^{eq}$ an algebraically closed parameter set. Suppose $M$ is strongly $|A|^+$-homogeneous, and fix a cross-section on $M$ (which exists by saturation). Let $H$ be the divisible $\varnothing$-type-definable group $\bigcap_{N>0} (k^*)^N$. Suppose for all $r \in \mathrm{ac}(K(A))$, there exists $a \in \mathcal{O}_A$ such that* $\mathrm{res}(a) \bmod H = r \bmod H$*. Then every $A$-definable subset of $k$ is actually* $\mathrm{res}(K(A))$*-definable in the pure field $k(M)$.*

*Proof.* Let $X$ be an $A$-definable subset of $k$. Any element of $\mathrm{div}(\mathrm{val}(K(A)))$ is a $\mathbb{Q}$-linear combination of elements from $\mathrm{val}(K(A))$, so there exists $(\gamma_i)_i$ a family of values from $\mathrm{val}(K(A))$ that is a $\mathbb{Q}$-basis of $\mathrm{div}(\mathrm{val}(K(A)))$. In other words, $\mathrm{div}(\mathrm{val}(K(A)))$ is the direct sum of the Abelian groups $(\mathbb{Q}\gamma_i)_i$. Let $a_i \in K(A)$ such that $\mathrm{val}(a_i) = \gamma_i$, and $a_i' \in \mathcal{O}_A$ such that $\mathrm{res}(a_i') \bmod H = \mathrm{ac}(a_i) \bmod H$. Let $r_i = \mathrm{ac}\left(\frac{a_i}{a_i'}\right) \in H$.

Write for short $k' = k(M)$, $G = \Gamma(M)$, and let us identify $M$ with a $\mathcal{L}_{cs}$-substructure of $k'((t^G))$ which is $\mathcal{L}_{Pas}$-elementary by using proposition 6.3.11. Note that the Hahn series in $k'((t^G))$ that have finite support belong to $M$ by the conditions of proposition 6.3.11.

By lemma 6.1.4, and by strong-homogeneity, there exists $\sigma \in \mathrm{Aut}(M)$ such that $\sigma\left(\frac{a_i}{a_i'}\right) = r_i t^{\gamma_i}$. By lemma 6.3.2 applied to $H$, as $\mathbb{Q}\gamma_i$ is isomorphic to $\mathbb{Q}$, choose $\tau_i \in \mathrm{Hom}(\mathbb{Q}\gamma_i, k'^*)$ such that $\tau_i(\gamma_i) = r_i^{-1}$. By the universal property of the direct sum, we can find $\tau' \in \mathrm{Hom}\left(\mathrm{div}(\mathrm{val}(K(A))), k'^*\right)$ extending each of the $\tau_i$. Now, choose $B$ a $\mathbb{Q}$-vector subspace of $\mathrm{div}(G)$ such that $\mathrm{div}(G) = \mathrm{div}(\mathrm{val}(K(A))) \oplus B$, and extend $\tau'$ to $\tau'' \in \mathrm{Hom}\left(\mathrm{div}(G), k'^*\right)$ such that $\tau''_{|B} = 1$. Let $\tau = \mathrm{aut}^2_{k',G}(\tau''_{|G})$. Then $\tau$ is an automorphism[2] of $k'((t^G))$ leaving $k'$, $G$ pointwise-invariant, and sending each $r_i t^{\gamma_i}$ to $t^{\gamma_i}$. Let $L$ be the subfield of $K(A)$ generated by the $\frac{a_i}{a_i'}$. Then we clearly have $\mathrm{res}(\tau \circ \sigma(L)) =$

[2]The map $\tau$ is an automorphism of valued fields, but it is not an $\mathcal{L}_{Pas}$-automorphism.

$\mathrm{ac}(\tau\circ\sigma(L))$. By lemma 6.3.14, we have $\mathrm{res}(\tau\circ\sigma(K(A))) = \mathrm{ac}(\tau\circ\sigma(K(A)))$. By fact 6.1.5, $\tau\circ\sigma(X)$ is $\mathrm{res}(\tau\circ\sigma(K(A)))$-definable in the pure field $k'$. Now, as $M$ is an elementary substructure, $\sigma$ and hence $\tau\circ\sigma$ is an elementary map, so $X$ is $\mathrm{res}(K(A))$-definable in the pure field $k(M)$. □

**Lemma 6.3.16.** *Let $M \vDash HF_{0,0}$, and $A \subseteq M^{eq}$ such that $A = \mathrm{acl}^{eq}(K(A))$. Suppose $[k^* : (k^*)^N]$ is finite for each $N > 0$. Suppose for each $N > 0$, $\alpha \in {k^*}/{(k^*)^N}$, there exists $a \in \mathcal{O}_A$ such that $\mathrm{res}(a) \bmod (k^*)^N = \alpha$. Then any $A$-definable subset of $k$ is $\mathrm{res}(K(A))$-definable in the pure field $k(M)$.*

*Proof.* By the finiteness hypothesis, we can freely replace $M$ by an $\aleph_1$-saturated, strongly $|A|^+$-homogeneous elementary extension. Let $X$ be an $A$-definable subset of $k$, and $H = \bigcap_{N>0} (k^*)^N$. For each $N > 0$, let $\pi_N$ be the canonical projection ${k^*}/{H} \longrightarrow {k^*}/{(k^*)^N}$. Let $\Sigma$ be the partial type

$$\left\{\mathrm{res}(x_\alpha) \bmod (k^*)^N = \pi_N(\alpha) | \alpha \in {k^*(M)}/{H(M)}, N > 0\right\}$$

The partial type $\Sigma$ does not necessarily have a realization in $K(A)$, but it is finitely satisfiable in $K(A)$ by hypothesis on $A$.

Let $E_N$ be the set of every formula with parameters in $K(A)$ and variables in $(x_\alpha)_{\alpha\in k^*(M)/H(M)}$, $y_1 \dots y_N$. For each $\phi \in E_N$ and $m > 0$, let $\bar{\phi}_m$ be the formula in the variables $(x_\alpha)_{\alpha\in k^*/H}$ stating that the set:

$$\{y_1 \dots y_N \in K | \phi((x_\alpha)_\alpha, y_1 \dots y_N)\}$$

has exactly $m$ elements. Let $F_N$ be the set of every formula (in the language of the pure field $k$) without parameters, with variables $y_1 \dots y_N \in k$, $z \in k$. Let $\Pi$ be the following partial type in the variables $(x_\alpha)_\alpha$:

$$\left\{ \begin{array}{l|l} \bar{\phi}_m \Longrightarrow [ & N, m > 0 \\ \forall y_1 \dots y_N \in \mathcal{O}\ (\phi((x_\alpha)_\alpha, y_1 \dots y_N) \Longrightarrow & \phi \in E_N \\ \neg[\forall z \in k\ (z \in X \Longleftrightarrow \psi(\mathrm{res}(y_1) \dots \mathrm{res}(y_N), z))])] & \psi \in F_N \end{array} \right\}$$

The realizations of $\Pi$ are the tuples $(a_\alpha)_\alpha$ for which $X$ is not definable over $\mathrm{res}(K(\mathrm{acl}(A(a_\alpha)_\alpha)))$. The definition of $\Pi$ can be read "for all $y_i$, if $y_i$ is $A(x_\alpha)_\alpha$-algebraic via the formula $\phi$, then the formula $\psi(\mathrm{res}(y_i), z)$ does not define $X$".

Now lemma 6.3.15 implies that the type $\Sigma \cup \Pi$ is inconsistent. By compactness, as $\Sigma$ is finitely satisfiable in $K(A)$, $K(A)$ has a tuple that is not a realization of $\Pi$. This concludes the proof. □

**Corollary 6.3.17.** *Let $M \vDash HF_{0,0}$, and $A \subseteq M$, such that $A = \mathrm{acl}(K(A))$. Suppose $[k^* : (k^*)^N]$ is finite for every $N > 0$, and $A$ satisfies the conditions of lemma 6.3.16. Then, for all $r$, $s \in k(M)$, if $r \equiv_{k(A)} s$ in the field $k(M)$, then we have $r \equiv_A s$ in the valued field $M$.*

Lemma 6.3.16 also allows us to prove an analogue of corollary 6.2.4. It will be useful for one of our main results (theorem 6.4.9).

**Corollary 6.3.18.** *With the same hypothesis as in lemma 6.3.16, suppose additionally:*

- *For every $\gamma \in \Gamma(\mathrm{dcl}(\varnothing))$, there exists $a \in K(A)$ such that $\mathrm{val}(a) = \gamma$.*
- *The field $k(M)$ is algebraically bounded.*

*then, for any $\alpha \in \mathrm{RV}(A)$, there exists $a \in K(A)$ such that $\mathrm{rv}(a) = \alpha$.*

*Proof.* Let $\alpha \in \mathrm{RV}(A)$. By lemma 6.2.9, any value of $\Gamma(A)$ can be pulled back to $K(A)$, so we can assume by scaling everything that $\mathrm{val}(\alpha) = 0$, i.e. $\alpha \in k^*$. Now, $\alpha \in \mathrm{acl}^{eq}(K(A))$, so by lemma 6.3.16 we actually have $\alpha \in \mathrm{acl}(\mathrm{res}(K(A)))$ in the pure field structure of $k(M)$. By hypothesis, $\alpha \in \mathrm{res}(K(A))^{\mathbf{alg}} \cap k(M)$. Therefore, we have a polynomial $P$ with non-zero coefficients in $\mathcal{O}_M \cap K(A)$ such that $\alpha$ is a root of $\mathrm{res}(P)$. If we choose such a $P$ as a pullback of the minimal polynomial of $\alpha$ over $\mathrm{res}(K(A))$, then, by Henselianity and as we are in residue characteristic zero, there exists a root of $P$ in $K(M)$ having residue $\alpha$. As $P$ is chosen with coefficients in $K(A)$, and $A = \mathrm{acl}^{eq}(K(A))$, this root is in $K(A)$, which concludes the proof. □

**Corollary 6.3.19.** *With the same hypothesis as corollary 6.3.18, every $A$-ball is pointed.*

*Proof.* Let $Y$ be an $A$-ball. If $Y$ is not open, then $\Gamma$ is dense, so we can apply corollary 6.2.4 to conclude. Else, by using proposition 6.1.1 with a reasoning similar to corollary 6.2.4, we can prove that $X$, the minimal closed ball strictly containing $Y$, is pointed. Let $a \in K(A) \cap X$. We can use corollary 6.3.18 to find $a' \in K(A)$ such that $\mathrm{rv}(a') = \mathrm{rv}(Y - a)$, which concludes the proof as $a' + a \in Y$. □

## 6.4 Finding extension bases

In this section, we fix $\kappa > \aleph_0$, $M$ a $\kappa$-saturated and strongly $\kappa$-homogeneous Henselian valued field of residue characteristic zero.

In order to study extension bases, we use the argument of corollary 6.0.5 to reduce to the case of unary types. We start this section by giving sufficient conditions for a global unary type to be non-forking over some parameter set $A$. We have three distinct cases that correspond to the three following lemmas. The cases depend on the nature of the chain of $A$-balls of which our unary type is $A$-generic, as well as whether $A$ is generated by parameters from $K$.

The following lemma deals with the case of a non-residual chain, with $A$ not necessarily generated by field parameters.

**Lemma 6.4.1.** *Let $A$ be a small subset of $M^{eq}$, $\mathfrak{B}$ a chain of $A$-balls, $B$ a $(2^{|A|})^{++}$-saturated, strongly $(2^{|A|})^{++}$-homogeneous small elementary substructure of $M^{eq}$ that contains $A$, and $c \in K(M)$. Assume the following conditions hold:*

- $A = \mathrm{acl}^{eq}(A)$.
- $\mathfrak{B}$ *is not residual.*
- $c$ *is $B$-generic of* $\mathfrak{B}$.
- *There exists $b \in K(B) \cap (\cap \mathfrak{B})$ such that* $\mathrm{val}(c-b) \underset{\Gamma(A)}{\downarrow^{\mathbf{f}}} \Gamma(B)$ *in the ordered group structure of* $\Gamma(M)$.
- *For every $N < \omega$, $[k^* : (k^*)^N]$ is finite.*

*then we have* $c \underset{A}{\downarrow^{\mathbf{f}}} B$.

*Proof.* First of all, let us notice that the smallest closed ball that contains $b$ and $c$ is contained in $\cap\mathfrak{B}$, and this inclusion is strict, otherwise this closed ball would be the least element of $\mathfrak{B}$, contradicting the hypothesis. By genericity of $c$, this ball does not belong to $B$, but it has a point $b$ in $B$, so its radius $\mathrm{val}(c-b)$ must not be in $\Gamma(B)$.

Suppose by contradiction that $c \underset{A}{\not\downarrow^{\mathbf{f}}} B$. Then, by fact 1.1.19, we have $c \underset{A}{\not\downarrow^{\mathbf{d}}} B$. By proposition 6.2.6 (since $B = \mathrm{acl}^{eq}(K(B))$), the type of $c$ over $B$ is entirely controlled by that of $\mathrm{rv}(c-b)$. As a result, a $B$-definable set containing $c$ which divides over $A$ may be written $X' = \{x \in K | \mathrm{rv}(x-b) \in X\}$, with $X$ some unary $B$-definable subset of RV containing $\mathrm{rv}(c-b)$. By lemma 1.1.9, there exists $b' \in B$ such that $b \equiv_A b'$ and $X'$ divides over $Ab'$. We have $\mathrm{val}(c-b) > \mathrm{val}(b-b')$, so $\mathrm{rv}(c-b) = \mathrm{rv}(c-b')$. Thus, by replacing $X$ by $X \cap \mathrm{val}^{-1}(]-\mathrm{val}(b-b'), +\infty[) \ni \mathrm{rv}(c-b)$, we have $\mathrm{rv}(x-b) \in X$ if and only if $\mathrm{rv}(x-b') \in X$ for all $x \in K$, and we can suppose $b' = b$ without loss of generality. One can notice that a witness for division of $X'$ over $Ab'$ is also a witness for division of the translate $X' - b'$ over $A$, so we can suppose $b = b' = 0$ without loss of generality (so $Ab' = A$). By proposition 1.1.7, we can deduce that $X$ divides over $A$, as $X'$ is the preimage of $X$ under the $A$-definable function rv. Now, let $(\gamma_i)_i$, $(X_i)_i$, $N < \omega$ be a witness of lemma 6.3.13 applied to $X$. The family $(\gamma_i)_i$ is finite, so we can assume that the $\gamma_i$ are algebraic over $B$, i.e. they belong to $\Gamma(B)$, so they are all distinct from $\mathrm{val}(c)$. As a result, we can replace $X$ by $X \setminus \left(\bigcup_i X_i\right)$ without loss of generality. As $X$ is invariant under multiplication by $(k^*)^N$, $X$ is the preimage of $\bar{X} = X \bmod (k^*)^N$, so $\bar{X}$ divides over $A$. The $A$-definable group homomorphism $\mathrm{val} : {\mathrm{RV}}/{(k^*)^N} \longrightarrow \Gamma$ has a finite kernel ${k^*}/{(k^*)^N}$. We can apply lemma 1.1.8 to the relation:

$$R = \left\{(x,y) \in \Gamma \times {\mathrm{RV}}/{(k^*)^N} | \mathrm{val}(y) = x\right\}$$

to show that $\mathrm{val}(\bar{X})$ divides over $A$. Now we are done: take a witness for division $(\sigma_n)_n$ of $\mathrm{val}(\bar{X})$ over $A$, $(\sigma_{n|\Gamma(M)})_n$ is a witness for division of $\mathrm{val}(\bar{X})$ over $\Gamma(A)$ in the ordered group $\Gamma(M)$ ; we know that $\mathrm{val}(\bar{X})$ is $\Gamma(B)$-definable by fact 6.1.5, and $\mathrm{val}(c) = \mathrm{val}(c-b) \in \mathrm{val}(\bar{X})$, so $\mathrm{val}(c-b) \not\mathop{\smile}^{\mathbf{d}}_{\Gamma(A)} \Gamma(B)$, a contradiction. □

*Remark* 6.4.2. In lemma 6.4.1, the notions of $B$-genericity and $B$-weak genericity are the same: $B$ is a model, so every $B$-ball is pointed.

The following lemma deals with the case of a residual chain, with $A$ not necessarily generated by field parameters.

**Lemma 6.4.3.** *Let $A$ be a small subset of $M^{eq}$, $Z$ a closed $A$-ball, $B$ a $(2^{|A|})^{++}$-saturated, strongly $(2^{|A|})^{++}$-homogeneous small elementary substructure of $M^{eq}$ that contains $A$, and $c \in K(M)$. Assume the following conditions hold:*

- $A = \mathrm{acl}^{eq}(A)$.
- *$c$ is $B$-generic of $\{Z\}$.*
- *In the field $k(M)$, we have $\mathop{\smile\!\!\!\!\!\mid}^{\mathbf{f}} = \mathop{\smile\!\!\!\!\!\mid}^{\mathbf{alg}}$.*

*then we have $c \mathop{\smile\!\!\!\!\!\mid}^{\mathbf{f}}_{A} B$.*

*Proof.* Suppose towards contradiction that $c \not\mathop{\smile\!\!\!\!\!\mid}^{\mathbf{d}}_{A} B$. Let $X$, $X'$ be as in the proof of lemma 6.4.1. We can assume $X' \subseteq Z$. By lemma 1.1.9, there exist $b$, $b' \in B$ such that $b \in Z$, $\mathrm{val}(b') = \mathrm{rad}(Z)$, and $X'$ divides over $Abb'$. Just like in the previous proof, a witness for division of $X'$ over $Abb'$ is a witness for division of $\frac{X'-b}{b'}$ over $A$, so we can assume $b = 0$, $b' = 1$, and $Z = \mathcal{O}$. By genericity of $c$, and by proposition 6.2.6, we know that, for all $c' \in K(M)$, if $\mathrm{res}(c') = \mathrm{res}(c)$, then $c \equiv_B c'$. As a result, the set:

$$\{x \in X' | \forall y \in \mathcal{O}\ \ (\mathrm{res}(x) = \mathrm{res}(y) \implies y \in X')\}$$

is a $B$-definable subset of $X'$ containing $c$, so it still divides over $A$, and we can assume that $X'$ coincides with this set. Now $X'$ is the preimage of $\mathrm{res}(X')$, so by proposition 1.1.7 $\mathrm{res}(X')$ divides over $A$. By fact 6.1.5 we have $\mathrm{res}(c) \not\mathop{\smile\!\!\!\!\!\mid}^{\mathbf{f}}_{k(A)} k(B)$ in the field $k$, so by hypothesis $\mathrm{res}(c) \not\mathop{\smile\!\!\!\!\!\mid}^{\mathbf{alg}}_{k(A)} k(B)$, and $\mathrm{res}(c)$ is not transcendental over $k(B)$. This contradicts the genericity of $c$. □

In the next lemma, we also deal with the residual case, but without the assumption $\mathop{\smile\!\!\!\!\!\mid}^{\mathbf{f}} = \mathop{\smile\!\!\!\!\!\mid}^{\mathbf{alg}}$ in $k(M)$. However, we have to strengthen our hypothesis on $A$ to make it work, and suppose that $A$ is generated by field elements.

**Lemma 6.4.4.** *Let $A$ be a small subset of $M^{eq}$, $Z$ a closed $A$-ball, $B$ a $(2^{|A|})^{++}$-saturated, strongly $(2^{|A|})^{++}$-homogeneous small elementary substructure of $M^{eq}$ that contains $A$, and $c \in K(M)$. Assume the following conditions hold:*

- $A = \mathrm{acl}^{eq}(K(A))$.
- *$c$ is $B$-generic of $\{Z\}$.*
- *There exist $b \in Z(A)$ and $b' \in K(A)$ such that* $\mathrm{val}(b') = \mathrm{rad}(Z)$*, and, in the pure field $k(M)$, we have* $\mathrm{res}\left(\frac{c-b}{b'}\right) \mathop{\smile\!\!\!\!\!\mid}^{\mathbf{f}}_{k(A)} k(B)$.

*then we have* $c \underset{A}{\downarrow}^{\mathbf{f}} B$.

*Proof.* Let $c' = \frac{c-b}{b'}$. It is clear that $c'$ is $B$-generic of $\{\mathcal{O}\}$. As we have $\operatorname{res}\left(\frac{c-b}{b'}\right) \underset{k(A)}{\downarrow}^{\mathbf{f}} k(B)$ in the field $k(M)$, we have $c' \underset{A}{\downarrow}^{\mathbf{f}} B$ with a reasoning similar to the previous lemma. Now, this time $b$, $b' \in K(A)$, so $c \in \operatorname{acl}(Ac')$, therefore $c \underset{A}{\downarrow}^{\mathbf{f}} B$ by corollary 1.1.22. □

Now that we have seen several sufficient conditions for a global unary type to be non-forking over $A$, we need to be able to build these types, i.e. realize these conditions. It will be easy in some cases, and the construction will be described quickly in the proofs of our theorems. The two following lemmas deal with cases where a more technical approach is required. The next lemma will be used to build a non-forking global extension of the type of an $A$-generic point of an $A$-immediate chain. The second lemma does the same for the generic of an $A$-ramified chain, but there we will have to assume that $A$ is generated by field elements.

*Remark* 6.4.5. In the residual case, the non-forking global extension is easy to build when $\downarrow^{\mathbf{f}} = \downarrow^{\mathbf{alg}}$ in $k(M)$: lemma 6.4.3 merely requires $c$ to be $B$-generic of the correct chain. One can note that this lemma holds with a weaker hypothesis: instead of requiring $\downarrow^{\mathbf{f}} = \downarrow^{\mathbf{alg}}$ in $k(M)$, we may just assume that for all elements $b$, $b' \in K(B)$, if $b \in Z$ and $\operatorname{val}(b') = \operatorname{rad}(Z)$, then we have $\operatorname{res}\left(\frac{c-b}{b'}\right) \underset{k(A)}{\downarrow}^{\mathbf{f}} k(B)$ in the field structure of $k(M)$. However, without the hypothesis $\downarrow^{\mathbf{f}} = \downarrow^{\mathbf{alg}}$ in $k(M)$, it is not clear whether the existence of the non-forking global extension can be established, because a random $B$-generic point of $\{Z\}$ might not satisfy these new conditions. Maybe such a point still exists, this is one of the interesting questions that naturally follow from this paper.

**Lemma 6.4.6.** *Let $A$ be a small subset of $M^{eq}$. Let $\mathfrak{B}$ be an non-residual chain of $A$-balls. Let $B$ be a $(2^{|A|})^{++}$-saturated, strongly $(2^{|A|})^{++}$-homogeneous small elementary substructure of $M^{eq}$ that contains $A$, and $b \in K(B)$ such that $b \in \cap\mathfrak{B}$. Suppose the following conditions hold:*

- $A = \operatorname{acl}^{eq}(A)$.
- *In the ordered group $\Gamma(M)$, $\Gamma(A)$ is an extension base.*

*then there exists* $c \in K(M)$ *such that* $c$ *is* $B$*-generic of* $\mathfrak{B}$*, and such that* $\mathrm{val}(c-b) \underset{\Gamma(A)}{\downarrow^{\mathbf{f}}} \Gamma(B)$ *in the ordered group structure.*

*Proof.* Define:

$$\bar{A} = \{\mathrm{rad}(X) | X \in \mathfrak{B}\} \subseteq \Gamma(A)$$

$$\bar{B} = \left\{\delta \in \Gamma(B) | \forall \gamma \in \bar{A}\ \gamma > \delta\right\}$$

Let $X = \bigcap_{\gamma \in \bar{A}, \delta \in \bar{B}} ]\delta, \gamma[$, $(\sigma_n)_n \in \mathrm{Aut}(\Gamma(M)/\Gamma(A))^\omega$, and $\bar{B}' = \bigcup_n \sigma_n(\bar{B})$. Then, one may easily see by compactness that the type-definable set $Y = \bigcap_{\gamma \in \bar{A}, \delta \in \bar{B}'} ]\delta, \gamma[$ is non-empty. However, we have $Y = \bigcap_n \sigma_n(X)$, so $(\sigma_n)_n$ cannot be a witness for division of (any $\Gamma(B)$-definable set containing) $X$ over $\Gamma(A)$. Therefore, $X$ does not divide over $\Gamma(A)$ in the ordered group $\Gamma(M)$. In this ordered group, the induced theory is dependent ([GS84]), and $\Gamma(A)$ is an extension base by hypothesis, so forking over $\Gamma(A)$ coincides with dividing over $\Gamma(A)$ ([CK12], Theorem 1.2). As a result, $X$ does not fork over $\Gamma(A)$ in $\Gamma(M)$. By corollary 1.1.16, there exists $\gamma \in \Gamma(M)$ such that $\gamma \in X$, and $\gamma \underset{\Gamma(A)}{\downarrow^{\mathbf{f}}} \Gamma(B)$.

Now we are done: a witness of the statement is any $c \in K(M)$ such that $\mathrm{val}(c-b) = \gamma$. Such a point $c$ must be $B$-generic of $\mathfrak{B}$: $\mathrm{val}(c-b) < \bar{A}$ implies $c \in \cap\mathfrak{B}$, and $\mathrm{val}(c-b) > \bar{B}$ implies that $c$ does not belong to any smaller $B$-ball. □

**Lemma 6.4.7.** *Let* $A$ *be a small subset of* $M^{eq}$*. Let* $\mathfrak{B}$ *be an* $A$*-ramified chain of* $A$*-balls. Let* $B$ *be a* $(2^{|A|})^{++}$*-saturated, strongly* $(2^{|A|})^{++}$*-homogeneous small elementary substructure of* $M^{eq}$ *that contains* $A$*,* $c \in K(M)$*, and* $b \in K(A)$ *such that* $b \in \cap\mathfrak{B}$*. Suppose the following conditions hold:*

- $A = \mathrm{acl}^{eq}(K(A))$.
- *$c$ is $A$-generic of* $\mathfrak{B}$.
- *In the ordered group* $\Gamma(M)$*,* $\Gamma(A)$ *is an extension base.*

*then there exists* $c' \equiv_A c$ *such that* $c'$ *is* $B$*-generic of* $\mathfrak{B}$*, and such that, in the reduct to the value group, we have* $\mathrm{val}(c'-b) \underset{\Gamma(A)}{\downarrow^{\mathbf{f}}} \Gamma(B)$.

*Proof.* Let $\bar{A}$, $\bar{B}$, $X$, $(\sigma_n)_n$, $\bar{B}'$, $Y$ be as in the proof of the previous lemma. Let $B'$ be a small $|B|^+$-saturated elementary substructure of $M$ containing $B \cup \bar{B}'$. By proposition 6.2.8, there exists $d \equiv_A c$ such that $d$ is $B'$-generic of $\mathfrak{B}$. The value $\mathrm{val}(d-b)$ not only is a point of $Y$, put it is also a $\Gamma(A)$-conjugate of $\mathrm{val}(c-b)$. As a result, the partial type (of the ordered group) $\mathrm{tp}(\mathrm{val}(c-b)/\Gamma(A)) \cup \{x \in X\}$ does not divide, and hence does not fork over $\Gamma(A)$ by the third hypothesis of the list, with the same argument as in lemma 6.4.6. Let $\gamma \in \Gamma(M)$ be a realization of this type such that we have $\gamma \mathop{\downarrow}^{\mathbf{f}}_{\Gamma(A)} \Gamma(B)$. Then $\gamma \equiv_{\Gamma(A)} \mathrm{val}(c-b)$, so $\gamma \equiv_A \mathrm{val}(c-b)$ by corollary 6.1.6. Let $\sigma \in \mathrm{Aut}(M/A)$ such that $\sigma(\mathrm{val}(c-b)) = \gamma$. As $b \in K(A)$, $\sigma(b) = b$, so $c' = \sigma(c)$ witnesses the statement. □

We can finally prove the main results of this chapter:

**Theorem 6.4.8.** *Let $A$ be a small subset of $K(M)$. Suppose the following conditions hold:*

- *Every subset of $\Gamma(M)$ is an extension base with respect to the ordered group structure.*
- *$\Gamma(M)$ is dense.*
- *Every subset of $k(M)$ is an extension base with respect to the field structure.*
- *For every $N < \omega$, $[k^* : (k^*)^N]$ is finite.*
- *For every $\gamma \in \Gamma(\mathrm{dcl}(\varnothing))$, there exists an $a \in K(\mathrm{acl}(A))$ such that $\mathrm{val}(a) = \gamma$.*
- *For every $N > 0$, $\alpha \in {}^{k^*}\!/_{(k^*)^N}(M)$, there exists $a \in \mathcal{O}^*_{\mathrm{acl}(A)}$ such that $\mathrm{res}(a) \bmod (k^*)^N = \alpha$.*

*then $A$ is an extension base.*

*Proof.* Let $\mathfrak{C}$ be the class of every small subsets $A'$ of $K(M)$ satisfying the two last points of the theorem. Clearly, for all $A' \in \mathfrak{C}$, for all $c \in K(M)$, we have $A'c \in \mathfrak{C}$.[3] By corollary 6.0.5, we just have to show that $c \mathop{\downfootline}^{\mathbf{f}}_{A} A$ for every singleton $c \in K(M)$.

Let $c \in K(M)$. As $c \mathop{\downfootline}^{\mathbf{f}}_{A} A$ if and only if $c \mathop{\downfootline}^{\mathbf{f}}_{\mathrm{acl}^{eq}(A)} \mathrm{acl}^{eq}(A)$, we can replace $A$ by $\mathrm{acl}^{eq}(A)$. Let $B$ be a $(2^{|A|})^{++}$-saturated, strongly $(2^{|A|})^{++}$-homogeneous small elementary substructure of $M^{eq}$ that contains $A$. We have to find $c' \equiv_A c$ for which $c' \mathop{\downfootline}^{\mathbf{f}}_{A} B$. Let $\mathfrak{B}$ be the set of every $A$-ball containing $c$.

- Suppose $\mathfrak{B}$ is $A$-immediate. Let $c'$ witness lemma 6.4.6. We have $c' \mathop{\downfootline}^{\mathbf{f}}_{A} B$ by lemma 6.4.1. Let us show $c \equiv_A c'$. Let $\mathfrak{B}'$ be the chain of every pointed $A$-ball that contains $c$. If $\cap\mathfrak{B}'$ has no point in $K(A)$, then we have $c \equiv_A c'$ by proposition 6.2.6, as $c'$ is clearly in $\cap\mathfrak{B}'$. Else, if $a \in K(A)$ is in $\cap\mathfrak{B}'$, then $\mathfrak{B}'$ is not $A$-immediate, so $\cap\mathfrak{B}' \neq \cap\mathfrak{B}$, and $c$ is not $A$-generic of $\mathfrak{B}'$. By proposition 6.2.3, $\mathfrak{B}$ has a least element $Y$, an open ball which does not contain $a$. As $c$ and $c'$ are both in $Y$, we clearly have $\mathrm{rv}(c-a) = \mathrm{rv}(c'-a)$, which also implies $c \equiv_A c'$ by proposition 6.2.6.

- Suppose $\mathfrak{B}$ is $A$-ramified, i.e. $\cap\mathfrak{B}$ has a point $b \in K(A)$. Let $c'$ witness lemma 6.4.7 applied to $b$. We have $c' \mathop{\downfootline}^{\mathbf{f}}_{A} B$ by lemma 6.4.1.

- Suppose $\mathfrak{B}$ is residual, let $X = \min(\mathfrak{B})$. Then $X$ is pointed by corollary 6.2.4, as $\Gamma$ is dense. Let $b \in K(A)$ such that $b \in X$. By lemma 6.2.9, there also exists $b' \in K(A)$ such that $\mathrm{val}(b') = \mathrm{rad}(X)$. As $k(A)$ is an extension base in the residue field, we can use corollary 1.1.16 to find $\alpha \equiv_{k(A)} \mathrm{res}\left(\frac{c-b}{b'}\right)$ such that $\alpha \mathop{\downfootline}^{\mathbf{f}}_{k(A)} k(B)$. By corollary 6.3.17, we have $\alpha \equiv_A \mathrm{res}\left(\frac{c-b}{b'}\right)$, so there exists $c' \equiv_A c$ for which $\mathrm{res}\left(\frac{c'-b}{b'}\right) \mathop{\downfootline}^{\mathbf{f}}_{k(A)} k(B)$. By lemma 6.4.4, we have $c' \mathop{\downfootline}^{\mathbf{f}}_{A} B$. □

**Theorem 6.4.9.** *Theorem 6.4.8 also holds if we replace the hypothesis on the density of* $\Gamma$ *by the the fact that the theory of the field* $k$ *is algebraically bounded.*

[3]There is a subtlety here. The $k^* / (k^*)^N$ are fixed. If they were infinite (hence unbounded), maybe adding $c$ to $A'$ would add new classes of $(k^*)^N$ that cannot be pulled back to $K(A'c)$! So the finiteness hypothesis is crucial for the induction argument of corollary 6.0.5.

*Proof.* The only difference in the proof is in the residual case, where we have to use corollary 6.3.19 instead of corollary 6.2.4 to show that $X$ is pointed. □

**Theorem 6.4.10.** *Let $A$ be a small subset of $K(M)$. Suppose the following conditions hold:*

- *Every subset of $\Gamma(M)$ is an extension base with respect to the ordered group structure.*
- $\mathop{\smile\!\!\!\!\!\!\mid}^{\mathbf{f}} = \mathop{\smile\!\!\!\!\!\!\mid}^{\mathbf{alg}}$ *in the field $k(M)$.*
- *For every $N < \omega$, $[k^* : (k^*)^N]$ is finite.*

*then $A$ is an extension base.*

*Proof.* We do the same proof as above with $\mathfrak{C}$ the class of every small subset of $K(M)$. The only case where we have to act differently is the residual case. Here, we use proposition 6.2.8 to find $c' \equiv_A c$ which is $B$-generic of $\mathfrak{B}$. Then $c' \mathop{\smile\!\!\!\!\!\!\mid}^{\mathbf{f}}_A B$ by lemma 6.4.3. □

**Theorem 6.4.11.** *Let $A$ be a small subset of $M^{eq}$. Suppose the following conditions hold:*

- $\Gamma(M) \equiv \mathbb{Z}$.
- $\mathop{\smile\!\!\!\!\!\!\mid}^{\mathbf{f}} = \mathop{\smile\!\!\!\!\!\!\mid}^{\mathbf{alg}}$ *in the field $k(M)$.*
- *For every $N < \omega$, $[k^* : (k^*)^N]$ is finite.*

*then $A$ is an extension base.*

*Proof.* This time $\mathfrak{C}$ is the class of every small subset of $M^{eq}$. We have the same disjunction of three cases:

- $\mathfrak{B}$ is $A$-immediate. We proceed just like the proof of theorem 6.4.8, except that we have to use proposition 6.2.7 instead of proposition 6.2.6 (because $A$ is not necessarily generated by field elements) to prove $c' \equiv_A c$.

- $\mathfrak{B}$ is $A$-ramified, i.e. $\cap\mathfrak{B}$ contains an $A$-ball $X$. By proposition 6.2.8, let $c' \equiv_A c$ a $B$-generic point of $\mathfrak{B}$. By genericity of $c'$, we can apply theorem 4.4.1 to show that $\mathrm{val}(c - X) \mathop{\downarrow}^{\mathbf{f}}_{\Gamma(A)} \Gamma(B)$ in the ordered group $\Gamma(M)$. Then we can apply lemma 6.4.1 (with $b$ any point in $K(B) \cap X$) to get $c \mathop{\downarrow}^{\mathbf{f}}_{A} B$.

- $\mathfrak{B}$ is residual. This time, we proceed like the proof of theorem 6.4.10.

□

**Corollary 6.4.12.** *Forking coincides with dividing in (any structure interpretable in)* $\mathrm{PL}_0$.

*Proof.* By theorem 6.4.11, every set in any structure interpretable in $\mathrm{PL}_0$ is an extension base. As the completions of $\mathrm{PL}_0$ are $\mathrm{NTP}_2$, it follows from ([CK12], Theorem 1.2) that forking coincides with dividing. □

# Chapter 7

# Separated extensions

Most results in the literature regarding forking in valued fields take place in separated extensions (see definition 2.2.11), as this allows minimal interference from immediate extensions. Haskell-Hrushovski-Macpherson initially proved domination results in ACVF for those extensions in ([HHM05]), with applications to stable domination, and some applications to forking. More recently, Ealy-Haskell-Simon ([EHS23]) and Vicaria ([Vic21]) proved domination results for general Henselian valued fields of residue characteristic zero. We mostly focus on the work of Ealy-Haskell-Simon. The main difference between their results and the one of Vicaria is that they ask that the quotients $\mathrm{RV}^* \diagup (\mathrm{RV}^*)^N$ are finite and lift to field elements in the parameter sets for any $N > 0$, while Vicaria does not make that assumption, but works instead in a language where there are sorts for the cosets of those quotients. Our results may work in the context of Vicaria's paper, though one would have to adapt the contents of the previous chapter to this new language to make things work. Our goal in this chapter is to discuss the aforementioned results, we notably state equivalent hypothesis and conclusions which are easier to comprehend. Additionally, we make a contribution of our own in section 7.3, where we relate forking to invariant Keisler measures.

**Fact 7.0.1.** *For every $N < \omega$, we have $[K^* : (K^*)^N] \geqslant [\mathrm{RV}^* : (\mathrm{RV}^*)^N] \geqslant [k^* : (k^*)^N]$.*

**Assumptions 7.0.2.** Let $M \vDash \mathrm{HVF}_{0,0}$, and $B \geqslant A \leqslant C$ valued subfields of $K(M)$, such that $M$ is $|BC|^+$-saturated, strongly $|BC|^+$-homogeneous. We assume that the index $[K^* : (K^*)^N]$ is finite for every $N > 0$. In particular, by fact 7.0.1, $[\mathrm{RV}^* : (\mathrm{RV}^*)^N]$ and $[k^* : (k^*)^N]$ are finite. Furthermore, we name boundedly many (at most $2^{\aleph_0}$) constants:

- Assume that for every $N > 0$, for every coset $X$ of $(\mathrm{RV}^*)^N$ in $\mathrm{RV}^*$, there exists some $a \in A^*$ such that $\mathrm{rv}(a) \in X$.
- Likewise, assume that for every coset $X$ of $(k^*)^N$ in $k^*$, there exists some $a \in A \cap \mathcal{O}^*$ such that $\mathrm{res}(a) \in X$.
- Lastly, assume that for every $\gamma \in \Gamma(\mathrm{dcl}^{eq}(\varnothing))$, there exists some $a \in A$ such that $\mathrm{val}(a) = \gamma$.

**Proposition 7.0.3.** *Under assumptions 7.0.2, the following hold for any field $D \leqslant K(M)$ extending $A$:*

1. *Every $D$-definable subset of $k$ is $\mathrm{res}(D)$-definable in the field $k(M)$.*
2. *Every $D$-definable subset of $\Gamma$ is $\mathrm{val}(D)$-definable in $\Gamma(M)$.*
3. *For every $\gamma \in \Gamma(\mathrm{acl}^{eq}(D))$, there exists some $d \in K(\mathrm{acl}(D))$ such that $\mathrm{val}(d) = \gamma$.*
4. *If $k$ is algebraically bounded, then, for every $r \in k(\mathrm{acl}^{eq}(D))$, there exists $d \in K(\mathrm{acl}(D))$ such that $\mathrm{res}(D) = r$.*
5. *If $k$ is algebraically bounded, then every $\mathrm{acl}^{eq}(D)$-ball has a point in $K(\mathrm{acl}(D))$.*

Note that lemma 6.3.16 is not enough to prove this proposition, as it relies on fact 6.1.5, which requires the parameter set $A$ to be algebraically closed. It turns out that this requirement is not necessary for several of the conditions.

*Proof.* First of all, assumptions 7.0.2 hold for any extension of $A$, with the exception of saturation and strong homogeneity. As the conditions of the statement are insensitive to the choice of the ambient model, those conditions of saturation and strong homogeneity do not matter, hence we may assume $D = A$.

Condition 3 follows from lemma 6.2.9, condition 4 follows from corollary 6.3.18, and condition 5 follows from corollary 6.3.19. As for conditions 1 and 2, we refer the reader to ([EHS23], Remark 1.4). □

*Remark* 7.0.4. The book ([vdD14], Corollary 5.24) shows the well-known fact that $k$ and $\Gamma$ are *orthogonal* sorts, i.e. for $F \subseteq k$, $G \subseteq \Gamma$, $\mathrm{tp}(F/A)$ and $\mathrm{tp}(G/A)$ are weakly orthogonal.

Note that in pseudo-local fields of residue characteristic zero ($\mathrm{PL}_0$), the finiteness assumptions from assumptions 7.0.2 hold, and $k$ is algebraically bounded.

We make one last observation which follows from proposition 7.0.3. By algebraic boundedness of $K$, the extension $A \leqslant K(\mathrm{acl}(A))$ is algebraic, hence $\mathrm{val}(A) \leqslant \mathrm{val}(K(\mathrm{acl}(A))) = \Gamma(\mathrm{acl}^{eq}(A))$ is pure under assumptions 7.0.2. Moreover, if $k$ is algebraically bounded, then $\mathrm{res}(A) \leqslant \mathrm{res}(K(\mathrm{acl}(A))) = k(\mathrm{acl}^{eq}(A))$ is algebraic under 7.0.2.

We use in this chapter the following domination result from Ealy-Haskell-Simon:

**Theorem 7.0.5** ([EHS23], Corollary 3.2)**.** *Under hypothesis slightly weaker than assumptions 7.0.2 (see their Subsection 1.2), suppose* $\mathrm{res}(C)$ *and* $\mathrm{res}(B)$ *are linearly disjoint over* $\mathrm{res}(A)$*, and* $\mathrm{val}(C) \cap \mathrm{val}(B) = \mathrm{val}(A)$*. Suppose the extension $C \geqslant A$ is separated. Then* $\mathrm{tp}(B/A, \mathrm{res}(C)\mathrm{val}(C)) \vDash \mathrm{tp}(B/C)$*.*

They relate this result to forking (in fact, to dividing) in their Corollary 3.4. We can use an abstract argument to relate type implication and dividing:

**Lemma 7.0.6.** *Let $A, B, C, D$ be parameter sets in any first-order structure, such that* $\mathrm{tp}(B/AC) \vDash \mathrm{tp}(B/ACD)$*. Then we have $CD \mathop{\downarrow}^{\mathbf{d}}_A B$ if and only if $C \mathop{\downarrow}^{\mathbf{d}}_A B$.*

*Proof.* The left-to-right direction is trivial.

Conversely, suppose $C \mathop{\downarrow}^{\mathbf{d}}_A B$. Let $(B_i)_i$ be an $A$-indiscernible sequence to which $B$ belongs. By fact 1.1.5, there exists $(B'_i)_i \equiv_{AB} (B_i)_i$ which is $AC$-indiscernible. Let $X_B$ be an $AB$-definable set to which $CD$ belongs. Then $\mathrm{tp}(B/AC) \vDash CD \in X_B$, therefore $CD$ belongs to $\bigcap_i X_{B'_i}$. It follows that $\bigcap_i X_{B_i}$ is consistent, concluding the proof. □

**Corollary 7.0.7.** *With the hypothesis of theorem 7.0.5, we have $C \mathop{\downarrow}^{\mathbf{d}}_A B$ if and only if* $\mathrm{res}(C)\mathrm{val}(C) \mathop{\downarrow}^{\mathbf{d}}_A B$*.*

Just like in our theorem 5.4.3, the problem is reduced to forking/dividing in the full structure $M$. In the next two sections, we show how to reduce the problem to forking/dividing in the reducts $k(M)$ and $\Gamma(M)$. The part about forking will require algebraic boundedness of the residue field.

## 7.1 Dividing

We adopt assumptions 7.0.2.

**Lemma 7.1.1.** *The following conditions are equivalent:*

- $\mathrm{res}(C)\mathrm{val}(C) \mathop{\downarrow}^{\mathbf{d}}_{A} B$
- $\mathrm{res}(C) \mathop{\downarrow}^{\mathbf{d}}_{A} \mathrm{res}(B)$ *and* $\mathrm{val}(C) \mathop{\downarrow}^{\mathbf{d}}_{A} \mathrm{val}(B)$

*Proof.* The top-to-bottom direction follows from the definition of dividing.

Conversely, suppose we have $\mathrm{res}(C) \mathop{\downarrow}^{\mathbf{d}}_{A} \mathrm{res}(B)$ and $\mathrm{val}(C) \mathop{\downarrow}^{\mathbf{d}}_{A} \mathrm{val}(B)$. Let $X$ be a $B$-definable set containing (some enumeration of) $\mathrm{res}(C)\mathrm{val}(C)$. By orthogonality, we may assume that $X$ can be written $X_k \times X_\Gamma$, with $X_k$ (resp. $X_\Gamma$) a $B$-definable set containing $\mathrm{res}(C)$ (resp. $\mathrm{val}(C)$). By the first two items of proposition 7.0.3 applied to $B$, $X_k$ (resp. $X_\Gamma$) is $\mathrm{res}(B)$ (resp. $\mathrm{val}(B)$)-definable in the reduct $k(M)$ (resp. $\Gamma(M)$). We may write $X_k = \phi_k(x_k, \mathrm{res}(B))$, $X_\Gamma = \phi_\Gamma(x_\Gamma, \mathrm{val}(B))$, with $\phi_k(x_k, y)$ and $\phi_\Gamma(x_\Gamma, z)$ parameter-free formulas in the reducts.

Let $(F_i G_i)_{i<\omega}$ be an $A$-indiscernible sequence starting at $\mathrm{res}(B)\mathrm{val}(B)$. It suffices to show that the partial type:

$$q(x): \ \{\phi_k(x_k, F_i) \wedge \phi_\Gamma(x_\Gamma, G_i) | i < \omega\}$$

is consistent. By hypothesis, both the partial types:

$$q_k(x_k): \ \{\phi_k(x_k, F_i) | i < \omega\}$$

$$q_\Gamma(x_\Gamma): \ \{\phi_\Gamma(x_\Gamma, G_i) | i < \omega\}$$

are consistent. Choose $c_k$, $c_\Gamma$ realizations of $q_k$, $q_\Gamma$. Then $c_k c_\Gamma$ realizes $q$ by its definition, which concludes the proof. $\square$

*Remark* 7.1.2. The condition from corollary 7.0.7 always sits in-between the two conditions of lemma 7.1.1, i.e. the condition:

$$\mathrm{res}(C)\mathrm{val}(C) \mathop{\downarrow}^{\mathbf{d}}_{A} B$$

implies the condition:

$$\mathrm{res}(C)\mathrm{val}(C) \mathop{\downarrow}^{\mathbf{d}}_{A} \mathrm{res}(B)\mathrm{val}(B)$$

which implies the condition:

$$\mathrm{val}(C) \mathop{\downarrow}^{\mathbf{d}}_{A} \mathrm{val}(B), \text{ and } \mathrm{res}(C) \mathop{\downarrow}^{\mathbf{d}}_{A} \mathrm{res}(B)$$

By lemma 7.1.1, those three conditions are equivalent.

Our lemma 7.1.1 is stated more or less explicitly in the original paper of Ealy-Haskell-Simon, and a similar statement is proved explicitly in the paper of Vicaria. This equivalence can be further refined:

**Lemma 7.1.3.** *We have* $\mathrm{res}(C) \mathop{\downarrow}^{\mathbf{d}}_{A} \mathrm{res}(B)$ *if and only if, in the reduct* $k(M)$, *we have* $\mathrm{res}(C) \mathop{\downarrow}^{\mathbf{d}}_{\mathrm{res}(A)} \mathrm{res}(B)$.

*Likewise, we have* $\mathrm{val}(C) \mathop{\downarrow}^{\mathbf{d}}_{A} \mathrm{val}(B)$ *if and only if, in the reduct* $\Gamma(M)$, *we have* $\mathrm{val}(C) \mathop{\downarrow}^{\mathbf{d}}_{\mathrm{val}(A)} \mathrm{val}(B)$.

*Proof.* Suppose we have $\mathrm{res}(C) \mathop{\downarrow}^{\mathbf{d}}_{A} \mathrm{res}(B)$. Let $I$ be a $\mathrm{res}(A)$-indiscernible sequence in the reduct starting at $\mathrm{res}(B)$. By the first condition of proposition 7.0.3, $I$ is $A$-indiscernible in $M$. By hypothesis, and fact 1.1.5, there exists $J \equiv_{A \cup \mathrm{res}(B)} I$ such that $J$ is $A \cup \mathrm{res}(C)$-indiscernible. In particular, $J \equiv_{\mathrm{res}(B)} I$ in $k(M)$ and $J$ is $\mathrm{res}(C)$-indiscernible in $k(M)$, therefore $\mathrm{res}(C) \mathop{\downarrow}^{\mathbf{d}}_{\mathrm{res}(A)} \mathrm{res}(B)$ in $k(M)$.

Conversely, suppose $\mathrm{res}(C) \mathop{\downarrow}^{\mathbf{d}}_{\mathrm{res}(A)} \mathrm{res}(B)$ in $k(M)$. Let $I$ be some $A$-indiscernible sequence in $M$ starting at $\mathrm{res}(B)$. Then $I$ is $\mathrm{res}(A)$-indiscernible in $M$, hence in $k(M)$ by proposition 7.0.3. Let $J$ be a sequence which is $\mathrm{res}(B)$-conjugate of $I$ in the reduct, and $\mathrm{res}(C)$-indiscernible in the reduct. Then $J$ is $C$ (and hence $\mathrm{res}(C)$)-indiscernible in $M$, and $B$ (and hence $\mathrm{res}(B)$)-conjugate of $I$. This concludes the equivalence.

The same proof goes through for $\Gamma$. □

**Corollary 7.1.4.** *With assumptions 7.0.2, and the hypothesis from theorem 7.0.5, we have* $C \mathop{\downarrow}^{\mathbf{d}}_{A} B$ *if and only if we have* $\mathrm{res}(C) \mathop{\downarrow}^{\mathbf{d}}_{\mathrm{res}(A)} \mathrm{res}(B)$ *and* $\mathrm{val}(C) \mathop{\downarrow}^{\mathbf{d}}_{\mathrm{val}(A)} \mathrm{val}(B)$ *in the respective reducts* $k(M)$ *and* $\Gamma(M)$.

This statement simplifies and refines the hypothesis and the conclusion from ([EHS23], Corollary 3.4).

## 7.2 Forking

We adopt assumptions 7.0.2. Let us adapt the result from last section to forking. We start by proving analogues of lemma 7.1.1 and lemma 7.1.3 for forking.

**Lemma 7.2.1.** *The following conditions are equivalent:*

- $\mathrm{res}(C)\mathrm{val}(C) \mathop{\downarrow}^{\mathbf{f}}_{A} B$
- $\mathrm{res}(C) \mathop{\downarrow}^{\mathbf{f}}_{A} \mathrm{res}(B)$ *and* $\mathrm{val}(C) \mathop{\downarrow}^{\mathbf{f}}_{A} \mathrm{val}(B)$

*Proof.* The top-to-bottom direction is trivial.

Conversely, suppose $\mathrm{res}(C) \mathop{\downarrow}^{\mathbf{f}}_{A} \mathrm{res}(B)$ and $\mathrm{val}(C) \mathop{\downarrow}^{\mathbf{f}}_{A} \mathrm{val}(B)$. Then there exist $F \equiv_{A\cup\mathrm{res}(B)} \mathrm{res}(C)$, $G \equiv_{A\cup\mathrm{val}(B)} \mathrm{val}(C)$ such that $F \mathop{\downarrow}^{\mathbf{d}}_{A} M$ and $G \mathop{\downarrow}^{\mathbf{d}}_{A} M$. By the first two conditions of proposition 7.0.3 applied to $B$, we have $F \equiv_B \mathrm{res}(C)$ and $G \equiv_B \mathrm{val}(C)$. By orthogonality, we have $FG \equiv_B \mathrm{res}(C)\mathrm{val}(C)$. In some sufficiently strongly homogeneous elementary extension $N$ of $M$, let $\sigma \in \mathrm{Aut}(N/B)$ be such that $\sigma(F) = \mathrm{res}(C)$, $\sigma(G) = \mathrm{val}(C)$. Let $D = \sigma^{-1}(C)$. Then we have $D \equiv_B C$, $\mathrm{res}(D) = F$ and $\mathrm{val}(D) = G$. By lemma 7.1.1 applied to $D$ and $M$, we have $FG \mathop{\downarrow}^{\mathbf{d}}_{A} M$, which implies $\mathrm{res}(C)\mathrm{val}(C) \mathop{\downarrow}^{\mathbf{f}}_{A} B$, concluding the proof. □

**Lemma 7.2.2.** *We have* $\mathrm{res}(C) \mathop{\downarrow}^{\mathbf{f}}_{A} \mathrm{res}(B)$ *if and only if, in the reduct* $k(M)$, *we have* $\mathrm{res}(C) \mathop{\downarrow}^{\mathbf{f}}_{\mathrm{res}(A)} \mathrm{res}(B)$.

*Likewise, we have* $\mathrm{val}(C) \mathop{\downarrow}^{\mathbf{f}}_{A} \mathrm{val}(B)$ *if and only if, in the reduct* $\Gamma(M)$, *we have* $\mathrm{val}(C) \mathop{\downarrow}^{\mathbf{f}}_{\mathrm{val}(A)} \mathrm{val}(B)$.

*Proof.* Suppose $\mathrm{res}(C) \mathop{\downarrow}^{\mathbf{f}}_{A} \mathrm{res}(B)$. In some elementary extension of $M$, let $\sigma$ be an automorphism leaving $A \cup \mathrm{res}(B)$ pointwise-invariant, such that $\sigma(\mathrm{res}(C)) \mathop{\downarrow}^{\mathbf{d}}_{A} M$. By lemma 7.1.3 applied to $\sigma(C)$ and $M$, we have in the reduct: $\mathrm{res}(\sigma(C)) \mathop{\downarrow}^{\mathbf{d}}_{\mathrm{res}(A)} k(M)$. As $\sigma$ restricts to an automorphism of the reduct (of the elementary extension of $M$) leaving $\mathrm{res}(B)$ pointwise-invariant, we have $\mathrm{res}(C) \mathop{\downarrow}^{\mathbf{f}}_{\mathrm{res}(A)} \mathrm{res}(B)$ in the reduct $k(M)$.

Conversely, suppose $\operatorname{res}(C) \mathop{\downfree}^{\mathbf{f}}_{\operatorname{res}(A)} \operatorname{res}(B)$ in the reduct $k(M)$. Let $F$ be a $\operatorname{res}(B)$-conjugate of $\operatorname{res}(C)$ in the reduct such that $F \mathop{\downfree}^{\mathbf{d}}_{\operatorname{res}(A)} k(M)$ in the reduct. In some elementary extension, let $D \equiv_B C$ be such that $\operatorname{res}(D) = F$. By lemma 7.1.3 applied to $D$ and $M$, we have $F \mathop{\downfree}^{\mathbf{d}}_{A} k(M)$ in $M$. By stable embeddedness of $k$, the type of $F$ over $k(M)$ implies its type over $M$, hence $F \mathop{\downfree}^{\mathbf{d}}_{A} M$. As $F \equiv_B \operatorname{res}(C)$, we have $\operatorname{res}(C) \mathop{\downfree}^{\mathbf{f}}_{A} B$, hence $\operatorname{res}(C) \mathop{\downfree}^{\mathbf{f}}_{A} \operatorname{res}(B)$, which concludes the equivalence.

The same proof goes through for $\Gamma$. □

We may now prove a refinement of ([EHS23], Corollary 3.4) without the assumption that forking equals dividing:

**Proposition 7.2.3.** *Under assumptions 7.0.2, suppose that $k$ is algebraically bounded, that* $\operatorname{res}(C)$ *is a regular extension of* $\operatorname{res}(A)$*, and that* $\operatorname{val}(C)$ *is a pure extension of* $\operatorname{val}(A)$*. With the additional assumptions from theorem 7.0.5, we have $C \mathop{\downfree}^{\mathbf{f}}_{A} B$ if and only if, in the reducts $k(M)$ and $\Gamma(M)$, we have* $\operatorname{res}(C) \mathop{\downfree}^{\mathbf{f}}_{\operatorname{res}(A)} \operatorname{res}(B)$ *and* $\operatorname{val}(C) \mathop{\downfree}^{\mathbf{f}}_{\operatorname{val}(A)} \operatorname{val}(B)$*.*

*Proof.* Suppose $C \mathop{\downfree}^{\mathbf{f}}_{A} B$. Then we have $\operatorname{res}(C)\operatorname{val}(C) \mathop{\downfree}^{\mathbf{f}}_{A} B$, and we can apply lemma 7.2.1 and lemma 7.2.2 to show that we have $\operatorname{res}(C) \mathop{\downfree}^{\mathbf{f}}_{\operatorname{res}(A)} \operatorname{res}(B)$ and $\operatorname{val}(C) \mathop{\downfree}^{\mathbf{f}}_{\operatorname{val}(A)} \operatorname{val}(B)$ in the reducts.

Conversely, suppose $\operatorname{res}(C) \mathop{\downfree}^{\mathbf{f}}_{\operatorname{res}(A)} \operatorname{res}(B)$ and $\operatorname{val}(C) \mathop{\downfree}^{\mathbf{f}}_{\operatorname{val}(A)} \operatorname{val}(B)$ in the reducts. Once again, by lemma 7.2.1 and lemma 7.2.2, we have:

$$\operatorname{res}(C)\operatorname{val}(C) \mathop{\downfree}^{\mathbf{f}}_{A} B$$

therefore there exists $D \equiv_B C$ in some elementary extension of $M$ such that $\operatorname{res}(D)\operatorname{val}(D) \mathop{\downfree}^{\mathbf{d}}_{A} M$. As $M$ is an expansion of the infinite field $k(M)$, we may use proposition 1.2.34 to show that:

$$\operatorname{res}(D) \mathop{\downfree}^{\mathbf{alg}}_{k(\operatorname{dcl}^{eq}(A))} k(M)$$

Now, $k(\operatorname{acl}^{eq}(A))$ is an algebraic extension of $\operatorname{res}(A)$ by remark 7.0.4, therefore we have:

$$\operatorname{res}(D) \mathop{\downfree}^{\mathbf{alg}}_{\operatorname{res}(A)} k(M)$$

As $D$ and $C$ are $B$-conjugates, $\mathrm{res}(D)$ is also a regular extension of $\mathrm{res}(A)$, hence it is linearly disjoint from $k(M) = \mathrm{res}(K(M))$ over $\mathrm{res}(A)$. Likewise, we have $\mathrm{val}(D) \cap \mathrm{val}(K(M)) = \mathrm{val}(A)$. As the extension $D \geqslant A$ is separated, we apply theorem 7.0.5 to $K(M)$ to get the type implication $\mathrm{tp}(K(M)/A, \mathrm{res}(D)\mathrm{val}(D)) \vDash \mathrm{tp}(K(M)/D)$, hence $D \mathop{\downarrow}^{\mathbf{d}}_{A} M$. We conclude that $C \mathop{\downarrow}^{\mathbf{f}}_{A} B$, as $C$ and $D$ are $B$-conjugates. □

Note that we need to assume algebraic boundedness of the residue field, as well as the fact that $\mathrm{res}(C)$ is a regular extension of $\mathrm{res}(A)$ to establish linear disjointness, which is a necessary condition of theorem 7.0.5. In order to adapt corollary 7.1.4 to forking, one needs to work in a setting where forking coincides with dividing over $A$, just as in the original papers of Ealy-Haskell-Simon and Vicaria. In case the residue field is $\mathrm{NTP}_2$, we point out that our work in chapter 6 (theorem 6.4.8, theorem 6.4.9, theorem 6.4.10, theorem 6.4.11) provides tools to recognize such a setting. Alternatively, the hypothesis of proposition 7.2.3 could be weakened if one could prove a variant of lemma 7.0.6 for forking. Whether such a variant holds for forking is an open question. We conjecture that it fails.

## 7.3 Invariant measures for Abhyankar extensions

In the case of pseudo-local fields, we know (see the beginning of section 1.3) that non-forking in the residue field is always witnessed by global invariant Keisler measures. This is also the case in the value group (the measure at hand is a type). This motivates the question of whether one can use an invariant Keisler measure in the residue field to build an invariant Keisler measure in the valued field, hereby witnessing non-forking via Keisler measures. In this section, we give a positive answer to this question for a special class of separated extensions of valued fields, those that are *Abhyankar*, basically those that do not involve immediate transcendental extensions (see assumptions 7.3.1 below).

**Assumptions 7.3.1.** On top of assumptions 7.0.2, suppose $k$ is algebraically bounded. Fix $c_\Gamma = (c_{\Gamma,i})_i$ a tuple of $C$ such that $(\mathrm{val}(c_{\Gamma,i}))_i$ is a $\mathbb{Q}$-basis of $\mathrm{val}(C)$ over $\Gamma(\mathrm{acl}^{eq}(A))$. Likewise, fix $c_k = (c_{k,j})_j$ a tuple of $\mathcal{O} \cap C$ such that $(\mathrm{res}(c_{k,j}))_j$ is a transcendence basis of $\mathrm{res}(C)$ over $k(\mathrm{acl}^{eq}(A))$. Recall from remark 7.0.4 that $\mathrm{res}(\mathrm{acl}(A)) = k(\mathrm{acl}^{eq}(A))$, and likewise $\mathrm{val}(\mathrm{acl}(A)) = \Gamma(\mathrm{acl}^{eq}(A))$, thus the tuple $c_k c_\Gamma$ is algebraically independent over $K(\mathrm{acl}(A))$, by fact 2.2.15.

We assume that $C \geqslant A$ is an *Abhyankar* extension, i.e. the tuple $c_\Gamma c_k$ is a transcendence basis of $C$ over $A$. More generally, in relation with fact 2.2.15, we assume that for every intermediate extension $A \leqslant D \leqslant C$ of finite transcendence degree, the transcendence degree of $D$ over $A$ equals the sum of the $\mathbb{Q}$-dimension of $\mathrm{val}(D) / \mathrm{val}(A)$ and the transcendence degree of $\mathrm{res}(D)$ over $\mathrm{res}(A)$.

*Remark* 7.3.2. Note that $(\mathrm{val}(c_{\Gamma,i}))_i$ is also a $\mathbb{Q}$-basis of $\mathrm{val}(C)$ over $\mathrm{val}(A)$, and $(\mathrm{res}(c_{k,j}))_j$ is a transcendence basis of $\mathrm{res}(C)$ over $\mathrm{res}(A)$.

**Lemma 7.3.3.** *Let $F$ be a subfield of $A$ such that for every $N > 0$, for every coset $X$ of $(\mathrm{RV}^*)^N$ in $\mathrm{RV}^*$, there exists some $a \in F^*$ such that $\mathrm{rv}(a) \in X$. Let $q$ be the type of $(\mathrm{val}(c_{\Gamma,i}))_i$ over $\mathrm{val}(A)$ in the reduct $\Gamma(M)$. Then the partial type:*

$$q((\mathrm{val}(x_i))_i) \cup \mathrm{tp}(c_\Gamma / F)$$

*implies the full type of $c_\Gamma$ over $A$.*

Before we dive into the proof, we recall technical definitions from chapter 5. Given the expansion of a pure short exact sequence of Abelian structures:

$$0 \longrightarrow A \underset{\iota}{\longrightarrow} B \underset{\nu}{\longrightarrow} C \longrightarrow 0$$

we choose an adequate set of formulas $\Phi$, a *fundamental system*, which generates all p.p.-formulas (see definition 5.1.1). For each $\phi \in \Phi$, we define $\rho_\phi : B \longrightarrow A / \phi(A)$, which is zero outside of $\phi(B) + \iota(A)$, and extends the group homomorphism:

$$\phi(B) + \iota(A) \longrightarrow (\phi(B) + \iota(A)) / \phi(B) \simeq A / \phi(A)$$

Lastly, we define $A/C$-formulas and $A/C$-types in definition 5.1.2. They basically encode the algebraic relations satisfied by the variables in the reducts $A$ and $C$.

In the setting of valued fields, the short exact sequence that we are interested in is:

$$0 \longrightarrow k^* \underset{\iota}{\longrightarrow} \mathrm{RV} \underset{\mathrm{val}}{\longrightarrow} \Gamma \longrightarrow 0$$

where the group $k^*$ is expanded to the field $k$, and the group $\Gamma^*$ is expanded to the ordered group $\Gamma$ with the constant $-\infty$. The fundamental system we consider is the set of formulas:

$$\phi_N: \ \exists y \ x = N \cdot y$$

for every $N < \omega$.

Let us now prove lemma 7.3.3.

*Proof.* We may assume $c_\Gamma$ is finite.

Let $d = (d_i)_i \equiv_F c_\Gamma$ be such that $(\mathrm{val}(d_i))_i$ realizes $q$. It suffices to show that $d \equiv_A c_\Gamma$. By corollary 6.1.6, and strong homogeneity, we may assume $\mathrm{val}(d_i) = \mathrm{val}(c_{\Gamma,i})$ for every $i$. By the same argument as in lemma 6.1.4, it suffices to prove that $(\mathrm{rv}(d_i))_i \equiv_A (\mathrm{rv}(c_{\Gamma,i}))_i$. By proposition 6.1.1, it suffices to prove that $(\mathrm{rv}(d_i))_i \equiv_{\mathrm{rv}(A^*)} (\mathrm{rv}(c_{\Gamma,i}))_i$ in the expansion of a pure short exact sequence:

$$0 \longrightarrow k \longrightarrow \mathrm{RV}^* \longrightarrow \Gamma \longrightarrow 0$$

as this first-order structure naturally interprets $(\mathrm{RV}, \mathrm{rv}(0), \cdot, \oplus)$. As $\mathrm{val}(d_i) = \mathrm{val}(c_{\Gamma,i})$ for every $i$, we have:

$$\mathrm{tp}_C((\mathrm{rv}(d_i))_i/\mathrm{rv}(A^*)) = \mathrm{tp}_C((\mathrm{rv}(c_{\Gamma,i}))_i/\mathrm{rv}(A^*))$$

By theorem 5.1.3, it suffices to prove that they have the same $A$-type. By definition of $A$-formulas (definition 5.1.2), it suffices to prove that for every $\lambda_i \in \mathbb{Z}$, for every $a \in A^*$, we have:

$$\rho_{\phi_N}\left(\mathrm{rv}\left(\frac{\prod_i d_i^{\lambda_i}}{a}\right)\right) = \rho_{\phi_N}\left(\mathrm{rv}\left(\frac{\prod_i c_{\Gamma,i}^{\lambda_i}}{a}\right)\right)$$

We have three cases:

- Suppose the $\lambda_i$ are all zero. Then the above equation is trivial.
- Suppose $N = 0$, and the $\lambda_i$ are not all zero. Then, as $(\mathrm{val}(c_{\Gamma,i}))_i$ is $\mathbb{Q}$-free over $\mathrm{val}(A)$, we have $\mathrm{val}\left(\prod_i c_{\Gamma,i}^{\lambda_i}\right) \neq \mathrm{val}(a)$, hence:

$$M \vDash \neg\phi_0\left(\mathrm{val}\left(\frac{\prod_i c_{\Gamma,i}^{\lambda_i}}{a}\right)\right)$$

therefore:

$$\rho_{\phi_0}\left(\mathrm{rv}\left(\frac{\prod_i c_{\Gamma,i}^{\lambda_i}}{a}\right)\right) = 1$$

and likewise for $d$.

- Suppose $N > 0$. Then, as the fibers of $\rho_{\phi_N}$ are unions of cosets of $\phi_N(\mathrm{RV}^*) = (\mathrm{RV}^*)^N$, it suffices to show that:

$$\mathrm{rv}\left(\prod_i \left(\frac{d_i}{c_{\Gamma,i}}\right)^{\lambda_i}\right) \in (\mathrm{RV}^*)^N$$

however, all the cosets of $(\mathrm{RV}^*)^N$ are definable over $F$. We conclude, as $d \equiv_F c_\Gamma$.

This concludes the proof. □

We use proposition 6.1.1 to reduce to the short exact sequence, but we should give credit to the original result by Basarab ([Bas91], Theorem A), where such a reduction is more explicit.

Note that lemma 7.3.3 also holds when replacing $F$ with $\mathrm{acl}^{eq}(\varnothing)$, which also contains all the cosets of $(\mathrm{RV}^*)^N$. Note also that the proof only requires that $A$ lifts enough $\varnothing$-algebraic imaginaries in $K$, and that $(\mathrm{val}(c_{\Gamma,i}))_i$ is $\mathbb{Q}$-free over the value group of $A$. This lemma can (and will) be applied to extensions of $A$ over the value group of which $(\mathrm{val}(c_{\Gamma,i}))_i$ is $\mathbb{Q}$-free.

**Proposition 7.3.4.** *If* $(\mathrm{val}(c_{\Gamma,i}))_i \mathop{\downarrow}^{\mathbf{inv}}_{\mathrm{val}(A)} \mathrm{val}(B)$ *in the reduct* $\Gamma(M)$*, then* $c_\Gamma \mathop{\downarrow}^{\mathbf{inv}}_{A} B$ *in* $M$.

*Proof.* Let $q$ be an $\mathrm{Aut}(\Gamma(M)/\mathrm{val}(A))$-invariant type in $\Gamma(M)$, with parameters in $\Gamma(M)$, such that $q$ extends $\mathrm{tp}((\mathrm{val}(c_{\Gamma,i}))_i/\mathrm{val}(B))$. As $q$ does not divide over $\mathrm{val}(A)$, its realizations must be $\mathbb{Q}$-free over $\Gamma(M)$. By lemma 7.3.3 (with $A$ replaced by $K(M)$, and $F$ replaced by $A$), the partial type:

$$r(x) = q((\mathrm{val}(x_i))_i) \cup \mathrm{tp}(c_\Gamma/A)$$

is a complete type in $M$. As $q$ is $\mathrm{Aut}(\Gamma(M)/\mathrm{val}(A))$-invariant (and hence $\mathrm{Aut}(M/A)$-invariant by proposition 7.0.3 and strong homogeneity), $r$ is similarly $\mathrm{Aut}(M/A)$-invariant. As $(\mathrm{val}(c_{\Gamma,i}))_i \mathop{\downarrow}^{\mathbf{inv}}_{\mathrm{val}(A)} \mathrm{val}(B)$, $(\mathrm{val}(c_{\Gamma,i}))_i$ must be $\mathbb{Q}$-free over $\mathrm{val}(B)$, hence, by lemma 7.3.3 applied to $B$, $r$ extends $\mathrm{tp}(c/B)$. □

**Lemma 7.3.5.** *Let $F$ be a subfield of $K(M)$ which extends $A$. In some elementary extension of $M$, let $c = (c_l)_l$ and $d = (d_l)_l$ be two tuples of $\mathcal{O}$ such that* $\mathrm{res}(c_l) = \mathrm{res}(d_l)$ *for each $l$, and such that* $(\mathrm{res}(c_l))_l$ *is algebraically independent over* $\mathrm{res}(F)$*. Then $c \equiv_F d$.*

*Proof.* Choose $l$, and suppose by induction that $c_{<l} \equiv_F d_{<l}$. We may assume $c_{<l} = d_{<l}$. Let $L = K(\mathrm{acl}(Fc_{<l}))$. It suffices to show that $c_l \equiv_L d_l$. As $c_l$ and $d_l$ are weakly $L$-generic of $\{\mathcal{O}\}$, we use proposition 6.2.6 to conclude. □

**Proposition 7.3.6.** *Suppose* $(\mathrm{res}(c_{k,j}))_j \downarrow^{\mathbf{inv}\text{-KM}}_{\mathrm{res}(A)} \mathrm{res}(B)$ *in $k(M)$. Then we have $c_k \downarrow^{\mathbf{inv}\text{-KM}}_{A} B$ in $M$.*

*Proof.* Let $\mu(y)$ be a Keisler measure in the reduct $k(M)$ which witnesses $(\mathrm{res}(c_{k,j}))_j \downarrow^{\mathbf{inv}\text{-KM}}_{\mathrm{res}(A)} \mathrm{res}(B)$ in $k(M)$. We recall that by proposition 1.3.19, every definable set which divides over $A$ must have measure zero according to any $A$-Lascar-invariant (hence in any $\mathrm{Aut}(M/A)$-invariant) Keisler measure. We also recall that a union of definable sets of measure zero also has measure zero. As a result, and by proposition 1.2.34, and the fact that the extension $\mathrm{res}(A) \leqslant k(\mathrm{dcl}^{eq}(A))$ is algebraic, with $\mu(y)$-probability one, $y$ is $k(M)$-algebraically independent. By proposition 1.3.3, $\mu$ restricts to an $\mathrm{Aut}(M/A)$-invariant, finitely additive probability measure $\nu$ on the Boolean algebra associated to the Stone space $S$ of $|c_k|$-ary $k$-types with parameters in $k(M)$ of $k(M)$-algebraically independent tuples. Let $S'$ be the Stone space of $|c_k|$-ary $K$-types with parameters in $M$ of tuples $x$ such that $\mathrm{res}(x)$ is $k(M)$-algebraically independent. By lemma 7.3.5 applied to $K(M)$, the map $\mathrm{tp}(x/M) \longmapsto \mathrm{tp}(\mathrm{res}(x)/k(M))$ is an $\mathrm{Aut}(M)$-equivariant homeomorphism $S' \longrightarrow S$. Via this map, $\nu$ induces a finitely additive probability measure $\nu'$ on the Boolean algebra associated to $S'$. By proposition 7.0.3, and as $\nu$ is $\mathrm{Aut}(k(M)/\mathrm{res}(A))$-invariant and the homeomorphism is $\mathrm{Aut}(M)$-equivariant, $\nu'$ is $\mathrm{Aut}(M/A)$-invariant. By proposition 1.3.3, $\nu'$ induces an $\mathrm{Aut}(M/A)$-invariant Keisler measure $\mu'(x)$ in $M$ giving $S'$ probability one. Let $q = \mathrm{tp}((\mathrm{res}(c_{k,j}))_j/\mathrm{res}(B))$. By lemma 7.3.5 applied to $B$, the partial type:

$$\{\phi(\mathrm{res}(x)) | q(y) \vDash \phi(y)\}$$

is complete over $B$, hence coincides with $\mathrm{tp}(c_k/B)$. As $\mu$ is supported on $q$, $\mu'$ must be supported on $\mathrm{tp}(c_k/B)$, concluding that we have $c_k \downarrow^{\mathbf{inv}\text{-KM}}_{A} B$. □

So far, in case $(\mathrm{val}(c_{\Gamma,i}))_i \mathop{\downarrow}^{\mathbf{inv}}_{\mathrm{val}(A)} \mathrm{val}(B)$ and $(\mathrm{res}(c_{k,j}))_j \mathop{\downarrow}^{\mathbf{inv}\text{-KM}}_{\mathrm{res}(A)} \mathrm{res}(B)$, we managed to build a global invariant type and an invariant Keisler measure witnessing $c_\Gamma \mathop{\downarrow}^{\mathbf{inv}}_A B$ and $c_k \mathop{\downarrow}^{\mathbf{inv}\text{-KM}}_A B$. Now, we build tools in order to glue them into an invariant Keisler measure witnessing $C \mathop{\downarrow}^{\mathbf{inv}\text{-KM}}_A B$.

**Proposition 7.3.7.** *Suppose we have* $c_\Gamma \mathop{\downarrow}^{\mathbf{f}}_A B$ *and* $c_k \mathop{\downarrow}^{\mathbf{f}}_A B$. *Then* $\mathrm{tp}(c_\Gamma/B)$ *and* $\mathrm{tp}(c_k/B)$ *are weakly orthogonal.*

*Proof.* It suffices to prove that $\mathrm{tp}(c_\Gamma/B)$ implies $\mathrm{tp}(c_\Gamma/Bc_k)$. Recall that each $c_{k,j}$ is $\mathrm{acl}^{eq}(Bc_{k,<j})$(and hence $B(c_{k,<j})$)-generic of $\{\mathcal{O}\}$. As $\mathcal{O}$ is obviously pointed in any subfield, the third condition from fact 2.2.13 applied to the extension $B(c_{k,\leqslant j}) \geqslant B(c_{k,<j})$ fails, hence $\mathrm{val}(B(c_{k,\leqslant j})) = \mathrm{val}(B(c_{k,<j}))$. By induction, we have $\mathrm{val}(B(c_k)) = \mathrm{val}(B)$. Now, as $c_\Gamma \mathop{\downarrow}^{\mathbf{f}}_A B$, the values of (any $B$-conjugate of) $c_\Gamma$ are $\mathbb{Q}$-free over $\mathrm{val}(B)$. We conclude the proof by applying lemma 7.3.3 with $A$ replaced by $B(c_k)$, $F$ replaced by $B$. ☐

Let us explain some terminology that we use in the following proposition. Let $x$, $y$ be two disjoint tuples of variables. Given a definable set $X \subseteq M^x$, we identify $X$ with $X \times M^y \subseteq M^x \times M^y$. Given two Keisler measures $\mu(x)$, $\nu(y)$ and $\eta(xy)$ on the variables $x$, $y$, $xy$ respectively, we say that $\eta(xy)$ *extends* $\mu(x)$ if for every definable set $X \subseteq M^x$, we have $\eta(X) = \mu(X)$. We point out that when $\eta(xy)$ extends $\mu(x)$ and $\nu(y)$, this does not necessarily mean that $\eta(xy)$ agrees with the product measure $\mu \times \nu$.

**Proposition 7.3.8.** *In any theory, let* $\mu(x) \in \mathrm{KM}^x$ *be a Keisler measure, and let* $p(y) \in S^y(M)$ *be a global type. Suppose for every* $q \in \mathrm{supp}(\mu)$, *$p$ and $q$ are weakly orthogonal. Then there exists a unique measure* $\nu(xy) \in \mathrm{KM}^{xy}$ *which extends both* $\mu(x)$ *and* $p(y)$.

*Proof.* Let us prove uniqueness. Let $\nu$ be such a measure, and let $Z$ be a definable set on the variables $xy$. Let $E(Z)$ be the set of all $q \in \mathrm{supp}(\mu)$ such that $q(x) \cup p(y) \vDash xy \in Z$. Note that $E(Z)$ is not Borel in general. For each $q \in E(Z)$, there exists by weak orthogonality some clopen (i.e. some definable set) $X(q, Z)$, such that $q(x) \vDash x \in X(q, Z)$, and $p(y) \vDash \forall x \in X(q, Z)\ xy \in Z$. Note that, if the projection of $Z$ on $y$ is not implied by $p$, then $E(Z)$ is empty.

Let us show that $\nu(Z)$ is the measure given by $\mu$ to the open set $O(Z) = \bigcup_{q \in E(Z)} X(q, Z)$. Let $F(Z)$ be the closed set $\bigcap_{q \in E(Z)} (\forall z \in X(q, Z)\ zy \in Z)$. By definition of $X(q,Z)$, we have $p(y) \vDash y \in F(Z)$. As $\nu$ extends $p$, we have $\nu(F(Z)) = 1$. Moreover, we clearly have $O(Z) \times F(Z) \subseteq Z$, therefore:

$$\nu(Z) \geqslant \nu(O(Z) \times F(Z)) = \nu(O(Z)) = \mu(O(Z))$$

Let us show that $\nu(Z \setminus (O(Z) \times F(Z))) = 0$. If not, then, as $\nu(\mathrm{supp}(\nu)) = 1$, there must exist $r \in \mathrm{supp}(\nu)$ such that $r(xy) \vDash xy \in Z \setminus (O(Z) \times F(Z))$. Let $q$ be the projection of $r$ on the coordinates $x$, which must be in $\mathrm{supp}(\mu)$. Then $r$ coincides with the complete type $q \cup p$, in particular $q \in E(Z)$. However, $r \notin O(Z) \times F(Z)$, and $p \in F(Z)$, therefore $q \notin O(Z)$, which implies that $q \notin X(q,Z)$, a contradiction. It follows that $\nu(Z) = \mu(O(Z))$, and we proved uniqueness. This also shows that $\mu(O(Z))$ does not depend on the choice of $X(q,Z)$.

To prove existence, it suffices to prove that $\nu : \ Z \longmapsto \mu(O(Z))$ is a Keisler measure extending $\mu$ and $p$. If $Z$ is the true formula on the variables $xy$, then we can choose $X(q,Z)$ as the true formula on $x$, this choice is consistent with our hypothesis on $X(q,Z)$. It follows that $O(Z)$ is the true formula on $x$, therefore $\nu(Z) = 1$. If $Z$ and $Z'$ are disjoint, then, for all $q \in E(Z)$, $q' \in E(Z')$, $X(q,Z)$ and $X(q',Z')$ must be disjoint. It follows that $O(Z)$ and $O(Z')$ are disjoint, therefore $\nu(Z \cup Z') = \nu(Z) + \nu(Z')$. We proved that $\nu$ is a Keisler measure. If $Z$ only involves the variables $x$, say, $Z: \ x \in X$, then $E(Z) = X \cap \mathrm{supp}(\mu)$, and we can choose $X(q,Z) = X$, hence $O(Z) = X$, and we can conclude that $\nu$ extends $\mu$. It remains to show that $\nu$ extends $p$. Let $Y$ be a definable set in the variables $y$. If $p(y) \vDash y \notin Y$, then $E(Y) = \varnothing$, therefore $O(Y) = \varnothing$, and $\nu(Y) = 0$. Else, we have $E(Y) = \mathrm{supp}(\nu) \subseteq O(Y)$, and $\nu(Y) \geqslant \mu(\mathrm{supp}(\mu)) = 1$. This concludes the proof. □

**Definition 7.3.9.** We say that a measure and a global type are *weakly orthogonal* when they satisfy the hypothesis in proposition 7.3.8.

*Example* 7.3.10. We show a very basic example which illustrates that proposition 7.3.8 fails for two measures.

In INFSET, let $a$, $b$ be two distinct elements. Let $\mu(x)$ be the average $\frac{1}{2}(x = a) + \frac{1}{2}(x = b)$. Then, for every $p(x) \in \mathrm{supp}(\mu(x))$, $q(y) \in \mathrm{supp}(\mu(y))$, $p(x)$ and $q(y)$ are both realized, hence weakly orthogonal. However, the following are two distinct measures on $xy$ which both extend $\mu(x)$ and $\mu(y)$:

$$\nu_1(xy) = \frac{1}{2}(xy = aa) + \frac{1}{2}(xy = bb)$$

$$\nu_2(xy) = \frac{1}{2}(xy = ab) + \frac{1}{2}(xy = ba)$$

Recall that we make at the beginning of the chapter some finiteness hypothesis, and we assume that $A$ lifts some imaginaries to field elements. Recall also that we assume at the beginning of the section that the extension $C \geqslant A$ is Abhyankar.

**Theorem 7.3.11.** *Under assumptions 7.3.1, suppose* $(\mathrm{val}(c_{\Gamma,i}))_i \mathop{\downarrow}^{\mathbf{inv}}_{\mathrm{val}(A)} \mathrm{val}(B)$ *in the reduct* $\Gamma(M)$ *and* $(\mathrm{res}(c_{k,j}))_j \mathop{\downarrow}^{\mathbf{inv}\text{-KM}}_{\mathrm{res}(A)} \mathrm{res}(B)$ *in the reduct* $K(M)$. *Then we have* $C \mathop{\downarrow}^{\mathbf{inv}\text{-KM}}_{A} B$.

*Proof.* By proposition 7.3.4 and proposition 7.3.6, let $p$ be an $\mathrm{Aut}(M/A)$-invariant global extension of $\mathrm{tp}(c_\Gamma/B)$, and let $\mu$ be an $\mathrm{Aut}(M/A)$-invariant Keisler measure supported on $\mathrm{tp}(c_k/B)$. By the proof of proposition 7.3.6, we may assume that $\mu$ gives probability one to the space $S'$ of types $q(x)$ which imply that $\mathrm{res}(x)$ is $k(M)$-algebraically independent. By proposition 7.3.7 applied to $M$, $p$, and any type in $S'$, every type in $S'$ is weakly orthogonal to $p$, hence $\nu$ and $p$ are weakly orthogonal. Let $\eta(x_\Gamma, x_k)$ be the unique Keisler measure extending $\nu$ and $p$. As $\nu$ and $p$ are $\mathrm{Aut}(M/A)$-invariant, so is $\eta$. Moreover, $\eta$ is supported on $\mathrm{tp}(c_\Gamma/B)$ and $\mathrm{tp}(c_k/B)$, hence $\eta$ is supported on $\mathrm{tp}(c_\Gamma/B) \cup \mathrm{tp}(c_k/B)$. By proposition 7.3.7, $\eta$ is supported on $\mathrm{tp}(c_\Gamma c_k/B)$, therefore $c_k c_\Gamma \mathop{\downarrow}^{\mathbf{inv}\text{-KM}}_{A} B$. We conclude that $C \mathop{\downarrow}^{\mathbf{inv}\text{-KM}}_{A} B$ using proposition 1.3.18, since $C \subseteq \mathrm{acl}(A c_k c_\Gamma)$. □

**Corollary 7.3.12.** *Under assumptions 7.3.1, suppose* $k(M)$ *is pseudo-finite, and* $\downarrow^{\mathbf{f}} = \downarrow^{\mathbf{inv}}$ *over the base* $\mathrm{val}(A)$ *in the reduct* $\Gamma(M)$. *Then we have* $C \mathop{\downarrow}^{\mathbf{f}}_{A} B$ *if and only if* $C \mathop{\downarrow}^{\mathbf{inv}\text{-KM}}_{A} B$.

*Proof.* The right-to-left direction is trivial.

Conversely, suppose $C \mathop{\downarrow}^{\mathbf{f}}_{A} B$. Then, by lemma 7.2.1 and lemma 7.2.2, we have $(\mathrm{val}(c_{\Gamma,i}))_i \mathop{\downarrow}^{\mathbf{f}}_{\mathrm{val}(A)} \mathrm{val}(B)$ and $(\mathrm{res}(c_{k,j}))_j \mathop{\downarrow}^{\mathbf{f}}_{\mathrm{res}(A)} \mathrm{res}(B)$ in the respective reducts $k(M)$ and $\Gamma(M)$. By hypothesis, we have $(\mathrm{val}(c_{\Gamma,i}))_i \mathop{\downarrow}^{\mathbf{inv}}_{\mathrm{val}(A)} \mathrm{val}(B)$ in $\Gamma(M)$. Let $S$ be the space defined in the proof of proposition 7.3.6. By ([CvdDM92], Section 4), there exists an $\mathrm{Aut}(M)$-invariant Keisler measure $\mu$ such that $\mathrm{supp}(\mu) = S$. In particular, $\mu$ is supported on any partial type which admits a completion in $S$. As $\mathrm{res}(C) \underset{\mathrm{res}(A)}{\downfree^{\mathbf{f}}} \mathrm{res}(B)$, and $k(\mathrm{dcl}^{eq}(A))$ is an algebraic extension of $\mathrm{res}(A)$, we have $\mathrm{res}(C) \underset{\mathrm{res}(A)}{\downfree^{\mathbf{alg}}} \mathrm{res}(B)$, hence $\mu$ is supported on $\mathrm{tp}((\mathrm{res}(c_{k,j}))_j/\mathrm{res}(B))$. This concludes the proof by theorem 7.3.11. □

*Remark* 7.3.13. By theorem 7.0.5, the type of a separated extension is always controlled by the type of some intermediate Abhyankar extension, hence theorem 7.3.11 would hold for separated extensions if we had an analogue of lemma 7.0.6 for $\downfree^{\mathbf{inv}-\mathrm{KM}}$. Just as for forking, we do not know if such an analogue exists.

We point out that the assumption that $\downfree^{\mathbf{f}} = \downfree^{\mathbf{inv}}$ over $\mathrm{val}(A)$ is weaker than $\mathrm{val}(A)$ being an elementary substructure of $\Gamma(M)$, as ordered Abelian groups are NIP, therefore non-forking coincides with invariance over models (as explained in remark 3.1.1). We also point out (see the second part of corollary 4.4.3) that $\downfree^{\mathbf{f}} = \downfree^{\mathbf{inv}}$ over any parameter set in the ordered Abelian groups elementarily equivalent to $\mathbb{Q}$ or $\mathbb{Z}$. In particular, corollary 7.3.12 applies for free to $\mathrm{PL}_0$, and the only assumptions needed on $A$ are the ones from assumptions 7.0.2.

**Corollary 7.3.14.** *Under 7.3.1, suppose $M \vDash \mathrm{PL}_0$. Then we have $C \underset{A}{\downfree^{\mathbf{f}}} B$ if and only if $C \underset{A}{\downfree^{\mathbf{inv}-\mathrm{KM}}} B$.*

# Bibliography

[Abh56] Shreeram Abhyankar. On the valuations centered in a local domain. *American Journal of Mathematics*, 78(2):321–348, 1956.

[ACGZ22] Matthias Aschenbrenner, Artem Chernikov, Allen Gehret, and Martin Ziegler. Distality in valued fields and related structures. *Transactions of The American Mathematical Society - TRANS AMER MATH SOC*, 375(7):4641 – 4710, 2022.

[AK65a] James Ax and Simon Kochen. Diophantine problems over local fields i. *American Journal of Mathematics*, 87(3):605–630, 1965.

[AK65b] James Ax and Simon Kochen. Diophantine problems over local fields ii. a complete set of axioms for p-adic number theory.*. *American Journal of Mathematics*, 87:631, 1965.

[AK66] James Ax and Simon Kochen. Diophantine problems over local fields: Iii. decidable fields. *Annals of Mathematics*, 83:437, 1966.

[Bas91] Serban Basarab. Relative elimination of quantifiers for henselian valued fields. *Annals of Pure and Applied Logic*, 53(1):51–74, 1991.

[BYC14] Itaï Ben Yaacov and Artem Chernikov. An independence theorem for ntp2 theories. *The Journal of Symbolic Logic*, 79(1):135–153, 2014.

[Cas11] Enrique Casanovas. *Simple Theories and Hyperimaginaries*. Cambridge University Press, March 2011.

[CH11] Raf Cluckers and Immanuel Halupczok. Quantifier elimination in ordered abelian groups. *Confluentes Mathematici*, 03(4):587 – 615, October 2011.

[Che76] Gregory Cherlin. *Model theoretic algebra: Selected topics.* Springer Berlin, Heidelberg, 1976.

[Che14] Artem Chernikov. Theories without the tree property of the second kind. *Annals of Pure and Applied Logic*, 165(2):695 – 723, February 2014.

[CHK+23] Artem Chernikov, Ehud Hrushovski, Alex Kruckman, Krzysztof Krupiński, Slavko Moconja, Anand Pillay, and Nicholas Ramsey. Invariant measures in simple and in small theories. *Journal of Mathematical Logic*, 23(02):2250025, 2023.

[CK12] Artem Chernikov and Itay Kaplan. Forking and dividing in ntp2 theories. *J. Symb. Log.*, 77(1):1 – 20, March 2012.

[CK23] Gabriel Conant and Alex Kruckman. Three surprising instances of dividing. `https://arxiv.org/abs/2311.00609`, 2023.

[CvdDM92] Zoé Chatzidakis, Lou van den Dries, and Angus Macintyre. Definable sets over finite fields. *Journal für die reine und angewandte Mathematik (Crelles Journal)*, 1992:107 – 136, 1992.

[Del82] Francoise Delon. *Quelques propriétés des corps valués en théorie des modeles.* PhD thesis, Paris VII, 1982.

[Dol04] Alfred Dolich. Forking and independence in o-minimal theories. *Journal of Symbolic Logic*, 69, 03 2004.

[EHS23] Clifton Ealy, Deirdre Haskell, and Pierre Simon. Residue field domination in some henselian valued fields. *Model Theory*, 2(2):255–284, October 2023.

[Ers65] Yuri Ershov. On elementary theories of local fields. *Algebra i Logika*, 4:5–30, 01 1965.

[Far17] Rafel Farré. Strong ordered abelian groups and dp-rank. `https://arxiv.org/abs/1706.05471`, 06 2017.

[Fle11] Joseph Flenner. Relative decidability and definability in henselian valued fields. *Journal of Symbolic Logic*, 76(4):1240 – 1260, December 2011.

[Gra56] Karen Gravett. Ordered abelian groups. *Quarterly Journal of Mathematics*, 7:57–63, 1956.

[GS84] Yuri Gurevich and Peter Schmitt. The theory of ordered abelian groups does not have the independence property. *Transactions of The American Mathematical Society - TRANS AMER MATH SOC*, 284(1):171 – 182, July 1984.

[HHM05] Deirdre Haskell, Ehud Hrushovski, and Dugald Macpherson. *Stable domination and independence in algebraically closed valued fields*. Cambridge University Press, December 2005.

[HK06] Ehud Hrushovski and David Kazhdan. Integration in valued fields. In Victor Ginzburg, editor, *Algebraic Geometry and Number Theory: In Honor of Vladimir Drinfeld's 50th Birthday*, pages 261–405. Birkhäuser Boston, Boston, MA, 2006.

[HKP20] Ehud Hrushovski, Krzysztof Krupiński, and Anand Pillay. On first order amenability. *arXiv: Logic*, 2020.

[HL16] Ehud Hrushovski and François Loeser. *Non-Archimedean Tame Topology and Stably Dominated Types*. Princeton University Press, 2016.

[HMRC18] Ehud Hrushovski, Ben Martin, Silvain Rideau, and Raf Cluckers. Definable equivalence relations and zeta functions of groups. *Journal of the European Mathematical Society*, 20(10):2467 – 2537, July 2018.

[Hol95] Jan Holly. Canonical forms for definable subsets of algebraically closed and real closed valued fields. *Journal of Symbolic Logic*, 60(3):843–860, 1995.

[Hos23a] Akash Hossain. Extension bases in henselian valued fields. `https://arxiv.org/abs/2210.01567`, 2023.

[Hos23b] Akash Hossain. Forking and invariant types in regular ordered abelian groups. `https://arxiv.org/abs/2312.12279`, 2023.

[HP07] Ehud Hrushovski and Anand Pillay. On nip and invariant measures. *Journal of the European Mathematical Society*, 13, 11 2007.

[Joh20] Will Johnson. Forking and dividing in fields with several orderings and valuations. *J. Math. Log.*, 22:2150025:1–2150025:43, 2020.

[Kei87] Howard Keisler. Measures and forking. *Ann. Pure Appl. Log.*, 34:119–169, 1987.

[KP97] Byunghan Kim and Anand Pillay. Simple theories. *Annals of Pure and Applied Logic*, 88(2):149–164, 1997. Joint AILA-KGS Model Theory Meeting.

[KP98] Byunghan Kim and Anand Pillay. From stability to simplicity. *Bulletin of Symbolic Logic*, 4(1):17 – 36, March 1998.

[KS22] Salma Kuhlmann and Michele Serra. The automorphism group of a valued field of generalised formal power series. *Journal of Algebra*, 605:339 – 376, September 2022.

[Lan05] Serge Lang. *Algebra*. Graduate Texts in Mathematics. Springer New York, 2005.

[Men20] Rosario Mennuni. The domination monoid in o-minimal theories. *J. Math. Log.*, 22:2150030:1–2150030:36, 2020.

[PS23] Anand Pillay and Atticus Stonestrom. Forking and invariant measures in nip theories. `https://arxiv.org/abs/2307.11037`, 2023.

[Rob56] Abraham Robinson. *Complete Theories*. Studies in logic and the foundations of mathematics. North-Holland, 1956.

[RZ60] Abraham Robinson and Elias Zakon. Elementary properties of ordered abelian groups. *Transactions of the American Mathematical Society*, 96:222–236, 1960.

[She82] Saharon Shelah. Classification theory and the number of non-isomorphic models. *Journal of Symbolic Logic*, 47(3):694–696, 1982.

[Sim11] Pierre Simon. On dp-minimal ordered structures. *The Journal of Symbolic Logic*, 76:448 – 460, 2011.

[Sim13] Pierre Simon. Distal and non-distal nip theories. *Annals of Pure and Applied Logic*, 164(3):294–318, 2013.

[Sim15] Pierre Simon. *A Guide to NIP Theories*. Lecture Notes in Logic. Cambridge University Press, 2015.

[TZ12] Katrin Tent and Martin Ziegler. *A Course in Model Theory*. Cambridge University Press, June 2012.

[vdD14] Lou van den Dries. Lectures on the model theory of valued fields. In *Model Theory in Algebra, Analysis and Arithmetic: Cetraro, Italy 2012, Editors: H. Dugald Macpherson, Carlo Toffalori*, pages 55 – 157. Springer Berlin Heidelberg, January 2014.

[Vic21] Mariana Vicaría. Residue field domination in henselian valued fields of equicharacteristic zero. `https://arxiv.org/abs/2109.08243`, 2021.

[Wei81] Volker Weispfenning. Elimination of quantifiers for certain ordered and lattice-ordered abelian groups. *Bull. Soc. Math. Belg.*, 33:131 – 156, 1981.

[Wei84] Volker Weispfenning. Quantifier elimination and decision procedures for valued fields. In Gert H. Müller and Michael M. Richter, editors, *Models and Sets*, pages 419–472, Berlin, Heidelberg, 1984. Springer Berlin Heidelberg.

[Wei95] André Weil. *Basic number theory*. Classics in Mathematics. Springer-Verlag, Berlin, 1995. Reprint of the second (1973) edition.

[Zak61] Elias Zakon. Generalized archimedean groups. *Transactions of the American Mathematical Society*, 99:21–40, 1961.